\documentclass[leqno]{article}

\def\supp{\mathop{\rm supp}}
\def\re{\mathop{\rm Re}}
\def\im{\mathop{\rm Im}}
\def\con{\mathop{\rm Con}}
\def\pol{\mathop{\rm Pol}}
\def\hol{\mathop{\rm Hol}}
\def\Hom{\mathop{\rm Hom}}

\newtheorem{theorem}{Theorem}
\newtheorem{lemma}[theorem]{Lemma}
\newtheorem{proposition}[theorem]{Proposition}
\newtheorem{definition}[theorem]{Definition}
\newtheorem{corollary}[theorem]{Corollary}

\newcommand{\begintheorem}{\addtocounter{equation}{1}\begin{theorem}}
\newcommand{\beginlemma}{\addtocounter{equation}{1}\begin{lemma}}
\newcommand{\beginproposition}{\addtocounter{equation}{1}\begin{proposition}}
\newcommand{\begindefinition}{\addtocounter{equation}{1}\begin{definition}}
\newcommand{\begincorollary}{\addtocounter{equation}{1}\begin{corollary}}

\begin{document}

\title{Notes on algebras and vector spaces of functions}

\author{Stephen Semmes \\
        Rice University}

\date{}

\maketitle

\begin{abstract}
These informal notes are concerned with spaces of functions in various
situations, including continuous functions on topological spaces,
holomorphic functions of one or more complex variables, and so on.
\end{abstract}

\tableofcontents

\part{Elements of functional analysis}

\section{Norms and seminorms}
\label{norms, seminorms}
\setcounter{equation}{0}

        Let $V$ be a vector space over the real numbers ${\bf R}$ or
complex numbers ${\bf C}$.  A nonnegative real-valued function $N(v)$
on $V$ is said to be a \emph{seminorm} on $V$ if
\begin{equation}
        N(t \, v) = |t| \, N(v)
\end{equation}
for every $v \in V$ and $t \in {\bf R}$ or ${\bf C}$, as appropriate,
and
\begin{equation}
        N(v + w) \le N(v) + N(w)
\end{equation}
for every $v, w \in V$.  Here $|t|$ denotes the absolute value of $t$
when $t$ is a real number, and the usual modulus of $t$ when $t$ is a
complex number.  A seminorm $N(v)$ on $V$ is said to be a \emph{norm}
if $N(v) > 0$ for every $v \in V$.  Of course, the absolute value
defines a norm on ${\bf R}$, and the modulus defines a norm on ${\bf
C}$.

        As a basic class of examples, let $E$ be a nonempty set, and
let $V$ be the vector space of real or complex-valued functions on
$E$, with respect to pointwise addition and scalar multiplication.
If $x \in E$ and $f \in V$, then
\begin{equation}
        N_x(f) = |f(x)|
\end{equation}
defines a seminorm on $V$.  Let $\ell^\infty(E)$ be the linear
subspace of $V$ consisting of bounded functions on $E$, which may be
denoted $\ell^\infty(E, {\bf R})$ or $\ell^\infty(E, {\bf C})$ to
indicate whether the functions are real or complex-valued.  It is easy
to see that
\begin{equation}
\label{||f||_infty = sup_{x in E} |f(x)|}
        \|f\|_\infty = \sup_{x \in E} |f(x)|
\end{equation}
defines a norm on $\ell^\infty(E)$.

\section{Norms and metrics}
\label{norms, metrics}
\setcounter{equation}{0}

        Let $V$ be a vector space over the real or complex numbers,
and let $\|v\|$ be a norm on $V$.  It is easy to see that
\begin{equation}
        d(v, w) = \|v - w\|
\end{equation}
defines a metric on $V$, using the corresponding properties of a norm.
More precisely, $d(v, w)$ is a nonnegative real-valued function
defined for $v, w \in V$ which is equal to $0$ if and only if $v = w$,
$d(v, w)$ is symmetric in $v$ and $w$, and
\begin{equation}
        d(v, z) \le d(v, w) + d(w, z)
\end{equation}
for every $v, w, z \in V$.  Thus open and closed subsets of $V$,
convergence of sequences, and so on may be defined as in the context
of metric spaces.

        Moreover, one can check that the topology on $V$ determined by
the metric associated to the norm is compatible with the algebraic
structure corresponding to the vector space operations.  This means
that addition of vectors is continuous as a mapping from the Cartesian
product of $V$ with itself into $V$, and that scalar multiplication is
continuous as a mapping from the Cartesian product of ${\bf R}$ or
${\bf C}$ with $V$ into $V$.  This can also be described in terms of
the convergence of a sum of two convergent sequences in $V$, and the
convergence of a product of a convergent sequence in ${\bf R}$ or
${\bf C}$ with a convergent sequence in $V$.

\section{Seminorms and topologies}
\label{seminorms, topologies}
\setcounter{equation}{0}

        Let $V$ be a real or complex vector space, and let
$\mathcal{N}$ be a collection of seminorms on $V$.  A set $U \subseteq
V$ is said to be open with respect to $\mathcal{N}$ if for each $u \in
U$ there are finitely many seminorms $N_1, \ldots, N_l \in
\mathcal{N}$ and positive real numbers $r_1, \ldots, r_l$ such that
\begin{equation}
        \{v \in V : N_j(u - v) < r_j, \, j = 1, \ldots, l\} \subseteq U.
\end{equation}
It is easy to see that this defines a topology on $V$.  If $u \in V$,
$N \in \mathcal{N}$, and $r > 0$, then one can check that the
corresponding ball
\begin{equation}
        \{v \in V : N(u - v) < r\}
\end{equation}
is an open set in $V$, using the triangle inequality.  By
construction, the collection of these open balls is a subbase for the
topology on $V$ associated to $\mathcal{N}$.

       Let us say that $\mathcal{N}$ is \emph{nice} if for every $v
\in V$ with $v \ne 0$ there is an $N \in \mathcal{N}$ such that $N(v)
> 0$.  This is equivalent to the condition that $\{0\}$ be a closed
set in $V$ with respect to the topology associated to $\mathcal{N}$,
which is to say that $V \backslash \{0\}$ is an open set in this
topology.  If $\mathcal{N}$ is nice, then the topology on $V$
associated to $\mathcal{N}$ is Hausdorff.  If $\|v\|$ is a norm on
$V$, then the collection of seminorms on $V$ consisting only of
$\|v\|$ is nice, and the corresponding topology on $V$ is the same as
the one determined by the metric associated to $\|v\|$, as in the
previous section.

        If $\mathcal{N}$ is any collection of seminorms on $V$, then
addition of vectors defines a continuous mapping from $V \times V$
into $V$, and scalar multiplication defines a continuous mapping from
${\bf R} \times V$ or ${\bf C} \times V$, as appropriate, into $V$.
Thus $V$ is a \emph{topological vector space}, at least when
$\mathcal{N}$ is nice, since it is customary to ask that $\{0\}$ be a
closed set in a topological vector space.  In particular, a vector
space with a norm is a topological vector space, with respect to the
topology determined by the metric associated to the norm, as in the
previous section.  If $V$ is the space of real or complex-valued
functions on a nonempty set $E$, and if $\mathcal{N}$ is the
collection of seminorms of the form $N_x(f) = |f(x)|$, $x \in E$, as
in Section \ref{norms, seminorms}, then $\mathcal{N}$ is a nice
collection of seminorms on $V$.  In this case, $V$ can be identified
with a Cartesian product of copies of ${\bf R}$ or ${\bf C}$, indexed
by $E$, and the topology on $V$ associated to $\mathcal{N}$ is the
same as the product topology.

\section{Convergent sequences}
\label{convergent sequences}
\setcounter{equation}{0}

        Remember that a sequence of elements $\{x_j\}_{j = 1}^\infty$ of
a topological space $X$ is said to converge to an element $x$ of $X$ if
for every open set $U$ in $X$ with $x \in U$ there is an $L \ge 1$ such that
\begin{equation}
        x_j \in U
\end{equation}
for each $j \ge L$.  If the topology on $X$ is determined by a metric
$d(x, y)$, then this is equivalent to the condition that
\begin{equation}
        \lim_{j \to \infty} d(x_j, x) = 0.
\end{equation}
Similarly, if $V$ is a real or complex vector space with a norm
$\|\cdot \|$, and if $\{v_j\}_{j = 1}^\infty$ is a sequence of
elements of $V$, then $\{v_j\}_{j = 1}^\infty$ converges to another
element $v$ of $V$ when
\begin{equation}
        \lim_{j \to \infty} \|v_j - v\| = 0.
\end{equation}
If instead the topology on $V$ is determined by a collection
$\mathcal{N}$ of seminorms on $V$, then $\{v_j\}_{j = 1}^\infty$
converges to $v$ when
\begin{equation}
        \lim_{j \to \infty} N(v_j - v) = 0
\end{equation}
for every $N \in \mathcal{N}$.  In these last two cases, $\{v_j\}_{j =
1}^\infty$ converges to $v$ if and only if $\{v_j - v\}_{j =
1}^\infty$ converges to $0$.

        A topological space $X$ has a countable local base for the
topology at $x \in X$ if there is a sequence $U_1(x), U_2(x), \ldots$
of open subsets of $X$ such that $x \in U_l(x)$ for each $l$, and for
each open set $U \subseteq X$ with $x \in U$ there is an $l \ge 1$
such that $U_l(x) \subseteq U$.  In this case, one can also ask that
$U_{l + 1}(x) \subseteq U_l(x)$ for each $l$, by replacing $U_l(x)$
with the intersection of $U_1(x), \ldots, U_l(x)$ if necessary.  Under
this condition, if $x$ is in the closure of a set $E \subseteq X$,
then there is a sequence of elements of $E$ that converges to $x$.
Otherwise, one may have to use nets or filters instead of sequences.
Of course, the limit of a convergent sequence of elements of a set $E
\subseteq X$ is in the closure of $E$ in any topological space $X$.

        If $X$ has a countable local base for the topology at each
point, then the closed subsets of $X$ can be characterized in terms of
convergent sequences, as in the previous paragraph.  Equivalently, the
topology on $X$ is determined by convergence of sequences.  If the
topology on $X$ is defined by a metric, then $X$ automatically
satisfies this condition, with $U_l(x)$ equal to the open ball
centered at $x$ with radius $1/l$.  In particular, this applies to a
real or complex vector space $V$ with a norm.

        Suppose that the topology on $V$ is given by a nice collection
$\mathcal{N}$ of seminorms.  If $\mathcal{N}$ consists of only finitely
many seminorms $N_1, \ldots, N_l$, then
\begin{equation}
        \|v\| = \max_{1 \le j \le l} N_j(v)
\end{equation}
is a norm on $V$, and the topology on $V$ associated to $\mathcal{N}$
is the same as the one associated to $\|v\|$.  If $\mathcal{N}$
consists of an infinite sequence $N_1, N_2, \ldots$ of seminorms and
$v \in V$, then
\begin{equation}
        U_l(v) = \{w \in V : N_1(v - w), \ldots, N_l(v - w) \le 1/l\}
\end{equation}
is a countable local base for the topology of $V$ at $v$.  Conversely,
suppose that $U_1, U_2, \ldots$ is a sequence of open subsets of $V$
such that $0 \in U_l$ for each $l$, and for each open set $U$ in $V$
with $0 \in U$ there is an $l \ge 1$ such that $U_l \subseteq U$.
By the definition of the topology on $V$ associated to $\mathcal{N}$,
for each $l \ge 1$ there are finitely many seminorms $N_{l, 1},
\ldots, N_{l, n_l} \in \mathcal{N}$ and positive real numbers $r_{l,
1}, \ldots, r_{l, n_l}$ such that
\begin{equation}
 \{v \in V : N_{l, j}(v) < r_{l, j}, \ j = 1, \ldots, n_l\} \subseteq U_l.
\end{equation}
If $\mathcal{N}'$ is the collection of seminorms of the form $N_{l,
j}$, $1 \le j \le n_l$, $l \ge 1$, then $\mathcal{N}'$ is a subset of
$\mathcal{N}$ with only finitely or countably many elements.  One can
also check that the topology on $V$ determined by $\mathcal{N}'$ is
the same as the topology on $V$ determined by $\mathcal{N}$.

\section{Metrizability}
\label{metrizability}
\setcounter{equation}{0}

        Let $X$ be a set, and let $\rho(x, y)$ be a nonnegative
real-valued function defined for $x, y \in X$.  We say that $\rho(x, y)$
is a \emph{semimetric} on $X$ if it satisfies the same conditions as a
metric, except that $\rho(x, y)$ may be equal to $0$ even when $x \ne y$.
Thus $\rho(x, y)$ is a semimetric if $\rho(x, x) = 0$ for each $x \in X$,
\begin{equation}
        \rho(x, y) = \rho(y, x)
\end{equation}
for every $x, y \in X$, and
\begin{equation}
        \rho(x, z) \le \rho(x, y) + \rho(y, z)
\end{equation}
for every $x, y, z \in X$.  If $N$ is a seminorm on a real or complex
vector space $V$, then
\begin{equation}
\label{rho(v, w) = N(v - w)}
        \rho(v, w) = N(v - w)
\end{equation}
defines a semimetric on $V$.

        If $\rho(x, y)$ is a semimetric on a set $X$ and $t$ is a
positive real number, then
\begin{equation}
\label{rho_t(x, y) = min (rho(x, y), t)}
        \rho_t(x, y) = \min (\rho(x, y), t)
\end{equation}
is also a semimetric on $X$.  The main point is that $\rho_t(x, y)$
also satisfies the triangle inequality, since $\rho(x, y)$ does.  If
$\rho(x, y)$ is a metric on $X$, then $\rho_t(x, y)$ is too, and they
determine the same topology on $X$.

        Let $V$ be a real or complex vector space, and let
$\mathcal{N}$ be a nice collection of seminorms on $V$.  If
$\mathcal{N}$ consists of only finitely many seminorms, then their
maximum is a norm on $V$ which determines the same topology on $V$ as
$\mathcal{N}$, as in the preceding section.  If $\mathcal{N}$ consists
of an infinite sequence of seminorms $N_1, N_2, \ldots$, then
\begin{equation}
        d(v, w) = \max_{l \ge 1} \min(N_l(v - w), 1/l)
\end{equation}
defines a metric on $V$ that determines the same topology on $V$ as
$\mathcal{N}$.  More precisely, if $v = w$, then $N_l(v - w) = 0$ for
each $l$, and so $d(v, w) = 0$.  If $v \ne w$, then $N_j(v - w) > 0$
for some $j$, because $\mathcal{N}$ is nice, and
\begin{equation}
        \min(N_l(v - w), 1/l) \le 1/l < N_j(v - w)
\end{equation}
for all but finitely many $l$, so that the maximum in the definition
of $d(v, w)$ always exists.  This also shows that $d(v, w) > 0$ when
$v \ne w$, and $d(v, w)$ is obviously symmetric in $v$ and $w$.  
It is not difficult to check that $d(v, w)$ satisfies the triangle
inequality, using the fact that
\begin{equation}
        \min (N_l(v - w), 1/l)
\end{equation}
satisfies the triangle inequality for each $l$, as in the previous
paragraphs.  If $r$ is a positive real number, then $d(v, w) < r$ if
and only if $N_l(v - w) < r$ when $l \le 1/r$, and one can use this to
show that $d(v, w)$ determines the same topology on $V$ as
$\mathcal{N}$.

        Suppose now that $\mathcal{N}$ is a nice collection of
seminorms on $V$, and that there is a countable local base for the
topology on $V$ associated to $\mathcal{N}$ at $0$.  This implies that
there is a subset $\mathcal{N}'$ of $\mathcal{N}$ with only finitely
or countably many elements that determines the same topology on $V$,
as in the preceding section.  It follows that there is a metric on $V$
that determines the same topology on $V$, as in the previous
paragraph.  Note that this metric is invariant under translations on
$V$, since it depends only on $v - w$.

\section{Comparing topologies}
\label{comparing topologies}
\setcounter{equation}{0}

        Let $V$ be a real or complex vector space, and let
$\mathcal{N}$, $\mathcal{N}'$ be collections of seminorms on $V$.
Suppose that every open set in $V$ with respect to $\mathcal{N}'$ is
also an open set with respect to $\mathcal{N}$.  If $N' \in
\mathcal{N}'$, then it follows that the open unit ball with respect to
$N'$ is an open set with respect to $\mathcal{N}$.  This implies that
there are finitely many seminorms $N_1, \ldots, N_l \in \mathcal{N}$
and positive real numbers $r_1, \ldots, r_l$ such that
\begin{equation}
\label{{v in V : N_j(v) < r_j, j = 1, ldots, l} subseteq {v in V : N'(v) < 1}}
        \{v \in V : N_j(v) < r_j, \ j = 1, \ldots, l\}
         \subseteq \{v \in V : N'(v) < 1\},
\end{equation}
since $0$ is an element of the open unit ball corresponding to $N'$.
Equivalently,
\begin{equation}
 N'(v) < 1 \quad\hbox{when}\quad \max_{1 \le j \le l} r_j^{-1} \, N_j(v) < 1,
\end{equation}
and so
\begin{equation}
        N'(v) \le \max_{1 \le j \le l} r_j^{-1} \, N_j(v)
\end{equation}
for every $v \in V$.  This implies in turn that
\begin{equation}
\label{N'(v) le C max_{1 le j le l} N_j(v)}
        N'(v) \le C \max_{1 \le j \le l} N_j(v)
\end{equation}
for every $v \in V$, where $C$ is the maximum of $r_1^{-1}, \ldots,
r_l^{-1}$.  Conversely, if for every $N' \in \mathcal{N}'$ there are
finitely many seminorms $N_1, \ldots, N_l \in \mathcal{N}$ such that
(\ref{N'(v) le C max_{1 le j le l} N_j(v)}) holds for some $C \ge 0$,
then every open set with respect to $\mathcal{N}'$ is also open with
respect to $\mathcal{N}$.  Of course, one can interchange the roles of
$\mathcal{N}$ and $\mathcal{N}'$, so that they determine the same
topology on $V$ if and only if $\mathcal{N}$ and $\mathcal{N}'$ both
satisfy this condition relative to the other.

        Let us apply this to the case where $\mathcal{N}'$ consists of
a single norm $\|v\|$.  If every open set in $V$ with respect to this
norm is also an open set with respect to $\mathcal{N}$, then there are
finitely many seminorms $N_1, \ldots, N_l \in \mathcal{N}$ such that
\begin{equation}
        \|v\| \le C \max_{1 \le j \le l} N_j(v)
\end{equation}
for some $C > 0$ and every $v \in V$.  In particular,
\begin{equation}
        \|v\|' = \max_{1 \le j \le l} N_j(v)
\end{equation}
is also a norm on $V$ in this case.  Similarly, if every open set in
$V$ with respect to $\mathcal{N}$ is also an open set with respect to
$\|v\|$, then for each $N \in \mathcal{N}$ there is a $C(N) \ge 0$ such that
\begin{equation}
        N(v) \le C(N) \, \|v\|
\end{equation}
for every $v \in V$.  If the topologies on $V$ associated to
$\mathcal{N}$ and $\|v\|$ are the same, then $\|v\|'$ also determines
the same topology on $V$.

        As a basic class of examples, let $V$ be the vector space of
real or complex-valued functions on a nonempty set $E$, and let
$\mathcal{N}$ be the collection of seminorms on $V$ of the form
$N_x(f) = |f(x)|$, $x \in E$.  If there is a norm $\|v\|$ on $V$ such
that the open unit ball in $V$ with respect to $\|v\|$ is an open set
with respect to $\mathcal{N}$, then it follows that the maximum of
finitely many elements of $\mathcal{N}$ is a norm on $V$, as in the
previous paragraphs.  This implies that $E$ has only finitely many
elements.  Conversely, if $E$ has only finitely many elements, then
the maximum of $N_x(f)$, $x \in E$, is a norm on $V$ that determines
the same topology.  Note that the topology on $V$ is metrizable if and
only if $E$ has only finitely or countably many elements, as in the
preceding section.

        Now let $E$ be the set ${\bf Z}_+$ of positive integers, and
let $V$ be the vector space of real or complex-valued functions on
${\bf Z}_+$ that are rapidly decreasing in the sense that $f(j)$ is
bounded by a constant multiple of $j^{-k}$ for each nonnegative
integer $k$.  Put
\begin{equation}
\label{N_k(f) = sup_{j ge 1} j^k |f(j)|}
        N_k(f) = \sup_{j \ge 1} j^k \, |f(j)|
\end{equation}
for each $k \ge 0$, which is a norm on $V$ that reduces to the
$\ell^\infty$ norm when $k = 0$ and is monotone increasing in $k$.  It
is easy to see that the topology on $V$ associated to this collection
of norms is not determined by finitely many of these norms.  Hence the
topology on $V$ associated to this collection of norms is not
determined by any single norm at all.  However, this topology is
metrizable, as in the preceding section.

\section{Continuous linear functionals}
\label{continuous linear functionals}
\setcounter{equation}{0}

        Let $V$ be a real or complex vector space with a nice
collection of seminorms $\mathcal{N}$.  As usual, a linear functional
on $V$ is a linear mapping from $V$ into the real or complex numbers,
as appropriate.  Let $V^*$ be the space of linear functionals on $V$
that are continuous with respect to the topology on $V$ determined by
$\mathcal{N}$.  This may be described as the topological dual of $V$,
to distinguish it from the algebraic dual of all linear functionals on
$V$.  These dual spaces are also vector spaces over the real or
complex numbers, as appropriate, using pointwise addition and scalar
multiplication of functions.

        If $\lambda \in V^*$, then the set of $v \in V$ such that
$|\lambda(v)| < 1$ is open, because $\lambda$ is continuous.  Of
course, $0$ is an element of this set, because $\lambda(0) = 0$.  It
follows that there are finitely many seminorms $N_1, \ldots, N_l \in
\mathcal{N}$ and positive real numbers $r_1, \ldots, r_l$ such that
\begin{equation}
        \{v \in V : N_j(v) < r_j, \ j = 1, \ldots, l\}
         \subseteq \{v \in V : |\lambda(v)| < 1\}.
\end{equation}
As in the previous section, this implies that
\begin{equation}
        |\lambda(v)| \le \max_{1 \le j \le l} r_j^{-1} \, N_j(v)
\end{equation}
for every $v \in V$.  In particular, if $C$ is the maximum of
$r_1^{-1}, \ldots, r_l^{-1}$, then
\begin{equation}
\label{|lambda(v)| le C max_{1 le j le l} N_j(v)}
        |\lambda(v)| \le C \max_{1 \le j \le l} N_j(v)
\end{equation}
for every $v \in V$.

        Conversely, suppose that $\lambda$ is a linear functional on
$V$ for which there are finitely many seminorms $N_1, \ldots, N_l \in
\mathcal{N}$ and a nonnegative real number $C$ such that
(\ref{|lambda(v)| le C max_{1 le j le l} N_j(v)}) holds.  In this case,
\begin{equation}
\label{|lambda(v) - lambda(w)| = ... le C max_{1 le j le l} N_j(v - w)}
        |\lambda(v) - \lambda(w)| = |\lambda(v - w)|
                                   \le C \max_{1 \le j \le l} N_j(v - w)
\end{equation}
for every $v, w \in V$, because $\lambda$ is linear.  It is easy to
see that $\lambda$ is continuous on $V$ with respect to the topology
associated to $\mathcal{N}$ under these conditions.  More precisely,
for each $v \in V$ and $\epsilon > 0$, we have that
\begin{equation}
        |\lambda(v) - \lambda(w)| < \epsilon
\end{equation}
for every $w \in V$ such that $N_j(v - w) < C^{-1} \, \epsilon$ for $j
= 1, \ldots, l$.  Remember that open balls defined in terms of
seminorms in $\mathcal{N}$ are automatically open sets with respect to
$\mathcal{N}$, as in Section \ref{seminorms, topologies}.

        If the topology on $V$ is determined by a single norm $\|v\|$,
then the previous discussion can be simplified.  If $\lambda$ is a
continuous linear functional on $V$, then $|\lambda(v)| < 1$ on an
open ball around $0$ in $V$.  As before, this implies that there is
a nonnegative real number $C$ such that
\begin{equation}
\label{|lambda(v)| le C ||v||}
        |\lambda(v)| \le C \, \|v\|
\end{equation}
for every $v \in V$.  Conversely, if $\lambda$ is a linear functional
on $V$ that satisfies (\ref{|lambda(v)| le C ||v||}) for some $C \ge
0$, then
\begin{equation}
        |\lambda(v) - \lambda(w)| = |\lambda(v - w)| \le C \, \|v - w\|
\end{equation}
for every $v, w \in V$, because of linearity.  This clearly implies
that $\lambda$ is continuous with respect to the metric $d(v, w) = \|v
- w\|$ associated to $V$, as in Section \ref{norms, metrics}.

\section{${\bf R}^n$ and ${\bf C}^n$}
\label{R^n, C^n}
\setcounter{equation}{0}

        Let $n$ be a positive integer, and let ${\bf R}^n$, ${\bf
C}^n$ be the space of $n$-tuples of real and complex numbers,
respectively.  As usual, these are vector spaces with respect to
coordinatewise addition and scalar multiplication.  Put
\begin{equation}
        \|v\|_\infty = \max_{1 \le j \le n} |v_j|
\end{equation}
for each $v = (v_1, \ldots, v_n) \in {\bf R}^n$ or ${\bf C}^n$.  It is
easy to see that this defines a norm on ${\bf R}^n$, ${\bf C}^n$, for
which the corresponding topology is the standard topology.  The latter
is the same as the product topology on ${\bf R}^n$, ${\bf C}^n$ as the
Cartesian product of $n$ copies of ${\bf R}$, ${\bf C}$, with their
standard topologies.

        Another simple norm on ${\bf R}^n$, ${\bf C}^n$ is given by
\begin{equation}
        \|v\|_1 = \sum_{j = 1}^n |v_j|.
\end{equation}
Note that
\begin{equation}
        \|v\|_\infty \le \|v\|_1
\end{equation}
for every $v \in {\bf R}^n$ or ${\bf C}^n$.  Similarly,
\begin{equation}
        \|v\|_1 \le n \, \|v\|_\infty
\end{equation}
for each $v \in {\bf R}^n$, ${\bf C}^n$.  It follows that $\|v\|_1$
also determines the standard topology on ${\bf R}^n$, ${\bf C}^n$.

        If $a_1, \ldots, a_n$ are real or complex numbers, then
\begin{equation}
\label{lambda(v) = sum_{j = 1}^n a_j v_j}
        \lambda(v) = \sum_{j = 1}^n a_j \, v_j
\end{equation}
defines a linear functional on ${\bf R}^n$ or ${\bf C}^n$, as
appropriate.  It is easy to see that $\lambda$ is continuous with
respect to the standard topology on ${\bf R}^n$ or ${\bf C}^n$.  Of
course, every linear functional on ${\bf R}^n$, ${\bf C}^n$ is of this
form.  More precisely, if $\lambda$ is any linear functional on ${\bf
R}^n$ or ${\bf C}^n$, then $\lambda$ can be expressed as in
(\ref{lambda(v) = sum_{j = 1}^n a_j v_j}), with
\begin{equation}
\label{a_j = lambda(e_j)}
        a_j = \lambda(e_j)
\end{equation}
for each $j$, where $e_1, \ldots, e_n$ are the standard basis vectors
in ${\bf R}^n$, ${\bf C}^n$.  These are defined by taking the $l$th
component of $e_j$ equal to $1$ when $j = l$ and $0$ otherwise,
so that
\begin{equation}
        v = \sum_{j = 1}^n v_j \, e_j
\end{equation}
for each $v \in {\bf R}^n$ or ${\bf C}^n$.

        If $N$ is any seminorm on ${\bf R}^n$ or ${\bf C}^n$, then
\begin{equation}
        N(v) = N\Big(\sum_{j = 1}^n v_j \, e_j\Big)
             \le \sum_{j = 1}^n N(e_j) \, |v_j|.
\end{equation}
This implies that
\begin{equation}
\label{N(v) le (sum_{j = 1}^n N(e_j)) ||v||_infty}
        N(v) \le \Big(\sum_{j = 1}^n N(e_j)\Big) \, \|v\|_\infty
\end{equation}
and
\begin{equation}
\label{N(v) le (max_{1 le j le n} N(e_j)) ||v||_1}
        N(v) \le \Big(\max_{1 \le j \le n} N(e_j)\Big) \, \|v\|_1
\end{equation}
for every $v \in {\bf R}^n$ or ${\bf C}^n$.  Thus $N$ is automatically
bounded by constant multiples of the basic norms $\|v\|_\infty$, $\|v\|_1$.

        Using the triangle inequality, we get that
\begin{equation}
        N(v) - N(w) \le N(v - w)
\end{equation}
and
\begin{equation}
        N(w) - N(v) \le N(v - w)
\end{equation}
for every $v, w \in {\bf R}^n$ or ${\bf C}^n$, as appropriate.  It
follows that
\begin{equation}
\label{|N(v) - N(w)| le N(v - w)}
        |N(v) - N(w)| \le N(v - w)
\end{equation}
for every $v$, $w$.  Combining this with the estimates in the previous
paragraph, we get that $N$ is continuous as a real-valued function on
${\bf R}^n$ or ${\bf C}^n$.

        Suppose now that $N$ is a norm on ${\bf R}^n$ or ${\bf C}^n$.
The set of $v \in {\bf R}^n$ or ${\bf C}^n$ with $\|v\|_\infty = 1$ is
closed and bounded, and hence compact, with respect to the standard
topology.  Because $N$ is continuous, it attains its minimum on this
set, which is therefore positive.  Hence there is a positive real
number $c$ such that
\begin{equation}
        N(v) \ge c
\end{equation}
when $\|v\|_\infty = 1$, which implies that
\begin{equation}
\label{N(v) ge c ||v||_infty}
        N(v) \ge c \, \|v\|_\infty
\end{equation}
for every $v \in {\bf R}^n$ or ${\bf C}^n$, as appropriate, by
homogeneity.  We already know from (\ref{N(v) le (sum_{j = 1}^n
N(e_j)) ||v||_infty}) that $N(v)$ is bounded from above by a constant
multiple of $\|v\|_\infty$, and we may now conclude that the topology
on ${\bf R}^n$ or ${\bf C}^n$ determined by $N$ is the same as the
standard topology.

        Let $\mathcal{N}$ be any nice collection of seminorms on ${\bf
R}^n$ or ${\bf C}^n$, and let us check that the topology on ${\bf
R}^n$ or ${\bf C}^n$ associated to $\mathcal{N}$ is the same as the
standard topology.  Let $N_1$ be an element of $\mathcal{N}$ that is
not identically zero.  If $N_1$ is a norm, then we stop, and otherwise
we choose $N_2 \in \mathcal{N}$ such that $N_2(v) > 0$ for some $v \in
{\bf R}^n$ or ${\bf C}^n$ with $v \ne 0$ and $N_1(v) = 0$.  Note that
the set of $v \in {\bf R}^n$ or ${\bf C}^n$ such that $N_1(v) = 0$ is
a proper linear subspace of ${\bf R}^n$ or ${\bf C}^n$.  If this
linear subspace contains a nonzero element, then the set of $v \in
{\bf R}^n$ or ${\bf C}^n$ such that $N_1(v) = N_2(v) = 0$ is a proper
linear subspace of it.  By repeating the process, we get finitely many
seminorms $N_1, \ldots, N_l \in \mathcal{N}$ with $l \le n$ whose
maximum defines a norm on ${\bf R}^n$ or ${\bf C}^n$, as appropriate.
The topology on ${\bf R}^n$ or ${\bf C}^n$ associated to this norm is
the same as the standard topology, as before.  It follows that the
topology on ${\bf R}^n$ or ${\bf C}^n$ associated to $\mathcal{N}$ is
the same as the standard topology, since every seminorm on ${\bf
R}^n$, ${\bf C}^n$ is bounded by a constant multiple of the usual
norms $\|v\|_\infty$,$\|v\|_1$.

\section{Weak topologies}
\label{weak topologies}
\setcounter{equation}{0}

        Let $V$ be a real or complex vector space.  If $\lambda$ is
any linear functional on $V$, then
\begin{equation}
\label{N_lambda(v) = |lambda(v)|}
        N_\lambda(v) = |\lambda(v)|
\end{equation}
defines a seminorm on $V$.  Let $\Lambda$ be a collection of linear
functionals on $V$, and let $\mathcal{N}(\Lambda)$ be the
corresponding collection of seminorms $N_\lambda$, $\lambda \in
\Lambda$.  If $\Lambda$ is nice in the sense that for each $v \in V$
with $v \ne 0$ there is a $\lambda \in \Lambda$ such that $\lambda(v)
\ne 0$, then $\mathcal{N}(\Lambda)$ is a nice collection of seminorms
on $V$.  This leads to a topology on $V$, as in Section
\ref{seminorms, topologies}, which is the weak topology associated to
$\Lambda$.

        Under these conditions, each element of $\Lambda$ is a
continuous linear functional on $V$ with respect to the weak topology
associated to $\Lambda$.  This implies that any finite linear
combination of elements of $\Lambda$ is also continuous with respect
to this topology.  Conversely, if $\lambda$ is a continuous linear
functional on $V$ with respect to the weak topology associated to
$\Lambda$, then there are finitely many elements $\lambda_1, \ldots,
\lambda_n$ of $\Lambda$ and a nonnegative real number $C$ such that
\begin{equation}
\label{|lambda(v)| le C max_{1 le j le n} |lambda_j(v)|}
        |\lambda(v)| \le C \max_{1 \le j \le n} |\lambda_j(v)|
\end{equation}
for every $v \in V$.  In particular, $\lambda(v) = 0$ when
$\lambda_j(v) = 0$ for $j = 1, \ldots, n$, and an elementary argument
in linear algebra shows that $\lambda$ can be expressed as a linear
combination of the $\lambda_j$'s.  One may wish to reduce first to the
case where the $\lambda_j$'s are linearly independent, by discarding
any that are linear combinations of the rest.

        Let $E$ be a nonempty set, and let $V$ be the vector space of
real or complex-valued functions on $E$.  Note that $\lambda_x(f) =
f(x)$ is linear functional on $V$ for each $x \in E$.  This defines a
nice collection of linear functionals on $V$, for which the
corresponding collection of seminorms has been mentioned previously.
It follows from the discussion in the previous paragraph that a linear
functional $\lambda$ on $V$ is continuous with respect to the topology
associated to this collection of seminorms if and only if it is a
finite linear combination of $\lambda_x$'s, $x \in E$.

        Let $V$ be any real or complex vector space, and let
$\mathcal{N}$ be a nice collection of seminorms on $V$.  This leads to
the corresponding dual space $V^*$ of continuous linear functionals on
$V$.  If $v \in V$ and $v \ne 0$, then there is a $\lambda \in V^*$
such that $\lambda(v) \ne 0$.  This follows from the Hahn--Banach
theorem, as in the next section.  Thus $V^*$ is itself a nice
collection of linear functionals on $V$, which determines a weak
topology on $V$ as before, also known as the weak topology associated
to $\mathcal{N}$.  Note that every open set in $V$ with respect to
this weak topology is also an open set with respect to the topology
associated to $\mathcal{N}$, because the elements of $V^*$ are
continuous with respect to the topology associated to $\mathcal{N}$.
Every element of $V^*$ is automatically continuous with respect to the
weak topology on $V$, and conversely every continuous linear
functional on $V$ with respect to the weak topology is continuous with
respect to the topology associated to $\mathcal{N}$.  Hence $V^*$ is
also the space of continuous linear functionals on $V$ with respect to
the weak topology, which follows from the earlier discussion for the
weak topology associated to any collection of linear functionals on
$V$ as well.

\section{The Hahn--Banach theorem}
\label{hahn--banach theorem}
\setcounter{equation}{0}

        Let $V$ be a real or complex vector space, and let $N$ be a
seminorm on $V$.  Also let $\lambda$ be a linear functional on a
linear subspace $W$ of $V$ such that
\begin{equation}
\label{|lambda(v)| le C N(v)}
        |\lambda(v)| \le C \, N(v)
\end{equation}
for some $C \ge 0$ and every $v \in W$.  The Hahn--Banach theorem
states that there is an extension of $\lambda$ to a linear functional
on $V$ that satisfies (\ref{|lambda(v)| le C N(v)}) for every $v \in
V$, with the same constant $C$.  We shall not go through the proof
here, but we would like to mention some aspects of it, and some
important consequences.

        Sometimes the Hahn--Banach theorem is stated only in the case
where $N$ is a norm on $V$.  This does not really matter, because
essentially the same proof works for seminorms.  Alternatively, if $N$
is a seminorm on $V$, then
\begin{equation}
\label{Z = {v in V : N(v) = 0}}
        Z = \{v \in V : N(v) = 0\}
\end{equation}
is a linear subspace of $V$.  One can begin by extending $\lambda$ to
the linear span of $W$ and $Z$ by setting
\begin{equation}
        \lambda(w + z) = \lambda(w)
\end{equation}
for every $w \in W$ and $z \in Z$, which makes sense because
$\lambda(v) = 0$ when $v$ is in $W \cap Z$, by (\ref{|lambda(v)| le C
N(v)}).  One can then reduce to the case of norms by passing to the
quotient of $V$ by $Z$.

        To prove the Hahn--Banach theorem in the real case, one first
shows that $\lambda$ can be extended to the linear span of $W$ and any
element of $V$, while maintaining (\ref{|lambda(v)| le C N(v)}).  If
$W$ has finite codimension in $V$, then one can apply this repeatedly
to extend $\lambda$ to $V$.  If $N$ is a norm on $V$ and $V$ has a
countable dense set, then one can apply this repeatedly to extend
$\lambda$ to a dense linear subspace of $V$, and then extend $\lambda$
to all of $V$ using continuity.  Otherwise, the extension of $\lambda$
to $V$ is obtained using the axiom of choice, through Zorn's lemma or
the Hausdorff maximality principle.  The complex case can be reduced
to the real case, by treating the real part of $\lambda$ as a linear
functional on $W$ as a real vector space, and then complexifying the
extension to $V$ afterwards.

        As an application, let $\mathcal{N}$ be a nice collection of
seminorms on $V$, let $u \in V$ with $u \ne 0$ be given, and choose $N
\in \mathcal{N}$ such that $N(u) > 0$.  We can define $\lambda$ on the
one-dimensional subspace $W$ of $V$ spanned by $u$ by
\begin{equation}
\label{lambda(t u) = t N(u)}
        \lambda(t \, u) = t \, N(u)
\end{equation}
for each $t \in {\bf R}$ or ${\bf C}$, as appropriate.  This satisfies
(\ref{|lambda(v)| le C N(v)}) with $C = 1$, and the Hahn--Banach
theorem implies that there is an extension of $\lambda$ to $V$ that
also satisfies (\ref{|lambda(v)| le C N(v)}) with $C = 1$.  In
particular, this extension is a continuous linear functional on $V$
with respect to the topology associated to $\mathcal{N}$ such that
$\lambda(u) \ne 0$, as in the previous section.  One can also use the
Hahn--Banach theorem to show that a closed linear subspace of $V$ with
respect to the topology associated to $\mathcal{N}$ is also closed
with respect to the weak topology.

\section{Dual norms}
\label{dual norms}
\setcounter{equation}{0}

        Let $V$ be a real or complex vector space with a norm $\|v\|$.
Remember that a linear functional $\lambda$ on $V$ is continuous with
respect to the topology associated to $\|v\|$ if and only if there is
a nonnegative real number $C$ such that
\begin{equation}
        |\lambda(v)| \le C \, \|v\|
\end{equation}
for every $v \in V$.  In this case, the dual norm $\|\lambda\|_*$ of
$\lambda$ is defined by
\begin{equation}
        \|\lambda\|_* = \sup \{|\lambda(v)| : v \in V, \ \|v\| \le 1\}.
\end{equation}
This is the same as the smallest value of $C$ for which the previous
inequality holds.  It is not difficult to check that $\|\lambda\|_*$
defines a norm on the dual space $V^*$ of continuous linear
functionals on $V$.

        If $v \in V$ and $v \ne 0$, then there is a $\lambda \in V^*$
such that $\|\lambda\|_* = 1$ and
\begin{equation}
        \lambda(v) = \|v\|.
\end{equation}
This uses the Hahn--Banach theorem, as in the previous section.  More
precisely, the argument in the previous section shows that
$\|\lambda\|_* \le 1$, and equality holds because of the value of
$\lambda(v)$.

        Suppose that $V = {\bf R}^n$ or ${\bf C}^n$ for some positive
integer $n$.  If $a \in V$, then
\begin{equation}
\label{lambda_a(v) = sum_{j = 1}^n a_j v_j}
        \lambda_a(v) = \sum_{j = 1}^n a_j \, v_j
\end{equation}
defines a linear functional on $V$, and every linear functional on $V$
is of this form.  Note that
\begin{equation}
 |\lambda_a(v)| \le \Big(\sum_{j = 1}^n |a_j|\Big) \max_{1 \le j \le n} |v_j|
                 = \|a\|_1 \, \|v\|_\infty
\end{equation}
for every $a, v \in V$, where $\|a\|_1$, $\|v\|_\infty$ are as in
Section \ref{R^n, C^n}.  This shows that the dual norm of $\lambda_a$
on $V$ with respect to $\|v\|_\infty$ is less than or equal to
$\|a\|_1$.  If one chooses $v \in V$ such that $\|v\|_\infty = 1$ and
$a_j \, v_j = |a_j|$ for each $j$, then one gets that
\begin{equation}
        |\lambda_a(v)| = \|a\|_1,
\end{equation}
and hence the dual norm of $\lambda_a$ with respect to $\|v\|_\infty$
is equal to $\|a\|_1$.  Similarly,
\begin{equation}
        |\lambda_a(v)| \le \|a\|_\infty \, \|v\|_1
\end{equation}
for every $a, v \in V$.  This shows that the dual norm of $\lambda_a$
on $V$ with respect to $\|v\|_1$ is less than or equal to
$\|a\|_\infty$, and one can check that the dual norm is equal to
$\|a\|_\infty$ using standard basis vectors for $v$ to get equality in
the previous inequality.

\section{Topological vector spaces}
\label{topological vector spaces}
\setcounter{equation}{0}

        A \emph{topological vector space} is basically a vector space
with a topology that is compatible with the vector space operations.
More precisely, let $V$ be a vector space over the real or complex
numbers, and suppose that $V$ is also equipped with a topological
structure.  In order for $V$ to be a topological vector space, the
vector space operations of addition and scalar multiplication ought to
be continuous.  Addition of vectors corresponds to a mapping from the
Cartesian product $V \times V$ of $V$ with itself into $V$, and
continuity of addition means that this mapping should be continuous,
where $V \times V$ is equipped with the product topology associated to
the given topology on $V$.  Similarly, scalar multiplication
corresponds to a mapping from ${\bf R} \times V$ or ${\bf C} \times V$
into $V$, depending on whether $V$ is a real or complex vector space.
Continuity of scalar multiplication means that this mapping is
continuous when ${\bf R} \times V$ or ${\bf C} \times V$ is equipped
with the product topology associated to the standard topology on ${\bf
R}$ or ${\bf C}$ and the given topology on $V$.  It is customary to
ask that topological vector spaces also satisfy a separation
condition, which will be mentioned in a moment.

        Note that continuity of addition implies that the translation
mapping
\begin{equation}
\label{tau_a(v) = a + v}
        \tau_a(v) = a + v
\end{equation}
is continuous as a mapping from $V$ into itself for every $a \in V$.
This implies that $\tau_a$ is actually a homeomorphism from $V$ onto
itself for each $a \in V$, since $\tau_a$ is a one-to-one mapping from
$V$ onto itself whose inverse is $\tau_{-a}$, which is also continuous
for the same reason.  In the same way, the dilation mapping
\begin{equation}
\label{delta_t(v) = t cdot v}
        \delta_t(v) = t \cdot v
\end{equation}
is a continuous mapping on $V$ for every $t \in {\bf R}$ or ${\bf C}$,
as appropriate, because of continuity of scalar multiplication.  If $t
\ne 0$, then $\delta_t$ is a one-to-one mapping from $V$ onto itself,
with inverse equal to $\delta_{1/t}$, and hence a homeomorphism.

        The additional separation condition for $V$ to be a
topological vector space is that the set $\{0\}$ consisting of the
additive identity element $0$ in $V$ be a closed set in $V$.  This
implies that every subset of $V$ with exactly one element is closed,
because of the continuity of the translation mappings.  One can also
use continuity of addition at $0$ to show that $V$ is Hausdorff under
these conditions.

        It is easy to see that ${\bf R}^n$ and ${\bf C}^n$ are
topological vector spaces with respect to their standard topologies.
If a real or complex vector space $V$ is equipped with a norm $N$,
then $V$ is a topological vector space with respect to the topology
determined by the metric determined by $N$ as in Section \ref{norms,
metrics}.  If instead $V$ is equipped with a nice collection
$\mathcal{N}$ of seminorms, then $V$ is a topological vector space
with respect to the topology defined in Section \ref{seminorms,
topologies}.  In particular, the requirement that $\mathcal{N}$ be
nice corresponds exactly to the separation condition discussed in the
previous paragraph.

        In linear algebra, one is often interested in linear mappings
between vector spaces.  Similarly, in topology, one is often
interested in continuous mappings between topological spaces.  In the
context of topological vector spaces, one is often interested in
continuous linear mappings between topological vector spaces.  This
includes continuous linear functionals from a topological vector space
into the real or complex numbers, as appropriate.  Thus the
topological dual $V^*$ of a topological vector space $V$ may be
defined as the space of continuous linear functionals on $V$, as in
Section \ref{continuous linear functionals}.

        Let $V$ and $W$ be topological vector spaces, both real or
both complex.  If $\phi$ is a one-to-one linear mapping from $V$ onto
$W$, then the inverse mapping $\phi^{-1}$ is a one-to-one linear
mapping from $W$ onto $V$ as well.  If $\phi : V \to W$ and $\phi^{-1}
: W \to V$ are also continuous, so that $\phi$ is a homeomorphism from
$V$ onto $W$, then $\phi$ is said to be an \emph{isomorphism} between
$V$ and $W$ as topological vector spaces.  It can be shown that a
finite-dimensional real or complex topological vector space of
dimension $n$ is isomorphic to ${\bf R}^n$ or ${\bf C}^n$, as
appropriate, with its standard topology.

        A topological vector space $V$ is said to be \emph{locally
convex} if there is a local base for the topology of $V$ at $0$
consisting of convex open subsets of $V$.  If the topology on $V$ is
determined by a nice collection of seminorms, then it is easy to see
that $V$ is locally convex.  Conversely, if $V$ is locally convex,
then one can show that the topology on $V$ may be described by a nice
collection of seminorms.

        In any topological space, a necessary condition for the
existence of a metric that describes the same topology is that there
be a countable local base for the topology at each point.  If a
topological vector space $V$ has a counatble local base for the
topology at $0$, then it has a countable local base for the topology
at every point, because the topology is invariant under translations.
In this case, it can be shown that there is a metric on $V$ that
describes the same topology and which is invariant under translations.
If $V$ has a countable local base for the topology at $0$ and the
topology on $V$ is determined by a nice collection of seminorms, then
only finitely or countably many seminorms are necessary to describe
the topology, as in Section \ref{convergent sequences}, and one can
get a translation-invariant metric as in Section \ref{metrizability}.

        The definition of a Cauchy sequence can be extended to
topological vector spaces, as follows.  A sequence $\{v_j\}_{j =
1}^\infty$ of elements of a topological vector space $V$ is said to be
a \emph{Cauchy sequence} if for every open set $U$ in $V$ with $0 \in
U$ there is a positive integer $L$ such that
\begin{equation}
\label{v_j - v_l in U}
        v_j - v_l \in U
\end{equation}
for every $j, l \ge L$.  If $d(v, w)$ is a metric on $V$ that
determines the given topology on $V$, and if $d(v, w)$ is invariant
under translations on $V$ in the sense that
\begin{equation}
        d(v - z, w - z) = d(v, w)
\end{equation}
for every $v, w, z \in V$, then it is easy to see that the usual
definition of a Cauchy sequence in $V$ with respect to $d(v, w)$ is
equivalent to the preceding condition using the topological vector
space structure.

        Remember that a metric space $X$ is said to be \emph{complete}
if every Cauchy sequence of elements of $X$ converges to another
element of $X$.  Similarly, let us say that a topological vector space
$V$ is \emph{sequentially complete} if every Cauchy sequence of
elements of $V$ as in the previous paragraph converges to an element
of $V$.  If there is a countable local base for the topology of $V$ at
$0$, then this is equivalent to completeness of $V$ with respect to
any translation-invariant metric that determines the same topology on
$V$.  Otherwise, one can also consider Cauchy conditions for nets or
filters on $V$.

\section{Summable functions}
\label{summable functions}
\setcounter{equation}{0}

        Let $E$ be a nonempty set, and let $f(x)$ be a real or complex
valued function on $E$.  We say that $f$ is \emph{summable} on $E$ if
the sums
\begin{equation}
\label{sum_{x in A} |f(x)|}
        \sum_{x \in A} |f(x)|
\end{equation}
over finite subsets $A$ of $E$ are uniformly bounded.  Of course, this
holds trivially when $E$ has only finitely many elements, since we can
take $A = E$.  If $E$ is the set ${\bf Z}_+$ of positive integers,
then this is equivalent to saying that $\sum_{j = 1}^\infty |f(j)|$
converges, which means that $\sum_{j = 1}^\infty f(j)$ converges
absolutely.

        We would like to define the sum
\begin{equation}
        \sum_{x \in E} f(x)
\end{equation}
when $f$ is a summable function on $E$.  Again this is trivial when
$E$ has only finitely many elements.  If $E = {\bf Z}_+$, then the sum
may be considered as a convergent infinite series, since it converges
absolutely.  If $E$ is a countably infinite set, then one can reduce
to the case where $E = {\bf Z}_+$ using an enumeration of $E$.
Different enumerations lead to the same value of the sum, because the
sum of an absolutely convergent series is invariant under
rearrangements.  If $f$ is a summable function on any infinite set
$E$, then one can check that the set of $x \in E$ such that $|f(x)|
\ge \epsilon$ has only finitely many elements for each $\epsilon > 0$.
This implies that the set of $x \in E$ such that $f(x) \ne 0$ has only
finitely or countably many elements, so that the definition of the sum
can be reduced to the previous case.

        Alternatively, if $f$ is a nonnegative real-valued summable
function on $E$, then one can define the sum over $E$ to be the
supremum of the subsums (\ref{sum_{x in A} |f(x)|}) over all finite
subsets of $E$.  If $f$ is any summable function on $E$, then $f$ can
be expressed as a linear combination of nonnegative real-valued
summable functions, so that the definition of the sum can be reduced
to that case.  It is easy to see that this approach is compatible with
the one in the previous paragraph.  The space of summable functions on
$E$ is denoted $\ell^1(E)$, or more precisely $\ell^1(E, {\bf R})$ or
$\ell^1(E, {\bf C})$ to indicate whether real or complex-valued
functions on $E$ are being used.  One can check that this is a vector
space with respect to pointwise addition and scalar multiplication,
and that the sum over $E$ defines a linear functional on $\ell^1(E)$.

        If $f \in \ell^1(E)$, then put
\begin{equation}
        \|f\|_1 = \sum_{x \in E} |f(x)|.
\end{equation}
One can check that this is a norm on $\ell^1(E)$, and that
\begin{equation}
        \biggl|\sum_{x \in E} f(x)\biggr| \le \|f\|_1.
\end{equation}
Let us say that a function $f$ on $E$ has \emph{finite support} if
$f(x) = 0$ for all but finitely many $x \in E$.  It is not difficult
to show that these functions are dense in $\ell^1$, by considering
finite sets $A \subseteq E$ for which $\sum_{x \in A} |f(x)|$
approximates $\|f\|_1$.  Of course, $\sum_{x \in E} f(x)$ reduces to a
finite sum when $f$ has finite support on $E$.  This gives another way
to look at the sum of an arbitrary summable function on $E$.  Namely,
it is the unique continuous linear functional on $\ell^1(E)$ that is
equal to the ordinary finite sum on the dense linear subspace of
functions with finite support.

\section{$c_0(E)$}
\label{c_0(E)}
\setcounter{equation}{0}

        Let $E$ be a nonempty set, and let us say that a real or
complex-valued function $f(x)$ on $E$ \emph{vanishes at infinity} if
for each $\epsilon > 0$ the set of $x \in E$ such that $|f(x)| \ge
\epsilon$ has only finitely many elements.  The space of these
functions is denoted $c_0(E)$, or $c_0(E, {\bf R})$, $c_0(E, {\bf C})$
to indicate whether real or complex-valued functions are being used.
It is easy to see that these functions are bounded, and that they form
a closed linear subspace of $\ell^\infty(E)$ with respect to the
$\ell^\infty$ norm.  Moreover, functions on $E$ with finite support
are dense in $c_0(E)$, and $c_0(E)$ is the closure of the linear
subspace of functions on $E$ with finite support in $\ell^\infty(E)$.

        If $f$ is a bounded function on $E$ and $g$ is a summable
function on $E$, then $f \, g$ is a summable function on $E$, and
\begin{equation}
        \|f \, g\|_1 \le \|f\|_\infty \, \|g\|_1.
\end{equation}
In particular,
\begin{equation}
        \lambda_g(f) = \sum_{x \in E} f(x) \, g(x)
\end{equation}
is well-defined and satisfies
\begin{equation}
        |\lambda_g(f)| \le \|f\|_\infty \, \|g\|_1.
\end{equation}
This shows that $\lambda_g$ defines a continuous linear functional on
$\ell^\infty(E)$, whose dual norm with respect to the $\ell^\infty$
norm is less than or equal to $\|g\|_1$.  The dual norm of $\lambda_g$
with respect to the $\ell^\infty$ norm is actually equal to $\|g\|_1$,
because one can choose $f \in \ell^\infty(E)$ so that $\|f\|_\infty =
1$ and $f(x) \, g(x) = |g(x)|$ for each $x \in E$.

        We can also restrict $\lambda_g$ to $c_0(E)$, to get a
continuous linear functional on $c_0(E)$ whose dual norm is less than
or equal to $\|g\|_1$.  The dual norm of $\lambda_g$ on $c_0(E)$ with
respect to the $\ell^\infty$ norm is still equal to $\|g\|_1$, but we
have to do a bit more to show that.  The problem is that the function
$f$ mentioned at the end of the previous paragraph may not vanish at
infinity on $E$.  To fix that, we can choose for each nonempty finite
set $A \subseteq E$ a function $f_A(x)$ such that $f_A(x) = 0$ when $x
\in E \backslash A$, $f_A(x) \, g(x) = |g(x)|$ when $x \in A$, and
$\|f_A\|_\infty = 1$.  Thus $f_A$ has finite support on $E$, and hence
vanishes at infinity.  By construction,
\begin{equation}
        \lambda_g(f_A) = \sum_{x \in A} |g(x)|.
\end{equation}
This shows that the dual norm of $\lambda_g$ on $c_0(E)$ is greater
than or equal to $\sum_{x \in A} |g(x)|$ for every nonempty finite set
$A \subseteq E$, and it follows that the dual norm is equal to
$\|g\|_1$, by taking the supremum over $A$.

        Suppose now that $\lambda$ is any continuous linear functional
on $c_0(E)$.  If $x \in E$, then let $\delta_x(y)$ be the function on
$E$ equal to $1$ when $x = y$ and to $0$ otherwise.  Thus $\delta_x
\in c_0(E)$, and we can put
\begin{equation}
\label{g(x) = lambda(delta_x)}
        g(x) = \lambda(\delta_x)
\end{equation}
for each $x \in E$.  If $f$ is a function on $E$ with finite support,
then $f$ can be expressed as a linear combination of $\delta$'s, and
we get that
\begin{equation}
        \lambda(f) = \sum_{x \in E} f(x) \, g(x),
\end{equation}
by linearity.  Using functions like $f_A$ in the previous paragraph,
we get that
\begin{equation}
        \sum_{x \in A} |g(x)| \le \|\lambda\|_*
\end{equation}
for every nonempty finite set $A \subseteq E$, where $\|\lambda\|_*$
is the dual norm of $\lambda$ on $c_0(E)$.  Hence $g \in \ell^1(E)$,
and $\|g\|_1 \le \|\lambda\|_*$.  We have already seen that
$\lambda(f) = \lambda_g(f)$ when $f$ has finite support on $E$, and it
follows that this holds for every $f \in c_0(E)$, since functions with
finite support are dense in $c_0(E)$, and both $\lambda$ and
$\lambda_g$ are continuous on $c_0(E)$.

\section{The dual of $\ell^1$}
\label{dual of ell^1}
\setcounter{equation}{0}

        Let $E$ be a nonempty set, and suppose that $f$ is a summable
function on $E$, and that $g$ is a bounded function on $E$.  As in the
previous section, $f \, g$ is a summable function on $E$, and
\begin{equation}
        \|f \, g\|_1 \le \|f\|_1 \, \|g\|_\infty.
\end{equation}
Hence
\begin{equation}
        \lambda_g(f) = \sum_{x \in E} f(x) \, g(x)
\end{equation}
is well-defined and satisfies
\begin{equation}
        |\lambda_g(f)| \le \|f\|_1 \, \|g\|_\infty.
\end{equation}
Thus $\lambda_g$ defines a continuous linear functional on
$\ell^1(E)$, with dual norm less than or equal to $\|g\|_\infty$.  One
can check that the dual norm of $\lambda_g$ is actually equal to
$\|g\|_\infty$, using functions $f$ on $E$ that are equal to $1$ at
one point and $0$ elsewhere.

        Conversely, suppose that $\lambda$ is a bounded linear
functional on $\ell^1(E)$.  Let $\delta_x$ be as in the previous
section, and put $g(x) = \lambda(\delta_x)$ for each $x \in E$, as
before.  Thus
\begin{equation}
        |g(x)| = |\lambda(\delta_x)| \le \|\lambda\|_* \, \|\delta_x\|_1
                                        = \|\lambda\|_*
\end{equation}
for each $x \in E$, where $\|\lambda\|_*$ is the dual norm of
$\lambda$ on $\ell^1(E)$.  This shows that $g \in \ell^\infty(E)$, and
that $\|g\|_\infty \le \|\lambda\|_*$.  In particular, $\lambda_g$ is
a continuous linear functional on $\ell^1(E)$, as in the preceding
paragraph.  If $f$ has finite support on $E$, then $f$ can be
expressed as a linear combination of $\delta$'s, and hence $\lambda(f)
= \lambda_g(f)$.  It follows that this holds for every $f \in
\ell^1(E)$, because functions with finite support are dense in
$\ell^1(E)$, and $\lambda$, $\lambda_g$ are continuous on $\ell^1$.

        Suppose now that $E$ is an infinite set, and let $c(E)$ be the
space of real or complex-valued functions $f(x)$ on $E$ that have a
limit at infinity.  This means that there is a real or complex number
$a$, as appropriate, such that for each $\epsilon > 0$,
\begin{equation}
\label{|f(x) - a| < epsilon}
        |f(x) - a| < \epsilon
\end{equation}
for all but finitely many $x \in E$.  Equivalently, $f \in c(E)$ if
there is an $a \in {\bf R}$ or ${\bf C}$ such that $f(x) - a \in
c_0(E)$.  As usual, this space may also be denoted $c(E, {\bf R})$ or
$c(E, {\bf C})$, to indicate whether real or complex-valued functions
are being used.  It is easy to see that these functions are bounded,
and that they form a closed linear subspace of $\ell^\infty(E)$ with
respect to the $\ell^\infty$ norm.

        If $f$, $a$ are as in the previous paragraph, then put
\begin{equation}
        \lim_{x \to \infty \atop x \in E} f(x) = a.
\end{equation}
It is easy to see that this limit is unique when it exists, and that
\begin{equation}
        \biggl|\lim_{x \to \infty \atop x \in E} f(x)\biggr| \le \|f\|_\infty.
\end{equation}
Thus the limit defines a continuous linear functional on $c(E)$.  The
Hahn--Banach theorem implies that there is a continuous linear
functional $L$ on $\ell^\infty(E)$ with dual norm equal to $1$ such
that $L(f)$ is equal to this limit when $f \in c(E)$.  However, one
can also check that there is no $g \in \ell^1(E)$ such that
$\lambda_g(f) = L(f)$ for every $f \in c(E)$, where $\lambda_g$ is as
in the previous section.

\section{Filters}
\label{filters}
\setcounter{equation}{0}

        A nonempty collection $\mathcal{F}$ of nonempty subsets of a
set $X$ is said to be a \emph{filter} if
\begin{equation}
        A \cap B \in \mathcal{F} \hbox{ for every } A, B \in \mathcal{F},
\end{equation}
and
\begin{eqnarray}
 && E \in \mathcal{F} \hbox{ for every } E \subseteq X \hbox{ for which} \\
 && \hbox{there is an } A \in \mathcal{F} \hbox{ such that } A \subseteq E.
                                                                \nonumber
\end{eqnarray}
Suppose that $X$ is a topological space, and that $p$ is an element of
$X$.  A filter $\mathcal{F}$ on $X$ is said to \emph{converge} to $p$ if
\begin{equation}
        U \in \mathcal{F}
\end{equation}
for every open set $U$ in $X$ with $p \in U$.  If $X$ is Hausdorff,
then the limit of a convergent filter on $X$ is unique.

        A filter $\mathcal{F}'$ on a set $X$ is said to be a
\emph{refinement} of another filter $\mathcal{F}$ on $X$ if
$\mathcal{F} \subseteq \mathcal{F}'$, as collections of subsets of
$X$.  Suppose that $X$ is a topological space, and let $p$ be an
element of $X$.  Remember that $\overline{A}$ denotes the closure in
$X$ of a subset $A$ of $X$.  If $\mathcal{F}$ is a filter on $X$ and
\begin{equation}
        p \in \bigcap_{A \in \mathcal{F}} \overline{A},
\end{equation}
then there is a refinement of $\mathcal{F}$ that converges to $p$.  To
see this, let $\mathcal{F}'$ be the collection of subsets $E$ of $X$
such that $A \cap U \subseteq E$ for some $A \in \mathcal{F}$ and open
set $U \subseteq X$ with $p \in U$.  By hypothesis, $p \in
\overline{A}$, and hence $A \cap U \ne \emptyset$ under these
conditions.  One can check that the intersection of two elements of
$\mathcal{F}'$ is also an element of $\mathcal{F}$, because of the
corresponding properties of $\mathcal{F}$ and open neighborhoods of
$p$.  We also have that $E \in \mathcal{F}'$ for every $E \subseteq X$
for which there is a $B \in \mathcal{F}'$ such that $B \subseteq E$,
by construction.  Thus $\mathcal{F}'$ is a filter on $X$, which is
clearly a refinement of $\mathcal{F}$, since we can take $U = X$.  If
$U$ is any open set in $X$ that contains $p$, then $U \in
\mathcal{F}'$, since $A \cap U \subseteq U$ for every $A \in
\mathcal{F}$.  This shows that $\mathcal{F}'$ is a filter which is a
refinement of $\mathcal{F}$ that converges to $p$, as desired.

        Conversely, suppose that $\mathcal{F}'$ is a refinement of
$\mathcal{F}$ that converges to $p$.  If $U$ is an open set in $X$
that contains $p$, then $U \in \mathcal{F}'$, and so $A \cap U \in
\mathcal{F}'$ for every $A \in \mathcal{F}'$.  Hence $A \cap U \ne
\emptyset$, and this holds in particular for every $A \in
\mathcal{F}$, because $\mathcal{F}'$ is a refinement of $\mathcal{F}$.
It follows that $p \in \overline{A}$, since this works for every open
neighborhood $U$ of $p$ in $X$.  Thus $p \in \overline{A}$ for every
$A \in \mathcal{F}$, as before.

        Now let $V$ be a real or complex topological vector space.  A
filter $\mathcal{F}$ on $V$ satisfies the \emph{Cauchy condition} if
for every open set $U$ in $V$ with $0 \in U$ there is an $A \in
\mathcal{F}$ such that
\begin{equation}
        A - A \subseteq U,
\end{equation}
where
\begin{equation}
        A - A = \{v - w : v, w \in A\}.
\end{equation}
It is easy to see that convergent filters on $V$ satisfy the Cauchy
condition, using the continuity of
\begin{equation}
\label{(v, w) mapsto v - w}
        (v, w) \mapsto v - w
\end{equation}
as a mapping from $V \times V$ into $V$.  One can say that $V$ is
complete if every Cauchy filter on $V$ converges, as in Section
\ref{topological vector spaces}.

\section{Compactness}
\label{compactness}
\setcounter{equation}{0}

        Remember that a topological space $X$ is \emph{compact} if for
every collection $\{U_\alpha\}_{\alpha \in A}$ of open subsets of $X$
such that
\begin{equation}
        X = \bigcup_{\alpha \in A} U_\alpha,
\end{equation}
there are finitely many indices $\alpha_1, \ldots, \alpha_n \in A$
such that
\begin{equation}
        X = U_{\alpha_1} \cup \cdots \cup U_{\alpha_n}.
\end{equation}
A collection $\{E_i\}_{i \in I}$ of closed subsets of $X$ is said to
have the \emph{finite intersection property} if
\begin{equation}
        E_{i_1} \cap \cdots \cap E_{i_n} \ne \emptyset
\end{equation}
for every collection $i_1, \ldots, i_n$ of finitely many indices in
$I$.  If $X$ is compact and $\{E_i\}_{i \in I}$ is a collection of
closed subsets of $X$ with the finite intersection property, then
\begin{equation}
        \bigcap_{i \in I} E_i \ne \emptyset.
\end{equation}
Otherwise, if $\bigcap_{i \in I} E_i = \emptyset$, then $U_i = X
\backslash E_i$ would be an open covering of $X$ with no finite
subcovering.

        Conversely, if $X$ is not compact, then there is an open
covering $\{U_\alpha\}_{\alpha \in A}$ of $X$ for which there is no
finite subcovering.  If $E_\alpha = X \backslash U_\alpha$ for each
$\alpha \in A$, then it is easy to see that $\{E_\alpha\}_{\alpha \in
A}$ is a collection of closed subsets of $X$ with the finite
intersection property.  However, the intersection of all of the
$E_\alpha$'s is empty, because $\{U_\alpha\}_{\alpha \in A}$ is an
open covering of $X$.  This shows that $X$ is compact when the
intersection of any collection of closed subsets of $X$ with the
finite intersection property is nonempty.

        Let $\mathcal{F}$ be a filter on $X$.  As in the previous
section, $\mathcal{F}$ has a refinement that converges to an element
of $X$ if and only if
\begin{equation}
\label{bigcap_{A in mathcal{F}} overline{A} ne emptyset}
        \bigcap_{A \in \mathcal{F}} \overline{A} \ne \emptyset.
\end{equation}
Note that $\{\overline{A} : A \in \mathcal{F}\}$ has the finite
intersection property, because of the definition of a filter and the
elementary fact that
\begin{equation}
        \overline{A \cap B} \subseteq \overline{A} \cap \overline{B}
\end{equation}
for every $A, B \subseteq X$.  If $X$ is compact, then it follows that
every filter on $X$ has a refinement that converges to an element of
$X$.

        Conversely, let $\{E_i\}_{i \in I}$ be a collection of closed
subsets of $X$ with the finite intersection property.  Let
$\mathcal{F}$ be the set of all $A \subseteq X$ such that
\begin{equation}
\label{E_{i_1} cap cdots cap E_{i_n} subseteq A}
        E_{i_1} \cap \cdots \cap E_{i_n} \subseteq A
\end{equation}
for some finite collection of indices $i_1, \ldots, i_n \in I$.  It is
easy to see that $\mathcal{F}$ is a filter on $X$.  If there is a
refinement of $\mathcal{F}$ that converges to an element of $X$, then
it follows that $\bigcap_{i \in I} E_i \ne \emptyset$.  Thus $X$ is
compact when every filter on $X$ has a refinement that converges to an
element of $X$.

\section{Ultrafilters}
\label{ultrafilters}
\setcounter{equation}{0}

        A maximal filter on a set $X$ is said to be an
\emph{ultrafilter}.  More precisely, a filter $\mathcal{F}$ on $X$ is
an ultrafilter if the only filter on $X$ that is a refinement of
$\mathcal{F}$ is itself.  If $p \in X$ and $\mathcal{F}_p$ is the
collection of subsets $A$ of $X$ such that $p \in A$, then it is easy
to see that $\mathcal{F}_p$ is an ultrafilter on $X$.  One can show
that every filter has a refinement that is an ultrafilter, using the
axiom of choice through Zorn's lemma or the Hausdorff maximality
principle.  If $X$ is a compact topological space and $\mathcal{F}$ is
an ultrafilter on $X$, then it follows that $\mathcal{F}$ converges to
an element of $X$.  More precisely, $\mathcal{F}$ has a refinement
that converges, as in the previous section, and this refinement is the
same as $\mathcal{F}$, since $\mathcal{F}$ is an ultrafilter.
Conversely, if every ultrafilter on a topological space $X$ converges,
then $X$ is compact.  This is because every filter on $X$ has a
refinement which is an ultrafilter, and hence converges.

        Suppose that $\mathcal{F}$ is a filter on a set $X$, and that
$B$ is a subset of $X$ such that $A \cap B \ne \emptyset$ for every $A
\in \mathcal{F}$.  Let $\mathcal{F}_B$ be the collection of subsets
$E$ of $X$ for which there is an $A \in \mathcal{F}$ such that
\begin{equation}
        A \cap B \subseteq E.
\end{equation}
It is easy to see that $\mathcal{F}_B$ is a filter on $X$ that is a
refinement of $\mathcal{F}$.  If $\mathcal{F}$ is an ultrafilter on
$X$, then it follows that $\mathcal{F}_B = \mathcal{F}$, and hence
that $B \in \mathcal{F}$.  Conversely, suppose that $\mathcal{F}$ is a
filter on $X$ such that $B \in \mathcal{F}$ for every $B \subseteq X$
such that $A \cap B \ne \emptyset$ for every $A \in \mathcal{F}$.  If
$\mathcal{F}'$ is a filter on $X$ that is a refinement of
$\mathcal{F}$, and if $B \in \mathcal{F}'$, then $A \cap B \in
\mathcal{F}'$ for every $A \in \mathcal{F} \subseteq \mathcal{F}'$,
and hence $A \cap B \ne \emptyset$ for every $A \in \mathcal{F}$.  It
follows that every $B \in \mathcal{F}'$ is also in $\mathcal{F}$,
which means that $\mathcal{F}' = \mathcal{F}$.  Thus $\mathcal{F}$ is
an ultrafilter under these conditions.

        Let $\mathcal{F}$ be an ultrafilter on a set $X$, and let $B$
be a subset of $X$.  If $A \cap B = \emptyset$ for some $A \in
\mathcal{F}$, then $A \subseteq X \backslash B$, and hence $X
\backslash B \in \mathcal{F}$.  Otherwise, if $A \cap B \ne \emptyset$
for every $A \in \mathcal{F}$, then $B \in \mathcal{F}$, as in the
previous paragraph.  This shows that for every $B \subseteq X$, either
\begin{equation}
        B \in \mathcal{F} \quad\hbox{or}\quad X \backslash B \in \mathcal{F}
\end{equation}
when $\mathcal{F}$ is an ultrafilter.  Conversely, if $\mathcal{F}$ is
a filter on $X$ with this property, then $\mathcal{F}$ is an
ultrafilter.  To see this, let $\mathcal{F}'$ be a filter on $X$ that
is a refinement of $\mathcal{F}$.  If $B \in \mathcal{F}'$ and $X
\backslash B \in \mathcal{F} \subseteq \mathcal{F}'$, then we get a
contradiction, since $B \cap (X \backslash B) = \emptyset$.  Thus each
$B \in \mathcal{F}'$ is an element of $\mathcal{F}$, which implies
that $\mathcal{F}' = \mathcal{F}$, as desired.

        Let $X$, $Y$ be sets, and let $f$ be a mapping from $X$ into $Y$.
If $\mathcal{F}$ is a filter on $X$, then one can check that
\begin{equation}
        f_*(\mathcal{F}) = \{A \subseteq Y : f^{-1}(A) \in \mathcal{F}\}
\end{equation}
is a filter on $Y$.  If $\mathcal{F}$ is an ultrafilter on $X$, then
$f_*(\mathcal{F})$ is an ultrafilter on $Y$.  To see this, let $B$ be
a subset of $Y$, and note that $f^{-1}(Y \backslash B) = X \backslash
f^{-1}(B)$, so that $f^{-1}(B)$ or $f^{-1}(Y \backslash B)$ is an
element of $\mathcal{F}$.  Thus $B$ or $Y \backslash B$ is an element
of $f_*(\mathcal{F})$, which implies that $\mathcal{F}$ is an
ultrafilter on $Y$, as in the previous paragraph.

\section{Tychonoff's theorem}
\label{tychonoff's theorem}
\setcounter{equation}{0}

        Let $\{X_i\}_{i \in I}$ be a collection of compact topological
spaces, and let $X = \prod_{i \in I} X_i$ be their Cartesian product.
A famous theorem of Tychonoff states that $X$ is also compact with
respect to the product topology.  There is a well-known proof of this
using ultrafilters, as follows.  Let $\mathcal{F}$ be an ultrafilter
on $X$, and let us show that $\mathcal{F}$ converges.  Let $p_i$ be
the standard coordinate projection from $X$ onto $X_i$ for each $i \in
I$.  As before, $(p_i)_*(\mathcal{F})$ is an ultrafilter on $X_i$ for
each $i \in I$, and hence converges to an element $x_i$ of $X_i$, by
compactness.  If $x \in X$ satisfies $p_i(x) = x_i$ for each $i$, then
we would like to check that $\mathcal{F}$ converges to $x$.  Let $U$
be an open set in $X$ such that $x \in U$.  By the definition of the
product topology, there are open sets $U_i \subseteq X_i$ for each $i
\in I$ such that $x_i \in U_i$ for each $i$, $U_i = X_i$ for all but
finitely many $i$, and
\begin{equation}
        \prod_{i \in I} U_i \subseteq U.
\end{equation}
Because $(p_i)_*(\mathcal{F})$ converges to $x_i$ for each $i$, we get
that $U_i \in (p_i)_*(\mathcal{F}_i)$ for each $i$, which means that
$p_i^{-1}(U_i) \in \mathcal{F}$ for each $i$.  Of course,
\begin{equation}
        \prod_{i \in I} U_i = \bigcap_{i \in I} p_i^{-1}(U_i).
\end{equation}
This is the same as the intersection of $p_i^{-1}(U_i)$ over finitely
many $i \in I$, since $U_i = X_i$ and hence $p_i^{-1}(U_i) = X$ for
all but finitely many $i$.  It follows that the intersection is
contained in $\mathcal{F}$, which implies that $U$ is contained in
$\mathcal{F}$, as desired.

\section{The weak$^*$ topology}
\label{weak^* topology}
\setcounter{equation}{0}

        Let $V$ be a real or complex topological vector space, and let
$V^*$ be the dual space of continuous linear functionals on $V$.  If
$v \in V$, then
\begin{equation}
        L_v(\lambda) = \lambda(v)
\end{equation}
defines a linear functional on $V^*$.  This is automatically a nice
collection of linear functionals on $V^*$ in the sense of Section
\ref{weak topologies}, since $\lambda = 0$ in $V^*$ when $\lambda(v) =
0$ for each $v \in V$.  The weak topology on $V^*$ corresponding to
this collection of linear functionals is known as the \emph{weak$^*$
topology}.

        Suppose now that $V$ is equipped with a norm $\|v\|$ that
determines the given topology on $V$, and let $\|\lambda\|_*$ be the
corresponding dual norm on $V^*$.  Observe that $L_v$ is a continuous
linear functional on $V^*$ with respect to $\|\lambda\|_*$ for each $v
\in V$.  More precisely,
\begin{equation}
        |L_v(\lambda)| \le \|\lambda\|_* \, \|v\|
\end{equation}
for every $v \in V$ and $\lambda \in V^*$, by definition of
$\|\lambda\|_*$.  This shows that the dual norm of $L_v$ as a
continuous linear functional on $V^*$ with respect to $\|\lambda\|_*$
is less than or equal to $\|v\|$ for each $v \in V$.  The dual norm of
$L_v$ on $V^*$ is actually equal to $\|v\|$, since for each $v \in V$
there is a $\lambda \in V^*$ such that $\|\lambda\|_* = 1$ and
$\lambda(v) = \|v\|$, by the Hahn--Banach theorem.

        Note that every open set in $V^*$ with respect to the weak$^*$
topology is also an open set with respect to the dual norm.  This
follows from the fact that $L_v$ is continuous on $V^*$ with respect
to the dual norm for each $v \in V$, as in the previous paragraph.

        Consider the closed unit ball $B^*$ in $V^*$ with respect to
the dual norm, which consists of all $\lambda \in V^*$ with
$\|\lambda\|_* \le 1$.  This is the same as the set of $\lambda \in
V^*$ such that $|\lambda(v)| \le 1$ for every $v \in V$ with $\|v\|
\le 1$.  It follows easily from this description that $B^*$ is a
closed set in the weak$^*$ topology.  The Banach--Alaoglu theorem
states that $B^*$ is a compact set with respect to the weak$^*$
topology.  The basic idea is to show that $B^*$ is homeomorphic to a
closed subset of a product of closed intervals in the real case, or a
product of closed disks in the complex case, and then use Tychonoff's
theorem.

\section{Filters on subsets}
\label{filters, subsets}
\setcounter{equation}{0}

        Let $X$ be a nonempty set, and let $E$ be a nonempty subset of
$X$.  If $\mathcal{F}_0$ is a filter on $E$, then there is a natural
filter $\mathcal{F}_1$ on $X$ associated to it, given by
\begin{equation}
        \mathcal{F}_1 = \{B \subseteq X : B \cap E \in \mathcal{F}_0\}.
\end{equation}
In particular, $\mathcal{F}_0 \subseteq \mathcal{F}_1$.  Equivalently,
if $f : E \to X$ is the inclusion mapping that sends every element of
$E$ to itself as an element of $X$, then $\mathcal{F}_1 =
f_*(\mathcal{F}_0)$.

        Conversely, if $\mathcal{F}_1$ is a filter on $X$ such that $E
\in \mathcal{F}_1$, then
\begin{equation}
\label{mathcal{F}_0 = {A subseteq E : A in mathcal{F}_1}}
        \mathcal{F}_0 = \{A \subseteq E : A \in \mathcal{F}_1\}
\end{equation}
is a filter on $E$.  It is easy to see that this transformation
between filters is the inverse of the one described in the previous
paragraph.  Thus we get a one-to-one correspondence between filters on
$E$ and filters on $X$ that contain $E$ as an element.  Moreover,
refinements of filters on $E$ correspond exactly to refinements of
filters on $X$ that contain $E$ as an element in this way.  Of course,
any refinement of a filter on $X$ that contains $E$ as an element also
contains $E$ as an element, and hence corresponds to a filter on $E$
too.

        Suppose now that $X$ is a topological space.  It is easy to
check that a filter $\mathcal{F}_1$ on $X$ that contains $E$ as an
element converges to a point $p \in E$ if and only if the
corresponding filter $\mathcal{F}_0$ on $E$ converges to $p$ with
respect to the induced topology on $E$.  Using the remarks about
refinements in the previous paragraph, it follows that $\mathcal{F}_1$
has a refinement on $X$ that converges to an element of $E$ if and
only if $\mathcal{F}_0$ has a refinement on $E$ that converges to an
element of $E$ with respect to the induced topology.  Hence $E$ is
compact if and only if every filter $\mathcal{F}_1$ on $X$ that
contains $E$ as an element has a refinement that converges to an
element of $E$.

        In the same way, ultrafilters on $E$ correspond exactly to
ultrafilters on $X$ that contain $E$ as an element, and $E$ is compact
if and only if every ultrafilter on $X$ that contains $E$ as an
element converges to an element of $E$.

\section{Bounded linear mappings}
\label{bounded linear mappings}
\setcounter{equation}{0}

        Let $V$ and $W$ be vector spaces, both real or both complex,
and equipped with norms $\|v\|_V$, $\|w\|_W$, respectively.  A linear
mapping $T : V \to W$ is said to be \emph{bounded} if there is a
nonnegative real number $A$ such that
\begin{equation}
        \|T(v)\|_W \le A \, \|v\|_V
\end{equation}
for every $v \in V$.  Because of linearity, this implies that
\begin{equation}
\label{||T(v) - T(v')||_W le A ||v - v'||_V}
        \|T(v) - T(v')\|_W \le A \, \|v - v'\|_V
\end{equation}
for every $v, v' \in V$, and hence that $T$ is uniformly continuous
with respect to the metrics on $V$ and $W$ associated to their norms.
Conversely, if $T : V \to W$ is continuous at $0$, then there is a
$\delta > 0$ such that
\begin{equation}
\label{||T(v)||_W < 1}
        \|T(v)\|_W < 1
\end{equation}
for every $v \in V$ with $\|v\|_V < \delta$.  It is easy to see that
this implies that $T$ is bounded, with $A = 1/\delta$.

        If $T$ is a bounded linear mapping from $V$ into $W$, then the
\emph{operator norm} $\|T\|_{op}$ of $T$ is defined by
\begin{equation}
        \|T\|_{op} = \sup \{\|T(v)\|_W : v \in V, \, \|v\|_V \le 1\}.
\end{equation}
The boundedness of $T$ says exactly that the supremum is finite, and
is less than or equal to the nonnegative real number $A$ mentioned in
the previous paragraph.  Equivalently, $T$ satisfies the boundedness
condition in the previous paragraph with $A = \|T\|_{op}$, and this is
the smallest value of $A$ with this property.  Note that $\|T\|_{op} =
0$ if and only if $T = 0$.

        Let $T$ be a bounded linear mapping from $V$ into $W$, and let
$a$ be a real or complex number, as appropriate.  Of course, the
product $a \, T$ of $a$ and $T$ is the linear mapping that sends $v
\in V$ to $a \, T(v)$ in $W$.  It is easy to see that $a \, T$ is also
a bounded linear mapping, and that
\begin{equation}
        \|a \, T\|_{op} = |a| \, \|T\|_{op}.
\end{equation}
Similarly, if $R$ is another bounded linear mapping from $V$ into
$W$, then the sum $R + T$ is defined as the linear mapping that sends
$v \in V$ to $R(v) + T(v)$ in $W$.  It is easy to see that $R + T$
is also a bounded linear mapping from $V$ into $W$, and that
\begin{equation}
        \|R + T\|_{op} \le \|R\|_{op} + \|T\|_{op}.
\end{equation}

        Let $\mathcal{BL}(V, W)$ be the space of bounded linear
mappings from $V$ into $W$.  It follows from the remarks in the
preceding paragraph that $\mathcal{BL}(V, W)$ is a real or complex
vector space, as appropriate, with respect to pointwise addition and
scalar multiplication of linear mappings, and that the operator norm
defines a norm on this vector space.  If $W$ is the one-dimensional
vector space of real or complex numbers, as appropriate, then
$\mathcal{BL}(V, W)$ is the same as the dual space $V^*$ of bounded
linear functionals on $V$, and the operator norm is the same as the
dual norm on $V^*$.

        Suppose that $W$ is complete as a metric space with respect to
the metric associated to the norm, so that $W$ is a \emph{Banach space}.
In this case,  the space $\mathcal{BL}(V, W)$ of bounded linear
mappings from $V$ into $W$ is also complete with respect to the
operator norm, and thus a Banach space.  To see this, let $\{T_j\}_{j
= 1}^\infty$ be a Cauchy sequence in $\mathcal{BL}(V, W)$.  This means
that for each $\epsilon > 0$ there is an $L(\epsilon) \ge 1$ such that
\begin{equation}
\label{||T_j - T_l||_{op} le epsilon}
        \|T_j - T_l\|_{op} \le \epsilon
\end{equation}
for every $j, l \ge L(\epsilon)$.  Equivalently,
\begin{equation}
\label{||T_j(v) - T_l(v)||_W le epsilon ||v||_V}
        \|T_j(v) - T_l(v)\|_W \le \epsilon \, \|v\|_V
\end{equation}
for every $j, l \ge L(\epsilon)$ and $v \in V$, so that $\{T_j(v)\}_{j
= 1}^\infty$ is a Cauchy sequence in $W$ for every $v \in V$.  Because
$W$ is complete, it follows that $\{T_j(v)\}_{j = 1}^\infty$ converges
in $W$ for every $v \in V$.  Put
\begin{equation}
        T(v) = \lim_{j \to \infty} T_j(v)
\end{equation}
for every $v \in V$.  It is easy to see that $T$ is a linear mapping
from $V$ into $W$, because of the linearity of the $T_j$'s.

        Observe that
\begin{equation}
\label{||T_j(v) - T(v)||_W le epsilon ||v||}
        \|T_j(v) - T(v)\|_W \le \epsilon \, \|v\|
\end{equation}
for every $j \ge L(\epsilon)$ and $v \in V$, by taking the limit as $l
\to \infty$ in (\ref{||T_j(v) - T_l(v)||_W le epsilon ||v||_V}).  In
particular,
\begin{equation}
        \|T(v)\|_W \le \|T_j(v)\|_W + \epsilon \, \|v\|
\end{equation}
when $j \ge L(\epsilon)$.  Applying this to $\epsilon = 1$ and $j =
L(1)$, we get that
\begin{equation}
\label{||T(v)||_W le ... le ||T_{L(1)}||_{op} ||v|| + ||v||}
        \|T(v)\|_W \le \|T_{L(1)}(v)\|_W + \|v\|
                    \le \|T_{L(1)}\|_{op} \, \|v\| + \|v\|.
\end{equation}
This implies that $T$ is bounded, with $\|T\|_{op} \le
\|T_{L(1)}\|_{op} + 1$.

        Using (\ref{||T_j(v) - T(v)||_W le epsilon ||v||}), we also get that
\begin{equation}
\label{||T_j - T||_{op} le epsilon}
        \|T_j - T\|_{op} \le \epsilon
\end{equation}
when $j \ge L(\epsilon)$.  Thus $\{T_j\}_{j = 1}^\infty$ converges to
$T$ with respect to the operator norm.  This shows that every Cauchy
sequence in $\mathcal{BL}(V, W)$ converges to an element of
$\mathcal{BL}(V, W)$ when $W$ is complete, as desired.  In particular,
we can apply this to $W = {\bf R}$ or ${\bf C}$, as appropriate, to
get that the dual space $V^*$ of bounded linear functionals on $V$ is
complete with respect to the dual norm, since the real and complex
numbers are complete with respect to their standard metrics.

        Suppose now that $V_1$, $V_2$, and $V_3$ are vector spaces,
all real or all complex, and equipped with norms.  Let $T_1$ be a
bounded linear mapping from $V_1$ into $V_2$, and let $T_2$ be a
bounded linear mapping from $V_2$ into $V_3$.  As usual, the
composition $T_2 \circ T_1$ is the linear mapping from $V_1$ into
$V_3$ that sends $v \in V_1$ to $T_2(T_1(v))$.  It is easy to see that
$T_2 \circ T_1$ is also bounded, and that
\begin{equation}
        \|T_2 \circ T_1\|_{op, 13} \le \|T_1\|_{op, 12} \, \|T_2\|_{op, 23}.
\end{equation}
Here the subscripts in the operator norms are included to indicate the
vector spaces and norms being used.

\section{Topological vector spaces, continued}
\label{topological vector spaces, continued}
\setcounter{equation}{0}

        Let $V$ be a topological vector space over the real or complex
numbers.  If $v \in V$ and $A \subseteq V$, then put
\begin{equation}
        v + A = \{v + a : a \in A\}.
\end{equation}
If $A$ is an open or closed set in $V$, then $v + A$ has the same
property, because translations determine homeomorphisms on $V$.
Similarly, if $A, B \subseteq V$, then put
\begin{equation}
        A + B = \{a + b : a \in A, b \in B\}.
\end{equation}
Equivalently,
\begin{equation}
\label{A + B = bigcup_{a in A} (a + B) = bigcup_{b in B} (b + A)}
        A + B = \bigcup_{a \in A} (a + B) = \bigcup_{b \in B} (b + A),
\end{equation}
which shows that $A + B$ is an open set in $V$ as soon as either $A$
or $B$ is open, since it is the union of a collection of open sets.

        If $A \subseteq V$ and $t$ is a real or complex number, as
appropriate, then we put
\begin{equation}
        t \, A = \{t \, a : a \in A\}.
\end{equation}
If $t \ne 0$ and $A$ is an open or closed set in $V$, then $t \, A$
has the same property, because multiplication by $t$ defines a
homeomorphism on $V$.  Of course, $t \, A = \{0\}$ when $t = 0$ and $A
\ne \emptyset$.  If $t = -1$, then $t \, A$ may be expressed simply as
$- A$.

        Suppose that $U$ is an open set in $V$ that contains $0$.
Continuity of addition at $0$ in $V$ implies that there are open
sets $U_1, U_2 \subseteq V$ that contain $0$ and satisfy
\begin{equation}
        U_1 + U_2 \subseteq U.
\end{equation}
If $v \in V$ and $v \ne 0$, then one can apply this to $U = V
\backslash \{v\}$, to get that
\begin{equation}
        U_1 \cap (v - U_2) = \emptyset.
\end{equation}
Here $v - U_2 = v + (- U_2)$, which is an open set in $V$ that
contains $v$.  This shows that $V$ is Hausdorff, using also
translation-invariance of the topology on $V$.

        Let $U$ be an open set in $V$ that contains $0$ again.
Continuity of scalar multiplication at $0$ implies that there is an
open set $U_0 \subseteq V$ that contains $0$ and a positive real
number $\delta > 0$ such that
\begin{equation}
        t \, U_0 \subseteq U
\end{equation}
for every $t \in {\bf R}$ or ${\bf C}$, as appropriate, with $|t| < \delta$.
Consider
\begin{equation}
\label{widetilde{U}_0 = bigcup_{|t| < delta} t U_0}
        \widetilde{U}_0 = \bigcup_{|t| < \delta} t \, U_0,
\end{equation}
where more precisely the union is taken over all real or complex
numbers $t$ such that $|t| < \delta$, as appropriate.  Equivalently,
\begin{equation}
        \widetilde{U}_0 = \bigcup_{0 < |t| < \delta} t \, U_0,
\end{equation}
since $0 \in U_0$, and hence $0 \in \widetilde{U}_0$.  This shows that
$\widetilde{U}_0$ is an open set in $V$, because it is a union of open
sets.

        A set $E \subseteq V$ is said to be \emph{balanced} if
\begin{equation}
        r \, E \subseteq E
\end{equation}
for every $r \in {\bf R}$ or ${\bf C}$, as appropriate, with $|r| \le 1$.
Thus every nonempty balanced set contains $0$ automatically.  It is
easy to see that the set $\widetilde{U}_0$ described in the previous
paragraph is balanced by construction.  This shows that for every open
set $U \subseteq V$ with $0 \in U$ there is a nonempty balanced open
set $\widetilde{U}_0 \subseteq U$.  To put it another way, the nonempty
balanced open sets in $V$ form a local base for the topology of $V$ at $0$.

        Let $U$ be an open set in $V$ that contains $0$ again, and let
$v \in V$ be given.  Because $0 \, v = 0 \in U$, continuity of scalar
multiplication at $v$ implies that there is a $\delta(v, U) > 0$ such
that
\begin{equation}
\label{t v in U}
        t \, v \in U
\end{equation}
for every $t \in {\bf R}$ or ${\bf C}$, as appropriate, with $|t| <
\delta(v, U)$.

        Let $A$ be any subset of $V$, and let $U \subseteq V$ be an
open set that contains $0$.  If $v$ is an element of the closure
$\overline{A}$ of $A$ in $V$, then
\begin{equation}
\label{(v - U) cap A ne emptyset}
        (v - U) \cap A \ne \emptyset,
\end{equation}
since $v - U$ is an open set in $V$ that contains $v$.  Equivalently,
\begin{equation}
        v \in A + U.
\end{equation}
It follows that
\begin{equation}
\label{overline{A} subseteq A + U}
        \overline{A} \subseteq A + U.
\end{equation}

\section{Bounded sets}
\label{bounded sets}
\setcounter{equation}{0}

        Let $V$ be a topological vector space over the real or complex
numbers.  A set $E \subseteq V$ is said to be \emph{bounded} if for
every open set $U \subseteq V$ with $0 \in U$ there is a real or
complex number $t_1$, as appropriate, such that
\begin{equation}
\label{E subseteq t_1 U}
        E \subseteq t_1 \, U.
\end{equation}
If $U$ is balanced, then it follows that
\begin{equation}
\label{E subseteq t U}
        E \subseteq t \, U
\end{equation}
for every $t \in {\bf R}$ or ${\bf C}$, as appropriate, such that $|t|
\ge |t_1|$.

        If $U$ is any open set in $V$ that contains $0$, then there is
a nonempty balanced open set $U' \subseteq U$, as in the previous
section.  In order to check that a set $E \subseteq V$ is bounded, it
is therefore enough to consider nonempty balanced open sets in $V$,
instead of arbitrary neighborhoods of $0$.  If $U$ is an arbitrary
neighborhood of $0$ in $V$, then we also get that (\ref{E subseteq t
U}) holds for all real or complex numbers $t$, as appropriate, for
which $|t|$ is sufficiently large.  This follows by applying the
stronger form of boundedness to a nonempty balanced open subset of
$U$.

        If $E \subseteq V$ has only finitely many elements, then it is
easy to see that $E$ is bounded, using the property (\ref{t v in U})
of neighborhoods of $0$ in $V$.  It is also easy to see that the union
of finitely many bounded subsets of $V$ is bounded, using the stronger
form of boundedness described in the previous paragraphs.  If the
topology on $V$ is determined by a norm $\|v\|$, then a set $E
\subseteq V$ is bounded if and only if $\|v\|$ is bounded on $E$.
Similarly, if the topology on $V$ is determined by a nice collection
of seminorms $\mathcal{N}$, then one can check that a set $E \subseteq
V$ is bounded if and only if each seminorm $N \in \mathcal{N}$ is
bounded on $V$.  Of course, subsets of bounded sets are also bounded.

        If $U$ is a neighborhood of $0$ in $V$, then
\begin{equation}
        \bigcup_{n = 1}^\infty n \, U = V,
\end{equation}
because of (\ref{t v in U}).  If $K \subseteq V$ is compact, then it
follows that
\begin{equation}
\label{K subseteq n_1 U cup cdots cup n_l U}
        K \subseteq n_1 \, U \cup \cdots \cup n_l \, U
\end{equation}
for some finite collection $n_1, \ldots, n_l$ of positive integers.
If $U$ is also balanced, then (\ref{K subseteq n_1 U cup cdots cup n_l
U}) implies that
\begin{equation}
\label{K subseteq n U}
        K \subseteq n \, U,
\end{equation}
where $n$ is the maximum of $n_1, \ldots, n_l$.  This shows that
compact subsets of $V$ are bounded, since it suffices to check
boundedness with respect to nonempty balanced open subsets of $V$, as
before.

        Suppose that $E_1, E_2 \subseteq V$ are bounded, and let us
check that $E_1 + E_2$ is bounded as well.  Let $U$ be a neighborhood
of $0$ in $V$, and let $U_1$, $U_2$ be neighborhoods of $0$ such that
$U_1 + U_2 \subseteq U$, as in the previous section.  Thus
\begin{equation}
\label{E_j subseteq t U_j}
        E_j \subseteq t \, U_j
\end{equation}
when $t$ is a real or complex number, as appropriate, for which $|t|$
is sufficiently large, and $j = 1, 2$.  This implies that
\begin{equation}
        E_1 + E_2 \subseteq t \, U_1 + t \, U_2 \subseteq t \, U
\end{equation}
when $t \in {\bf R}$ or ${\bf C}$ is sufficiently large, as desired.
In particular, it follows that translations of bounded sets are
bounded, since sets with only one element are bounded.

        Let us show now that the closure $\overline{E}$ of a bounded
set $E \subseteq V$ is bounded.  Let $U$ be a neighborhood of $0$ in
$V$ again, and let $U_1$, $U_2$ be neighborhoods of $0$ such that $U_1
+ U_2 \subseteq V$.  Because $E$ is bounded, there is a nonzero real
or complex number $t$, as appropriate, such that
\begin{equation}
        E \subseteq t \, U_1.
\end{equation}
We also have that
\begin{equation}
\label{overline{E} subseteq E + t U_2}
        \overline{E} \subseteq E + t \, U_2,
\end{equation}
as in (\ref{overline{A} subseteq A + U}), since $t \, U_2$ is a
neighborhood of $0$ in $V$.  Hence
\begin{equation}
\label{overline{E} subseteq E + t U_2 subseteq t U_1 + t U_2 subseteq t U}
        \overline{E} \subseteq E + t \, U_2 \subseteq t \, U_1 + t \, U_2
                                             \subseteq t \, U,
\end{equation}
as desired.

        If $E \subseteq V$ is bounded and $r$ is a real or complex
number, then it is easy to see that $r \, E$ is bounded too.  More
generally, suppose that $W$ is another topological vector space over
the real or complex numbers, depending on whether $V$ is real or
complex.  If $T$ is a continuous linear mapping from $V$ into $W$ and
$E \subseteq V$ is bounded, then it is easy to see that $T(E)$ is
bounded in $W$ as well.

\section{Uniform boundedness}
\label{uniform boundedness}
\setcounter{equation}{0}

        Let $M$ be a complete metric space, and let $\mathcal{E}$ be a
nonempty collection of continuous nonnegative real-valued functions on
$M$.  Suppose that $\mathcal{E}$ is bounded pointwise on $M$, in the
sense that
\begin{equation}
\label{mathcal{E}(x) = {f(x) : x in M}}
        \mathcal{E}(x) = \{f(x) : x \in M\}
\end{equation}
is a bounded set of real numbers for each $x \in M$.  Put
\begin{equation}
 \mathcal{E}_n = \{x \in M : f(x) \le n \hbox{ for every } f \in \mathcal{E}\}
\end{equation}
for each positive integer $n$.  Thus $\mathcal{E}_n$ is a closed set
in $M$ for every $n$, because the elements of $\mathcal{E}$ are
supposed to be continuous functions on $M$, and
\begin{equation}
        \bigcup_{n = 1}^\infty \mathcal{E}_n = M,
\end{equation}
by the hypothesis that $\mathcal{E}$ be bounded pointwise on $M$.  The
Baire category theorem implies that $\mathcal{E}_n$ has nonempty
interior for some $n$, so that $\mathcal{E}$ is uniformly bounded on a
nonempty open set in $M$.

        Now let $V$ be a real or complex vector space with a norm
$\|v\|$, and let $\Lambda$ be a nonempty collection of continuous
linear functionals on $V$.  Suppose that $\Lambda$ is bounded
pointwise on $V$, so that
\begin{equation}
        \Lambda(v) = \{\lambda(v) : \lambda \in \Lambda\}
\end{equation}
is a bounded set of real or complex numbers, as appropriate, for every
$v \in V$.  If $V$ is also complete, then it follows from the argument
in the previous paragraph that $\Lambda$ is uniformly bounded on a
nonempty open set in $V$.  Using the linearity of the elements of $V$,
one can show that $\Lambda$ is actually bounded on the unit ball in
$V$, which means that the dual norms of the elements of $\Lambda$ are
uniformly bounded.  This is a version of the Banach--Steinhaus
theorem.

        Let $V^*$ be the dual space of continuous linear functionals
on $V$, as usual.  Thus $V^*$ is equipped with the dual norm
$\|\lambda\|_*$, as in Section \ref{dual norms}, and also the weak$^*$
topology, as in Section \ref{weak^* topology}.  It is easy to see that
every bounded set in $V^*$ with respect to the dual norm is also
bounded with respect to the weak$^*$ topology, in the sense described
in the previous section.  Conversely, if $V$ is complete, then every
bounded set in $V^*$ with respect to the weak$^*$ topology is also
bounded with respect to the dual norm, by the principle of uniform
boundedness described in the previous paragraph.

        Similarly, we can consider $V$ equipped with the weak topology
associated to the collection of all continuous linear functionals on
$V$ with respect to the norm, as in Section \ref{weak topologies}.  If
$E \subseteq V$ is bounded with respect to the norm, then it is easy
to see that $E$ is also bounded with respect to the weak topology on
$V$.  Conversely, suppose that $E$ is bounded with respect to the weak
topology on $V$.  This means that
\begin{equation}
\label{E(lambda) = {lambda(v) : v in E}}
        E(\lambda) = \{\lambda(v) : v \in E\}
\end{equation}
is a bounded set of real or complex numbers, as appropriate, for each
$\lambda \in V^*$.  As in Section \ref{weak^* topology},
\begin{equation}
\label{L_v(lambda) = lambda(v)}
        L_v(\lambda) = \lambda(v)
\end{equation}
defines a continuous linear functional on $V^*$ with respect to the
dual norm $\|\lambda\|_*$ for every $v \in V$.  Let $V^{**} = (V^*)^*$
be the space of continuous linear functionals on $V^*$ with respect to
the dual norm on $V^*$.  Thus $V^{**}$ is also equipped with a
weak$^*$ topology, as the dual of $V^*$.  Consider
\begin{equation}
        \mathcal{L} = \{L_v : v \in E\},
\end{equation}
as a subset of $V^{**}$.  It is easy to see that $\mathcal{L}$ is
bounded with respect to the weak$^*$ topology on $V^{**}$, because $E$
is bounded with respect to the weak topology on $V$.  We also know
that $V^*$ is complete with respect to the dual norm, as in Section
\ref{bounded linear mappings}.  It follows from the discussion in the
previous paragraph that $\mathcal{L}$ is bounded with respect to the
dual norm on $V^{**}$ associated to the dual norm on $V^*$.  As in
Section \ref{weak^* topology}, the dual norm of $L_v$ as a continuous
linear functional on $V^*$ is equal to the norm of $v$ as an element
of $V$ for every $v \in V$, by the Hahn--Banach theorem.  This implies
that $E$ is bounded with respect to the norm on $V$.

\section{Bounded linear mappings, continued}
\label{bounded linear mappings, continued}
\setcounter{equation}{0}

        Let $V$, $W$ be topological vector spaces, both real or both
complex.  A linear mapping $T : V \to W$ is said to be \emph{bounded}
if for every bounded set $E \subseteq V$, $T(E)$ is a bounded set in
$W$.  It is easy to see that continuous linear mappings are bounded in
this sense, as mentioned at the end of Section \ref{bounded sets}.
Conversely, if the topology on $V$ is determined by a norm and $T : V
\to W$ is bounded, then $T$ is continuous.  More precisely, if there
is an open set $U \subseteq V$ that contains $0$ such that $T(U)$ is
bounded in $W$, then it is not difficult to check that $T$ is
continuous.  In particular, this condition holds when $T : V \to W$ is
bounded and there is a bounded neighborhood $U$ of $0$ in $V$.  If the
topology on $V$ is determined by a norm, then one can simply take $U$
to be the open unit ball in $V$.

        Let $V$ be a vector space over the real or complex numbers
equipped with a norm.  As in the previous section, the uniform
boundedness principle implies that every bounded set in $V$ with
respect to the weak topology is also bounded with respect to the norm.
Equivalently, the identity mapping on $V$ is bounded as a mapping from
$V$ with the weak topology into $V$ with the norm topology.  However,
the identity mapping on $V$ is not continuous as a mapping from $V$
with the weak topology into $V$ with the norm topology, unless $V$ is
finite-dimensional.  This is because the open unit ball in $V$ with
respect to the norm is not an open set with respect to the weak
topology when $V$ is infinite-dimensional, since an open set in $V$
with respect to the weak topology that contains $0$ also contains a
linear subspace of $V$ of finite codimension.

        Similarly, if $V$ is complete, then every bounded set in $V^*$
with respect to the weak$^*$ topology is bounded with respect to the
dual norm, as in the previous section.  This implies that the identity
mapping on $V^*$ is bounded as a mapping from $V^*$ with the weak$^*$
topology into $V^*$ with the topology determined by the dual norm.  As
in the preceding paragraph, this mapping is not continuous when $V$ is
infinite-dimensional.  Note that $V^*$ is infinite-dimensional when
$V$ is, by the Hahn--Banach theorem.

        Let $V$, $W$ be topological vector spaces again, both real or
both complex.  If $T : V \to W$ is a bounded linear mapping and $a$ is
a real or complex number, as appropriate, then $a \, T$ is also a
bounded linear mapping from $V$ into $W$.  This follows from the fact
that a scalar multiple of a bounded set in a topological vector space
is bounded as well.  Similarly, if $R : V \to W$ is another bounded
linear mapping, then the sum $R + T$ is bounded too.  This uses the
boundedness of the sum of two bounded subsets of a topological vector
space.  Now let $V_1$, $V_2$, and $V_3$ be topological vector spaces,
all real or all complex.  If $T_1 : V_1 \to V_2$ and $T_2 : V_2 \to
V_3$ are bounded linear mappings, then it is easy to see that their
composition $T_2 \circ T_1$ is a bounded linear mapping from $V_1$
into $V_3$, directly from the definition of a bounded linear mapping.

\section{Bounded sequences}
\label{bounded sequences}
\setcounter{equation}{0}

        Let $V$ be a topological vector space over the real or complex
numbers.  A sequence $\{v_j\}_{j = 1}^\infty$ of elements of $V$ is
said to be bounded if the set of $v_j$'s is bounded in $V$.  If
$\{v_j\}_{j = 1}^\infty$ converges to an element $v$ of $V$, then it
is easy to see that the set $K$ consisting of the $v_j$'s and $v$ is
compact, which works as well in any topological space.  This implies
that convergent sequences are bounded, since compact sets are bounded.
One can also show this more directly from the definitions, which is
especially simple when $\{v_j\}_{j = 1}^\infty$ converges to $0$.
Similarly, one can check that Cauchy sequences are bounded in $V$.  If
$\{v_j\}_{j = 1}^\infty$ is a bounded sequence in $V$, and $\{t_j\}_{j
= 1}^\infty$ is a sequence of real or complex numbers, as appropriate,
that converges to $0$, then it is easy to see that $\{t_j \, v_j\}_{j
= 1}^\infty$ converges to $0$ in $V$.

        Suppose now that there is a countable local base for the
topology of $V$ at $0$.  This means that there is a sequence $U_1,
U_2, \ldots$ of open subsets of $V$ that contain $0$ with the property
that if $U$ is any other open set in $V$ containing $0$, then $U_l
\subseteq U$ for some $l$.  As in Section \ref{topological vector
spaces, continued}, we can also take the $U_l$'s to be balanced
subsets of $V$.  We may as well ask that $U_{l + 1} \subseteq U_l$ for
each $l$ too, since otherwise we can replace $U_l$ with $U_1 \cap
\cdots \cap U_l$ for each $l$.  Let $\{v_j\}_{j = 1}^\infty$ be a
sequence of elements of $V$ that converges to $0$, and let us show
that there is a sequence of positive real numbers $\{r_j\}_{j =
1}^\infty$ such that $r_j \to \infty$ as $j \to \infty$ and $\{r_j \,
v_j\}_{j = 1}^\infty$ converges to $0$ in $V$.

        Because $\{v_j\}_{j = 1}^\infty$ converges to $0$ and
$l^{-1} \, U_l$ is an open set in $V$ that contains $0$ for each $l$,
there is a positive integer $N_l$ for each $l$ such that
\begin{equation}
        v_j \in l^{-1} \, U_l
\end{equation}
when $j \ge N_l$.  We may as well ask that $N_{l + 1} > N_l$ for every
$l$ too, by increasing the $N_l$'s if necessary.  Put
\begin{equation}
        r_j = l \quad\hbox{when } N_l \le j < N_{l + 1},
\end{equation}
and $r_j = 1$ when $1 \le j < N_1$ if $N_1 > 1$.  Thus $r_j \to
\infty$ as $j \to \infty$, and
\begin{equation}
        r_j \, v_j \in U_l \quad\hbox{when } N_l \le j < N_{l + 1}.
\end{equation}
This implies that $r_j \, v_j \in U_l$ when $j \ge N_l$, since $U_{l +
1} \subseteq U_l$ for each $l$.  It follows that $\{r_j \, v_j\}_{j =
1}^\infty$ converges to $0$ in $V$, as desired.  In particular, $\{r_j
\, v_j\}_{j = 1}^\infty$ is a bounded sequence in $V$.

        Let $W$ be another topological vector space, which is real if
$V$ is real and complex if $V$ is complex.  If $T$ is a bounded linear
mapping from $V$ into $W$ and $V$ has a countable local base for its
topology at $0$, then a well known theorem states that $T$ is
continuous.  To see this, it suffices to show that if $\{v_j\}_{j =
1}^\infty$ is a sequence of elements of $V$ that converges to $0$,
then $\{T(v_j)\}_{j = 1}^\infty$ converges to $0$ in $W$.  Let
$\{r_j\}_{j = 1}^\infty$ be a sequence of positive real numbers such
that $r_j \to +\infty$ as $j \to \infty$ and $\{r_j \, v_j\}_{j =
1}^\infty$ converges to $0$ in $V$, as in the previous paragraphs.
Thus $\{r_j \, v_j\}_{j = 1}^\infty$ is bounded in $V$, which implies
that $\{T(r_j \, v_j)\}_{j = 1}^\infty$ is bounded in $W$, since $T :
V \to W$ is bounded by hypothesis.  It follows that $T(v_j) = r_j^{-1}
\, T(r_j \, v_j)$ converges to $0$ as $j \to \infty$ in $W$, because
$\{r_j^{-1}\}_{j = 1}^\infty$ converges to $0$ in the real line.  This
shows that $T$ is sequentially continuous at $0$, and hence that $T$
is continuous at $0$, since $V$ has a countable local base for its
topology at $0$.  Of course, a linear mapping between topological
vector spaces is continuous at every point as soon as it is continuous
at $0$.

\section{Bounded linear functionals}
\label{bounded linear functionals}
\setcounter{equation}{0}

        If $V$ is a topological vector space over the real or complex numbers,
then we can restrict our attention in the previous section to the case where
$W$ is the one-dimensional vector space of real or complex numbers, as
appropriate.  Thus a bounded linear functional on $V$ is a linear functional
on $V$ that is bounded as a linear mapping into ${\bf R}$ or ${\bf C}$.

        Suppose now that $V$ is equipped with a norm $\|v\|$, so that
a linear functional on $V$ is bounded if and only if it is continuous,
as in the previous section.  Let $V^*$ be the dual space of bounded
linear functionals on $V$, which is equipped with the dual norm
$\|\lambda\|_*$, as in Section \ref{dual norms}.  Let $V^{**}$ be the
space of bounded linear functionals on $V^*$, which is equipped with a
dual norm $\|L\|_{**}$ associated to the dual norm $\|\lambda\|_*$ on
$V^*$.  As in Section \ref{weak^* topology}, each $v \in V$ determines
a bounded linear functional $L_v$ on $V^*$, defined by
\begin{equation}
        L_v(\lambda) = \lambda(v),
\end{equation}
and we also have that
\begin{equation}
        \|L_v\|_{**} = \|v\|.
\end{equation}
This defines an isometric linear embedding $v \mapsto L_v$ of $V$ into
$V^{**}$.

        A Banach space $V$ is said to be \emph{reflexive} if every
bounded linear functional on $V^*$ is of the form $L_v$ for some $v
\in V$.  It is easy to see that finite-dimensional Banach spaces are
automatically reflexive.  If $E$ is a nonempty set, then we have seen
in Section \ref{c_0(E)} that the dual of $c_0(E)$ may be identified
with $\ell^1(E)$, and we have seen in Section \ref{dual of ell^1} that
the dual of $\ell^1(E)$ may be identified with $\ell^\infty(E)$.  In
this case, the natural embedding of $c_0(E)$ into $c_0(E)^{**}$
described in the previous paragraph corresponds exactly to the
standard inclusion of $c_0(E)$ in $\ell^\infty(E)$ as a linear
subspace.  If $E$ has infinitely many elements, then $c_0(E)$ is a
proper linear subspace of $\ell^\infty(E)$, and it follows that
$c_0(E)$ is not reflexive.

        If $V$ is a real or complex vector space equipped with a norm
$\|v\|$, then every subset of $V^*$ that is bounded with respect to
the dual norm is also bounded with respect to the weak$^*$ topology.
This implies that every bounded linear functional on $V^*$ with
respect to the weak$^*$ topology is also bounded with respect to the
dual norm.  Conversely, if $V$ is also complete with respect to the
norm, then every bounded subset of $V^*$ with respect to the weak$^*$
topology is also bounded with respect to the dual norm, as in Section
\ref{uniform boundedness}.  This implies that every bounded linear
functional on $V^*$ with respect to the dual norm is also bounded with
respect to the weak$^*$ topology.  However, a linear functional on
$V^*$ is continuous with respect to the weak$^*$ topology if and only
if it is of the form $L_v$ for some $v \in V$, as in Section \ref{weak
topologies}.

\section{Uniform boundedness, continued}
\label{uniform boundedness, continued}
\setcounter{equation}{0}

        Let $V$ be a topological vector space over the real or complex
numbers, and let $\Lambda$ be a nonempty collection of continuous
linear functionals on $V$.  Suppose that $\Lambda$ is bounded
pointwise on $V$, in the sense that
\begin{equation}
        \Lambda(v) = \{\lambda(v) : \lambda \in \Lambda\}
\end{equation}
is a bounded set of real or complex numbers, as appropriate, for each
$v \in V$.  This is equivalent to asking that $\Lambda$ be bounded
with respect to the weak$^*$ topology on the dual space $V^*$ of
continuous linear functionals on $V$.  If the topology on $V$ is
determined by a norm, and if $V$ is complete with respect to this
norm, then $\Lambda$ is bounded with respect to the dual norm on
$V^*$, as in Section \ref{uniform boundedness}.

        Suppose now that $V$ is metrizable and complete, even if the
topology on $V$ may not be determined by a norm.  If $\Lambda
\subseteq V^*$ is bounded with respect to the weak$^*$ topology on
$V$, and hence bounded pointwise on $V$, then it follows from the
Baire category theorem that there is a nonempty open set $U_1
\subseteq V$ on which $\Lambda$ is uniformly bounded, as before.  If
$u_1 \in U_1$, then $U = U_1 - u_1$ is an open set in $V$ that
contains $0$, and $\Lambda$ is also uniformly bounded on $U$, because
the elements of $\Lambda$ are linear.  This is another version of the
theorem of Banach and Steinhaus.

        Let us restrict our attention now to the case where the
topology on $V$ is determined by a nice collection of seminorms
$\mathcal{N}$.  More precisely, we ask that $\mathcal{N}$ have only
finitely or countably many elements, so that $V$ is metrizable, and we
still ask that $V$ be complete.  If $\Lambda \subseteq V^*$ is bounded
pointwise on $V$, then $\Lambda$ is uniformly bounded on a
neighborhood of $0$, as in the previous paragraph.  In this case, this
implies that there are finitely many seminorms $N_1, \ldots, N_l \in
\mathcal{N}$ and a nonnegative real number $C$ such that
\begin{equation}
        |\lambda(v)| \le C \, \max_{1 \le j \le l} N_j(v)
\end{equation}
for every $v \in V$.  This is analogous to the discussion in
Section \ref{continuous linear functionals}.

        Of course, if there are finitely many seminorms $N_1, \ldots,
N_l \in \mathcal{N}$ and a $C \ge 0$ such that the preceding condition
holds for every $\lambda \in \Lambda$ and $v \in V$, then $\Lambda$ is
bounded pointwise on $V$.  In this situation, the choice of $N_1,
\ldots, N_l$ is part of the uniform boundedness condition.

\section{Another example}
\label{another example}
\setcounter{equation}{0}

        Let $V$ be the vector space of real or complex-valued
functions on the set ${\bf Z}_+$ of positive integers.  If $f \in V$
and $\rho$ is a positive real-valued function on ${\bf Z}_+$, then put
\begin{equation}
\label{B_rho(f) = {g in V : |f(l) - g(l)| < rho(l) for every l in {bf Z}_+}}
        B_\rho(f) = \{g \in V : |f(l) - g(l)| < \rho(l)
                                 \hbox{ for every } l \in {\bf Z}_+\}.
\end{equation}
Let us say that a set $U \subseteq V$ is an open set if for every $f \in U$
there is a positive real-valued function $\rho$ on ${\bf Z}_+$ such that
\begin{equation}
        B_\rho(f) \subseteq U.
\end{equation}
It is easy to see that this defines a topology on $V$, and that
$B_\rho(f)$ is an open set in $V$ with respect to this topology for
every $f \in V$ and positive function $\rho$ on ${\bf Z}_+$.
Equivalently, $V$ is the same as the Cartesian product of a sequence
of copies of the real or complex numbers, and this topology on $V$
corresponds to the ``strong product topology'', generated by arbitrary
products of open subsets of ${\bf R}$ or ${\bf C}$.

        One can also check that
\begin{equation}
        (f, g) \mapsto f + g
\end{equation}
defines a continuous mapping from $V \times V$ into $V$, using the
product topology on $V \times V$ determined by the topology just
described on $V$.  Similarly,
\begin{equation}
\label{f mapsto t f}
        f \mapsto t \, f
\end{equation}
is continuous as a mapping from $V$ into itself for each $t \in {\bf
R}$ or ${\bf C}$, as appropriate.  However, if $f(l) \ne 0$ for
infinitely many $l \in {\bf Z}_+$, then
\begin{equation}
\label{t mapsto t f}
        t \mapsto t \, f
\end{equation}
is not continuous as a mapping from the real or complex numbers with
the standard topology into $V$, and so $V$ is not a topological vector
space.

        Let $V_0$ be the linear subspace of $V$ consisting of
functions $f$ such that $f(l) = 0$ for all but finitely many $l \in
{\bf Z}_+$.  It is not difficult to verify that $V_0$ is a topological
vector space with respect to the topology induced by the one just
defined on $V$.  In particular, if $f \in V_0$, then (\ref{t mapsto t
f}) is continuous as a mapping from the real or complex numbers into
$V$.  

        Let us check that $V_0$ is a closed set in $V$.  Let $f \in V
\backslash V_0$ be given, and let $\rho$ be defined on ${\bf Z}_+$ by
\begin{equation}
\label{rho(l) = |f(l)|}
        \rho(l) = |f(l)|
\end{equation}
when $f(l) \ne 0$, and $\rho(l) = 1$ otherwise.  Thus $\rho(l) > 0$ for
every $l \in {\bf Z}_+$.  If $g \in B_\rho(f)$, then
\begin{equation}
\label{|f(l) - g(l)| < rho(l) = |f(l)|}
        |f(l) - g(l)| < \rho(l) = |f(l)|
\end{equation}
when $f(l) \ne 0$, which implies that $g(l) \ne 0$ for infinitely many
$l \in {\bf Z}_+$.  This shows that
\begin{equation}
        B_\rho(f) \subseteq V \backslash V_0,
\end{equation}
and hence that $V \backslash V_0$ is an open set in $V$, as desired.

        Let $U_1, U_2, \ldots$, be a sequence of relatively open sets
in $V_0$ containing $0$.  By construction, there is a sequence
$\rho_1, \rho_2, \ldots$ of positive functions on ${\bf Z}_+$ such
that
\begin{equation}
        B_{\rho_j}(0) \cap V_0 \subseteq U_j
\end{equation}
for each $j$.  Put
\begin{equation}
\label{rho(j) = frac{rho_j(j)}{2}}
        \rho(j) = \frac{\rho_j(j)}{2}
\end{equation}
for each $j \in {\bf Z}_+$, so that $\rho(j)$ is another positive
function on ${\bf Z}_+$.  Thus $B_\rho(0) \cap V_0$ is another
relatively open set in $V_0$ that contains $0$, and
\begin{equation}
\label{B_{rho_j}(0) cap V_0 not subseteq B_rho(0) cap V_0}
        B_{\rho_j}(0) \cap V_0 \not\subseteq B_\rho(0) \cap V_0
\end{equation}
for each $j$, because $\rho(j) < \rho_j(j)$ for each $j$.  This
implies that
\begin{equation}
        U_j \not\subseteq B_\rho(0) \cap V_0
\end{equation}
for each $j$, and it follows that $V_0$ does not have a countable
local base for its topology at $0$.

        Let $E$ be a subset of $V_0$, and let $L(E)$ be the set of $l
\in {\bf Z}_+$ for which there is an $f \in E$ such that $f(l) \ne 0$.
Also let $\rho$ be a positive function on ${\bf Z}_+$ such that
\begin{equation}
        \rho(l) = \frac{|f_l(l)|}{l}
\end{equation}
for some $f_l \in E$ with $f_l(l) \ne 0$ when $l \in L(E)$.  Thus
\begin{equation}
        f_l \not\in t \, B_\rho(0)
\end{equation}
when $l \in L(E)$ and $t \in {\bf R}$ or ${\bf C}$ satisfies $|t| \le l$.
If $L(E)$ has infinitely many elements, then it follows that
\begin{equation}
        E \not\subseteq t \, B_\rho(0)
\end{equation}
for any real or complex number $t$, as appropriate.  This shows that
$E$ can have only finitely or countably many elements when $E$ is
bounded in $V$.

        Let $V_{0, n}$ be the $n$-dimensional linear subspace of $V_0$
consisting of functions $f$ on ${\bf Z}_+$ such that $f(l) = 0$ when
$l > n$, for each positive integer $n$.  Note that
\begin{equation}
        \bigcup_{n = 1}^\infty V_{0, n} = V_0.
\end{equation}
If $E \subseteq V_0$ is bounded, then $E \subseteq V_{0, n}$ for some
$n$, as in the previous paragraph.  In this case, $E$ is also bounded
as a subset of $V_{0, n}$ in the usual sense, which is to say that
\begin{equation}
        E_j = \{f(j) : f \in E\}
\end{equation}
is bounded in ${\bf R}$ or ${\bf C}$, as appropriate, for each $j \le
n$.  Conversely, if $E \subseteq V_{0, n}$ and $E_j$ is bounded for
each $j \le n$, then $E$ is bounded in $V_{0, n}$, and hence in $V_0$.

        Let $\tau$ be a positive real-valued function on ${\bf Z}_+$,
and consider the norm $N_\tau$ on $V_0$ defined by
\begin{equation}
        N_\tau(f) = \max_{j \ge 1} |f(j)| \, \tau(j).
\end{equation}
If $\mathcal{N}$ is the collection of all of these norms $N_\tau$ on
$V_0$, then it is not difficult to check that the topology on $V_0$
associated to $\mathcal{N}$ is the same as the topology on $V_0$
induced from the one on $V$ as before.  To see this, observe that the
open unit ball in $V_0$ with respect to $N_\tau$,
\begin{equation}
        \{f \in V_0 : N_\tau(f) < 1\},
\end{equation}
is the same as the set of $f \in V_0$ for which there is a positive
real number $r < 1$ such that
\begin{equation}
        |f(j)| < r \, \tau(j)^{-1}
\end{equation}
for each $j \in {\bf Z}_+$.  This is contained in $B_\rho(0)$ with
$\rho = 1/\tau$, and more precisely it is equal to
\begin{equation}
\label{bigcup_{0 < r < 1} B_{r rho}(0)}
        \bigcup_{0 < r < 1} B_{r \, \rho}(0),
\end{equation}
which is close enough to show that the topologies are the same.

        Similarly, if $\sigma$ is a positive real-valued function on
${\bf Z}_+$, then
\begin{equation}
\label{N'_sigma(f) = sum_{j = 1}^infty |f(j)| sigma(j)}
        N'_\sigma(f) = \sum_{j = 1}^\infty |f(j)| \, \sigma(j)
\end{equation}
defines a norm on $V_0$.  Clearly
\begin{equation}
        N_\sigma(f) \le N'_\sigma(f)
\end{equation}
for every $f \in V_0$.  In the other direction, if we put
\begin{equation}
        \tau(j) = j^2 \, \sigma(j)
\end{equation}
for each $j \in {\bf Z}_+$, then
\begin{equation}
 N'_\sigma(f) \le \Big(\sum_{j = 1}^\infty \frac{1}{j^2}\Big) \, N_\tau(f)
\end{equation}
for every $f \in V_0$.  If $\mathcal{N}'$ is the collection of all of
these norms $N'_\sigma$ on $V_0$, then it follows that $\mathcal{N}'$
determines the same topology on $V_0$ as $\mathcal{N}$ does.  Hence
the topology on $V_0$ associated to $\mathcal{N}'$ is also the same as
the one induced on $V_0$ by the topology on $V$ defined at the
beginning of the section.

        Let $N$ be any seminorm on $V_0$, and let $\delta_j(l)$ be
the function on ${\bf Z}_+$ equal to $1$ when $j = l$ and to $0$
otherwise.  If
\begin{equation}
        N(\delta_j) \le \sigma(j)
\end{equation}
for each $j \in {\bf Z}_+$, then we get that
\begin{equation}
        N(f) \le N'_\sigma(f)
\end{equation}
for every $f \in V_0$.  More precisely, if $f \in V_{0, n}$, then
$f = \sum_{j = 1}^n f(j) \, \delta_j$, and hence
\begin{equation}
        N(f) \le \sum_{j = 1}^n |f(j)| \, N(\delta_j) \le N'_\sigma(f).
\end{equation}
This implies that open balls with respect to $N$ are also open sets in
$V_0$.

        Let $h$ be a real or complex-valued function on ${\bf Z}_+$,
as appropriate, and consider
\begin{equation}
\label{lambda_h(f) = sum_{j = 1}^infty f(j) h(j)}
        \lambda_h(f) = \sum_{j = 1}^\infty f(j) \, h(j)
\end{equation}
for $f \in V_0$.  This defines a linear functional on $V_0$, and every
linear functional on $V_0$ is of this form.  If
\begin{equation}
        |h(j)| \le \sigma(j)
\end{equation}
for each $j \in {\bf Z}_+$, then it follows that
\begin{equation}
        |\lambda_h(f)| \le N'_\sigma(f)
\end{equation}
for every $f \in V_0$.  Thus $\lambda_h$ is continuous on $V_0$, and
hence every linear functional on $V_0$ is continuous.

\part{Algebras of functions}

\section{Homomorphisms}
\label{homomorphisms}
\setcounter{equation}{0}

        Let $X$ be a set, and let $\mathcal{F}$ be an ultrafilter on
$X$.  If $f$ is a real or complex-valued function on $X$, then
$f_*(\mathcal{F})$ is an ultrafilter on ${\bf R}$ or ${\bf C}$, as
appropriate, as in Section \ref{ultrafilters}.  If $f$ is bounded on
$X$, then one can check that $f_*(\mathcal{F})$ converges to an
element of ${\bf R}$ or ${\bf C}$.  In this case, it is a bit simpler
to think of $f$ as taking values in a compact subset $K$ of ${\bf R}$
or ${\bf C}$, so that $f$ maps $\mathcal{F}$ to an ultrafilter on $K$,
which therefore converges.  Although this is not quite the same as
$f_*(\mathcal{F})$ as an ultrafilter on ${\bf R}$ or ${\bf C}$, they
are almost the same, and converge to the same limit, as in Section
\ref{filters, subsets}.

        Let $L_\mathcal{F}(f)$ denote the limit of $f_*(\mathcal{F})$,
which may also be described as the limit of $f$ along $\mathcal{F}$.
If $p \in X$ and $\mathcal{F}$ is the ultrafilter $\mathcal{F}_p$
based at $p$ as in the previous section, then $L_\mathcal{F}(f) =
f(p)$ for every bounded function $f$ on $X$.  It is easy to see that
every ultrafilter on $X$ is of this type when $X$ has only finitely
many elements.  Otherwise, if $X$ is an infinite set, then the
collection of subsets $A$ of $X$ such that $X \backslash A$ has only
finitely many elements is a filter on $X$.  Any ultrafilter on $X$
which is a refinement of this filter is not the same as
$\mathcal{F}_p$ for any $p \in X$.

        Observe that
\begin{equation}
        L_\mathcal{F}(f) \in \overline{f(X)}
\end{equation}
for every $f \in \ell^\infty(X)$.  In particular,
\begin{equation}
        |L_\mathcal{F}(f)| \le \|f\|_\infty.
\end{equation}
If $f$ is a constant function on $X$, then $L_\mathcal{F}(f)$ is equal
to this constant value.  One can also check that
\begin{equation}
        L_\mathcal{F}(f + g) = L_\mathcal{F}(f) + L_\mathcal{F}(f)
\end{equation}
and
\begin{equation}
        L_\mathcal{F}(f \, g) = L_\mathcal{F}(f) \, L_\mathcal{F}(g)
\end{equation}
for every $f, g \in \ell^\infty(X)$.  This is analogous to standard
facts about the limits of a sum and product being equal to the
corresponding sum or product of limits.

        If $X$ is an infinite set and $\mathcal{F} \ne \mathcal{F}_p$
for any $p \in X$, then one can check that $L_\mathcal{F}(f) = 0$ for
every $f \in c_0(X)$.  Similarly, $L_\mathcal{F}(f)$ is the same as
the limit of $f(x)$ at infinity when $f \in c(X)$, as in Section
\ref{dual of ell^1}.  Remember that the limit of $f(x)$ at infinity
defines a continuous linear functional on $c(X)$, with dual norm equal
to $1$ with respect to the $\ell^\infty$ norm.  Thus
$L_\mathcal{F}(f)$ is an extension of this linear functional on $c(X)$
to a continuous linear functional on $\ell^\infty(X)$, also with dual
norm equal to $1$.  The existence of such an extension was mentioned
before, as a consequence of the Hahn--Banach theorem.

\section{Homomorphisms, continued}
\label{homomorphisms, continued}
\setcounter{equation}{0}

        Let $X$ be a nonempty set, and note that the product of two
bounded real or complex-valued functions on $X$ is bounded as well.
Suppose that $L$ is a linear functional on $\ell^\infty(X)$ which is a
homomorphism with respect to multiplication of functions, in the sense
that
\begin{equation}
        L(f \, g) = L(f) \, L(g)
\end{equation}
for every $f, g \in \ell^\infty(X)$.  If $L(f) = 0$ for every $f \in
\ell^\infty(X)$, then $L$ satisfies these conditions trivially, and so
we suppose that $L(f) \ne 0$ for at least one $f \in \ell^\infty(X)$.
This implies that
\begin{equation}
        L({\bf 1}_X) = 1,
\end{equation}
where ${\bf 1}_X$ is the constant function equal to $1$ on $X$, since
${\bf 1}_X \, f = f$ and hence
\begin{equation}
        L(f) = L({\bf 1}_X \, f) = L({\bf 1}_X) \, L(f).
\end{equation}
We would like to show that $L$ is associated to an ultrafilter on $X$,
as in the previous section.

        If $A \subseteq X$, then let ${\bf 1}_A(x)$ be the indicator
function on $X$ associated to $A$, equal to $1$ when $x \in A$ and to
$0$ when $x \in X \backslash A$.  Thus ${\bf 1}_A^2 = {\bf 1}_A$,
which implies that
\begin{equation}
        L({\bf 1}_A) = L({\bf 1}_A^2) = L({\bf 1}_A)^2,
\end{equation}
and hence $L({\bf 1}_A) = 0$ or $1$.  Because ${\bf 1}_A + {\bf 1}_{X
\backslash A} = {\bf 1}_X$,
\begin{equation}
        L({\bf 1}_A) + L({\bf 1}_{X \backslash A}) = L({\bf 1}_X) = 1,
\end{equation}
so that exactly one of $L({\bf 1}_A)$ and $L({\bf 1}_{X \backslash
A})$ is equal to $1$.  If $A, B \subseteq X$, then ${\bf 1}_A \, {\bf
1}_B = {\bf 1}_{A \cap B}$, and so
\begin{equation}
        L({\bf 1}_{A \cap B}) = L({\bf 1}_A) \, L({\bf 1}_B).
\end{equation}
This shows that $L({\bf 1}_{A \cap B}) = 1$ when $L({\bf 1}_A) =
L({\bf 1}_B) = 1$.  Similarly, if $A \subseteq B$ and $L({\bf 1}_A) =
1$, then $A \cap B = A$, and we get that $L({\bf 1}_B) = 1$.  Of
course, ${\bf 1}_A = 0$ when $A = \emptyset$, so that $L({\bf 1}_A) =
0$.  If
\begin{equation}
        \mathcal{F}_L = \{A \subseteq X : L({\bf 1}_A) = 1\},
\end{equation}
then it follows that $\mathcal{F}_L$ is a filter on $X$.  More
precisely, $\mathcal{F}_L$ is an ultrafilter on $X$, since $A$ or $X
\backslash A$ is in $\mathcal{F}_L$ for each $A \subseteq X$.

        It is easy to see that $L({\bf 1}_A)$ is the same as the limit
of ${\bf 1}_A$ along $\mathcal{F}_L$ as in the previous section.  This
implies that $L(f)$ is equal to the limit of $f$ along $\mathcal{F}_L$
when $f$ is a finite linear combination of indicator functions of
subsets of $X$, by linearity.  One can also check that finite linear
combinations of indicator functions of subsets of $X$ are dense in
$\ell^\infty(X)$.  We already know that the limit along an ultrafilter
defines a continuous linear functional on $\ell^\infty(X)$, as in the
previous section, and we would like to check that $L$ is also a
continuous linear functional on $\ell^\infty(X)$.  This would imply
that $L(f)$ is equal to the limit of $f$ along $\mathcal{F}$ for every
$f \in \ell^\infty(X)$, by continuity and density.

        Suppose that $f$ is a bounded function on $X$ such that $f(x)
\ne 0$ for every $x \in X$ and $1/f$ is also bounded.  Thus
\begin{equation}
        L(f) \, L(1/f) = L({\bf 1}_X) = 1,
\end{equation}
and hence $L(f) \ne 0$ in particular.  Equivalently, $0 \in
\overline{f(X)}$ when $L(f) = 0$.  This implies that
\begin{equation}
        L(f) \in \overline{f(X)}
\end{equation}
for every $f \in \ell^\infty(X)$, since one can reduce to the case
where $L(f) = 0$ by subtracting $L(f) \, {\bf 1}_X$ from $f$, using
the fact that $L({\bf 1}_X) = 1$.  In particular,
\begin{equation}
        |L(f)| \le \|f\|_\infty
\end{equation}
for every $f \in \ell^\infty$, which implies that $L$ is a continuous
linear functional on $\ell^\infty(X)$ with dual norm equal to $1$, as
desired.

\section{Bounded continuous functions}
\label{bounded continuous functions}
\setcounter{equation}{0}

        Let $X$ be a topological space, and let $C_b(X)$ be the space
of bounded real or complex-valued continuous functions on $X$.  As
usual, this may also be denoted $C_b(X, {\bf R})$ or $C_b(X, {\bf
C})$, to indicate whether real or complex-valued functions are being
used.  Of course, $C_b(X)$ is the same as $\ell^\infty(X)$ when $X$ is
equipped with the discrete topology.  If $X$ is compact, then
continuous functions are automatically bounded on $X$.  Constant
functions on $X$ are always continuous, and the existence of
nonconstant functions on $X$ depends on the behavior of $X$.

        Remember that sums and products of continuous functions are
continuous.  Similarly, sums and products of bounded functions are
bounded, so that sums and products of bounded continuous functions are
bounded and continuous.  It follows that $C_b(X)$ is a vector space
with respect to pointwise addition and scalar multiplication, and a
commutative algebra with respect to multiplication of functions.  The
supremum norm on $C_b(X)$ is defined by
\begin{equation}
\label{||f||_{sup} = sup_{x in X} |f(x)|}
        \|f\|_{sup} = \sup_{x \in X} |f(x)|,
\end{equation}
and it is easy to see that this is indeed a norm.  Moreover,
\begin{equation}
        \|f \, g\|_{sup} \le \|f\|_{sup} \, \|g\|_{sup}
\end{equation}
for every $f, g \in C_b(X)$.

        Suppose that $\phi$ is a linear functional on $C_b(X)$ which
is also a homomorphism with respect to multiplication of functions,
in the sense that
\begin{equation}
        \phi(f \, g) = \phi(f) \, \phi(g)
\end{equation}
for every $f, g \in C_b(X)$.  If $\phi(f) = 0$ for each $f \in
C_b(X)$, then $\phi$ satisfies these conditions trivially, and so we
also ask that $\phi(f) \ne 0$ for some $f \in C_b(X)$.  As before,
this implies that
\begin{equation}
        \phi({\bf 1}_X) = 1,
\end{equation}
where ${\bf 1}_X$ is the constant function equal to $1$ at every point
in $X$.  Of course,
\begin{equation}
        \phi_p(f) = f(p)
\end{equation}
has these properties for every $p \in X$.

        If $f$ is a bounded continuous function on $X$ such that $f(x)
\ne 0$ for every $x \in X$, then $1/f$ is also a continuous function
on $X$.  If $1/f$ is bounded as well, then
\begin{equation}
        \phi(f) \, \phi(1/f) = \phi({\bf 1}_X) = 1,
\end{equation}
which implies that $\phi(f) \ne 0$.  If $f$ is any bounded continuous
function on $X$ such that $\phi(f) = 0$, then it follows that $0 \in
\overline{f(X)}$, since otherwise $1/f \in C_b(X)$.  This implies that
\begin{equation}
\label{phi(f) in overline{f(X)}}
        \phi(f) \in \overline{f(X)}
\end{equation}
for every $f \in C_b(X)$, by applying the previous statement to $f -
\phi(f) \, {\bf 1}_X$.

        In particular,
\begin{equation}
\label{|phi(f)| le ||f||_{sup}}
        |\phi(f)| \le \|f\|_{sup}
\end{equation}
for every $f \in C_b(X)$, so that $\phi$ is a continuous linear
functional on $C_b(X)$.  The dual norm of $\phi$ with respect to the
sumpremum norm is equal to $1$, since $\phi({\bf 1}_X) = 1$.  In the
complex case, (\ref{phi(f) in overline{f(X)}}) implies that $\phi(f)
\in {\bf R}$ when $f$ is real-valued.  In both the real and complex
cases, we get that
\begin{equation}
\label{phi(f) ge 0}
        \phi(f) \ge 0
\end{equation}
for every bounded nonnegative real-valued function $f$ on $X$.  If $A
\subseteq X$ is both open and closed, then the corresponding indicator
function ${\bf 1}_A$ is continuous on $X$, and $\phi({\bf 1}_A)$ is
either $0$ or $1$.

        Let $B^*$ be the closed unit ball in the dual of $C_b(X)$,
with respect to the dual norm associated to the supremum norm on
$C_b(X)$.  Thus multiplicative homomorphisms on $C_b(X)$ are elements
of $B^*$, because of (\ref{|phi(f)| le ||f||_{sup}}).  It is easy to
see that the set of multiplicative homomorphisms on $C_b(X)$ is closed
with respect to the weak$^*$ topology, since $\phi(f)$, $\phi(g)$, and
$\phi(f \, g)$ are continuous functions of $\phi \in C_b(X)^*$ with
respect to the weak$^*$ topology for every $f, g \in C_b(X)$.  The set
of nonzero multiplicative homomorphisms on $C_b(X)$ is also closed in
the weak$^*$ topology, since it can be described by the additional
condition $\phi({\bf 1}_X) = 1$, and $\phi({\bf 1}_X)$ is a continuous
function of $\phi$ with respect to the weak$^*$ topology.  Hence the
set of nonzero multiplicative homomorphisms on $C_b(X)$ is compact
with respect to the weak$^*$ topology, because it is a closed subset
of $B^*$, which is compact by the Banach--Alaoglu theorem.

        If $p \in X$, then $\phi_p(f) = f(p)$ is a nonzero
multiplicative homomorphism on $C_b(X)$, as before.  Thus $p \mapsto
\phi_p$ defines a mapping from $X$ into $B^*$.  It is easy to see that
this mapping is continuous with respect to the weak$^*$ topology on $B^*$,
since $\phi_p(f) = f(p)$ is continuous on $X$ for every $f \in C_b(X)$.

        If $X$ is equipped with the discrete topology, then $p \mapsto
\phi_p$ is a one-to-one mapping of $X$ into $B^*$, and the topology
induced on the set
\begin{equation}
        \{\phi_p : p \in X\}
\end{equation}
by the weak$^*$ topology is the same as the discrete topology.  If $X$
is infinite, then of course this set is not compact.  Let $\phi \in
B^*$ be a limit point of this set with respect to the weak$^*$
topology, which is therefore not in the set.  If $f \in c(X)$, then
one can check that $\phi(f)$ is equal to the limit of $f(x)$ at
infinity on $X$.  This is another way to get homomorphisms on
$\ell^\infty(X)$ extending the limit at infinity on $c(X)$.

\section{Compact spaces}
\label{compact spaces}
\setcounter{equation}{0}

        Let $X$ be a compact topological space, and let $C(X)$ be the
space of continuous real or complex-valued functions on $X$.  This may
also be denoted $C(X, {\bf R})$ or $C(X, {\bf C})$, to indicate
whether real or complex-valued functions are being used.  As before,
continuous functions on compact spaces are automatically bounded, so
that $C(X) = C_b(X)$.  Let $\phi$ be a nonzero multiplicative
homomorphism on $C(X)$, as in the previous section.  We would like to
show that there is a $p \in X$ such that $\phi(f) = f(p)$ for every $f
\in C(X)$.

        Suppose for the sake of a contradiction that for each $p \in
X$ there is a continuous function $f_p$ on $X$ such that $\phi(f_p)
\ne f_p(p)$.  We may as well ask also that $\phi(f_p) = 0$, since
otherwise we can replace $f_p$ with $f_p - \phi(f_p) \, {\bf 1}_X$,
using the fact that $\phi({\bf 1}_X) = 1$.  Thus $f_p(p) \ne 0$.

        Similarly, we may suppose that $f_p$ is a nonnegative
real-valued function on $X$ for each $p \in X$, by replacing $f_p$
with $|f_p|^2$ if necessary.  More precisely, in the real case,
\begin{equation}
        \phi(|f_p|^2) = \phi(f_p^2) = \phi(f_p)^2 = 0,
\end{equation}
while in the complex case,
\begin{equation}
        \phi(|f_p|^2) = \phi(f_p \overline{f_p})
                       = \phi(f_p) \, \phi(\overline{f_p}) = 0.
\end{equation}
Of course, we also get that $f_p(p) > 0$ after this substitution.

        Consider
\begin{equation}
        U(p) = \{x \in X : f_p(x) > 0\}.
\end{equation}
This is an open set in $X$ for each $p \in X$, because $f_p$ is
continuous, and $p \in U(p)$ by construction.  Thus $U(p)$, $p \in X$,
is an open covering of $X$, and so there are finitely many elements
$p_1, \ldots, p_n$ of $X$ such that
\begin{equation}
        X = \bigcup_{j = 1}^n U(p_j),
\end{equation}
by compactness.  If $f = \sum_{j = 1}^n f_{p_j}$, then $f$ is
continuous on $X$, $\phi(f) = 0$, and $f(x) > 0$ for every $x \in X$.
This is a contradiction, because $1/f$ is also a continuous function
on $X$, which implies that $\phi(f) \ne 0$, as in the previous
section.

\section{Closed ideals}
\label{closed ideals}
\setcounter{equation}{0}

        Let $X$ be a topological space, and let $C(X)$ be the space of
continuous real or complex-valued functions on $X$.  As usual, this is
a vector space with respect to pointwise addition and scalar
multiplication, and a commutative algebra with respect to pointwise
multiplication of functions.  A linear subspace $\mathcal{I}$ of
$C(X)$ is said to be an \emph{ideal} if for every $a \in C(X)$ and $f
\in \mathcal{I}$ we have that $a \, f \in C(X)$.  In this section, we
shall restrict our attention to compact Hausdorff spaces $X$, so that
continuous functions on $X$ are automatically bounded.  We shall also
be especially interested in ideals that are closed subsets of $C(X)$
with respect to the supremum norm.

        If $E \subseteq X$, then let $\mathcal{I}_E$ be the collection
of $f \in C(X)$ such that $f(x) = 0$ for every $x \in X$.  It is easy
to see that this is a closed ideal in $C(X)$, directly from the
definitions.  If $\overline{E}$ is the closure of $E$ in $X$, then
\begin{equation}
\label{mathcal{I}_{overline{E}} = mathcal{I}_E}
        \mathcal{I}_{\overline{E}} = \mathcal{I}_E,
\end{equation}
because any continuous function that vanishes on $E$ automatically
vanishes on the closure of $E$ as well.  Thus we may as well restrict
our attention to closed subsets $E$ of $X$.  We would like to show
that any closed ideal $\mathcal{I}$ in $C(X)$ is of the form
$\mathcal{I}_E$ for some closed set $E \subseteq X$.

        If $\mathcal{I}$ is any subset of $C(X)$, then
\begin{equation}
        E = \{x \in X : f(x) = 0\}
\end{equation}
is a closed set in $X$.  This is because the set where a continuous
function is equal to $0$ is a closed set, and $E$ is the same as the
intersection of the zero sets associated to the elements of
$\mathcal{I}$.  By construction,
\begin{equation}
        \mathcal{I} \subseteq \mathcal{I}_E.
\end{equation}
We would like to show that equality holds when $\mathcal{I}$ is a
closed ideal in $C(X)$.

        If $\phi$ is a real or complex-valued function on $X$, then
the \emph{support} of $\phi$ is denoted $\supp \phi$ and is defined to
be the closure of the set of $x \in X$ such that $\phi(x) \ne 0$.
Suppose that $\phi$ is a continuous function on $X$ whose support is
contained in the complement of $E$ in $X$.  If $p \in \supp \phi$, so
that $p \not\in E$, then there is an $f_p \in \mathcal{I}$ such that
$f_p(p) \ne 0$.  Note that
\begin{equation}
        U(p) = \{x \in X : f_p(x) \ne 0\}
\end{equation}
is an open set in $X$, because $f_p$ is continuous.  Thus $U(p)$, $p
\in \supp \phi$, is an open covering of $\supp \phi$ in $X$, since $p
\in U(p)$ for each $p$.  We also know that $\supp f$ is compact,
because it is a closed set in a compact space.  It follows that there
are finitely many elements $p_1, \ldots, p_n$ of $\supp \phi$ such that
\begin{equation}
        \supp \phi \subseteq \bigcup_{j = 1}^n U(p_j).
\end{equation}

        Observe that
\begin{equation}
        \sum_{l = 1}^n |f_{p_l}(x)|^2 > 0
\end{equation}
for every $x \in \supp \phi$.  Put
\begin{equation}
        \psi(x) = \phi(x) \, \Big(\sum_{l = 1}^n |f_{p_l}(x)|^2\Big)^{-1},
\end{equation}
which is interpreted as being $0$ when $x \not\in \supp \phi$.  This
is a continuous function on $X$, because it is equal to $0$ on a
neighborhood of every $x \in X \backslash \supp \phi$, and because it
is a quotient of continuous functions with nonzero denominator on a
neighborhood of every $x \in \supp \phi$.  In the real case, we have that
\begin{equation}
        \phi(x) = \sum_{j = 1}^n (\psi(x) \, f_j(x)) \, f_j(x)
\end{equation}
for every $x \in X$, and in the complex case we have that
\begin{equation}
        \phi(x) = \sum_{j = 1}^n (\psi(x) \, \overline{f_j(x)}) \, f_j(x),
\end{equation}
where $\overline{f_j(x)}$ is the complex conjugate of $f_j(x)$.  This
implies that $\phi \in \mathcal{I}$, since $f_j \in \mathcal{I}$ for
each $j$, $\psi \, f_j \in C(X)$ in the real case and $\psi \,
\overline{f_j} \in C(X)$ in the complex case, and $\mathcal{I}$ is an
ideal.

        Now let $f$ be a continuous function on $X$ such that $f(x) =
0$ for every $x \in E$, and let $\epsilon > 0$ be given.  Thus
\begin{equation}
        K(\epsilon) = \{x \in X : |f(x)| \ge \epsilon\}
\end{equation}
is a closed set in $X$ contained in the complement of $E$.  Let
$V(\epsilon)$ be an open set in $X$ such that $K(\epsilon) \subseteq
V(\epsilon)$ and $\overline{V(\epsilon)} \subseteq X \backslash E$.
By Urysohn's lemma, there is a continuous real-valued function
$\theta_\epsilon$ on $X$ such that $\theta(x) = 1$ for every $x \in
K(\epsilon)$, $\theta_\epsilon(x) = 0$ when $x \not\in V(\epsilon)$,
and $0 \le \theta_\epsilon(x) \le 1$ for every $x \in X$.  In
particular, the support of $\theta_\epsilon$ is contained in
$\overline{V(\epsilon)}$, which is contained in the complement of $E$.
Of course, the support of $\theta_\epsilon \, f$ is contained in the
support of $\theta_\epsilon$.  Hence
\begin{equation}
        \theta_\epsilon \, f \in \mathcal{I}
\end{equation}
for each $\epsilon > 0$, by the discussion in the previous paragraph.
Moreover,
\begin{equation}
 |\theta_\epsilon(x) \, f(x) - f(x)| = (1 - \theta_\epsilon(x)) \, |f(x)|
                                      < \epsilon
\end{equation}
for every $x \in X$, because $1 - \theta_\epsilon(x) = 0$ when $x \in
K(\epsilon)$, $|f(x)| < \epsilon$ when $x \in X \backslash
K(\epsilon)$, and $0 \le \theta_\epsilon(x) \le 1$ for every $x \in
X$.  This implies that $\theta_\epsilon \, f \to f$ uniformly on $X$
as $\epsilon \to 0$.  Thus $f \in \mathcal{I}$ when $\mathcal{I}$ is
closed with respect to the supremum norm, since $\theta_\epsilon \, f
\in \mathcal{I}$ for each $\epsilon > 0$.  This shows that
$\mathcal{I} = \mathcal{I}_E$ when $\mathcal{I}$ is a closed ideal in
$C(X)$ and $E$ is associated to $\mathcal{I}$ as before.

        Let $E$ be any closed set in $X$, and consider the
corresponding closed ideal $\mathcal{I}_E$.  In particular,
$\mathcal{I}_E$ is a linear subspace of $C(X)$, and the quotient space
$C(X) / \mathcal{I}_E$ can be defined as a real or complex vector
space, as appropriate.  By standard arguments in abstract algebra,
there is a natural operation of multiplication on the quotient, so
that the quotient mapping from $C(X)$ onto $C(X) / \mathcal{I}_E$ is a
multiplicative homomorphism, because $\mathcal{I}_E$ is an ideal in
$C(X)$.  If $E = \emptyset$, then $\mathcal{I}_E = C(X)$, and $C(X) /
\mathcal{I}_E = \{0\}$, and so we suppose from now on that $E \ne
\emptyset$.  We also have a homomorphism $R_E : C(X) \to C(E)$,
defined by sending a continuous function $f$ on $X$ to its restriction
$R_E(f)$ to $E$.  The kernel of this homomorphism is equal to
$\mathcal{I}_E$, which leads to a one-to-one homomorphism $r_E : C(X)
/ \mathcal{I}_E \to C(E)$.  By the Tietze extension theorem, every
continuous function on $E$ has an extension to a continuous function
on $X$.  This says exactly that $R_E(C(X)) = C(E)$, and hence that
$r_E$ maps $C(X) / \mathcal{I}_E$ onto $C(E)$.

\section{Locally compact spaces}
\label{locally compact spaces}
\setcounter{equation}{0}

        Let $X$ be a locally compact Hausdorff topological space, and
let $C(X)$ be the space of continuous real or complex-valued functions
on $X$, as usual.  If $K \subseteq X$ is nonempty and compact, then the
corresponding \emph{supremum seminorm} is defined on $C(X)$ by
\begin{equation}
        \|f\|_K = \sup_{x \in K} |f(x)|.
\end{equation}
Of course, every continuous function $f$ on $X$ is bounded on $K$,
because $f(K)$ is a compact set in ${\bf R}$ or ${\bf C}$, as appropriate.
It is easy to see that this is indeed a seminorm on $C(X)$, and that
\begin{equation}
        \|f \, g\|_K \le \|f\|_K \, \|g\|_K
\end{equation}
for every $f, g \in C(X)$.  It follows that multiplication of
functions is continuous as a mapping from $C(X) \times C(X)$ into
$C(X)$ with respect to the topology on $C(X)$ determined by the
collection of supremum seminorms associated to nonempty compact
subsets of $X$.

        If $X$ is compact, then we can take $X = K$, and simply use
the supremum norm on $X$.  Thus we shall focus on the case where $X$
is not compact in this section.  Suppose that $X$ is
\emph{$\sigma$-compact}, so that there is a sequence $K_1, K_2,
\ldots$ of compact subsets of $X$ such that $X = \bigcup_{l =
1}^\infty K_l$.  We may also ask that $K_l \ne \emptyset$ and $K_l
\subseteq K_{l + 1}$ for each $l$, by replacing $K_l$ with the union
of $K_1, \ldots, K_l$ if necessary.  Moreover, we can enlarge these
compact sets in such a way that $K_l$ is contained in the interior of
$K_{l + 1}$ for each $l$.  This uses the local compactness of $X$, to
get that any compact set in $X$ is contained in the interior of
another compact set.  In particular, it follows that the union of the
interiors of the $K_l$'s is all of $X$ under these conditions.  If $H$
is any compact set in $X$, then the interiors of the $K_l$'s form an
open covering of $H$, for which there is a finite subcovering.  This
implies that $H$ is contained in the interior of $K_l$ for some $l$,
since the $K_l$'s are increasing.  Hence $H \subseteq K_l$ for some
$l$, which implies that the supremum seminorms associated to the
$K_l$'s determine the same topology on $C(K)$ as the collection of
supremum seminorms corresponding to all nonempty compact subsets of
$X$.  Therefore this topology on $C(X)$ is metrizable in this case.

        If $E \subseteq X$ and $\mathcal{I}_E$ is the collection of $f
\in C(X)$ such that $f(x) = 0$ for every $x \in E$, then
$\mathcal{I}_E$ is a closed ideal in $C(X)$, as in the previous
section.  We also have that $\mathcal{I}_{\overline{E}} =
\mathcal{I}_E$, where $\overline{E}$ is the closure of $E$ in $X$.  If
$\mathcal{I}$ is any subset of $X$ and $E$ is the set of $x \in X$
such that $f(x) = 0$ for every $f \in \mathcal{I}$, then $E$ is a
closed set in $X$, and $\mathcal{I} \subseteq \mathcal{I}_E$.  We
would like to show that $\mathcal{I} = \mathcal{I}_E$ when
$\mathcal{I}$ is a closed ideal in $C(X)$, as before.  Let $f$ be a
continuous function on $X$ such that $f(x) = 0$ for every $x \in E$,
and let us check that $f \in \mathcal{I}$.

        If $f$ has compact support contained in the complement of $E$,
then one can show that $f \in \mathcal{I}$ in the same way as in the
previous section.  Otherwise, it suffices to show that $f$ can be
approximated by continuous functions with compact support contained in
$X \backslash E$ in the topology of $C(X)$.  Let $K$ be a nonempty
compact set in $X$, and let $\epsilon > 0$ be given.  Thus
\begin{equation}
\label{K(epsilon) = {x in K : |f(x)| ge epsilon}}
        K(\epsilon) = \{x \in K : |f(x)| \ge \epsilon\}
\end{equation}
is a compact set in $X$, since it is the intersection of the compact
set $K$ with the closed set where $|f(x)| \ge \epsilon$.  Also,
$K(\epsilon) \subseteq X \backslash E$, because $f = 0$ on $E$ by
hypothesis.  Let $V(\epsilon)$ be an open set in $X$ such that
$K(\epsilon) \subseteq V(\epsilon)$, $\overline{V(\epsilon)}$ is
compact, and $\overline{V(\epsilon)} \subseteq X \backslash E$.  This
is possible, because $X$ is locally compact and Hausdorff.  By
Urysohn's lemma, there is a continuous real-valued function
$\theta_\epsilon$ on $X$ which satisfies $\theta_\epsilon(x) = 1$ when
$x \in K(\epsilon)$, $\theta_\epsilon(x) = 0$ when $x \in X \backslash
V(\epsilon)$, and $0 \le \theta_\epsilon \le 1$ on all of $X$.  In
particular, the support of $\theta_\epsilon$ is contained in
$\overline{V(\epsilon)}$, which is a compact subset of $X \backslash
E$.  Hence $\theta_\epsilon \, f$ is a continuous function with
compact support in $X \backslash E$, which implies that
$\theta_\epsilon \, f \in \mathcal{I}$.  We also have that
\begin{equation}
 |\theta_\epsilon(x) \, f(x) - f(x)| = (1 - \theta_\epsilon(x)) \, |f(x)|
                                       < \epsilon
\end{equation}
for every $x \in K$, because $1 - \theta_\epsilon(x) = 0$ when $x \in
K(\epsilon)$, $|f(x)| < \epsilon$ when $x \in K \backslash
K(\epsilon)$, and $0 \le \theta_\epsilon(x) \le 1$ for every $x \in X$.
This shows that $f$ can be approximated by elements of $\mathcal{I}$
in the topology of $C(X)$, which implies that $f \in \mathcal{I}$,
as desired, since $\mathcal{I}$ is supposed to be closed in $C(X)$.

\section{Locally compact spaces, continued}
\label{locally compact spaces, continued}
\setcounter{equation}{0}

        Let $X$ be a locally compact Hausdorff space, and let
$C_{com}(X)$ be the space of continuous real or complex-valued
functions on $X$ with compact support.  As usual, this may be denoted
$C_{com}(X, {\bf R})$ or $C_{com}(X, {\bf C})$, to indicate whether
real or complex-valued functions are being used.  If $K \subseteq X$
is compact, then there is an open set $V$ in $X$ such that $K
\subseteq V$ and $\overline{V}$ is compact, because $X$ is locally
compact.  Urysohn's lemma implies that there is a continuous
real-valued function $\theta$ on $X$ such that $\theta(x) = 1$ when $x
\in K$, $\theta(x)= 0$ when $x \in X \backslash V$, and $0 \le
\theta(x) \le 1$ for every $x \in X$.  Thus the support of $\theta$ is
contained in $\overline{V}$, and hence is compact.  If $f$ is any
continuous function on $X$, then $\theta \, f$ is a continuous
function on $X$ with compact support that is equal to $f$ on $K$.  In
particular, this implies that $C_{com}(X)$ is dense in $C(X)$ with
respect to the topology determined by supremum seminorms associated to
nonempty compact subsets of $X$.

        Suppose that $\lambda$ is a continuous linear functional on
$C(X)$ with respect to this topology.  This implies that there is a nonempty
compact set $K \subseteq X$ and an nonnegative real number $A$ such that
\begin{equation}
\label{|lambda(f)| le A ||f||_K}
        |\lambda(f)| \le A \, \|f\|_K
\end{equation}
for every $f \in C(X)$.  In this context, it is not necessary to take
the maximum of finitely many seminorms on the right side of this
inequality, because the union of finitely many compact subsets of $X$
is also compact.  Note that $\lambda(f) = 0$ when $f(x) = 0$ for every
$x \in K$, so that $\lambda(f)$ depends only on the restriction of $f$
to $K$.

        Let $X^*$ be the one-point compactification of $X$.  Thus
$X^*$ is a compact Hausdorff space consisting of the elements of $X$
and an additional element ``at infinity'', for which the induced
topology on $X$ as a subset of $X^*$ is the same as its given
topology.  By construction, a set $K \subseteq X$ is closed as a
subset of $X^*$ if and only if it is compact in $X$.  In this case,
Tietze's extension theorem implies that every continuous function on
$K$ can be extended to a continuous function on $X^*$, and to a
continuous function on $X$ in particular.  If $X$ is already compact,
then one can simply use $X$ instead of $X^*$.

        If $\lambda$ is a continuous linear functional on $C(X)$ that
satisfies (\ref{|lambda(f)| le A ||f||_K}), then it follows that
$\lambda$ corresponds to a continuous linear functional $\lambda_K$ on
$C(K)$ in a natural way.  More precisely, if $g$ is a continuous
function on $K$, then there is a continuous function $f$ on $X$ such
that $f = g$ on $K$, and we put
\begin{equation}
        \lambda_K(g) = \lambda(f).
\end{equation}
This does not depend on the particular extension $f$ of $g$, by the
earlier remarks.  By construction, $\lambda_K$ satisfies the same
continuity condition on $C(K)$ as $\lambda$ does on $C(X)$, with the
same constant $A$.

        Let $C_0(X)$ be the space of continuous functions $f$ on $X$
that vanish at infinity.  This means that for every $\epsilon > 0$
there is a compact set $K_\epsilon \subseteq X$ such that
\begin{equation}
\label{|f(x)| < epsilon}
        |f(x)| < \epsilon
\end{equation}
for every $x \in X \backslash K_\epsilon$.  This space may also be
denoted $C_0(X, {\bf R})$ or $C_0(X, {\bf C})$, to indicate whether
real or complex-valued functions are being used.  If $X$ is compact,
then one can take $K_\epsilon = X$ for each $\epsilon$, and $C_0(X) =
C(X)$.  If $X$ is not compact, then $f \in C_0(X)$ if and only if $f$
has a continuous extension to the one-point compactification $X^*$ of
$X$ which is equal to $0$ at the point at infinity.

        Note that continuous functions on $X$ that vanish at infinity
are automatically bounded, so that $C_0(X) \subseteq C_b(X)$.  It is
not difficult to check that $C_0(X)$ is a closed linear subspace of
$C_b(X)$, with respect to the supremum norm.  Of course, continuous
functions with compact support automatically vanish at infinity, so
that $C_{com}(X) \subseteq C_0(X)$.  One can also check that $C_0(X)$
is the same as the closure of $C_{com}(X)$ in $C_b(X)$ with respect to
the supremum norm, using functions $\theta$ as before.  This is all
trivial when $X$ is compact, in which case these spaces are all the
same as $C(X)$.

        Suppose that $X$ is not compact, let $f$ be a real or
complex-valued continuous function on $X$, and let $a$ be a real or
complex number, as appropriate.  We say that $f(x) \to a$ as $x \to
\infty$ in $X$ if for each $\epsilon > 0$ there is a compact set
$K_\epsilon \subseteq X$ such that
\begin{equation}
        |f(x) - a| < \epsilon
\end{equation}
for every $x \in X \backslash K_\epsilon$.  Thus $f$ vanishes at
infinity if and only if this holds with $a = 0$.  If $a$ is any real
or complex number, then $f(x) \to a$ as $x \to \infty$ in $X$ if and
only if $f(x) - a$ vanishes at infinity.  It is easy to see that the
limit $a$ is unique when it exists.  Similarly, $f(x) \to a$ as $x \to
\infty$ in $X$ if and only if $f$ has a continuous extension to the
one-point compactification $X^*$ of $X$ which is equal to $a$ at the
point at infinity.  Note that $f$ is bounded when $f$ has a limit at
infinity.  One can also check that the collection of continuous
functions on $X$ which have a limit at infinity is a closed linear
subspace of $C_b(X)$ with respect to the supremum norm.

\section{$\sigma$-Compactness}
\label{sigma-compactness}
\setcounter{equation}{0}

        Let $X$ be a topological space, and let $\{U_\alpha\}_{\alpha
\in A}$ be a collection of open subsets of $X$ such that
$\bigcup_{\alpha \in A} U_\alpha = X$, which is to say an open
covering of $X$.  Suppose that $X$ is $\sigma$-compact, so that there
is a sequence $K_1, K_2, \ldots$ of compact subsets of $X$ such that
$X = \bigcup_{l = 1}^\infty K_l$.  Because $\{U_\alpha\}_{\alpha \in
A}$ is an open covering of $K_l$ for each $l$ and $K_l$ is compact,
there is a finite set of indices $A_l \subseteq A$ such that $K_l
\subseteq \bigcup_{\alpha \in A_l} U_\alpha$.  If $B = \bigcup_{l =
1}^\infty A_l$, then $B$ has only finitely or countably many elements,
and $\bigcup_{\alpha \in B} U_\alpha = X$.  Conversely, if $X$ is
locally compact and every open covering of $X$ can be reduced to a
subcovering with only finitely or countably many elements, then $X$ is
$\sigma$-compact.  This follows by using local compactness to cover
$X$ by open sets that are contained in compact sets.  In particular,
if $X$ is locally compact and there is a base for the topology of $X$
with only finitely or countably many elements, then $X$ is
$\sigma$-compact, since every open covering of $X$ can be reduced to a
subcovering with only finitely or countably many elements in this case.

        Suppose that the topology on $X$ is determined by a metric.
It is well known that there is a base for the topology of $X$ with
only finitely or countably many elements if and only if $X$ is
separable, in the sense that there is a dense set in $X$ with only
finitely or countably many elements.  Compact metric spaces are
separable, and it follows that $X$ is separable when $X$ is
$\sigma$-compact.  Urysohn's famous metrization theorem states that a
regular topological space is metrizable when there is a countable base
for its topology.  Note that locally compact Hausdorff spaces are
automatically regular.

        Suppose now that $X$ is a locally compact Hausdorff
topological space which is $\sigma$-compact.  As before, this implies
that there is a sequence $K_1, K_2, \ldots$ of compact subsets of $X$
such that $X = \bigcup_{l = 1}^\infty K_l$ and $K_l$ is contained in
the interior of $K_{l + 1}$ for each $l$.  By Urysohn's lemma, there
is a continuous real-valued function $\theta_l$ on $X$ for each
positive integer $l$ such that $\theta(x) > 0$ when $x \in K_l$, $0
\le \theta_l(x) \le 1$ for every $x \in X$, and the support of
$\theta_l$ is contained in $K_{l + 1}$.  Let $a_1, a_2, \ldots$ be a
sequence of positive real numbers such that $\sum_{l = 1}^\infty a_l$
converges, and consider
\begin{equation}
        f(x) = \sum_{l = 1}^\infty a_l \, \theta_l(x).
\end{equation}
This series converges everywhere on $X$, by the comparison test.  The
partial sums of this series converge uniformly on $X$, as in
Weierstrass' $M$-test.  Thus $f$ is a continuous function on $X$,
which also vanishes at infinity, because $\theta_l$ has compact
support for each $l$.  Moreover, $f(x) > 0$ for every $x \in X$, by
construction.

\section{Homomorphisms, revisited}
\label{homomorphisms, revisited}
\setcounter{equation}{0}

        Let $X$ be a locally compact Hausdorff topological space, and
let $C(X)$ be the algebra of real or complex-valued continuous
functions on $X$.  Also let $\phi$ be linear functional on $C(X)$
which is a homomorphism with respect to multiplication.  If $\phi(f)
\ne 0$ for some $f \in C(X)$, then it follows that $\phi({\bf 1}_X) =
1$, where ${\bf 1}_X$ is the constant function equal to $1$ on $X$, as
before.  Let us suppose from now on that this is the case.  If $f$ is
a continuous function on $X$ such that $f(x) \ne 0$ for every $x \in
X$, then $1/f$ is defines a continuous function on $X$ as well.  This
implies that $\phi(f) \ne 0$, since
\begin{equation}
        \phi(f) \, \phi(1/f) = \phi({\bf 1}_X) = 1.
\end{equation}
If $f$ is any continuous function on $X$ and $c$ is a real or complex
number, as appropriate, such that $c \not\in f(X)$, then $g = f - c \,
{\bf 1}_X$ is a continuous function on $X$ such that $g(x) \ne 0$ for
every $x \in X$, so that $\phi(g) \ne 0$.  Thus $\phi(f) \ne c$, and
hence
\begin{equation}
        \phi(f) \in f(X).
\end{equation}
In particular, if $C(X)$ is the algebra of complex-valued continuous
functions on $X$, and $f$ happens to be real-valued, then it follows
that $\phi(f) \in {\bf R}$.

        Suppose now that $\phi$ is continuous with respect to the
topology on $C(X)$ determined by the supremum seminorms corresponding
to nonempty compact subsets of $X$.  This means that there is a
nonempty compact set $K \subseteq X$ and a nonnegative real number $A$
such that
\begin{equation}
        |\phi(f)| \le A \, \|f\|_K
\end{equation}
for every $f \in C(X)$, as in Section \ref{locally compact spaces, continued}.
In particular, $\phi(f) = 0$ when $f(x) = 0$ for every $x \in K$, so that
$\phi(f)$ depends only on the restriction of $f$ to $K$.  As in Section
\ref{locally compact spaces, continued} again, every continuous real or
complex-valued function on $K$ has a continuous extension to $X$,
so that $\phi$ determines a continuous linear functional $\phi_K$ on $C(K)$.
It is easy to see that $\phi_K$ is also a homomorphism with respect to
multiplication on $C(K)$.  Hence there is a $p \in K$ such that
$\phi_K(f) = f(p)$ for every $f \in C(K)$, as in Section \ref{compact
spaces}.  This implies that
\begin{equation}
        \phi(f) = f(p)
\end{equation}
for every $f \in C(X)$.

        Alternatively, consider
\begin{equation}
        \mathcal{I}_\phi = \{f \in C(X) : \phi(f) = 0\}.
\end{equation}
It is easy to see that this is a closed ideal in $C(X)$ when $\phi$ is
a continuous homomorphism on $C(X)$.  As in Section \ref{locally
compact spaces}, there is a closed set $E \subseteq X$ such that
$\mathcal{I}_\phi = \mathcal{I}_E$, where $\mathcal{I}_E$ consists of
$f \in C(X)$ such that $f(x) = 0$ for every $x \in E$.  Note that
$\mathcal{I}_\phi$ has codimension $1$ as a linear subspace of $C(X)$,
since it is the same as the kernel of the nonzero linear functional
$\phi$.  Using this, one can check that $E$ has exactly one element,
which may be denoted $p$.  Thus $\phi(f) = 0$ for every $f \in C(X)$
such that $f(p) = 0$.  If $f$ is any continuous function on $X$, then
$f - f(p) \, {\bf 1}_X$ is equal to $0$ at $p$, and hence $\phi(f -
f(p) \, {\bf 1}_X) = 0$.  This implies that $\phi(f) = f(p)$ for every
$f \in C(X)$, since $\phi({\bf 1}_X) = 1$.

        Remember that the same conclusion holds for every nonzero
homomorphism $\phi$ on $C(X)$ when $X$ is compact, without the
additional hypothesis of continuity, as in Section \ref{compact
spaces}.  Suppose now that $X$ is a locally compact Hausdorff which is
not compact but $\sigma$-compact, and that $\phi$ is a nonzero
homomorphism on $C(X)$.  Let $X^*$ be the one-point compactification
of $X$, and note that the space $C(X^*)$ of continuous functions on
$X^*$ can be identified with the subalgebra of $C(X)$ consisting of
functions with a limit at infinity, as in Section \ref{locally compact
spaces, continued}.  The restriction of $\phi$ to this subalgebra
determines a homomorphism on $C(X^*)$, which is nonzero because it
sends constant functions to their constant values.  It follows that
there is a $p \in X^*$ such that $\phi(f) = f(p)$ when $f \in C(X)$
has a limit at infinity, as in Section \ref{compact spaces}.  If $p$
is the point at infinity in $X^*$, then $f(p)$ refers to the limit of
$f$ at infinity on $X$.  Let us check that $p$ cannot be the point at
infinity in $X^*$ when $X$ is $\sigma$-compact.  In this case, there
is a continuous real-valued function $f$ on $X$ that vanishes at
infinity such that $f(x) > 0$ for every $x \in X$, as in the previous
section.  Because $\phi$ is defined on all of $C(X)$, we also have
that $\phi(f) \ne 0$, as discussed at the beginning of the section.
If $p$ were the point at infinity, then we would have that $\phi(f) =
0$, since $f \in C_0(X)$.  Thus $p \in X^*$ is not the point at
infinity, which means that $p \in X$.  If $g$ is any bounded
continuous function on $X$, then $f \, g \in C_0(X)$, which implies
that
\begin{equation}
        \phi(f \, g) = f(p) \, g(p),
\end{equation}
and so
\begin{equation}
        \phi(f) \, \phi(g) = f(p) \, g(p),
\end{equation}
because $\phi$ is a homomorphism on $C(X)$.  This shows that $\phi(g)
= g(p)$ for every bounded continuous function $g$ on $X$.  If $h$ is
any continuous function on $X$ and $\epsilon > 0$, then
\begin{equation}
        h_\epsilon = \frac{h}{1 + \epsilon \, |h|^2}
\end{equation}
is a bounded continuous function on $X$, and so $\phi(h_\epsilon) =
h_\epsilon(p)$.  One can also check that
\begin{equation}
        \phi(h_\epsilon) = \frac{\phi(h)}{1 + \epsilon \, |\phi(h)|^2}
\end{equation}
for every $\epsilon > 0$, because $\phi$ is a homomorphism.  Hence
\begin{equation}
        \frac{\phi(h)}{1 + \epsilon \, |\phi(h)|^2}
           = \frac{h(p)}{1 + \epsilon \, |h(p)|^2}
\end{equation}
for every $\epsilon > 0$, which implies that $\phi(h) = h(p)$ for
every $h \in C(X)$.

\section{$\sigma$-Compactness, continued}
\label{sigma-compactness, continued}
\setcounter{equation}{0}

        Let $X$ be a locally compact Hausdorff topological space which
is $\sigma$-compact, and let $K_1, K_2, \ldots$ be a sequence of
compact subsets of $X$ such that $X = \bigcup_{l = 1}^\infty K_l$ and
$K_l$ is contained in the interior of $K_{l + 1}$ for each $l$.  By
Urysohn's lemma, there is a continuous real-valued function $\theta_l$
on $X$ for each positive integer $l$ such that $\theta_l(x) = 1$ for
every $x$ in a neighborhood of $K_l$, $0 \le \theta_l(x) \le 1$ for
every $x \in X$, and the support of $\theta_l$ is contained in $K_{l +
1}$.  In particular, $\theta_l(x) \le \theta_{l + 1}(x)$ for each $x
\in X$ and $l \ge 1$.  It will be convenient to also put $K_0 =
\emptyset$ and $\theta_0 = 0$.  Let $b_1, b_2, \ldots$ be a sequence
of nonnegative real numbers, and consider
\begin{equation}
\label{B(x) = ...}
        B(x) = b_1 \, \theta_1(x) + \sum_{l = 2}^\infty b_l \,
                                      (\theta_l(x) - \theta_{l - 2}(x)).
\end{equation}
Note that $\theta_l(x) - \theta_{l - 2}(x) = 0$ for every $x$ in a
neighborhood of $K_{l - 2}$, and when $x \in X \backslash K_{l + 1}$,
for $l \ge 2$.  This implies that at most three terms on the right
side of (\ref{B(x) = ...})  are different from $0$ for any $x \in X$,
and more precisely that every $x \in X$ has a neighborhood on which at
most three terms on the right side of (\ref{B(x) = ...}) are different
from $0$, so that $B(x)$ is continuous on $X$.  We also have that
\begin{equation}
        B(x) \ge b_1 \, \theta_1(x) \ge b_1
\end{equation}
when $x \in K_1$, and
\begin{equation}
        B(x) \ge b_l \, (\theta_l(x) - \theta_{l - 2}(x)) \ge b_l
\end{equation}
when $x \in K_l \backslash K_{l - 1}$, $l \ge 2$.

        Suppose that $E$ is a bounded subset of the space $C(X)$ of
continuous real or complex-valued continuous functions on $X$ with
respect to the collection of supremum seminorms associated to nonempty
compact subsets of $X$.  Thus the elements of $E$ are uniformly
bounded on compact subsets of $X$, and so for each positive integer
$l$ there is a nonnegative real number $b_l$ such that
\begin{equation}
        |f(x)| \le b_l
\end{equation}
for every $f \in E$ and $x \in K_l$.  This implies that
\begin{equation}
\label{|f(x)| le B(x)}
        |f(x)| \le B(x)
\end{equation}
for every $f \in E$ and $x \in X$, where $B$ is as in the previous
paragraph.

        Now let $\phi$ be a linear functional on $C(X)$ which is a
homomorphism with respect to multiplication, and which satisfies
$\phi(f) \ne 0$ for some $f \in C(X)$.  As in the previous section,
$\phi(f) \in f(X)$ for every $f \in C(X)$, and in particular $\phi(f)
\ge 0$ when $f$ is a nonnegative real-valued continuous function on
$X$.  If $f \in C(X)$ satisfies (\ref{|f(x)| le B(x)}), then it follows that
\begin{equation}
\label{|phi(f)| le phi(B)}
        |\phi(f)| \le \phi(B).
\end{equation}
More precisely, if $f$ is real-valued, then $B(x) \pm f(x) \ge 0$ for
each $x \in X$, and so
\begin{equation}
        \phi(B) \pm \phi(f) = \phi(B \pm f) \ge 0.
\end{equation}
Similarly, if $f$ is complex-valued, then one can use the fact that
$\re \alpha f(x) \le B$ for each $\alpha \in {\bf C}$ with $|\alpha| =
1$ to get that
\begin{equation}
        \re \alpha \, \phi(f) = \phi(\re \alpha \, f) \le \phi(B),
\end{equation}
which implies (\ref{|phi(f)| le phi(B)}).

         This shows that $\phi$ is uniformly bounded on every bounded
set $E \subseteq C(X)$ with respect to the collection of supremum
seminorms associated to nonempty compact subsets of $X$.  There is
also a countable local base for the topology at $0$ in $C(X)$ with
respect to this collection of seminorms, because $X$ is
$\sigma$-compact.  It follows that $\phi$ is continuous with respect
to this topology on $C(X)$, by the result discussed in Section
\ref{bounded sequences}.  This gives another way to show that there is
a point $p \in X$ such that $\phi(f) = f(p)$ for every $f \in C(X)$,
by reducing to the case of continuous homomorphisms, as in the
previous section.

\section{Holomorphic functions}
\label{holomorphic functions}
\setcounter{equation}{0}

        Let $U$ be a nonempty open set in the complex plane ${\bf C}$,
and let $C(U)$ be the algebra of continuous complex-valued functions
on $U$.  Of course, $U$ is locally compact with respect to the
topology inherited from the standard topology on ${\bf C}$, and it is
also $\sigma$-compact, because it is a separable metric space, and
hence has a countable base for its topology.  As usual, $C(U)$ gets a
nice topology from the collection of supremum seminorms associated to
nonempty compact subsets of $U$.

        Remember that a complex-valued function $f(z)$ on $U$ is said
to be complex-analytic or holomorphic if the complex derivative
\begin{equation}
\label{f'(z) = lim_{h to 0} frac{f(z + h) - f(z)}{h}}
        f'(z) = \lim_{h \to 0} \frac{f(z + h) - f(z)}{h}
\end{equation}
exists at every point $z$ in $U$.  In particular, the existence of the
limit implies that $f$ is continuous, so that the space
$\mathcal{H}(U)$ of holomorphic functions on $U$ is contained in
$C(U)$.  More precisely, $\mathcal{H}(U)$ is a linear subspace of
$C(U)$, which is actually a subalgebra, because the product of two
holomorphic functions is holomorphic as well.  Note that constant
functions on $U$ are automatically holomorphic, since they have
derivative equal to $0$ at every point.

        It is well known that $\mathcal{H}(U)$ is closed in $C(U)$,
with respect to the topology determined by the collection of supremum
seminorms associated to nonempty compact subsets of $U$.  This is
equivalent to the statement that if $\{f_j\}_{j = 1}^\infty$ is a
sequence of holomorphic functions on $U$ that converges uniformly on
compact subsets of $U$ to a function $f$ on $U$, then $f$ is also
holomorphic on $U$.  To see this, one can use the Cauchy integral
formula to show that the sequence of derivatives $\{f'_j\}_{j =
1}^\infty$ converges uniformly on compact subsets of $U$, and that the
limit is equal to the derivative $f'$ of $f$.

        Let $\phi$ be a linear functional on $\mathcal{H}(U)$ which is
a homomorphism with respect to multiplication.  As before, if $\phi(f)
\ne 0$ for some $f \in \mathcal{H}(U)$, then $\phi({\bf 1}_U) = 1$,
where ${\bf 1}_U$ is the constant function on $U$ equal to $1$.  Let
us suppose from now on that this is the case.  If $f$ is a holomorphic
function on $U$ such that $f(z) \ne 0$ for every $z \in U$, then it is
well known that $1/f$ is holomorphic on $U$ too.  This implies that
\begin{equation}
        \phi(f) \, \phi(1/f) = \phi({\bf 1}_U) = 1,
\end{equation}
and hence $\phi(f) \ne 0$.  If $c$ is a complex number such that $c
\not\in f(U)$, then we can apply this to $f - c \, {\bf 1}_U$ to get
that $\phi(f) \ne c$.  Thus $\phi(f) \in f(U)$, as in the context of
continuous functions.  In particular, this holds when $f(z) = z$ for
every $z \in U$, which is holomorphic with derivative equal to $1$ at
every point.  If $\phi(f)$ is denoted $p$ when $f(z) = z$ for every $z
\in U$, then it follows that $p \in U$.  We would like to show that
\begin{equation}
        \phi(g) = g(p)
\end{equation}
for every $g \in \mathcal{H}(U)$.  If $g(p) = 0$, then $g$ can be
expressed as
\begin{equation}
        g(z) = (z - p) \, h(z)
\end{equation}
for some $h \in \mathcal{H}(U)$, by standard results in complex
analysis.  This implies that $\phi(g) = 0$, by the definition of $p$
and the fact that $\phi$ is a homomorphism.  If $g(p) \ne 0$, then one
can reduce to the case where $g(p) = 0$ by subtracting a constant from
$g$.

\section{The disk algebra}
\label{disk algebra}
\setcounter{equation}{0}

        Let $U$ be the open unit disk in the complex plane ${\bf C}$,
\begin{equation}
        U = \{z \in {\bf C} : |z| < 1\}.
\end{equation}
Thus the closure $\overline{U}$ of $U$ is the closed unit disk,
\begin{equation}
        \overline{U} = \{z \in {\bf C} : |z| \le 1\},
\end{equation}
and the boundary $\partial U$ of $U$ is the same as the unit circle,
\begin{equation}
        \partial U = \{z \in {\bf C} : |z| = 1\}.
\end{equation}
Let $C(\overline{U})$ be the algebra of continuous complex-valued functions
on $\overline{U}$, equipped with the supremum norm.

        Let $\mathcal{A}$ be the collection of $f \in C(\overline{U})$
such that the restriction of $f$ to $U$ is holomorphic.  Thus
$\mathcal{A}$ is a subalgebra of $C(\overline{U})$, since sums and
products of holomorphic functions are also holomorphic, which is known
as the \emph{disk algebra}.  Note that constant functions on
$\overline{U}$ are elements of $\mathcal{A}$, and that $\mathcal{A}$
is a closed set in $C(\overline{U})$ with respect to the supremum
norm, for the same reasons as in the previous section.  If $f \in
\mathcal{A}$ and $f(z) \ne 0$ for every $z \in \overline{U}$, then
$1/f$ is continuous on $\overline{U}$ and holomorphic on $U$, and
hence is in $\mathcal{A}$ too.

        If $f \in C(\overline{U})$ and $0 \le r < 1$, then
\begin{equation}
        f_r(z) = f(r \, z)
\end{equation}
is an element of $C(\overline{U})$ as well.  Note that $f$ is
automatically uniformly continuous on $\overline{U}$, because $f$ is
continuous on $\overline{U}$ and $\overline{U}$ is a compact set in a
metric space.  Using this, it is easy to see that $f_r \to f$
uniformly on $\overline{U}$ as $r \to 1$.

        If $f$ is a holomorphic function on the open unit disk $U$, then
\begin{equation}
        f(z) = \sum_{j = 0}^\infty a_j \, z^j
\end{equation}
for some complex numbers $a_0, a_1, \ldots$ and every $z \in U$.  More
precisely, $z^j$ is interpreted as being equal to $1$ for every $z$
when $j = 0$, and the convergence of the series when $|z| < 1$ is part
of the conclusion.  The series actually converges absolutely for every
$z \in U$, and the partial sums converge uniformly on compact subsets
of $U$.

        If $0 \le r < 1$, then
\begin{equation}
        f_r(z) = f(r \, z) = \sum_{j = 0}^\infty a_j r^j \, z^j
\end{equation}
for every $z \in \overline{U}$.  Under these conditions, the series
converges absolutely when $|z| \le 1$, and the partial sums converge
uniformly on $\overline{U}$, by the remarks in the previous paragraph.
If $f \in \mathcal{A}$, then $f$ can be approximated uniformly by
$f_r$ as $r \to 1$, and $f_r$ is approximated uniformly by partial
sums of its series expansion for each $r < 1$.  It follows that $f$
can be approximated uniformly by polynomials in $z$ on $\overline{U}$
when $f \in \mathcal{A}$.

        Let $\phi$ be a linear functional on $\mathcal{A}$ which is a
homomorphism with respect to multiplication.  As usual, we suppose
that $\phi(f) \ne 0$ for some $f \in \mathcal{A}$, so that $\phi$
sends constant functions on $\overline{U}$ to their constant values.
If $f \in \mathcal{A}$ and $f(z) \ne 0$ for every $z \in
\overline{U}$, then $1/f \in \mathcal{A}$, and we get that $\phi(f)
\ne 0$.  This implies that
\begin{equation}
        \phi(f) \in f(\overline{U})
\end{equation}
for every $f \in \mathcal{A}$, as before.  In particular,
\begin{equation}
\label{|phi(f)| le sup_{|z| le 1} |f(z)|}
        |\phi(f)| \le \sup_{|z| \le 1} |f(z)|
\end{equation}
for every $f \in \mathcal{A}$, so that $\phi$ is continuous with
respect to the supremum norm on $\mathcal{A}$.  Of course, $f(z) = z$
defines an element of $\mathcal{A}$, and we can put $\phi(f) = p$ for
this choice of $f$.  Note that $p \in \overline{U}$, by the previous
remarks.  If $g$ is a polynomial in $z$, then
\begin{equation}
        \phi(g) = g(p),
\end{equation}
because $\phi$ is a homomorphism.  This also works for every $g \in
\mathcal{A}$, because polynomials are dense in $\mathcal{A}$ with
respect to the supremum norm, and because $\phi$ is continuous on
$\mathcal{A}$ with respect to the supremum norm.

        If $f \in \mathcal{A}$, then
\begin{equation}
        \sup_{|z| = 1} |f(z)| = \sup_{|z| \le 1} |f(z)|,
\end{equation}
by the maximum modulus principle.  In particular, if $f(z) = 0$ for
every $z \in \partial U$, then $f(z) = 0$ for every $z \in
\overline{U}$.  This implies that $f$ is determined on the closed disk
$\overline{U}$ by its restriction to the unit circle $\partial U$.
Using this, one can identify the disk algebra with a closed subalgebra
of the algebra of continuous complex-valued functions on the unit
circle.

\section{Bounded holomorphic functions}
\label{bounded holomorphic functions}
\setcounter{equation}{0}

        Let $U$ be the open unit disk in the complex plane again, and
let $C_b(U)$ be the algebra of bounded continuous complex-valued
functions on $U$, equipped with the supremum norm.  Also let
$\mathcal{B}$ be the collection of bounded holomorphic functions on
$U$, which is the same as the intersection of $C_b(U)$ with
$\mathcal{H}(U)$.  As usual, this is a closed subalgebra of $C_b(U)$
with respect to the supremum norm.

        Let $\phi$ be a linear functional on $\mathcal{B}$ which is a
homomorphism with respect to multiplication.  Suppose also that
$\phi(f) \ne 0$ for some $f \in \mathcal{B}$, which implies that
$\phi$ sends constant functions on $U$ to their constant values.  If
$f \in \mathcal{B}$ and $|f(z)| \ge \delta$ for some $\delta > 0$ and
every $z \in U$, then $1/f$ is also a bounded holomorphic function on
$U$, and it follows that $\phi(f) \ne 0$, because $\phi(f) \,
\phi(1/f) = 1$.  This implies that
\begin{equation}
        \phi(f) \in \overline{f(U)}
\end{equation}
for every $f \in \mathcal{B}$, as in the previous situations, and
hence that
\begin{equation}
        |\phi(f)| \le \sup_{|z| < 1} |f(z)|.
\end{equation}
Thus $\phi$ is a continuous linear functional on $\mathcal{B}$ with
respect to the supremum norm, with dual norm equal to $1$, since
$\phi$ sends constants to themselves.

        Each element $p$ of $U$ determines a nonzero homomorphism
$\phi_p$ on $\mathcal{B}$, given by evaluation at $p$, or
\begin{equation}
        \phi_p(f) = f(p).
\end{equation}
The collection of nonzero homomorphisms on $\mathcal{B}$ is contained
in the unit ball of the dual of $\mathcal{B}$ with respect to the
supremum norm, as in the previous paragraph, and it is also a closed
set with respect to the weak$^*$ topology, as in Section \ref{bounded
continuous functions}.  Hence the collection of nonzero homomorphisms
on $\mathcal{B}$ is compact with respect to the weak$^*$ topology on
the dual of $\mathcal{B}$, by the Banach--Alaoglu theorem.  Of course,
the restriction of any nonzero homomorphism on $C_b(U)$ is a nonzero
homomorphism on $\mathcal{B}$, which includes evaluation at elements
of $U$.

        Suppose that $z_1, z_2, \ldots$ is a sequence of elements of
$U$ such that $|z_j| \to 1$ as $j \to \infty$.  Also let $L$ be a
nonzero homomorphism on $\ell^\infty({\bf Z}_+)$ which is equal to $0$
on $c_0({\bf Z}_+)$.  This determines a nonzero homomorphism on
$C_b(U)$, by applying $L$ to $f(z_j)$ as a bounded function on ${\bf
Z}_+$ for each $f \in C_b(U)$.  If $w_1, w_2, \ldots$ is another
sequence of elements of $U$ such that $|w_j| \to 1$ as $j \to \infty$,
then we can apply $L$ to $f(w_j)$ to get another homomorphism on
$C_b(U)$.  If $z_j \ne w_l$ for every $j, l \ge 1$, then it is easy to
see that these are distinct homomorphisms on $C_b(U)$, because one can
choose a bounded continuous function $f$ on $U$ such that $f(z_j) = 0$
and $f(w_l) = 1$ for each $j$, $l$.

        If $f$ is a bounded holomorphic function on $U$, then one can
check that there is a $C \ge 0$ such that
\begin{equation}
        \sup_{|z| < 1} (1 - |z|) \, |f'(z)| \le C \, \sup_{|z| < 1} |f(z)|.
\end{equation}
This follows from the Cauchy integral formula for $f'(z)$ applied to
the disk centered at $z$ with radius $(1 - |z|) / 2$, for instance.

        Suppose that $z_j$, $w_l$ are as before, and satisfy the
additional property that
\begin{equation}
        \lim_{j \to \infty} \frac{|z_j - w_j|}{(1 - |z_j|)} = 0.
\end{equation}
If $f$ is a bounded holomorphic function on $U$, then
\begin{equation}
        \lim_{j \to \infty} (f(z_j) - f(w_j)) = 0.
\end{equation}
This follows from the fact that $(1 - |z|) |f'(z)|$ is bounded on $U$,
as in the previous paragraph.  If $L$ is a nonzero homomorphism on
$\ell^\infty({\bf Z}_+)$ that vanishes on $c_0({\bf Z}_+)$, then $L$
applied to $f(z_j) - f(w_j)$ is equal to $0$, so that $L$ applied to
$f(z_j)$ is the same as $L$ applied to $f(w_j)$.  This shows that
distinct homomorphisms on $C_b(U)$ may determine the same homomorphism
on $\mathcal{B}$.

        A sequence $\{z_j\}_{j = 1}^\infty$ of points in $U$ is said
to be an \emph{interpolating sequence} if for every bounded sequence
of complex numbers $\{a_j\}_{j = 1}^\infty$ there is a bounded
holomorphic function $f$ on $U$ such that $f(z_j) = a_j$ for each $j$.
Equivalently, $\{z_j\}_{j = 1}^\infty$ is an interpolating sequence in
$U$ if
\begin{equation}
\label{f mapsto {f(z_j)}_{j = 1}^infty}
        f \mapsto \{f(z_j)\}_{j = 1}^\infty
\end{equation}
maps $\mathcal{B}$ onto $\ell^\infty({\bf Z}_+)$.  Of course, (\ref{f
mapsto {f(z_j)}_{j = 1}^infty}) defines a bounded linear mapping from
$\mathcal{B}$ into $\ell^\infty({\bf Z}_+)$ for any sequence
$\{z_j\}_{j = 1}^\infty$ of elements of $U$, and is also a
homomorphism with respect to pointwise multiplication.  A famous
theorem of Carleson characterizes interpolating sequences in $U$.  In
particular, there are plenty of them.

\section{Density}
\label{density}
\setcounter{equation}{0}

        Let $X$ be a topological space, and let $\psi$ be a nonzero
homomorphism from $C_b(X)$ into the real or complex numbers, as
appropriate.  As in Section \ref{bounded continuous functions}, $\psi$
is automatically a bounded linear functional on $C_b(X)$, and thus an
element of the dual space $C_b(X)^*$.  If $p \in X$, then let
$\phi_p(f) = f(p)$ be the corresponding point evaluation homomorphism
on $C_b(X)$, as usual.  We would like to show that $\psi$ can be
approximated by point evealuations with respect to the weak$^*$
topology on $C_b(X)^*$, so that point evaluations are dense in the set
of nonzero homomorphisms on $C_b(X)$ with respect to the weak$^*$
topology on $C_b(X)^*$.

        More precisely, we would like to show that for any finite
collection of bounded continuous functions $f_1, \ldots, f_n$ on $X$
and any $\epsilon > 0$ there is a $p \in X$ such that
\begin{equation}
        |\psi(f_j) - \phi_p(f)| = |\psi(f_j) - f_j(p)| < \epsilon
\end{equation}
for $j = 1, \ldots, n$.  Otherwise, there are $f_1, \ldots, f_n \in
C_b(X)$ and $\epsilon > 0$ such that
\begin{equation}
\label{max_{1 le j le n} |psi(f_j) - f_j(p)| ge epsilon}
        \max_{1 \le j \le n} |\psi(f_j) - f_j(p)| \ge \epsilon
\end{equation}
for every $p \in X$.  We may as well ask also that $\psi(f_j) = 0$ for
each $j$, since this can always be arranged by subtracting $\psi(f_j)$
as a constant function on $U$ from $f_j$.  In this case, (\ref{max_{1
le j le n} |psi(f_j) - f_j(p)| ge epsilon}) reduces to
\begin{equation}
\label{max_{1 le j le n} |f_j(p)| ge epsilon}
        \max_{1 \le j \le n} |f_j(p)| \ge \epsilon
\end{equation}
for each $p \in X$.

        If
\begin{equation}
        g(p) = \sum_{j = 1}^n |f_j(p)|^2,
\end{equation}
then $g$ is a bounded continuous function on $X$, and $g(p) \ge
\epsilon^2$ for each $p \in U$, by (\ref{max_{1 le j le n} |f_j(p)| ge
epsilon}).  Thus $1/g$ is also a bounded continuous function on $X$,
which implies that $\psi(g) \ne 0$, as in Section \ref{bounded
continuous functions}.  Of course, $g$ can also be expressed as
\begin{equation}
        g = \sum_{j = 1}^n f_j^2
\end{equation}
in the real case, and as
\begin{equation}
        g = \sum_{j = 1}^n f_j \, \overline{f_j}
\end{equation}
in the complex case, where $\overline{f_j}$ is the complex conjugate of $f_j$.
In both cases, this implies that $\psi(g) = 0$, a contradiction, because 
$\psi(f_j) = 0$ for each $j$, and $\psi$ is a homomorphism.

        Let $\mathcal{B}$ be the algebra of bounded holomorphic
functions on the open unit disk $U$, as in the preceding section.
Carleson's corona theorem states that every nonzero homomorphism
$\psi$ on $\mathcal{B}$ can be approximated by point evaluations
$\phi_p(f) = f(p)$, $p \in U$, with respect to the weak$^*$ topology
on the dual of $\mathcal{B}$.  As before, if this were not the case,
then there would be bounded holomorphic functions $f_1, \ldots, f_n$
on $U$ and $\epsilon > 0$ such that $\psi(f_j) = 0$ for $j = 1,
\ldots, n$ and (\ref{max_{1 le j le n} |f_j(p)| ge epsilon}) holds.
However, the previous argument does not work, because $\overline{f_j}$
is not holomorphic on $U$ unless $f_j$ is constant.  Instead, one can
try to show that there are bounded holomorphic functions $g_1, \ldots,
g_n$ on $U$ such that
\begin{equation}
        \sum_{j = 1}^n f_j(p) \, g_j(p) = 1
\end{equation}
for every $p \in U$, which would give a contradiction as before.

\section{Mapping properties}
\label{mapping properties}
\setcounter{equation}{0}

        Let $X$ be a topological space, and let us use $\Hom(X)$ to
denote the set of nonzero homomorphisms from $C_b(X)$ into ${\bf R}$
or ${\bf C}$, as appropriate.  In situations in which other types of
algebras are considered as well, this may be denoted more precisely as
$\Hom(C_b(X))$, to avoid confusion.  As in Section \ref{bounded
continuous functions}, $\Hom(X)$ is a compact subset of $C_b(X)^*$
with respect to the weak$^*$ topology.

        If $p \in X$, then $\phi_p(f) = f(p)$ is an element of
$\Hom(X)$, and we let $\Hom_1(X)$ be the subset of $\Hom(X)$
consisting of homomorphisms on $C_b(X)$ of this form.  Thus $\Hom(X) =
\Hom_1(X)$ when $X$ is compact, as in Section \ref{compact spaces}.
Otherwise, $\Hom_1(X)$ is dense in $\Hom(X)$ with respect to the
weak$^*$ topology on $C_b(X)^*$ for any $X$, as in the previous
section.

        We have also seen in Section \ref{bounded continuous
functions} that $p \mapsto \phi_p$ is continuous as a mapping from $X$
into $C_b(X)^*$ with the weak$^*$ topology.  By definition, this
mapping sends $X$ onto $\Hom_1(X)$ in $C_b(X)^*$.  If $X$ is compact,
then it follows that $\Hom_1(X)$ is compact with respect to the
weak$^*$ topology on $C_b(X)^*$, and hence closed.  This gives another
way to show that $\Hom_1(X) = \Hom(X)$ when $X$ is compact, since
$\Hom(X)$ is the same as the closure of $\Hom_1(X)$ with respect to
the weak$^*$ topology on $C_b(X)^*$ for any $X$.

        Note that $p \mapsto \phi_p$ is a one-to-one mapping of $X$
into $C_b(X)^*$ exactly when continuous functions separate points on
$X$.  If $X$ is completely regular, then it is easy to see that $p
\mapsto \phi_p$ is a homeomorphism from $X$ onto $\Hom_1(X)$ with
respect to the topology on $\Hom_1(X)$ induced by the weak$^*$
topology on $C_b(X)^*$.  In particular, if $X$ is compact and
Hausdorff, then $p \mapsto \phi_p$ is a homeomorphism from $X$ onto
$\Hom(X)$ with respect to the topology on $\Hom(X)$ induced by the
weak$^*$ topology on $C_b(X)^*$.  Remember that compact Hausdorff
topological spaces are normal and hence completely regular.

        Now let $Y$ be another topological space, and let $\rho$ be a
continuous mapping from $X$ into $Y$.  This leads to a linear mapping
$T_\rho : C_b(Y) \to C_b(X)$, defined by
\begin{equation}
\label{T_rho(f) = f circ rho}
        T_\rho(f) = f \circ \rho
\end{equation}
for each $f \in C_b(Y)$.  Observe that
\begin{equation}
\label{||T_rho(f)||_{sup, X} le ||f||_{sup, Y}}
        \|T_\rho(f)\|_{sup, X} \le \|f\|_{sup, Y}
\end{equation}
for every $f \in C_b(Y)$, where the subscripts $X$, $Y$ indicate on
which space the supremum norm is taken.  This shows that $T_\rho$ is a
bounded linear mapping from $C_b(Y)$ into $C_b(X)$ with respect to the
supremum norm, with operator norm less than or equal to $1$, and the
operator norm is actually equal to $1$, because $T_\rho({\bf 1}_Y) =
{\bf 1}_X$.  If $\rho(X)$ is dense in $Y$, then $T_\rho$ is an
isometric embedding of $C_b(Y)$ into $C_b(X)$ with respect to their
supremum norms.

        Let $T_\rho^* : C_b(X)^* \to C_b(Y)^*$ be the dual mapping
associated to $T_\rho$.  This sends a bounded linear functional
$\lambda$ on $C_b(X)$ to the bounded linear functional $\mu =
T_\rho^*(\lambda)$ defined by
\begin{equation}
        \mu(f) = \lambda(T_\rho(f)) = \lambda(f \circ \rho)
\end{equation}
for each $f \in C_b(Y)$.  The fact that $\mu = T_\rho^*(\lambda)$ is a
bounded linear functional on $C_b(Y)$ uses the fact that $T_\rho$ is a
bounded linear mapping from $C_b(Y)$ into $C_b(X)$, as well as the
boundedness of $\lambda$ on $C_b(X)$.  Similarly, it is easy to see
that $T_\rho^*$ is bounded as a linear mapping from $C_b(X)^*$ into
$C_b(Y)^*$ with respect to the corresponding dual norms.  It is also
easy to see that $T_\rho^*$ is continuous as a mapping from $C_b(X)^*$
into $C_b(Y)^*$ with respect to their corresponding weak$^*$
topologies.

        Observe that $T_\rho$ is a homomorphism from $C_b(Y)$ into
$C_b(X)$, in the sense that
\begin{equation}
\label{T_rho(f g) = T_rho(f) T_rho(g}
        T_\rho(f \, g) = T_\rho(f) \, T_\rho(g)
\end{equation}
for every $f, g \in C_b(Y)$.  If $\lambda$ is a homomorphism from
$C_b(X)$ into the real or complex numbers, as appropriate, then it
follows that $T_\rho^*(\lambda)$ is a homomorphism on $C_b(Y)$ too.
If $\lambda$ is a nonzero homomorphism on $C_b(X)$, so that
$\lambda({\bf 1}_X) = 1$, then $T_\rho^*(\lambda)$ is nonzero on
$C_b(Y)$ too, because
\begin{equation}
        T_\rho^*(\lambda)({\bf 1}_Y) = \lambda(T_\rho({\bf 1}_Y)
           = \lambda({\bf 1}_Y \circ \rho) = \lambda({\bf 1}_X) = 1.
\end{equation}
Thus $T_\rho^*(\Hom(X)) \subseteq \Hom(Y)$.

        If $q \in Y$, then let $\psi_q(f) = f(q)$ be the corresponding
point evaluation on $C_b(Y)$.  Observe that
\begin{equation}
        T_\rho^*(\phi_p) = \psi_{\rho(p)}
\end{equation}
for each $p \in X$, since
\begin{equation}
        T_\rho^*(\phi_p)(f) = \phi_p(T_\rho(f)) = \phi_p(f \circ \rho) 
                             = f(\rho(p)) = \psi_{\rho(p)}(f)
\end{equation}
for every $f \in C_b(Y)$.  Thus $T_\rho^*(\Hom_1(X)) \subseteq
\Hom_1(Y)$.  If $\rho(X)$ is dense in $Y$, then it follows that
$T_\rho^*(\Hom_1(X))$ is dense in $\Hom_1(Y)$ with respect to the
weak$^*$ topology on $C_b(Y)^*$, because $q \mapsto \psi_q$ is a
continuous mapping from $Y$ into $C_b(Y)^*$ with respect to the
weak$^*$ topology on $C_b(Y)^*$.  This implies that
$T_\rho^*(\Hom_1(X))$ is dense in $\Hom(Y)$ with respect to the
weak$^*$ topology on $C_b(Y)^*$ when $\rho(X)$ is dense in $Y$, since
$\Hom_1(Y)$ is dense in $\Hom(Y)$ with respect to the weak$^*$
topology on $C_b(Y)^*$.

        If $\rho(X)$ is dense in $Y$, then we also get that
\begin{equation}
\label{T_rho^*(Hom(X)) = Hom(Y)}
        T_\rho^*(\Hom(X)) = \Hom(Y).
\end{equation}
Remember that $\Hom(X)$ is compact in $C_b(X)^*$ with respect to the
weak$^*$ topology, which implies that $T_\rho^*(\Hom(X))$ is compact
in $C_b(Y)^*$ with respect to the weak$^*$ topology, because
$T_\rho^*$ is a continuous mapping from $C_b(X)^*$ into $C_b(Y)^*$
with respect to their weak$^*$ topologies.  Hence $T_\rho^*(\Hom(X))$
is a closed set in $C_b(Y)^*$ with respect to the weak$^*$ topology.
This implies that $\Hom(Y)$ is contained in $T_\rho^*(\Hom(X))$,
because $T_\rho^*(\Hom_1(X)) \subseteq T_\rho^*(\Hom(X))$ is dense in
$\Hom(Y)$ with respect to the weak$^*$ topology on $C_b(Y)^*$ when
$\rho(X)$ is dense in $Y$, as in the previous paragraph.  Therefore
(\ref{T_rho^*(Hom(X)) = Hom(Y)}) holds, since $T_\rho^*(\Hom(X))$ is
contained in $\Hom(Y)$ automatically.

        Suppose now that $Y$ is compact and Hausdorff, so that
$\Hom_1(Y) = \Hom(Y)$, and $q \mapsto \psi_q$ defines a homeomorphism
from $Y$ onto $\Hom(Y)$ with respect to the topology on $\Hom(Y)$
induced by the weak$^*$ topology on $C_b(Y)^*$.  In this case, the
restriction of $T_\rho^*$ to $\Hom(X)$ can be identified with a
mapping into $Y$.  If $\rho(X)$ is dense in $Y$, then we get a mapping
from $\Hom(X)$ onto $Y$, as in the previous paragraph.  If $X$ is
completely regular, so that $p \mapsto \phi_p$ defines a homeomorphism
from $X$ onto $\Hom_1(X)$ with respect to the topology induced on
$\Hom_1(X)$ by the weak$^*$ topology on $C_b(X)^*$, then the
restriction of $T_\rho^*$ to $\Hom(X)$ is basically an extension of
$\rho$.  If $X$ is compact and Hausdorff, then the restriction of
$\rho$ to $\Hom(X)$ is essentially the same as $\rho$ itself.

\section{Discrete sets}
\label{discrete sets}
\setcounter{equation}{0}

        Let $X$ be a nonempty set, and let $\beta X$ be the set of all
untrafilters on $X$.  As in Sections \ref{homomorphisms} and
\ref{homomorphisms, continued}, there is a natural one-to-one
correspondence between $\beta X$ and the set of all nonzero
homomorphisms on $\ell^\infty(X)$.  If $X$ is equipped with the
discrete topology, then $\ell^\infty(X)$ is the same as $C_b(X)$, and
the set of nonzero homomorphisms on $\ell^\infty(X)$ is the same as
the set $\Hom(X)$ discussed in the previous section.  In this section,
we shall see how properties of $\Hom(X)$ can be described more
directly in terms of ultrafilters on $X$.

        If $A \subseteq X$, then let $\widehat{A} \subseteq \beta X$
be the set of ultrafilters $\mathcal{F}$ on $X$ such that $A \in
\mathcal{F}$.  Thus $\widehat{X} = \beta X$, and there is a natural
one-to-one correspondence between $\widehat{A}$ and $\beta A$ for any
$A$, in which an ultrafilter on $A$ is extended to an ultrafilter on
$X$ that contains $A$ as an element, as in Section \ref{filters,
subsets}.  It is easy to see that
\begin{equation}
        \widehat{A \cap B} = \widehat{A} \cap \widehat{B}
\end{equation}
for every $A, B \subseteq X$.  Moreover,
\begin{equation}
        \widehat{X \backslash A} = \widehat{X} \backslash \widehat{A}
                                  = \beta X \backslash \widehat{A}
\end{equation}
for every $A \subseteq X$, because any ultrafilter $\mathcal{F}$ on
$X$ contains exactly one of $A$ and $X \backslash A$ as an element.
It follows that
\begin{equation}
        \widehat{A \cup B} = \widehat{A} \cup \widehat{B}
\end{equation}
for every $A, B \subseteq X$.

        Let us define a topology on $\beta X$ by saying that a subset
of $\beta X$ is an open set if it can be expressed as a union of
subsets of the form $\widehat{A}$, $A \subseteq X$.  Equivalently,
$\widehat{A}$ is an open set in $\beta X$ for each $A \subseteq X$,
and these open subsets of $\beta X$ form a base for the topology of
$\beta X$.  It is easy to see that the intersection of two open
subsets of $\beta X$ is also open, so that this does define a topology
on $\beta X$, because of the fact about intersections mentioned in the
previous paragraph.  The fact about complements mentioned in the
previous paragraph implies that $\widehat{A}$ is both open and closed
for every $A \subseteq X$.

        If $\mathcal{F}$ is an ultrafilter on $X$, then let
$L_\mathcal{F}$ be the corresponding homomorphism on $\ell^\infty(X)$,
as in Section \ref{homomorphisms}.  Let $A$ be a subset of $X$, and
let ${\bf 1}_A$ be the indicator function on $X$ corresponding to $A$,
so that ${\bf 1}_A(x) = 1$ when $x \in A$ and ${\bf 1}_A(x) = 0$ when
$x \in X \backslash A$.  It is easy to check that
\begin{eqnarray}
\label{L_mathcal{F}({bf 1}_A) = ...}
 L_\mathcal{F}({\bf 1}_A) & = & 1 \hbox{ when } A \in \mathcal{F} \\
               & = & 0 \hbox{ when } X \backslash A \in \mathcal{F}, \nonumber
\end{eqnarray}
directly from the definition of $L_\mathcal{F}$.  Remember that
$\mathcal{F} \mapsto L_\mathcal{F}$ defines a one-to-one
correspondence between $\beta X$ and the set $\Hom(X)$ of nonzero
homomorphisms on $\ell^\infty(X) = C_b(X)$.  Using
(\ref{L_mathcal{F}({bf 1}_A) = ...}), one can check that $\widehat{A}$
corresponds to a relatively open subset of $\Hom(X)$ with respect to
the weak$^*$ topology on $\ell^\infty(X)^*$ for each $A \subseteq X$.
This implies that every open set in $\beta X$ with respect to the
topology described earlier corresponds to a relatively open set in
$\Hom(X)$ with respect to the weak$^*$ topology on $\ell^\infty(X)$.
Conversely, one can show that relatively open subsets of $\Hom(X)$
with respect to the weak$^*$ topology on $\ell^\infty(X)^*$ correspond
to open subsets of $\beta X$.  This uses the facts that finite linear
combinations of indicator functions of subsets of $X$ are dense in
$\ell^\infty(X)$, and that homomorphisms on $\ell^\infty(X)$ have
bounded dual norm.

        In particular, $\beta X$ should be compact and Hausdorff with
respect to the topology defined before, because of the corresponding
properties of $\Hom(X)$ with respect to the topology induced by
the weak$^*$ topology on $\ell^\infty(X)^*$.  Let us check these
properties directly from the definition of the topology on $\beta X$.
If $\mathcal{F}$, $\mathcal{F}'$ are distinct ultrafilters on $X$,
then there is a set $A \subseteq X$ such that $A \in \mathcal{F}$ and
$X \backslash A \in \widehat{F}'$.  Hence $\mathcal{F} \in
\widehat{A}$ and $\mathcal{F}' \in \widehat{X \backslash A}$, so that
$\mathcal{F}$, $\mathcal{F}'$ are contained in disjoint open subsets
of $\beta X$, which implies that $\beta X$ is Hausdorff.

        To show that $\beta X$ is compact, let $\mathcal{U}$ be an
arbitrary ultrafilter on $\beta X$, and let us show that $\mathcal{U}$
converges to an element of $\beta X$.  Let $\mathcal{F}$ be the
collection of subsets $A$ of $X$ such that $\widehat{A} \in
\mathcal{U}$.  It is easy to see that $\mathcal{F}$ is a filter on
$X$, because $\mathcal{U}$ is a filter on $\beta X$.  If $A \subseteq
X$, then either $\widehat{A}$ or $\widehat{X \backslash A} = \beta X
\backslash \widehat{A}$ is an element of $\mathcal{U}$, because
$\mathcal{U}$ is an ultrafilter on $\beta X$.  This implies that
either $A$ or $X \backslash A$ is an element of $\mathcal{F}$ for
every $A \subseteq X$, and hence that $\mathcal{F}$ is an ultrafilter
on $X$.  It remains to check that $\mathcal{U}$ converges to
$\mathcal{F}$ as an element of $\beta X$.  By definition, this means
that every neighborhood of $\mathcal{F}$ in $\beta X$ should be an
element of $\mathcal{U}$.  Because the sets $\widehat{A}$, $A
\subseteq X$, form a base for the topology of $\beta X$, it suffices
to have $\widehat{A} \in \mathcal{U}$ for every $A \subseteq X$ such
that $A \in \mathcal{F}$, which follows from the definition of
$\mathcal{F}$.

        If $p \in X$, then the collection $\mathcal{F}_p$ of $A
\subseteq X$ with $p \in A$ is an ultrafilter on $X$.  Thus $p \mapsto
\mathcal{F}_p$ defines a natural embedding of $X$ into $\beta X$.  It
is easy to see that the set of ultrafilters $\mathcal{F}_p$, $p \in
X$, is dense in $\beta X$ with respect to the topology defined
earlier.  One can also check that the homomorphism $L_{\mathcal{F}_p}$
on $\ell^\infty(X)$ corresponding to $\mathcal{F}_p$ is the same as
evaluation at $p$.

        Let $Y$ be a compact Hausdorff topological space, and let
$\rho$ be a mapping from $X$ into $Y$.  If $\mathcal{F}$ is an
ultrafilter on $X$, then we can define $\rho_*(\mathcal{F})$ as usual
as the collection of sets $E \subseteq Y$ such that $\rho^{-1}(E) \in
\mathcal{F}$.  In particular, we have seen that $\rho_*(\mathcal{F})$
is an ultrafilter on $Y$.  It follows that $\rho_*(\mathcal{F})$
converges to a unique element of $Y$, because $Y$ is compact and
Hausdorff.  Let $\widehat{\rho}(\mathcal{F})$ be the limit of
$\rho_*(\mathcal{F})$ in $Y$, which defines $\widehat{\rho}$ as a
mapping from $\beta X$ into $Y$.  If $p \in X$, then it is easy to see
that $\widehat{\rho}(\mathcal{F}_p) = \rho(p)$.  Thus $\widehat{\rho}$
is basically an extension of $\rho$ to a mapping from $\beta X$ into
$Y$.

        Let us check that $\widehat{\rho}$ is continuous as a mapping
from $\beta X$ into $Y$.  Let $\mathcal{F}$ be an ultrafilter on $X$,
and let $W$ be an open set in $Y$ that contains
$\widehat{\rho}(\mathcal{F})$ as an element.  Because $Y$ is compact
and Hausdorff, it is regular, which implies that there is an open set
$V$ in $Y$ such that $\widehat{\rho}(\mathcal{F}) \in V$ and the
closure $\overline{V}$ of $V$ in $Y$ is contained in $W$.  Remember
that $\rho_*(\mathcal{F})$ converges to $\widehat{\rho}(\mathcal{F})$
in $Y$, which implies that $V \in \rho_*(\mathcal{F})$.  This implies
in turn that $\rho^{-1}(V) \in \mathcal{F}$, by the definition of
$\rho_*(\mathcal{F})$.  Put $A = \rho^{-1}(V)$, so that $\widehat{A}$
is an open set in $\beta X$ that contains $\mathcal{F}$ as an element.
Let $\mathcal{F}'$ be any other ultrafilter on $X$ that is an element
of $\widehat{A}$.  This means that $\rho^{-1}(V) = A \in
\mathcal{F}'$, and hence that $A \in \rho_*(\mathcal{F}')$.  By
construction, $\rho_*(\mathcal{F}')$ converges to
$\widehat{\rho}(\mathcal{F}')$ in $Y$, which implies that
$\widehat{\rho}(\mathcal{F}') \in \overline{V}$.  This shows that
$\widehat{\rho}(\mathcal{F}') \in \overline{V} \subseteq W$ for every
$\mathcal{F}' \in \widehat{A}$, and hence that $\widehat{\rho}$ is
continuous at $\mathcal{F}$ for every $\mathcal{F} \in \beta X$, as
desired.

\section{Locally compact spaces, revisited}
\label{locally compact spaces, revisited}
\setcounter{equation}{0}

        Let $X$ be a locally compact Hausdorff topological space which
is not compact, and let $X^*$ be the one-point compactification of
$X$, as in Section \ref{locally compact spaces, continued}.  Also let
$C_{lim}(X)$ be the space of continuous real or complex-valued
functions on $X$ which have a limit at infinity, as in Section
\ref{locally compact spaces, continued}.  As usual, this may also be
denoted $C_{lim}(X, {\bf R})$ or $C_{lim}(X, {\bf C})$, to indicate
whether real or complex-valued functions are being used.  As in
Section \ref{locally compact spaces, continued}, $C_{lim}(X)$ is a
closed subalgebra of the algebra $C_b(X)$ of bounded continuous
functions on $X$ with respect to the supremum norm, and $C_{lim}(X)$
is the same as the linear span in $C_b(X)$ of the subspace $C_0(X)$ of
functions that vanish at infinity on $X$ and the constant functions on
$X$.  Equivalently, $C_{lim}(X)$ is the same as the space of
continuous functions on $X$ that have a continuous extension to $X^*$.

        Thus a nonzero homomorphism $\phi$ from $C_{lim}(X)$ into the
real or complex numbers, as appropriate, is basically the same as a
nonzero homomorphism on $C(X^*)$.  As in Section \ref{compact spaces},
every nonzero homomorphism on $C(X^*)$ can be represented by
evaluation at a point in $X^*$, because $X^*$ is compact.  This point
in $X^*$ is either an element of $X$, or the point at infinity in
$X^*$.  This implies that either there is a $p \in X$ such that
\begin{equation}
        \phi(f) = f(p)
\end{equation}
for every $f \in C_{lim}(X)$, or that
\begin{equation}
        f(x) \to \phi(f) \hbox{ as } x \to \infty
\end{equation}
for every $f \in C_{lim}(X)$.

        Suppose now that $\phi$ is a nonzero homomorphism on $C_b(X)$.
The restriction of $\phi$ to $C_{lim}(X)$ is a nonzero homomorphism on
$C_{\lim}(X)$, since $\phi({\bf 1}_X) = 1$.  If $\phi(f) = f(p)$ for
some $p \in X$ and every $f \in C_{lim}(X)$, then we would like to
check that this also holds for every $f \in C_b(X)$.  To see this, we
can use Urysohn's lemma to get a continuous function $\theta$ with
compact support on $X$ such that $\theta(p) = 1$.  Let $f$ be a
bounded continuous function on $X$, and observe that $\theta \, f \in
C_{lim}(X)$, because it has compact support on $X$.  This implies that
\begin{equation}
        \phi(\theta \, f) = (\theta \, f)(p) = \theta(p) \, f(p) = f(p),
\end{equation}
since $\theta(p) = 1$.  Similarly,
\begin{equation}
        \phi((1 - \theta) \, f) = \phi(1 - \theta) \, \phi(f)
                                = (1 - \theta(p)) \, \phi(f) = 0.
\end{equation}
More precisely, this uses the hypothesis that $\phi$ is a homomorphism
on $C_b(X)$ in the first step, and then the fact that $1 - \theta \in
C_{lim}(X)$ to get that $\phi(1 - \theta)$ is equal to $1 -
\theta(p)$.  Combining these two equations, we get that $\phi(f) =
f(p)$, as desired.

        Let $\rho$ be the standard embedding of $X$ into $X^*$, which
sends each $p \in X$ to itself as an element of $X^*$.  As in Section
\ref{mapping properties}, this leads to a mapping $T_\rho$ from
$C(X^*)$ into $C_b(X)$, which sends $C(X^*)$ onto $C_{lim}(X)$ in this
case.  The corresponding dual mapping $T_\rho^*$ sends the set
$\Hom(X)$ of nonzero homomorphisms on $C_b(X)$ into the analogous set
$\Hom(X^*)$ for $X^*$, which can be identified with $X^*$, because
$X^*$ is compact and Hausdorff.  Remember that $\Hom_1(X) \subseteq
\Hom(X)$ is the set of homomorphisms on $C_b(X)$ defined by evaluation
at elements of $X$, and that $T_\rho^*$ maps $\Hom_1(X)$ to the point
evaluations on $C(X^*)$ that correspond to elements of $X$.  The
discussion in the previous paragraph implies that $T_\rho^*$ sends
every other element of $\Hom(X)$ to the point evaluation on $C(X^*)$
that corresponds to the point at infinity in $X^*$.

\section{Mapping properties, continued}
\label{mapping properties, continued}
\setcounter{equation}{0}

        Let $U$ be the open unit disk in the complex plane, so that
$\overline{U}$ is the closed unit disk.  Also let $\rho$ be the
standard embedding of $U$ into $\overline{U}$, which sends each $z \in
U$ to itself as an element of $\overline{U}$.  This leads to a mapping
$T_\rho$ from $C(\overline{U})$ into $C_b(U)$, as in Section
\ref{mapping properties}, which sends a continuous function $f$ on
$\overline{U}$ to its restriction to $U$.  The dual mapping $T_\rho^*
: C_b(U)^* \to C(\overline{U})^*$ sends the set $\Hom(U)$ of
nonzero homomorphisms on $C_b(U)$ into the analogous set
$\Hom(\overline{U})$ for $\overline{U}$, as before.  If $\phi$
is a nonzero homomorphism on $C_b(U)$, then $T_\rho^*(\phi)$ is
basically the same as the restriction of $\phi$ to $C(\overline{U})$,
which is identified with a subalgebra of $C_b(U)$.  Each nonzero
homomorphism on $C(\overline{U})$ can be represented as a point
evaluation, as in Section \ref{compact spaces}.  If there is a $p \in
U$ such that $\phi(f) = f(p)$ for every $f \in C(\overline{U})$, then
the same relation holds for every $f \in C_b(U)$, as in the previous
section.  If $\phi \in \Hom(U)$ does not correspond to
evaluation at a point in $U$, then it follows that the restriction of
$\phi$ to $C(\overline{U})$ corresponds to evaluation at a point in
$\partial U$.

        Let $\mathcal{A}$ be the algebra of continuous complex-valued
functions on $\overline{U}$ that are holomorphic on $U$, as in Section
\ref{disk algebra}, and let $\mathcal{B}$ be the algebra of bounded
holomorphic functions on $U$, as in Section \ref{bounded holomorphic
functions}.  If $f \in \mathcal{A}$, then the restriction of $f$ to
$U$ is an element of $\mathcal{B}$, and $f$ is determined on
$\overline{U}$ by its restriction to $U$, by continuity.  Thus we can
identify $\mathcal{A}$ with a subalgebra of $\mathcal{B}$.

        Let $\Hom(\mathcal{A})$, $\Hom(\mathcal{B})$ denote the sets
of nonzero homomorphisms from $\mathcal{A}$, $\mathcal{B}$ into the
complex numbers, respectively.  As in Sections \ref{disk algebra} and
\ref{bounded holomorphic functions}, these are subsets of the duals of
$\mathcal{A}$, $\mathcal{B}$, and we are especially interested in the
topologies induced on $\Hom(\mathcal{A})$, $\Hom(\mathcal{B})$ by the
weak$^*$ topologies on the corresponding dual spaces.

        If $p \in \overline{U}$, then $\phi_p(f) = f(p)$ defines a
homomorphism on $\mathcal{A}$, and we have seen in Section \ref{disk
algebra} that every nonzero homomorphism on $\mathcal{A}$ is of this
form.  Of course, $\phi_p(f) = f(p)$ is a continuous function on
$\overline{U}$ for every $f \in \mathcal{A}$, by definition of
$\mathcal{A}$, which implies that $p \mapsto \phi_p$ is continuous as
a mapping from $\overline{U}$ into $\Hom(\mathcal{A})$ with respect to
the weak$^*$ topology on $\mathcal{A}$.  If $f_1(z)$ is the element of
$\mathcal{A}$ defined by $f_1(z) = z$ for each $z \in \overline{U}$,
then $\phi_p(f_1) = p$ for each $p \in \overline{U}$.  This shows that
$p \mapsto \phi_p$ is actually a homeomorphism from $\overline{U}$
onto $\Hom(\mathcal{A})$ with respect to the topology induced on
$\Hom(\mathcal{A})$ by the weak$^*$ topology on $\mathcal{A}^*$.

        Similarly, if $p \in U$, then $\phi_p(f) = f(p)$ defines a
nonzero homomorphism on $\mathcal{B}$, and $p \mapsto \phi_p$ defines
a continuous mapping from $U$ into $\Hom(\mathcal{B})$ with respect to
the weak$^*$ topology on $\mathcal{B}^*$.  Let $\Hom_1(\mathcal{B})$
be the set of homomorphisms on $\mathcal{B}$ of this form.  If
$f_1(z)$ is the element of $\mathcal{B}$ defined by $f_1(z) = z$ for
each $z \in U$, then $\phi_p(f_1) = p$ for each $p \in U$.  This
implies that $p \mapsto \phi_p$ is a homeomorphism from $U$ onto
$\Hom_1(\mathcal{B})$ with respect to the topology induced on
$\Hom_1(\mathcal{B})$ by the weak$^*$ topology on $\mathcal{B}^*$.

        If $\phi$ is a nonzero homomorphism on $\mathcal{B}$, then the
restriction of $\phi$ to $\mathcal{A}$ is a nonzero homomorphism on
$\mathcal{A}$.  This defines a natural mapping from
$\Hom(\mathcal{B})$ into $\Hom(\mathcal{A})$.  It is easy to see that
this mapping is continuous with respect to the topologies induced on
$\Hom(\mathcal{A})$, $\Hom(\mathcal{B})$ by the weak$^*$ topologies on
$\mathcal{A}^*$, $\mathcal{B}^*$, respectively.

        Let $f_1$ be the element of $\mathcal{B}$ defined by $f_1(z) =
z$ for each $z \in U$ again.  Also let $\phi$ be a nonzero
homomorphism on $\mathcal{B}$, and put $p = \phi(f_1)$.  Note that $p
\in \overline{U}$, since $\phi$ has dual norm equal to $1$ with
respect to the supremum norm on $\mathcal{B}$, as in Section
\ref{bounded holomorphic functions}.  If $f \in \mathcal{A}$, then
$\phi(f) = f(p)$, by the arguments in Section \ref{disk algebra}
applied to the restriction of $\phi$ to $\mathcal{A}$.

        Suppose that $p \in U$, and let us check that $\phi(f) = f(p)$
for every $f \in \mathcal{B}$.  Any holomorphic function $f$ on $U$
can be expressed as
\begin{equation}
        f(z) = f(p) + (z - p) \, g(z)
\end{equation}
for some holomorphic function $g$ on $U$, and $g$ is also bounded on
$U$ when $f$ is.  This implies that $\phi(f) = f(p)$ for every $f \in
\mathcal{B}$, because $\phi$ applied to $z - p$ is equal to $0$, by
definition of $p$.  If $\phi \in \Hom(\mathcal{B}) \backslash
\Hom_1(\mathcal{B})$, then it follows that $p \in \partial U$.  This
is analogous to the situation for bounded continuous functions on $U$
mentioned at the beginning of the section.

        Let us take $C_b(U)$ to be the algebra of bounded continuous
complex-valued functions on $U$, so that $\mathcal{B}$ is a subalgebra
of $C_b(U)$.  Let us also use $\Hom(C_b(U))$ to denote the set of
nonzero homomorphisms on $C_b(U)$, to be more consistent with the
notation for $\mathcal{B}$.  If $\phi$ is a nonzero homomorphism on
$C_b(U)$, then the restriction of $\phi$ to $\mathcal{B}$ is a nonzero
homomorphism on $\mathcal{B}$.  This defines a natural mapping $R$
from $\Hom(C_b(U))$ into $\Hom(\mathcal{B})$, which is easily seen to
be continuous with respect to the topologies induced by the weak$^*$
topologies on $C_b(U)^*$ and $\mathcal{B}^*$, respectively.

        By construction, $R$ sends $\Hom_1(C_b(U))$ onto
$\Hom_1(\mathcal{B})$.  We also know that $\Hom(C_b(U))$ is compact
with respect to the weak$^*$ topology on $C_b(U)^*$, which implies
that $R(\Hom(C_b(U)))$ is compact with respect to the weak$^*$
topology on $\mathcal{B}^*$.  In particular, $R(\Hom(C_b(U)))$ is
closed with respect to the weak$^*$ topology on $\mathcal{B}^*$.  As
in Section \ref{density}, Carleson's corona theorem states that
$\Hom_1(\mathcal{B})$ is dense in $\Hom(\mathcal{B})$ with respect to
the weak$^*$ topology on $\mathcal{B}^*$.  It follows that $R$ maps
$\Hom(C_b(U))$ onto $\Hom(\mathcal{B})$, so that every nonzero
homomorphism on $\mathcal{B}$ is the restriction to $\mathcal{B}$ of a
nonzero homomorphism on $C_b(U)$.

\section{Banach algebras}
\label{banach algebras}
\setcounter{equation}{0}

        A vector space $\mathcal{A}$ over the real or complex numbers
is said to be an (associative) algebra if every $a, b \in \mathcal{A}$
has a well-defined product $a \, b \in \mathcal{A}$ which is linear in
$a$ and $b$ separately and satisfies the associative law
\begin{equation}
 (a \, b) \, c = a \, (b \, c) \quad\hbox{for every } a, b, c \in \mathcal{A}.
\end{equation}
We shall be primarily concerned here with commutative algebras, so that
\begin{equation}
        a \, b = b \, a
\end{equation}
for each $a, b \in \mathcal{A}$.  We also ask that there be a nonzero
multiplicative identity element $e$ in $\mathcal{A}$, which means that
$e \ne 0$ and
\begin{equation}
        e \, a = a \, e = a
\end{equation}
for every $a \in \mathcal{A}$.  We have seen several examples of
algebras of functions in the previous sections, for which the
multiplicative identity element is the constant function equal $1$.

        Suppose that $\mathcal{A}$ is equipped with a norm $\|a\|$.
This norm should also be compatible with multiplication on
$\mathcal{A}$, in the sense that $\|e\| = 1$ and
\begin{equation}
\label{||a b|| le ||a|| ||b||}
        \|a \, b\| \le \|a\| \, \|b\|
\end{equation}
for every $a, b \in \mathcal{A}$.  We say that $\mathcal{A}$ is a
Banach algebra if it is also complete as a metric space with respect
to the metric $d(a, b) = \|a - b\|$ associated to the norm.  The
algebra of bounded continuous functions on any topological space is a
Banach algebra with respect to the supremum norm.  Closed subalgebras
of Banach algebras are also Banach algebras, such as the disk algebra
and the algebra of bounded holomorphic functions on the unit disk.

        Suppose that $\mathcal{A}$ is any Banach algebra, and let $a$
be an element of $\mathcal{A}$.  If $n$ is a positive integer, then
$a^n$ is the product $a \, a \, \cdots a$ of $n$ $a$'s in
$\mathcal{A}$, which can also be described by $a^n = a$ when $n = 1$,
and $a^{n + 1} = a \, a^n$ for every $n$.  This is interpreted as
being equal to the multiplicative identity element $e$ when $n = 0$.
Observe that
\begin{equation}
\label{||a^n|| le ||a||^n}
        \|a^n\| \le \|a\|^n
\end{equation}
for each $n \ge 0$, where again the right side is interpreted as being
equal to $1$ when $n = 0$.

        An element $a$ of $\mathcal{A}$ is said to be \emph{invertible}
if there is another element $a^{-1}$ of $\mathcal{A}$ such that
\begin{equation}
        a \, a^{-1} = a^{-1} \, a = e.
\end{equation}
It is easy to see that the inverse $a^{-1}$ of $a$ is unique when it
exists.  If $a$, $b$ are invertible elements of $\mathcal{A}$, then
their product $a \, b$ is also invertible, with
\begin{equation}
        (a \, b)^{-1} = b^{-1} \, a^{-1}.
\end{equation}
If $x$ is an invertible element of $\mathcal{A}$ and $y$ is another
element of $\mathcal{A}$ that commutes with $x$, so that $x \, y = y
\, x$, then $y$ also commutes with $x^{-1}$,
\begin{equation}
\label{y x^{-1} = x^{-1} y}
        y \, x^{-1} = x^{-1} \, y.
\end{equation}
If $a$, $b$ are commuting elements of $\mathcal{A}$ whose product $a
\, b$ is invertible, then $a$, $b$ are also invertible, with
\begin{equation}
        a^{-1} = b \, (a \, b)^{-1}, \quad b^{-1} = (a \, b)^{-1} \, a.
\end{equation}
This uses the fact that $a$, $b$ commute with $(a \, b)^{-1}$, since
they commute with $a \, b$.  Note that these statements do not involve
the norm on $\mathcal{A}$.

        If $a \in \mathcal{A}$ and $n$ is a positive integer, then
\begin{equation}
        (e - a) \, \Big(\sum_{j = 0}^n a^j\Big)
           = \Big(\sum_{j = 0}^n a^j\Big) \, (e - a) = e - a^{n + 1}.
\end{equation}
This is basically the same as for real or complex numbers.
If $\|a\| < 1$, then
\begin{equation}
        \lim_{n \to \infty} a^n = 0
\end{equation}
in $\mathcal{A}$, since $\|a^n\| \le \|a\|^n \to 0$ as $n \to \infty$.
Similarly,
\begin{equation}
        \sum_{j = 0}^\infty \|a^j\| \le \sum_{j = 0}^\infty \|a\|^j
                                     = \frac{1}{1 - \|a\|}.
\end{equation}
As in the context of real or complex numbers, the convergence of
$\sum_{j = 0}^\infty \|a^j\|$ means that $\sum_{j = 0}^\infty a^j$
converges absolutely.  More precisely, this implies that the partial
sums $\sum_{j = 0}^n a^j$ of $\sum_{j = 0}^\infty a^j$ form a Cauchy
sequence in $\mathcal{A}$, which converges when $\mathcal{A}$ is
complete.  It follows that
\begin{equation}
        (e - a) \, \Big(\sum_{j = 0}^\infty a^j\Big)
          = \Big(\sum_{j = 0}^\infty a^j\Big) \, (e - a) = e
\end{equation}
when $a \in \mathcal{A}$, $\|a\| < 1$, and $\mathcal{A}$ is a Banach
algebra.  Thus $e - a$ is invertible in $\mathcal{A}$ under these
conditions, with
\begin{equation}
        (e - a)^{-1} = \sum_{j = 0}^\infty a^j.
\end{equation}
We also get that
\begin{equation}
        \|(e - a)^{-1}\| \le \frac{1}{1 - \|a\|}.
\end{equation}
If $b$ is any invertible element of $\mathcal{A}$ and $\|a\| \,
\|b^{-1}\| < 1$, then $b - a$ is also invertible in $\mathcal{A}$,
because
\begin{equation}
        b - a = (e - a \, b^{-1}) \, b
\end{equation}
and $e - a \, b^{-1}$ is invertible by the previous argument.  This
shows that the invertible elements in a Banach algebra $\mathcal{A}$
form an open set in $\mathcal{A}$ with respect to the metric
associated to the norm.

        Let $\mathcal{A}$ be a real or complex algebra, and let $\phi$
be a linear functional on $\mathcal{A}$, which is to say a linear
mapping from $\mathcal{A}$ into the real or complex numbers, as
appropriate.  We say that $\phi$ is a \emph{homomorphism} on $\mathcal{A}$ if
\begin{equation}
        \phi(a \, b) = \phi(a) \, \phi(b)
\end{equation}
for every $a, b \in \mathcal{A}$.  Of course, $\phi$ satisfies this
condition trivially when $\phi(a) = 0$ for every $a \in \mathcal{A}$,
and we are primarily interested in the nonzero homomorphisms $\phi$,
which means that $\phi(a) \ne 0$ for some $a \in \mathcal{A}$.
This implies that
\begin{equation}
        \phi(e) = 1,
\end{equation}
because $\phi(a) = \phi(e) \, \phi(a)$, since $a = e \, a$.
If $b$ is any invertible element of $\mathcal{A}$, then we get that
\begin{equation}
        \phi(b) \, \phi(b^{-1}) = \phi(b \, b^{-1}) = \phi(e) = 1,
\end{equation}
and hence $\phi(b) \ne 0$.

        Suppose now that $\mathcal{A}$ is a Banach algebra again, and
let $\phi$ be a nonzero homomorphism on $\mathcal{A}$.  If $a \in
\mathcal{A}$ and $\|a\| < 1$, then $e - a$ is invertible, and so
\begin{equation}
        \phi(e - a) \ne 0,
\end{equation}
which means that $\phi(a) \ne 1$.  By the same argument, $\phi(t \, a)
= t \, \phi(a) \ne 1$ for every $t \in {\bf R}$ or ${\bf C}$, as
appropriate, such that $|t| < 1$.  This implies that $|\phi(a)| < 1$
when $a \in \mathcal{A}$ satisfies $\|a\| < 1$.  Hence
\begin{equation}
        |\phi(a)| \le \|a\|
\end{equation}
for every $a \in \mathcal{A}$, which shows that $\phi$ is a continuous
linear functional on $\mathcal{A}$ with dual norm less than or equal
to $1$.  The dual norm of $\phi$ is actually equal to $1$, because
$\phi(e) = 1$.  It is easy to see that the collection of nonzero
homomorphisms on $\mathcal{A}$ is closed with respect to the weak$^*$
topology on the dual of $\mathcal{A}$.  It follows that the collection
of nonzero homomorphisms on $\mathcal{A}$ is compact with respect to
the weak$^*$ topology, by the Banach--Alaoglu theorem.

        A linear subspace $\mathcal{I}$ of a real or complex algebra
$\mathcal{A}$ is said to be an \emph{ideal} in $\mathcal{A}$ if $a \,
x$ and $x \, a$ are contained in $\mathcal{I}$ for every $a \in
\mathcal{A}$ and $x \in \mathcal{I}$.  Of course, $\mathcal{A}$ itself
and the trivial subspace $\{0\}$ are ideals in $\mathcal{A}$, and an
ideal $\mathcal{I}$ in $\mathcal{A}$ is said to be proper if
$\mathcal{I} \ne \mathcal{A}$.  If $\mathcal{I}$ is an ideal in
$\mathcal{A}$ and $\mathcal{I}$ contains the identity element $e$, or
any invertible element $x$, then $\mathcal{I} = \mathcal{A}$.  If
$\mathcal{A}$ is a Banach algebra and $\mathcal{I}$ is an ideal in
$\mathcal{A}$, then it is easy to see that the closure
$\overline{\mathcal{I}}$ of $\mathcal{I}$ with respect to the norm on
$\mathcal{A}$ is also an ideal in $\mathcal{A}$.  If $\mathcal{I}$ is
a proper ideal in a Banach algebra $\mathcal{A}$, then $e \not\in
\overline{\mathcal{I}}$.  This is because elements of $\mathcal{A}$
sufficiently close to $e$ are invertible, as before.  Thus the closure
of a proper ideal in a Banach algebra is still proper.

        A proper ideal $\mathcal{I}$ in an algebra $\mathcal{A}$ is
said to be maximal if $\mathcal{A}$ and $\mathcal{I}$ are the only
ideals that contain $\mathcal{I}$.  It is easy to see that the kernel
of a nonzero homomorphism on $\mathcal{A}$ is maximal, since it has
codimension $1$.  A maximal ideal $\mathcal{I}$ in a Banach algebra
$\mathcal{A}$ is automatically closed, because its closure
$\overline{I}$ is a proper ideal that contains $\mathcal{I}$, and
hence is equal to $\mathcal{I}$.

        Using the axiom of choice, one can show that every proper
ideal in an algebra with nonzero multiplicative identity element is
contained in a maximal ideal.  More precisely, one can use Zorn's
lemma or the Hausdorff maximality principle, by checking that the
union of a chain of proper ideals is a proper ideal.  To get
properness, one uses the fact that the ideals do not contain the
identity element.

        If $\mathcal{A}$ is a commutative algebra and $a \in \mathcal{A}$, then
\begin{equation}
        \mathcal{I}_a = \{a \, b : b \in \mathcal{A}\}
\end{equation}
is an ideal in $\mathcal{A}$.  Moreover, $\mathcal{I}_a$ is a proper
ideal in $\mathcal{A}$ if and only if $a$ is not invertible in
$\mathcal{A}$.

        Suppose from now on that $\mathcal{A}$ is a complex Banach
algebra.  Let $a \in \mathcal{A}$ be given, and suppose that $t \, e -
a$ is invertible in $\mathcal{A}$ for every $t \in {\bf C}$.  If
$\lambda$ is a continuous linear functional on $\mathcal{A}$, then
one can show that
\begin{equation}
        f_\lambda(t) = \lambda((t \, e - a)^{-1})
\end{equation}
is a holomorphic function on the complex plane ${\bf C}$.  One can
also check that $(t \, e - a)^{-1} \to 0$ in $\mathcal{A}$ as $|t| \to
\infty$, so that $f_\lambda(t) \to 0$ as $|t| \to \infty$ for each
$\lambda$.  This implies that $f_\lambda(t) = 0$ for every $t \in {\bf
C}$ and continuous linear functional $\lambda$ on $\mathcal{A}$, by
standard results in complex analysis.  Using the Hahn--Banach theorem,
it follows that $(t \, e - a)^{-1} = 0$ for every $t \in {\bf C}$,
contradicting the fact that invertible elements of $\mathcal{A}$ are
not zero.  This is a brief sketch of the well-known fact that for each
$a \in \mathcal{A}$ there is a $t \in {\bf C}$ such that $t \, e - a$
is not invertible.

        Suppose that every nonzero element of $\mathcal{A}$ is
invertible.  If $a \in \mathcal{A}$, then there is a $t \in {\bf C}$
such that $t \, e - a$ is not invertible, as in the previous
paragraph.  In this case, it follows that $a = t \, e$, so that
$\mathcal{A}$ is isomorphically equivalent to the complex numbers.

        If $\mathcal{A}$ is an algebra $\mathcal{I}$ is an ideal in
$\mathcal{A}$, then the quotient $\mathcal{A} / \mathcal{I}$ defines
an algebra in a natural way, so that the corresponding quotient
mapping is a homomorphism from $\mathcal{A}$ onto $\mathcal{A} /
\mathcal{I}$ with kernel equal to $\mathcal{I}$.  If $\mathcal{A}$ has
a nonzero multiplicative identity element and $\mathcal{I}$ is proper,
then $\mathcal{A} / \mathcal{I}$ also has a nonzero multiplicative
identity element.  If $\mathcal{I}$ is a maximal ideal, then
$\mathcal{A} / \mathcal{I}$ contains no nontrivial proper ideals.  If
$\mathcal{A}$ is commutative and $\mathcal{I}$ is maximal, then it
follows that every nonzero element of $\mathcal{A} / \mathcal{I}$ is
invertible in the quotient.  If $\mathcal{A}$ is a Banach algebra and
$\mathcal{I}$ is a proper closed ideal in $\mathcal{A}$, then
$\mathcal{A} / \mathcal{I}$ is also a Banach algebra, with respect to
the usual quotient norm.  If $\mathcal{A}$ is a complex commutative
Banach algebra and $\mathcal{I}$ is a maximal ideal in $\mathcal{A}$,
then it follows that $\mathcal{A} / \mathcal{I}$ is isomorphic to the
complex numbers.  This implies that every maximal ideal in a
commutative complex Banach algebra $\mathcal{A}$ is the kernel of a
homomorphism from $\mathcal{A}$ onto the complex numbers.  If
$\mathcal{A}$ is a commutative complex Banach algebra and $a \in
\mathcal{A}$ is not invertible, then $a$ is contained in a maximal
ideal in $\mathcal{A}$, and hence there is a nonzero homomorphism
$\phi : \mathcal{A} \to {\bf C}$ such that $\phi(a) = 0$.

\section{Ideals and filters}
\label{ideals, filters}
\setcounter{equation}{0}

        Let $E$ be a nonempty set, and let $\mathcal{A}$ be the
algebra of all real or complex-valued functions on $E$.  Put
\begin{equation}
        Z(f) = \{x \in E : f(x) = 0\}
\end{equation}
for each $f \in \mathcal{A}$.  Thus
\begin{equation}
        Z(f) \cap Z(g) \subseteq Z(f + g)
\end{equation}
and
\begin{equation}
        Z(f) \cup Z(g) = Z(f \, g)
\end{equation}
for every $f, g \in \mathcal{A}$.  If $f$, $g$ are nonnegative
real-valued functions on $E$, then
\begin{equation}
        Z(f) \cap Z(g) = Z(f + g).
\end{equation}

        If $\mathcal{F}$ is a filter on $E$, then put
\begin{equation}
 \mathcal{I}(\mathcal{F}) = \{f \in \mathcal{A} : Z(f) \in \mathcal{F}\}.
\end{equation}
It is easy to see that this is an ideal in $\mathcal{A}$, using the
properties of the zero sets of sums and products of functions
mentioned in the previous paragraph.  More precisely,
$\mathcal{I}(\mathcal{F})$ is a proper ideal in $\mathcal{A}$, since
the elements of a filter are nonempty sets.  As a special case,
suppose that $A \subseteq E$ is not empty, and let $\mathcal{F}^A$ be
the collection of subsets $B$ of $E$ such that $A \subseteq B$.  
This is a filter on $E$, and the corresponding ideal
$\mathcal{I}(\mathcal{F}^A)$ is the same as
\begin{equation}
 \mathcal{I}_A = \{f \in \mathcal{A} : f(x) = 0 \hbox{ for every } x \in A\}.
\end{equation}
In this case, the quotient $\mathcal{A} / \mathcal{I}_A$ can be
identified with the algebra of real or complex-valued functions on
$A$, as appropriate.  In particular, if $A$ consists of a single
point, then the quotient is isomorphic to the real or complex numbers,
as appropriate.

        Conversely, if $\mathcal{I}$ is a proper ideal in
$\mathcal{A}$, then put
\begin{equation}
        \mathcal{F}(\mathcal{I}) = \{Z(f) : f \in \mathcal{I}\}.
\end{equation}
It is easy to check that this is a filter on $E$.  In connection with
this, note that
\begin{equation}
        Z(|f|) = Z(f)
\end{equation}
for every $f \in \mathcal{A}$, and that $|f| \in \mathcal{I}$ when $f
\in \mathcal{I}$.  This implies that $\mathcal{F}(\mathcal{I})$ is the
same as the collection of zero sets of nonnegative real-valued
functions on $E$ in $\mathcal{I}$.  Observe also that
\begin{equation}
        \mathcal{F}(\mathcal{I}(\mathcal{F})) = \mathcal{F}
\end{equation}
for every filter $\mathcal{F}$ on $E$, and that
\begin{equation}
        \mathcal{I}(\mathcal{F}(\mathcal{I})) = \mathcal{I}
\end{equation}
for every proper ideal $\mathcal{I}$ in $\mathcal{A}$.  This shows
that every proper ideal $\mathcal{I}$ in $\mathcal{A}$ is of the form
$\mathcal{I}(\mathcal{F})$ for some filter $\mathcal{F}$ on $E$.

        If $\mathcal{F}$, $\mathcal{F}'$ are filters on $E$, then it
is easy to see that
\begin{equation}
        \mathcal{I}(\mathcal{F}) \subseteq \mathcal{I}(\mathcal{F}')
\end{equation}
if and only if $\mathcal{F} \subseteq \mathcal{F}'$, which is to say
that $\mathcal{F}'$ is a refinement of $\mathcal{F}$.  It follows that
ultrafilters on $E$ correspond exactly to maximal ideals in
$\mathcal{A}$.  In particular, if $\mathcal{F}$ is an ultrafilter on
$E$, then $\mathcal{A} / \mathcal{I}(\mathcal{F})$ is a field.  One
can also see this more directly, as follows.  Suppose that $f \in
\mathcal{A}$ and $f \not\in \mathcal{I}(\mathcal{F})$, so that the
element of the quotient $\mathcal{A} / \mathcal{I}(\mathcal{F})$
corresponding to $f$ is not zero.  Thus $Z(f) \not\in \mathcal{F}$, by
definition of $\mathcal{I}(\mathcal{F})$, and so $E \backslash Z(f)
\in \mathcal{F}$, because $\mathcal{F}$ is an ultrafilter.  If $g \in
\mathcal{A}$ satisfies $f(x) \, g(x) = 1$ for every $x \in E
\backslash Z(f)$, then $f \, g - 1 \in \mathcal{I}(\mathcal{F})$,
which means that the product of the elements of the quotient
$\mathcal{A} / \mathcal{I}(\mathcal{F})$ corresponding to $f$, $g$ is
equal to the multiplicative identity element in the quotient, as
desired.

\section{Closure}
\label{closure}
\setcounter{equation}{0}

        Let $X$ be a topological space, and remember that $C_b(X)$ is
the algebra of bounded continuous real or complex-valued functions on
$X$, equipped with the supremum norm.  Put
\begin{equation}
\label{Z_epsilon(f) = {x in X : |f(x)| le epsilon}}
        Z_\epsilon(f) = \{x \in X : |f(x)| \le \epsilon\}
\end{equation}
for every $f \in C_b(X)$ and $\epsilon > 0$, which is a closed set in
$X$, since $f$ is continuous.  Note that $Z_\epsilon(f) = \emptyset$
for some $\epsilon > 0$ if and only if $f$ is invertible in $C_b(X)$.
If $\mathcal{F}$ is a filter on $X$, then let
$\mathcal{I}(\mathcal{F})$ be the collection of $f \in C_b(X)$ such
that $f_*(\mathcal{F})$ converges to $0$ in ${\bf R}$ or ${\bf C}$, as
appropriate.  Equivalently,
\begin{equation}
 \mathcal{I}(\mathcal{F}) = \{f \in C_b(X) : Z_\epsilon(f) \in \mathcal{F}
                                             \hbox{ for every } \epsilon > 0\}.
\end{equation}
This is analogous to, but different from, the definition in the
previous section.  It is not difficult to check that
$\mathcal{I}(\mathcal{F})$ is a proper closed ideal in $C_b(X)$ under
these conditions.  This uses the fact that
\begin{equation}
        Z_{\epsilon/2}(f) \cap Z_{\epsilon/2}(g) \subseteq Z_\epsilon(f + g)
\end{equation}
for every $f, g \in C_b(X)$ and $\epsilon > 0$, and that
\begin{equation}
        Z_{\epsilon/k}(f) \subseteq Z_\epsilon(f \, g)
\end{equation}
when $|g(x)| \le k$ for each $x \in X$ and $k > 0$.

        Let $\overline{\mathcal{F}}$ be the collection of subsets $B$
of $X$ for which there is an $A \in \mathcal{F}$ such that
$\overline{A} \subseteq B$.  One can check that
$\overline{\mathcal{F}}$ is also a filter on $X$, and that
$\mathcal{I}(\overline{\mathcal{F}}) = \mathcal{I}(\mathcal{F})$.
Thus one might as well restrict one's attention to filters on $X$
generated by closed subsets of $X$.  As a special case, if $A
\subseteq X$ is nonempty and $\mathcal{F}^A$ is the filter consisting
of $B \subseteq X$ such that $A \subseteq B$, then
$\overline{\mathcal{F}^A} = \mathcal{F}^{\overline{A}}$.

        Now let $\mathcal{I}$ be a proper ideal in $C_b(X)$, and put
\begin{equation}
 \mathcal{F}(\mathcal{I}) = \{A \subseteq X : Z_\epsilon(f) \subseteq A
               \hbox{ for some } f \in \mathcal{I} \hbox{ and } \epsilon > 0\}.
\end{equation}
Again this is analogous to, but different from, the corresponding
definition in the previous section.  One can also check that
$\mathcal{F}(\mathcal{I})$ is a filter on $X$ under these conditions.
This uses the fact that $Z_\epsilon(f) \ne \emptyset$ for each $f \in
\mathcal{I}$ and $\epsilon > 0$, because $\mathcal{I}$ is proper.  If
$f \in \mathcal{I}$, then $|f|^2 \in \mathcal{I}$, and
\begin{equation}
        Z_{\epsilon^2}(|f|^2) = Z_\epsilon(f),
\end{equation}
which means that one can restrict one's attention to nonnegative
real-valued functions in $\mathcal{I}$.  If $f$, $g$ are nonnegative
real-valued functions on $X$ and $\epsilon > 0$, then
\begin{equation}
        Z_\epsilon(f + g) \subseteq Z_\epsilon(f) \cap Z_\epsilon(g).
\end{equation}
This implies that $A \cap B \in \mathcal{F}(\mathcal{I})$ for every
$A, B \in \mathcal{F}(\mathcal{I})$.  Note that
$\mathcal{F}(\mathcal{I})$ is automatically generated by closed
subsets of $X$.  One can also check that $\mathcal{F}(\mathcal{I})$ is
the same as the filter associated to the closure of $\mathcal{I}$ in
$C_b(X)$, with respect to the supremum norm.  This uses the fact that
\begin{equation}
        Z_{\epsilon/2}(f) \subseteq Z_\epsilon(g)
\end{equation}
when $|f(x) - g(x)| \le \epsilon/2$ for every $x \in X$.

        By construction, $\mathcal{I} \subseteq
\mathcal{I}(\mathcal{F}(\mathcal{I}))$.  We have seen that
$\mathcal{I}(\mathcal{F})$ is closed in $C_b(X)$ for any filter
$\mathcal{F}$ on $X$, and so $\overline{\mathcal{I}} \subseteq
\mathcal{I}(\mathcal{F}(\mathcal{I}))$.  In order to show that
\begin{equation}
        \overline{\mathcal{I}} = \mathcal{I}(\mathcal{F}(\mathcal{I})),
\end{equation}
let $f \in \mathcal{I}(\mathcal{F}(\mathcal{I}))$ and $\epsilon > 0$ be given.
By definition of $\mathcal{I}(\mathcal{F}(\mathcal{I}))$, there are a
$g \in \mathcal{I}$ and a $\delta > 0$ such that
\begin{equation}
        Z_\delta(g) \subseteq Z_\epsilon(f).
\end{equation}
Put
\begin{equation}
        f_\eta = f \, \frac{|g|^2}{|g|^2 + \eta^2}
\end{equation}
for each $\eta > 0$.  Thus $f_\eta \in C_b(X)$ for each $\eta$, and in
fact $f_\eta \in \mathcal{I}$, because $g \in \mathcal{I}$.  We would
like to check that
\begin{equation}
        |f(x) - f_\eta(x)| = |f(x)| \, \frac{\eta^2}{|g(x)|^2 + \eta^2}
                           \le \epsilon
\end{equation}
for every $x \in X$ when $\eta$ is sufficiently small.  If $x \in
Z_\epsilon(f)$, then this holds for every $\eta > 0$, since $|f(x)|
\le \epsilon$ and $\eta^2/(|g(x)|^2 + \eta^2) \le 1$.  If $x \not\in
Z_\epsilon(f)$, then $x \not\in Z_\delta(g)$, so that $|g(x)| >
\delta$, and the desired estimate holds when $\eta$ is sufficiently
small, because $f$ is bounded.

\section{Regular topological spaces}
\label{regular spaces}
\setcounter{equation}{0}

        Remember that a topological space $X$ is said to be
\emph{regular}, or equivalently to satisfy the third separation
condition, if it has the following two properties.  First, $X$ should
satisfy the first separation condition, so that subsets of $X$ with
exactly one element are closed.  Second, for each $x \in X$ and closed
set $E \subseteq X$ with $x \not\in E$, there should be disjoint open
subsets $U$, $V$ of $X$ such that $x \in U$ and $E \subseteq V$.  In
particular, this implies that $X$ is Hausdorff, since one can take $E
= \{y\}$ when $y \in X$ and $y \ne x$.  Sometimes the term ``regular''
is used for topological spaces with the second property just
mentioned, and then the third separation condition is defined to be
the combination of regularity with the first separation condition.  We
shall include the first separation condition in the definition of
regularity here for the sake of simplicity.  As in Section
\ref{sigma-compactness}, it is well known that locally compact
Hausdorff topological spaces are regular.

        Equivalently, $X$ is regular if it satisfies the first
separation condition and for each $x \in X$ and open set $W \subseteq
X$ with $x \in W$ there is an open set $U \subseteq X$ such that $x
\in U$ and $\overline{U} \subseteq W$.  This corresponds to the
previous definition with $W = X \backslash E$.  Let $\mathcal{F}$ be a
filter on $X$, and let $\overline{\mathcal{F}}$ be the filter on $X$
generated by the closures of the elements of $\mathcal{F}$, as in the
previous section.  If $\mathcal{F}$ converges to a point $x \in X$ and
$X$ is regular, then it is easy to see that $\overline{\mathcal{F}}$
also converges to $x$ in $X$.  For if $W$ is an open set in $X$ that
contains $x$ and $U$ is an open set in $X$ that contains $x$ and
satisfies $\overline{U} \subseteq W$, then $U \in \mathcal{F}$,
because $\mathcal{F}$ converges to $x$, and hence $W \in
\overline{\mathcal{F}}$.

        Now let $x \in X$ be given, and let $\mathcal{F}(x)$ be the
collection of subsets $A$ of $X$ for which there is an open set $U
\subseteq X$ such that $x \in U$ and $U \subseteq A$.  This is a
filter on $X$ that converges to $x$, by construction.  The filter
$\overline{\mathcal{F}(x)}$ generated by the closed subsets of $X$ is
the same as the collection of subsets $B$ of $X$ for which there is an
open set $U \subseteq X$ such that $x \in U$ and $\overline{U}
\subseteq B$.  If $\overline{\mathcal{F}(x)}$ converges to $x$, then
for each open set $W \subseteq X$ with $x \in W$ there is an open set
$U \subseteq X$ such that $x \in U$ and $\overline{U} \subseteq W$.
It follows that $X$ is regular if it satisfies the first separation
condition and $\overline{\mathcal{F}(x)}$ converges to $x$ for every
$x \in X$.

        Of course, metric spaces are regular as topological spaces.
Real and complex topological vector spaces are also regular as
topological spaces.  To see this, remember that if $U$ is an open set
in a topological vector space $V$ that contains $0$, then there are
open subsets $U_1$, $U_2$ of $V$ that contain $0$ and satisfy
\begin{equation}
        U_1 + U_2 \subseteq U,
\end{equation}
as in Section \ref{topological vector spaces, continued}.  Moreover,
\begin{equation}
        \overline{U_1} \subseteq U_1 + U_2,
\end{equation}
as in (\ref{overline{A} subseteq A + U}).  Hence $\overline{U_1}
\subseteq U$, which implies that $V$ is regular, because of the
translation-invariance of the topology on $V$.

\section{Closed sets}
\label{closed sets}
\setcounter{equation}{0}

        Let $X$ be a topological space, and let us say that a nonempty
collection $\mathcal{E}$ of nonempty closed subsets of $X$ is a
\emph{C-filter} if $A \cap B \in \mathcal{E}$ for every $A, B \in
\mathcal{E}$, and if $E \in \mathcal{E}$ whenever $E \subseteq X$ is a
closed set such that $A \subseteq E$ for some $A \in \mathcal{E}$.
This is the same as a filter on $X$, except that we restrict our
attention to closed subsets of $X$.  If $\mathcal{F}$ is a filter on
$X$ and $\mathcal{E}(\mathcal{F})$ is the collection of closed subsets
of $X$ that are elements of $X$, then $\mathcal{E}(\mathcal{F})$ is a
C-filter.  This can also be described as the collection of closures of
elements of $\mathcal{F}$, since the closure of an element of
$\mathcal{F}$ is a closed set in $X$ that is contained in
$\mathcal{F}$.

        A C-filter $\mathcal{E}$ on $X$ also generates an ordinary
filter $\mathcal{F}(\mathcal{E})$ on $X$, consisting of the subsets
$B$ of $X$ that contain an element of $\mathcal{E}$ as a subset.  If
$\mathcal{F}$ is any filter on $X$, and $\mathcal{E}(\mathcal{F})$ is
the C-filter obtained from it as in the preceding paragraph, then the
filter generated by $\mathcal{E}(\mathcal{F})$ is the same as the
filter $\overline{\mathcal{F}}$ defined previously.  However, if
$\mathcal{E}$ is any C-filter on $X$, and $\mathcal{F}(\mathcal{E})$
is the ordinary filter generated by $\mathcal{E}$, then the C-filter
of closed sets in $\mathcal{F}(\mathcal{E})$ is the same as
$\mathcal{E}$.

        Let us say that a C-filter $\mathcal{E}$ on $X$ converges to a
point $x \in X$ if for every open set $U \subseteq X$ with $x \in U$
there is an $E \in \mathcal{E}$ such that $E \subseteq U$.  This is
equivalent to saying that $U \in \mathcal{F}(\mathcal{E})$ for every
open set $U \subseteq X$ with $x \in U$, so that $\mathcal{E}$
converges to $x$ if and only if $\mathcal{F}(\mathcal{E})$ converges
to $x$.  If $X$ is Hausdorff, then the limit of a convergent C-filter
on $X$ is unique, for the same reasons as for ordinary filters.  If
$\mathcal{F}$ is an ordinary filter on $X$ that converges to a point
$x \in X$ and $X$ is regular, then the corresponding C-filter
$\mathcal{E}(\mathcal{F})$ also converges to $x$, for the same reasons
as in the preceding section.

        Let $A$ be a nonempty subset of $X$, and let $\mathcal{E}^A$
be the collection of closed sets $B \subseteq X$ that contain $A$.
This is a C-filter on $X$, and it is easy to see that
$\mathcal{E}^{\overline{A}} = \mathcal{E}^A$ for every $A \subseteq
X$.  Note that $A \in \mathcal{E}^A$ if and only if $A$ is a closed
set in $X$.  If $A = \{p\}$ for some $p \in X$ and $X$ satisfies the
first separation condition, then $\{p\} \in \mathcal{E}^A$, and
$\mathcal{E}^A$ converges to $p$ in $X$.

        Suppose that $\mathcal{E}$ is a C-filter on $X$ that converges
to a point $p \in X$, and let $A \in \mathcal{E}$ be given.  If $U$ is
an open set in $X$ that contains $p$, then there is an $E \in
\mathcal{E}$ such that $E \subseteq U$, by definition of convergence.
This implies that $A \cap U \ne \emptyset$, because $A \cap E$ is
contained in $A \cap U$ and nonempty, since it is an element of
$\mathcal{E}$.  It follows that $p \in A$ for every $A \in
\mathcal{E}$, because every $A \in \mathcal{E}$ is a closed set in
$X$.

        Let $\mathcal{E}$ be a C-filter on $X$, and suppose that $B
\subseteq X$ satisfies $A \cap B \ne \emptyset$ for every $A \in
\mathcal{E}$.  Let $\mathcal{E}_B$ be the collection of closed subsets
$E$ of $X$ such that $A \cap B \subseteq E$ for some $A \in
\mathcal{E}$.  It is easy to see that this is also a C-filter on $X$,
which is a refinement of $\mathcal{E}$ in the sense that $\mathcal{E}
\subseteq \mathcal{E}_B$ as collections of subsets of $X$.  If $B$ is
a closed set in $X$, then $A \cap B \in \mathcal{E}_B$ for every $A
\in \mathcal{E}$.

        A C-filter $\mathcal{E}$ on $X$ may be described as a
\emph{C-ultrafilter} if it is maximal with respect to inclusion.  More
precisely, $\mathcal{E}$ is a C-ultrafilter if for every C-filter
$\mathcal{E}'$ such that $\mathcal{E} \subseteq \mathcal{E}'$, we have
that $\mathcal{E} = \mathcal{E}'$.  Using Zorn's lemma or the
Hausdorff maximality principle, one can show that every C-filter has a
refinement which is a C-ultrafilter, just as for ordinary
ultrafilters.

        For each $p \in X$, let $\mathcal{E}_p$ be the C-filter
consisting of all closed subsets of $X$ that contain $p$ as an
element.  This is the same as $\mathcal{E}^A$ with $A = \{p\}$, as
before.  If $X$ satisfies the first separation condition, then $\{p\}$
is a closed set in $X$, $\{p\} \in \mathcal{E}_p$, and it is easy to
see that $\mathcal{E}_p$ is a C-ultrafilter on $X$.  If $\mathcal{E}$
is any C-filter on $X$ and $p \in E$ for each $E \in \mathcal{E}$,
then $\mathcal{E} \subseteq \mathcal{E}_p$, and hence $\mathcal{E} =
\mathcal{E}_p$ when $\mathcal{E}$ is a $C$-ultrafilter.  In
particular, this holds when $\mathcal{E}$ converges to $p$.  If
$\mathcal{E}$ is a C-filter on $X$ and $X$ is compact, then
$\bigcap_{E \in \mathcal{E}} E \ne \emptyset$, because $\mathcal{E}$
has the finite intersection property.  If $\mathcal{E}$ is a
C-ultrafilter, then it follows that $\mathcal{E} = \mathcal{E}_p$ for
some $p \in X$.

        Let $\mathcal{E}$ be a C-filter on $X$, and suppose that $B$
is a closed set in $X$ such that $A \cap B \ne \emptyset$ for every $A
\in \mathcal{E}$.  This implies that $\mathcal{E} \subseteq
\mathcal{E}_B$, where $\mathcal{E}_B$ is the C-filter generated by the
intersections $A \cap B$ with $A \in \mathcal{E}$, as before.  If
$\mathcal{E}$ is a C-ultrafilter, then it follows that $\mathcal{E} =
\mathcal{E}_B$, and hence $B \in \mathcal{E}$.  Conversely, a C-filter
$\mathcal{E}$ is a C-ultrafilter when $B \in \mathcal{E}$ for every
closed set $B \subseteq X$ such that $A \cap B \ne \emptyset$ for
every $A \in \mathcal{E}$.  For if $\mathcal{E}'$ is a C-filter on $X$
such that $\mathcal{E} \subseteq \mathcal{E}'$, then $A \cap B$ is
contained in $\mathcal{E}'$ and is therefore nonempty for every $A \in
\mathcal{E}$ and $B \in \mathcal{E}'$.

        Let $X$, $Y$ be topological spaces, and let $f$ be a
continuous mapping from $X$ into $Y$.  Thus $f^{-1}(B)$ is a closed
set in $X$ for every closed set $B$ in $Y$.  Also let $\mathcal{E}$ be
a C-filter on $X$, and let $f_*(\mathcal{E})$ be the collection of
closed sets $B \subseteq Y$ such that $f^{-1}(B) \in \mathcal{E}$.  It
is easy to see that $f_*(\mathcal{E})$ is a C-filter on $Y$.  Note
that the closure of $f(A)$ in $Y$ is an element of $f_*(\mathcal{E})$
for each $A \in \mathcal{E}$.

        Suppose that $Y$ is compact, so that $\bigcap_{B \in
f_*(\mathcal{E})} B \ne \emptyset$, and let $q$ be an element of the
intersection.  Thus $q$ is contained in the closure of $f(A)$ in $Y$
for every $A \in \mathcal{E}$.  If $V$ is any open set in $Y$ that
contains $q$, then $f(A) \cap V \ne \emptyset$ for every $A \in
\mathcal{E}$, and hence $A \cap f^{-1}(V) \ne \emptyset$.  Let
$\mathcal{E}'$ be the collection of closed sets $E$ in $X$ such that
$A \cap f^{-1}(V) \subseteq E$ for some $A \in \mathcal{E}$ and open
set $V \subseteq Y$ with $q \in V$.  It is easy to see that
$\mathcal{E}'$ is a C-filter on $X$ that is a refinement of
$\mathcal{E}$, and that $A \cap f^{-1}(\overline{V}) \in \mathcal{E}'$
for every open set $V \subseteq Y$ with $q \in V$.  In particular,
$f^{-1}(\overline{V}) \in \mathcal{E}'$ under these conditions, which
means that $\overline{V} \in f_*(\mathcal{E}')$.  If $Y$ is also
Hausdorff, and hence regular, then it follows that $f_*(\mathcal{E}')$
converges to $q$ in $Y$.  If $\mathcal{E}$ is a C-ultrafilter on $X$,
then $\mathcal{E} = \mathcal{E}'$, and $f_*(\mathcal{E})$ converges to
$q$ in $Y$.

\section{Multi-indices}
\label{multi-indices}
\setcounter{equation}{0}

        Let $n$ be a positive integer, which will be kept fixed
throughout this section.  A \emph{multi-index} $\alpha = (\alpha_1,
\ldots, \alpha_n)$ is an $n$-tuple of nonnegative integers.  The sum
of two multi-indices is defined coordinatewise, and we put
\begin{equation}
        |\alpha| = \alpha_1 + \cdots + \alpha_n.
\end{equation}

        If $\alpha$ is a multi-index and $x = (x_1, \ldots, x_n) \in
{\bf R}^n$, then the corresponding monomial $x^\alpha$ is defined by
the product
\begin{equation}
\label{x^alpha}
        x^\alpha = x_1^{\alpha_1} \cdots x_n^{\alpha_n}.
\end{equation}
More precisely, $x_j^{\alpha_j}$ is interpreted as being equal to $1$
for every $x_j \in {\bf R}$ when $\alpha_j = 0$, so that $x^\alpha =
1$ for every $x \in {\bf R}^n$ when $\alpha = 0$.  Note that
$|\alpha|$ is the same as the degree of the monomial $x^\alpha$, and a
polynomial on ${\bf R}^n$ is the same as a linear combination of
finitely many monomials.  Moreover,
\begin{equation}
        x^{\alpha + \beta} = x^\alpha \, x^\beta
\end{equation}
for all multi-indices $\alpha$, $\beta$ and $x \in {\bf R}^n$.

        If $l$ is a positive integer, then $l!$ is ``$l$ factorial'',
the product of $1, \ldots, l$.  It is customary to include $l = 0$ by
setting $0! = 1$.  If $\alpha$ is a multi-index, then we put
\begin{equation}
        \alpha! = \alpha_1! \cdots \alpha_n!.
\end{equation}
If $\alpha$ is a multi-index and $x, y \in {\bf R}^n$, then
\begin{equation}
 (x + y)^\alpha = \sum_{\alpha = \beta + \gamma}
                   \frac{\alpha!}{\beta! \, \gamma!} \, x^\beta \, y^\gamma,
\end{equation}
where the sum is taken over all multi-indices $\beta$, $\gamma$ such
that $\alpha = \beta + \gamma$.  This follows from the binomial
theorem applied to $(x_j + y_j)^{\alpha_j}$ for $j = 1, \ldots, n$.

        Let $\partial_j = \partial / \partial x_j$ be the usual
partial derivative in $x_j$, $1 \le j \le n$.  If $\alpha$ is a
multi-index, then the corresponding differential operator
$\partial^\alpha$ is defined by
\begin{equation}
\label{partial^alpha}
        \partial^\alpha = \partial_1^{\alpha_1} \cdots \partial_n^{\alpha_n}.
\end{equation}
Here $\partial_j^{\alpha_j}$ is interpreted as being the identity
operator when $\alpha_j = 0$, so that $\partial^\alpha$ reduces to the
identity operator when $\alpha = 0$.  Observe that
\begin{equation}
        \partial^{\alpha + \beta} = \partial^\alpha \, \partial^\beta
\end{equation}
for all multi-indices $\alpha$, $\beta$.

\section{Smooth functions}
\label{smooth functions}
\setcounter{equation}{0}

        Let $U$ be a nonempty open set in ${\bf R}^n$ for some
positive integer $n$, and let $C^\infty(U)$ be the space of real or
complex-valued functions on $U$ that are smooth in the sense that they
are continuously-differentiable of all orders.  As usual, this may
also be denoted $C^\infty(U, {\bf R})$ or $C^\infty(U, {\bf C})$, to
indicate whether real or complex-valued functions are being used.  It
is well known that $C^\infty(U)$ is a commutative algebra with respect
to pointwise addition and multiplication, since sums and products of
smooth functions are smooth.

        If $\alpha$ is a multi-index and $K \subseteq U$ is a nonempty
compact set, then
\begin{equation}
\label{||f||_{alpha, K} = sup_{x in K} |partial^alpha f(x)|}
 \|f\|_{\alpha, K} = \sup_{x \in K} |\partial^\alpha f(x)|
\end{equation}
defines a seminorm on $C^\infty(U)$.  This is the same as the supremum
seminorm $\|f\|_K$ of $f$ over $K$ when $\alpha = 0$, and otherwise
this is the same as the supremum seminorm of $\partial^\alpha f$ over
$K$.  The collection of all of these seminorms defines a topology on
$C(U)$, as in Section \ref{seminorms, topologies}.  Of course, $U$ is
a locally compact Hausdorff topological space with respect to the
topology induced by the standard topology on ${\bf R}^n$, and one can
also check that $U$ is $\sigma$-compact.  As in Section \ref{locally
compact spaces}, there is a sequence of compact subsets $K_1, K_2,
\ldots$ of $U$ such that every compact set $H \subseteq U$ is
contained in $K_l$ for some $l$.  It follows that the seminorms
$\|f\|_{\alpha, K_l}$ are sufficient to determine the same topology on
$C^\infty(U)$ as the one that was just described, where $\alpha$ is a
multi-index and $l$ is a positive integer.  In particular, this
collection of seminorms on $C^\infty(U)$ is countable, since there are
only countably many multi-indices.

        If $f$, $g$ are smooth functions on $U$ and $\alpha$ is a
multi-index, then
\begin{equation}
\label{partial^alpha (f g) = ...}
 \partial^\alpha (f \, g) = \sum_{\alpha = \beta + \gamma}
\frac{\alpha!}{\beta! \, \gamma!} \, (\partial^\beta f) \, (\partial^\gamma g),
\end{equation}
where the sum is taken over all multi-indices $\beta$, $\gamma$ such
that $\alpha = \beta + \gamma$.  This can be derived from the usual
product rule for first derivatives, starting with the $n = 1$ case.
Using this identity, it is easy to check that multiplication of
functions is continuous as a mapping from $C^\infty(U) \times
C^\infty(U)$ into $C^\infty(U)$, with respect to the topology on
$C^\infty(U)$ defined in the previous paragraph.

        Let $\phi$ be a homomorphism from $C^\infty(U)$ into the real
or complex numbers, as appropriate.  As usual, we suppose also that
$\phi$ is nontrivial in the sense that $\phi(f) \ne 0$ for some $f \in
C^\infty(U)$.  This implies that $\phi({\bf 1}_U) = 1$, where ${\bf
1}_U$ is the constant function on $U$ equal to $1$ at every point.  If
$f$ is a smooth function on $U$ such that $f(x) \ne 0$ for every $x
\in U$, then $1/f(x)$ is also a smooth function on $U$, and it follows
that
\begin{equation}
        \phi(f) \, \phi(1/f) = \phi({\bf 1}_U) = 1.
\end{equation}
In particular, $\phi(f) \ne 0$ when $f(x) \ne 0$ for every $x \in U$.
Equivalently, if $f$ is any smooth function on $U$ and $\phi(f) = 0$,
then $f(x) = 0$ for some $x \in U$.  If $f$ is any smooth function on
$U$ and $\phi(f) = c$, then there is an $x \in U$ such that $f(x) =
c$, since one can apply the previous statement to $f - c \, {\bf
1}_U$.

        Let $f_j$ be the smooth function on $U$ defined by $f_j(x) =
x_j$, $j = 1, \ldots, n$, and put $p_j = \phi(f_j)$.  We would like to
check that
\begin{equation}
        p = (p_1, \ldots, p_n) \in U.
\end{equation}
Consider the smooth function on $U$ given by
\begin{equation}
        f(x) = \sum_{j = 1}^n (x_j - p_j)^2.
\end{equation}
Equivalently,
\begin{equation}
        f = \sum_{j = 1}^n (f_j - p_j \, {\bf 1}_U)^2,
\end{equation}
and so
\begin{equation}
        \phi(f) = \sum_{j = 1}^n (\phi(f_j) - p_j)^2 = 0,
\end{equation}
because $\phi$ is a homomorphism.  Hence $f(x) = 0$ for some $x \in
U$, as in the previous paragraph, which is only possible if $x = p$,
in which case $p \in U$.

        If $g$ is a smooth function on $U$ and $U$ is convex, then
\begin{eqnarray}
 g(x) - g(p) & = & 
        \int_0^1 (\partial / \partial t) g(t \, x + (1 - t) \, p) \, dt \\
             & = & \sum_{j = 1}^n (x_j - p_j) \,
              \int_0^1 (\partial_j g)(t \, x + (1 - t) \, p) \, dt. \nonumber
\end{eqnarray}
Hence there are smooth functions $g_1, \ldots, g_n$ on $U$ such that
\begin{equation}
\label{g(x) = g(p) + sum_{j = 1}^n (x_j - p_j) g_j(x)}
        g(x) = g(p) + \sum_{j = 1}^n (x_j - p_j) \, g_j(x).
\end{equation}
This also works when $g$ is the restriction to $U$ of a smooth
function on a convex open set that contains $U$, such as ${\bf R}^n$
itself.  In particular, this works when $g(x) = 0$ on the complement
of a closed ball contained in $U$.  Otherwise, if $g(x) = 0$ for every
$x$ in a neighborhood of $p$, then we can simply take $g_j(x)$ to be
$(x_j - p_j) / |x - p|^2$ times $g(x)$, where $|x - p|^2 = \sum_{l =
1}^n (x_l - p_l)^2$, and $g_j(p) = 0$.  Any smooth function on $U$ can
be expressed as the sum of a smooth function supported on a closed
ball in $U$ and a smooth function that vanishes on a neighborhood of
$p$, using standard cut-off functions.  It follows that every smooth
function $g$ on $U$ can be expressed as in (\ref{g(x) = g(p) + sum_{j
= 1}^n (x_j - p_j) g_j(x)}) for some smooth functions $g_1, \ldots,
g_n$ on $U$.

        Using this representation, we get that $\phi(g) = g(p)$ for
every $g \in C^\infty(U)$.  Of course, $\phi_p(g) = g(p)$ defines a
homomorphism on $C^\infty(U)$ for every $p \in U$.

\section{Polynomials}
\label{polynomials}
\setcounter{equation}{0}

        Let $\mathcal{P}({\bf R}^n)$ be the space of polynomials on
${\bf R}^n$ with real coefficients, which can be expressed as finite
linear combinations of the monomials $x^\alpha$, where $\alpha$ is a
multi-index.  This is an algebra in a natural way, corresponding to
pointwise addition and multiplication of functions.  If $p \in {\bf
R}^n$, then $\phi_p(f) = f(p)$ defines a homomorphism on
$\mathcal{P}({\bf R}^n)$, as usual.  Conversely, if $\phi$ is a
homomorphism on $\mathcal{P}({\bf R}^n)$ which is not identically $0$,
then $\phi = \phi_p$ for some $p \in {\bf R}^n$.  As in the previous
section, $p = (p_1, \ldots, p_n)$ is given by $p_j = \phi(f_j)$, where
$f_j(x) = x_j$.  In this case, the fact that $\phi(f) = f(p)$ for
every polynomial $f$ on ${\bf R}^n$ follows from simple algebra.
There are analogous statements for polynomials on ${\bf C}^n$ with
complex coefficients, which can be expressed as finite linear
combinations of monomials $z^\alpha = z_1^{\alpha_1} \cdots
z_n^{\alpha_n}$, $z = (z_1, \ldots, z_n) \in {\bf C}^n$.

\section{Continuously-differentiable functions}
\label{C^1 functions}
\setcounter{equation}{0}

        A real or complex-valued function $f$ on the closed unit
interval $[0, 1]$ is said to be \emph{continuously differentiable} if
it satisfies the following three conditions.  First, the derivative
$f'(x)$ of $f$ should exist at every $x$ in the open unit interval
$(0, 1)$.  Second, the appropriate one-sided derivatives should exist
at the endpoints $0$, $1$, which will also be denoted $f'(0)$, $f'(1)$
for simplicity.  Third, the resulting function $f'(x)$ should be
continuous on $[0, 1]$.  Of course, differentiability of $f$ implies
that $f$ is continuous on $[0, 1]$.

        Equivalently, a continuous function $f$ on $[0, 1]$ is
continuously differentiable if it is differentiable on $(0, 1)$, and
if the derivative can be extended to a continuous function on $[0,
1]$, also denoted $f'$.  More precisely, one can check that the
one-sided derivatives of $f$ exist at the endpoints, and are given by
the extension of $f'$ to $0$, $1$.  This follows from the fact that
\begin{equation}
        f(y) - f(x) = \int_x^y f'(t) \, dt
\end{equation}
when $0 \le x \le y \le 1$.

        The space of continuously-differentiable functions on $[0, 1]$
may be denoted $C^1([0, 1])$, or by $C^1([0, 1], {\bf R})$, $C^1([0,
1], {\bf C})$ to indicate whether real or complex-valued functions are
being used.  As usual, $C^1([0, 1])$ is an algebra with respect to
pointwise addition and scalar multiplication of functions.  If
\begin{equation}
        \|f\|_{sup} = \sup_{0 \le x \le 1} |f(x)|
\end{equation}
is the supremum norm of a bounded function on $[0, 1]$, then
\begin{equation}
        \|f\|_{C^1} = \|f\|_{C^1([0, 1])} = \|f\|_{sup} + \|f'\|_{sup}
\end{equation}
is a natural choice of norm on $C^1([0, 1])$.  In particular,
\begin{equation}
        \|f \, g\|_{C^1} \le \|f\|_{C^1} \, \|f\|_{C^1}
\end{equation}
for every $f, g \in C^1([0, 1])$.  To see this, remember that
\begin{equation}
        \|f \, g\|_{sup} \le \|f\|_{sup} \, \|g\|_{sup},
\end{equation}
so that
\begin{equation}
        \|f \, g\|_{C^1} = \|f \, g\|_{sup} + \|(f \, g)'\|_{sup}
            \le \|f\|_{sup} \, \|g\|_{sup} + \|(f \, g)'\|_{sup}.
\end{equation}
The product rule implies that
\begin{equation}
        \|(f \, g)'\|_{sup} = \|f' \, g + f \, g'\|_{sup}
                    \le \|f'\|_{sup} \|g\|_{sup} + \|f\|_{sup} \, \|g'\|_{sup},
\end{equation}
and hence
\begin{eqnarray}
 \|f \, g\|_{C^1} & \le & \|f\|_{sup} \, \|g\|_{sup}
            + \|f'\|_{sup} \, \|g\|_{sup} + \|f\|_{sup} \, \|g'\|_{sup} \\
& \le & (\|f\|_{sup} + \|f'\|_{sup}) \, (\|g\|_{sup} + \|g'\|_{sup})\nonumber\\
                    & = & \|f\|_{C^1} \, \|g\|_{C^1}. \nonumber
\end{eqnarray}
Note that the $C^1$ norm of a constant function is the same as the
absolute value or modulus of the corresponding real or complex number.
One can also check that $C^1([0, 1])$ is complete with respect to the
$C^1$ norm, so that $C^1([0, 1])$ is a Banach algebra.

        Remember that continuous functions on $[0, 1]$ can be
approximated uniformly by polynomials, by Weierstrass' approximation
theorem.  Using this, one can show that continuously-differentiable
functions on $[0, 1]$ can be approximated by polynomials in the $C^1$
norm.  More precisely, in order to approximate a
continuously-differentiable function $f$ on $[0, 1]$ by polynomials in
the $C^1$ norm, one can integrate polynomials that approximate $f'$
uniformly on $[0, 1]$.  One can choose the constant terms of these
approximations to $f$ to be equal to $f(0)$, so that the approximation
of $f$ follows from the approximation of $f'$.

        Let $\phi$ be a homomorphism from $C^1([0, 1])$ into the real
or complex numbers, as appropriate.  Suppose also that $\phi(f) \ne 0$
for some $f \in C^1([0, 1])$, so that $\phi$ takes the constant
function equal to $1$ on $[0, 1]$ to $1$, by the usual argument.  If
$f$ is a continuously-differentiable function on $[0, 1]$ such that
$f(x) \ne 0$ for every $x \in [0, 1]$, then $1/f$ is also a
continuously-differentiable function on $[0, 1]$, and hence $\phi(f)
\ne 0$.  This implies that $\phi(f) \in f([0, 1])$ for every $f \in
C^1([0, 1])$, as in previous situations.  In particular, it follows
that
\begin{equation}
\label{|phi(f)| le ||f||_{sup} le ||f||_{C^1}}
        |\phi(f)| \le \|f\|_{sup} \le \|f\|_{C^1}
\end{equation}
for every $f \in C^1([0, 1])$.

        Of course, $f_0(x) = x$ is a continuously-differentiable
function on $[0, 1]$.  Put $p = \phi(f_0)$, so that $p \in f_0([0, 1])
= [0, 1]$.  It follows that
\begin{equation}
        \phi(f) = f(p)
\end{equation}
when $f$ is a polynomial, by simple algebra.  The same relation holds
for every $f \in C^1([0, 1])$, because polynomials are dense in
$C^1([0, 1])$ with respect to the supremum norm.  We do not need the
stronger fact that polynomials are dense in $C^1([0, 1])$ with respect
to the $C^1$ norm here, because $\phi$ is continuous with respect to
the supremum norm, by (\ref{|phi(f)| le ||f||_{sup} le ||f||_{C^1}}).

        Alternatively, we can use the continuity of $\phi$ with
respect to the supremum norm to extend $\phi$ to a homomorphism on
$C([0, 1])$, since $C^1([0, 1])$ is dense in $C([0, 1])$ with respect
to the supremum norm.  This permits us to use the results about
homomorphisms on $C(X)$ when $X$ is compact, as in Section
\ref{compact spaces}.  This approach has the advantage of working in
more abstract situations, such as on compact manifolds.  The same type
of arguments as in Section \ref{compact spaces} can also be used
directly in these situations.

        At any rate, every nonzero homomorphism on $C^1([0, 1])$ can
be represented as $\phi(f) = f(p)$ for some $p \in [0, 1]$.  Of
course, $\phi_p(f) = f(p)$ is a homomorphism on $C^1([0, 1])$ for
every $p \in [0, 1]$.

\section{Spectral radius}
\label{spectral radius}
\setcounter{equation}{0}

        Let $(\mathcal{A}, \|\cdot \|)$ be a Banach algebra over the
real or complex numbers with nonzero multiplicative identity element
$e$.  If $x \in \mathcal{A}$ satisfies $\|x\| < 1$, then $e - x$ is
invertible in $\mathcal{A}$, as in Section \ref{banach algebras}.  The
same conclusion also holds when $\|x^n\| < 1$ for any positive integer
$n$.  One way to see this is to use the previous result to get that $e
- x^n$ is invertible, and then observe that
\begin{equation}
\label{(e - x) (sum_{j = 1}^{n - 1} x^j) = ... = e - x^n}
        (e - x) \Big(\sum_{j = 1}^{n - 1} x^j\Big)
          = \Big(\sum_{j = 1}^{n -1} x^j\Big) (e - x) = e - x^n.
\end{equation}
This shows that the product of $e - x$ with an element of
$\mathcal{A}$ that commutes with it is invertible, which implies that
$e - x$ is invertible too, as in Section \ref{banach algebras}.
Alternatively, one can check that $\sum_{j = 1}^\infty \|x^j\|$
converges when $\|x^n\| < 1$ for some $n$, and then argue as in
Section \ref{banach algebras} that $\sum_{j = 1}^\infty x^j$ converges
in $\mathcal{A}$, and that the sum is the inverse of $e - x$.  To do
this, note first that every positive integer $j$ can be represented as
$l \, n + r$ for some nonnegative integers $l$, $r$ with $r < n$.
This leads to the estimate
\begin{equation}
\label{||x^j|| le ||x^n||^l ||x||^r}
        \|x^j\| \le \|x^n\|^l \, \|x\|^r,
\end{equation}
which implies the convergence of $\sum_{j = 1}^\infty \|x^j\|$ when
$\|x^n\| < 1$.

        If $x$ is any element of $\mathcal{A}$, then put
\begin{equation}
        r(x) = \inf_{n \ge 1} \|x^n\|^{1/n},
\end{equation}
where more precisely the infimum is taken over all positive integers $n$.
Thus $e - x$ is invertible in $\mathcal{A}$ when $r(x) < 1$, as in the
previous paragraph.  Observe also that
\begin{equation}
        r(t \, x) = |t| \, r(x)
\end{equation}
for every real or complex number $t$, as appropriate.  It follows that
$e - t \, x$ is invertible in $\mathcal{A}$ when $|t| \, r(x) < 1$.
Equivalently, $t \, e - x$ is invertible in $\mathcal{A}$ when $|t| >
r(x)$.

        Let us check that
\begin{equation}
        \lim_{j \to \infty} \|x^j\|^{1/j} = r(x),
\end{equation}
where the existence of the limit is part of the conclusion.  Because
of the way that $r(x)$ is defined, it suffices to show that
\begin{equation}
        \limsup_{j \to \infty} \|x^j\|^{1/j} \le r(x),
\end{equation}
which is the same as saying that
\begin{equation}
\label{limsup_{j to infty} ||x^j||^{1/j} le ||x^n||^{1/n}}
        \limsup_{j \to \infty} \|x^j\|^{1/j} \le \|x^n\|^{1/n}
\end{equation}
for each $n \ge 1$.  As before, each positive integer $j$ can be
represented as $l \, n + r$ for some nonnegative integers $l$, $r$
with $r < n$, and (\ref{||x^j|| le ||x^n||^l ||x||^r}) implies that
\begin{equation}
 \|x^j\|^{1/j} \le (\|x^n\|^{1/n})^{l n / j} \, \|x\|^{r / j}
                =  (\|x^n\|^{1/n})^{1 - (r/j)} \, \|x\|^{r/j}.
\end{equation}
It is not difficult to derive (\ref{limsup_{j to infty} ||x^j||^{1/j}
le ||x^n||^{1/n}}) from this estimate, using the fact that $a^{1/j}
\to 1$ as $j \to \infty$ for every positive real number $a$.  This is
trivial when $x^n = 0$, since $x^j$ is then equal to $0$ for each $j
\ge n$.

        As a basic class of examples, suppose that $\mathcal{A}$ is
the algebra $C_b(X)$ of bounded continuous functions on a topological
space $X$, equipped with the supremum norm.  In this case, it is easy
to see that
\begin{equation}
        \|f^n\|_{sup} = \|f\|_{sup}^n
\end{equation}
for every $f \in C_b(X)$ and $n \ge 1$, and hence that
\begin{equation}
        r(f) = \|f\|_{sup}.
\end{equation}

        Suppose now that $\mathcal{A}$ is the algebra $C^1([0, 1])$ of
continuously-differentiable functions on the unit interval, as in the
previous section.  Thus $\|f\|_{C^1} \ge \|f\|_{sup}$, and hence
\begin{equation}
        r(f) \ge \|f\|_{sup}
\end{equation}
for every $f \in C^1([0, 1])$.  In the other direction,
\begin{eqnarray}
 \|f^n\|_{C^1} & = & \|f^n\|_{sup} + \|(f^n)'\|_{sup}          \\
 & = & \|f\|_{sup}^n + \|n \, f' \, f^{n - 1}\|_{sup} \nonumber \\
 & \le & \|f\|_{sup}^n + n \, \|f'\|_{sup} \, \|f\|_{sup}^{n - 1}.\nonumber
\end{eqnarray}
for each $n$.  Using this, it is not too difficult to show that
\begin{equation}
        r(f) = \lim_{n \to \infty} \|f^n\|_{C^1}^{1/n} = \|f\|_{sup}.
\end{equation}
This also uses the fact that $(a + b \, n)^{1/n} \to 1$ as $n \to
\infty$ for any two positive real numbers $a$, $b$.

        Let $\mathcal{A}$ be a complex Banach algebra, and put
\begin{equation}
\label{R(x) = sup {|t| : t in {bf C}  and t e - x is not invertible}}
        R(x) = \sup \{|t| : t \in {\bf C} \hbox{ and }
                            t \, e - x \hbox{ is not invertible}\}
\end{equation}
for every $x \in \mathcal{A}$.  We have already seen that $t \, e - x$
is invertible in $\mathcal{A}$ when $|t| > r(x)$, which works for both
real and complex Banach algebras.  If $\mathcal{A}$ is a complex
Banach algebra, then for each $x \in \mathcal{A}$ there is a $t \in
{\bf C}$ such that $t \, e - x$ is not invertible, as in Section
\ref{banach algebras}.  Thus the supremum in the definition of $R(x)$
makes sense, and $R(x) \le r(x)$.  A well-known theorem states that
$r(x) \le R(x)$ for every $x \in \mathcal{A}$ when $\mathcal{A}$ is a
complex Banach algebra, and hence $r(x) = R(x)$.

        To see this, note that $t \, e - x$ is invertible when $t \in
{\bf C}$ satisfies $|t| > R(x)$, which implies that $e - t \, x$ is
invertible when $|t| \, R(x) < 1$.  As in Section \ref{banach
algebras}, the basic idea is to look at
\begin{equation}
        f(t) = (e - t \, x)^{-1}
\end{equation}
as a holomorphic function on the disk where $|t| \, R(x) < 1$ with
values in $\mathcal{A}$.  In particular, the composition of $f$ with a
continuous linear functional on $\mathcal{A}$ defines a complex-valued
function on this disk which is holomorphic in the usual sense.  We
also know that $f(t)$ is given by the power series $\sum_{j =
0}^\infty t^j \, x^j$ when $|t|$ is sufficiently small, as in Section
\ref{banach algebras}.  By standard arguments in complex analysis, one
can estimate the size of the coefficients of this power series in $t$
in terms of the behavior of $f(t)$ on any circle $|t| = a$ with $a \,
R(x) < 1$.  Note that $f(t)$ is bounded on any circle of this type,
because the circle is compact and $f(t)$ is continuous on it.  More
precisely, one can show that for each positive real number $a$ with $a
\, R(x) < 1$, there is a $C(a) \ge 0$ such that
\begin{equation}
        a^j \, \|x^j\| \le C(a)
\end{equation}
for every $j \ge 1$.  Equivalently, $a \, \|x^j\|^{1/j} \le
C(a)^{1/j}$ for each $j$, which implies that $a \, r(x) \le 1$ when $a
\, R(x) < 1$, by taking the limit as $j \to \infty$.  Thus $r(x) \le
R(x)$, as desired.

\section{Topological algebras}
\label{topological algebras}
\setcounter{equation}{0}

        Let $\mathcal{A}$ be an associative algebra over the real or
complex numbers, as in Section \ref{banach algebras}.  Suppose that
$\mathcal{A}$ is also equipped with a topology which makes it into a
topological vector space, as in Section \ref{topological vector
spaces}.  In the same way, one can ask that multiplication in
$\mathcal{A}$ be continuous as a mapping from $\mathcal{A} \times
\mathcal{A}$ into $\mathcal{A}$, using the product topology on
$\mathcal{A} \times \mathcal{A}$ associated to the given topology on
$\mathcal{A}$.  Under these conditions, we can say that $\mathcal{A}$
is a \emph{topological algebra}.  As before, we are especially
interested here in the case where multiplication on $\mathcal{A}$ is
commutative.

        Of course, Banach algebras are topological algebras, with
respect to the topology associated to the norm.  If $X$ is a locally
compact Hausdorff topological space, then the algebra of continuous
functions on $X$ is a topological algebra with respect to the topology
determined by the collection of supremum seminorms corresponding to
nonempty compact subsets of $X$, as in Section \ref{locally compact spaces}.
If $U$ is a nonempty open set in ${\bf R}^n$, then the algebra of
smooth functions on $U$ is a topological algebra with respect to the
collection of supremum seminorms of derivatives of $f$ over nonempty
compact subsets of $U$, as in Section \ref{smooth functions}.

        As in the case of Banach algebras, one may wish to look at
topological algebras $\mathcal{A}$ that are complete as topological
vector spaces.  If $\mathcal{A}$ has a countable local base for its
topology at $0$, then this can be defined in terms of convergence of
Cauchy sequences, as usual.  Otherwise, one can consider more general
Cauchy conditions for nets or filters on $\mathcal{A}$.  It is not too
difficult to show that the examples of topological algebras of
continuous and smooth functions mentioned in the previous paragraph
are complete.

        If $U$ is a nonempty open set in the complex plane, then the
algebra $\mathcal{H}(U)$ of holomorphic functions on $U$ may be
considered as a subalgebra of the algebra $C(U)$ of continuous
complex-valued functions on $U$.  More precisely, we have seen that
$\mathcal{H}(U)$ is a closed subalgebra of $C(U)$ with respect to the
topology associated to the collection of supremum seminorms over
nonempty compact subsets of $U$.  Of course, $\mathcal{H}(U)$ is also
a topological algebra with respect to the topology determined by this
collection of seminorms, and it follows that $\mathcal{H}(U)$ is
complete as well, because $C(U)$ is complete.

\section{Fourier series}
\label{fourier series}
\setcounter{equation}{0}

        Let ${\bf T}$ be the unit circle in the complex plane,
consisting of the $z \in {\bf C}$ with $|z| = 1$.  It is well known that
\begin{equation}
\label{int_{bf T} z^j |dz| = 0}
        \int_{\bf T} z^j \, |dz| = 0
\end{equation}
for each nonzero integer $j$, where $|dz|$ is the element of arc
length along ${\bf T}$.  This integral is the same as $-i$ times the
line integral
\begin{equation}
\label{oint_{bf T} z^{j - 1} dz}
        \oint_{\bf T} z^{j - 1} \, dz,
\end{equation}
the vanishing of which when $j \ne 0$ is a basic fact in complex analysis.
More precisely, the relationship between these two integrals follows from
identifying $i \, z$ with the unit tangent vector to ${\bf T}$ at $z$
in the positive orientation.  Note that
\begin{equation}
\label{overline{(int_{bf T} z^j |dz|)} = int_{bf T} z^{-j} |dz|}
 \overline{\Big(\int_{\bf T} z^j \, |dz|\Big)} = \int_{\bf T} z^{-j} \, |dz|,
\end{equation}
since $\overline{z} = z^{-1}$ when $|z| = 1$, and so it suffices to
verify (\ref{int_{bf T} z^j |dz| = 0}) when $j$ is a positive integer.
If $j = 0$, then $z^j$ is interpreted as being equal to $1$ for each
$z$, so that the integral in (\ref{int_{bf T} z^j |dz| = 0}) is equal
to the length $2 \, \pi$ of ${\bf T}$.

        If $f$ is a continuous complex-valued function on ${\bf T}$
and $j$ is an integer, then the $j$th Fourier coefficient of $f$ is
defined by
\begin{equation}
\label{widehat{f}(j) = frac{1}{2 pi} int_{bf T} f(w) w^{-j} |dw|}
        \widehat{f}(j) = \frac{1}{2 \pi} \int_{\bf T} f(w) \, w^{-j} \, |dw|.
\end{equation}
The corresponding Fourier series is given by
\begin{equation}
\label{sum_{j = -infty}^infty widehat{f}(j) z^j}
        \sum_{j = -\infty}^\infty \widehat{f}(j) \, z^j.
\end{equation}
For the moment, this should be considered as a formal sum, without
regard to convergence.  If $f(z) = z^l$ for some integer $l$, then
$\widehat{f}(j)$ is equal to $1$ when $j = l$ and to $0$ when $j \ne
l$, as in the previous paragraph.  Thus the Fourier series
(\ref{sum_{j = -infty}^infty widehat{f}(j) z^j}) reduces to $f(z)$ in
this case, and also when $f(z)$ is a linear combination of $z^l$ for
finitely many integers $l$.

        Suppose that $f(z)$ is a continuous complex-valued function on
the closed unit disk in ${\bf C}$ which is holomorphic on the open
unit disk.  By standard results in complex analysis, $f(z)$ can be
represented by an absolutely convergent power series
\begin{equation}
        f(z) = \sum_{j = 0}^\infty a_j \, z^j
\end{equation}
on the open unit disk, which is to say for $z \in {\bf C}$ with $|z| <
1$.  In this case,
\begin{equation}
        a_j = \widehat{f}(j)
\end{equation}
for each $j \ge 0$, where $\widehat{f}(j)$ is the $j$th Fourier
coefficient of the restriction of $f$ to the unit circle.  This
follows from the usual Cauchy integral formulae, where one integrates
over the unit circle.  Normally one might integrate over circles of
radius $r < 1$ when dealing with holomorphic functions on the open
unit disk, but one can pass to the limit $r \to 1$ when $f$ extends to
a continuous function on the closed unit disk.

        Under these conditions, we also have that $\widehat{f}(j) = 0$
when $j < 0$.  This can be derived from Cauchy's theorem for line
integrals of holomorphic functions, starting with integrals over
circles of radius $r < 1$, and then passing to the limit $r \to 1$ as
in the previous paragraph.  Conversely, if $f$ is a continuous
function on the unit circle with $\widehat{f}(j) = 0$ when $j < 0$,
then it can be shown that $f$ has a continuous extension to the closed
unit disk which is holomorphic on the open unit disk.  More precisely,
the holomorphic function on the open unit disk is given by the power
series defined by the Fourier coefficients of $f$, as before.  The
remaining point is to show that the combination of this holomorphic
function on the open unit disk with the given function $f$ on the unit
circle is continuous on the closed unit disk, which will be discussed
in Section \ref{poisson kernel}.

\section{Absolute convergence}
\label{absolute convergence}
\setcounter{equation}{0}

        Let $\ell^1({\bf Z})$ be the space of doubly-infinite
sequences $a = \{a_j\}_{j = -\infty}^\infty$ of complex numbers such
that
\begin{equation}
        \|a\|_1 = \sum_{j = -\infty}^\infty |a_j|
\end{equation}
converges.  This is equivalent to the definition in Section
\ref{summable functions} with $E = {\bf Z}$, but in this case it is a
bit simpler to think of a sum over ${\bf Z}$ as a combination of two
ordinary infinite series, corresponding to sums over $j \ge 0$ and $j
< 0$.  In particular, if $a \in \ell^1({\bf Z})$, then $\sum_{j =
0}^\infty a_j$ and $\sum_{j = 1}^\infty a_{-j}$ converge absolutely,
so that their sum $\sum_{j = -\infty}^\infty a_j$ is well-defined, and
satisfies
\begin{equation}
        \biggl|\sum_{j = -\infty}^\infty a_j\biggr| \le \|a\|_1.
\end{equation}
As before, it is easy to see that $\|a\|_1$ defines a norm on
$\ell^1({\bf Z})$.

        If $a \in \ell^1({\bf Z})$, $z \in {\bf C}$, and $|z| = 1$, then put
\begin{equation}
        \widehat{a}(z) = \sum_{j = -\infty}^\infty a_j \, z^j,
\end{equation}
which is the Fourier transform of $a$.  This makes sense, because
\begin{equation}
 \sum_{j = -\infty}^\infty |a_j \, z^j| = \sum_{j = -\infty}^\infty |a_j|
\end{equation}
converges.  Moreover,
\begin{equation}
        \sup_{z \in {\bf T}} |\widehat{a}(z)| \le \|a\|_1.
\end{equation}
The partial sums $\sum_{j = -n}^n a_j \, z^j$ are continuous functions
that converge to $\widehat{a}(z)$ uniformly on the unit circle, by
Weierstrass' M-test, and so $\widehat{a}(z)$ is a continuous function
on ${\bf T}$.  It is easy to see that
\begin{equation}
        \widehat{(\widehat{a})}(j)
         = \frac{1}{2 \pi}\int_{\bf T} \widehat{a}(z) \, z^{-j} \, |dz| = a_j
\end{equation}
for each $j$, using the uniform convergence of the partial sums to
reduce to the identities discussed in the previous section.

        The convolution $a * b$ of $a, b \in \ell^1({\bf Z})$ is defined by
\begin{equation}
        (a * b)_j = \sum_{l = -\infty}^\infty a_{j - l} \, b_l.
\end{equation}
The sum on the right converges absolutely as soon as one of $a$, $b$
is summable and the other is bounded, and in particular when both $a$,
$b$ are summable.  We also have that
\begin{equation}
        |(a * b)_j| = \sum_{l = -\infty}^\infty |a_{j - l}| \, |b_l|,
\end{equation}
which implies that
\begin{equation}
        \sum_{j = -\infty}^\infty |(a * b)_j| 
 \le \sum_{j = -\infty}^\infty \sum_{l = -\infty}^\infty |a_{j - l}| \, |b_l|.
\end{equation}
Interchanging the order of summation, we get that
\begin{equation}
        \sum_{j = -\infty}^\infty |(a * b)_j|
 \le \sum_{l = -\infty}^\infty \sum_{j = -\infty}^\infty |a_{j - l}| \, |b_l|.
\end{equation}
Of course,
\begin{equation}
        \sum_{j = -\infty}^\infty |a_{j - l}| = \sum_{j = -\infty}^\infty |a_j|
\end{equation}
for each $l$, by making the change of variables $j \mapsto j + l$.  Thus
\begin{equation}
        \sum_{j = -\infty}^\infty |(a * b)_j|
         \le \Big(\sum_{j = -\infty}^\infty |a_j|\Big)
              \, \Big(\sum_{l = -\infty}^\infty |b_l|\Big),
\end{equation}
so that $a * b \in \ell^1({\bf Z})$ when $a, b \in \ell^1({\bf Z})$.
Equivalently,
\begin{equation}
        \|a * b\|_1 \le \|a\|_1 \, \|b\|_1.
\end{equation}

        If $a, b \in \ell^1({\bf Z})$, $z \in {\bf C}$, and $|z| = 1$, then
\begin{equation}
        \widehat{(a * b)}(z) = \sum_{j = -\infty}^\infty
                     \Big(\sum_{l = -\infty}^\infty a_{j - l} \, b_l\Big) z^j.
\end{equation}
This is the same as
\begin{equation}
 \sum_{j = -\infty}^\infty \sum_{l = -\infty}^\infty a_{j - l} \, z^{j - l}
                                                              \, b_l \, z^l,
\end{equation}
which is equal to
\begin{equation}
 \sum_{l = -\infty}^\infty \sum_{j = -\infty}^\infty a_{j - l} \, z^{j - l}
                                                              \, b_l \, z^l,
\end{equation}
by interchanging the order of summation.  This uses the absolute
summability shown in the previous paragraph.  As before, we can make
the change of variables $j \mapsto j + l$, to get that
\begin{equation}
        \sum_{j = -\infty}^\infty a_{j - l} \, z^{j - l}
         = \sum_{j = -\infty}^\infty a_j \, z^j = \widehat{a}(z)
\end{equation}
for each $l$.  Substituting this into the previous double sum, we get
that
\begin{equation}
\label{widehat{(a * b)}(z) = widehat{a}(z) widehat{b}(z)}
        \widehat{(a * b)}(z) = \widehat{a}(z) \, \widehat{b}(z)
\end{equation}
for every $z \in {\bf T}$.

        Let $\delta(n) = \{\delta_j(n)\}_{j = -\infty}^\infty$ be
defined for each integer $n$ by putting $\delta_j(n) = 1$ when $j = n$
and $\delta_j(n) = 0$ when $j \ne n$, so that $\|\delta(n)\|_1 = 1$
for each $n$.  It is easy to see that
\begin{equation}
        \delta(n) * \delta(r) = \delta(n + r)
\end{equation}
for every $n, r \in {\bf Z}$, and that
\begin{equation}
        \delta(0) * a = a * \delta(0) = a
\end{equation}
for every $a \in \ell^1({\bf Z})$.  One can also check that
\begin{equation}
        a * b = b * a
\end{equation}
and
\begin{equation}
        (a * b) * c = a * (b * c)
\end{equation}
for every $a, b, c \in \ell^1({\bf Z})$, directly from the definition
of convolution, or using the fact that linear combinations of the
$\delta(n)$'s are dense in $\ell^1({\bf Z})$.  It is well known and
not too difficult to show that $\ell^1({\bf Z})$ is complete with
respect to the $\ell^1$ norm $\|a\|_1$.  It follows that $\ell^1({\bf
Z})$ is a commutative Banach algebra, with convolution as
multiplication and $\delta(0)$ as the multiplicative identity element.

        Suppose that $\phi$ is a linear functional on $\ell^1({\bf
Z})$ that is also a homomorphism with respect to convolution, so that
$\phi(a * b) = \phi(a) \, \phi(b)$ for every $a, b \in \ell^1({\bf
Z})$.  If $\phi(a) \ne 0$ for some $a \in \ell^1({\bf Z})$, then
$\phi(\delta(0)) = 1$, and $\phi$ is a continuous linear functional on
$\ell^1({\bf Z})$ with dual norm $1$, as in Section \ref{banach
algebras}.  We would like to show that
\begin{equation}
        \phi(a) = \widehat{a}(z)
\end{equation}
for some $z \in {\bf T}$ and every $a \in \ell^1({\bf Z})$.  Of
course, we have already seen that $\phi_z(a) = \widehat{a}(z)$ defines
a homomorphism on $\ell^1({\bf Z})$ for every $z \in {\bf T}$.

        If $z = \phi(\delta(1))$, then $|z| \le 1$, because
$\|\delta(1)\|_1 = 1$ and $\phi$ has dual norm $1$.  We also know that
$\delta(-1) * \delta(1) = \delta(0)$, which implies that
$\phi(\delta(-1)) \, \phi(\delta(1)) = 1$.  Thus $z \ne 0$, $z^{-1} =
\phi(\delta(-1))$, and hence $|z^{-1}| \le 1$, because
$\|\delta(-1)\|_1 = 1$ and $\phi$ has dual norm $1$.  It follows that
$|z| = 1$, and that $\phi(\delta(n)) = z^n$ for each $n \in {\bf Z}$.
Equivalently, $\phi(a) = \widehat{a}(z)$ when $a = \delta(n)$ for some
$n$.  This also works when $a$ is a finite linear combination of
$\delta(n)$'s, by linearity.  Therefore $\phi(a) = \widehat{a}(z)$ for
every $a \in \ell^1({\bf Z})$, because linear combinations of the
$\delta(n)$'s are dense in $\ell^1({\bf Z})$.

\section{The Poisson kernel}
\label{poisson kernel}
\setcounter{equation}{0}

        Let $f(z)$ be a continuous complex-valued function on the unit
circle ${\bf T}$.  Note that the Fourier coefficients of $f$ are
bounded, with
\begin{equation}
 |\widehat{f}(j)| \le \frac{1}{2 \pi} \int_{\bf T} |f(w)| \, |dw|
                   \le \sup_{w \in {\bf T}} |f(w)|
\end{equation}
for each $j \in {\bf Z}$.  Put
\begin{equation}
 \phi(z) = \sum_{j = 0}^\infty \widehat{f}(j) \, z^j
             + \sum_{j = 1}^\infty \widehat{f}(-j) \, \overline{z}^j
\end{equation}
for each $z \in {\bf C}$ with $|z| < 1$, where $z^j$ is interpreted as
being equal to $1$ for each $z$ when $j = 0$, as usual.  These two
infinite series converge absolutely when $|z| < 1$, because
$\widehat{f}(j)$ is bounded.  If $|z| = 1$, then $\overline{z} =
z^{-1}$, and the sum of these two series is formally the same as the
Fourier series (\ref{sum_{j = -infty}^infty widehat{f}(j) z^j})
associated to $f$.

        Equivalently, $\phi = \phi_1 + \phi_2$, where
\begin{equation}
        \phi_1(z) = \sum_{j = 0}^\infty \widehat{f}(j) \, z^j, \quad
 \phi_2(z) = \sum_{j = 1}^\infty \widehat{f}(-j) \, \overline{z}^j.
\end{equation}
Of course, $\phi_1$ is a holomorphic function on the open unit disk,
and $\phi_2$ is the complex conjugate of a holomorphic function on the
open unit disk.  It is well known that a holomorphic function $h(z)$
is harmonic, meaning that it satisfies Laplace's equation
\begin{equation}
 \frac{\partial^2 h}{\partial x^2} + \frac{\partial^2 h}{\partial y^2} = 0
\end{equation}
when we identify the complex plane ${\bf C}$ with ${\bf R}^2$, and
where $x$, $y$ correspond to the real and imaginary parts of $z \in
{\bf C}$.  More precisely, Laplace's equation applies to the real and
imaginary parts of $h(z)$ separately, both of which are harmonic.
Thus the complex conjugate of a holomorphic function is also harmonic,
and hence $\phi$ is a harmonic function on the open unit disk.

        The Poisson kernel is defined by
\begin{equation}
\label{P(z, w) = ...}
 P(z, w) = \frac{1}{2 \pi} \Big(\sum_{j = 0}^\infty z^j \, \overline{w}^j
                              + \sum_{j = 1}^\infty \overline{z}^j \, w^j\Big)
\end{equation}
for $z, w \in {\bf C}$ with $|z| < 1$ and $|w| = 1$.  Of course, these
series converge absolutely under these conditions, and their partial
sums converge uniformly on the set where $|z| \le r$ and $|w| = 1$ for
every $r < 1$.  This implies that
\begin{equation}
\label{phi(z) = int_{bf T} P(z, w) f(w) |dw|}
        \phi(z) = \int_{\bf T} P(z, w) \, f(w) \, |dw|
\end{equation}
for every $z$ in the open unit disk, using uniform convergence for
$w \in {\bf T}$ to interchange the order of summation and integration.
In particular,
\begin{equation}
\label{int_{bf T} P(z, w) |dw| = 1}
        \int_{\bf T} P(z, w) \, |dw| = 1
\end{equation}
for every $z$ in the open unit disk, because $\phi(z) = 1$ for each
$z$ when $f$ is the constant function equal to $1$ on the unit circle.

        Observe that
\begin{equation}
        \sum_{j = 1}^\infty \overline{z}^j \, w^j
         = \overline{\Big(\sum_{j = 1}^\infty z^j \, \overline{w}^j\Big)},
\end{equation}
and hence
\begin{equation}
 P(z, w) = \frac{1}{2 \pi} \Big(2 \re \sum_{j = 0}^\infty z^j \, \overline{w}^j
                                                                    - 1\Big)
\end{equation}
for all $z$, $w$ as before.  Here $\re a$ denotes the real part of a
complex number $a$, and we are using the simple fact that $a +
\overline{a} = 2 \re a$.  Summing the geometric series, we get that
\begin{equation}
\label{sum_{j = 0}^infty z^j overline{w}^j = ...}
 \sum_{j = 0}^\infty z^j \, \overline{w}^j = \frac{1}{1 - z \, \overline{w}}
                      = \frac{1 - \overline{z} \, w}{|1 - z \, \overline{w}|^2}
\end{equation}
when $|z| < 1$ and $|w| = 1$.  Thus
\begin{equation}
        P(z, w) = \frac{1}{2 \pi} |1 - z \, \overline{w}|^{-2} (2 -
                      2 \re z \, \overline{w} - |1 - z \, \overline{w}|^2).
\end{equation}
We can expand $|1 - z \, \overline{w}|^2$ into $(1 - z \,
\overline{w}) (1 - \overline{z} \, w)$, which reduces to $1 - 2 \re z
\, \overline{w} - |z|^2$ when $|w| = 1$.  It follows that
\begin{equation}
\label{P(z, w) = ..., 2}
 P(z, w) = \frac{1}{2 \pi} \frac{1 - |z|^2}{|1 - z \, \overline{w}|^2}
         =  \frac{1}{2 \pi} \frac{1 - |z|^2}{|w - z|^2},
\end{equation}
using $|w| = 1$ again in the second step.  In particular, $P(z, w) > 0$.

        If $z_0, w \in {\bf T}$ and $z_0 \ne w$, then $P(z, w) \to 0$
as $z \to z_0$, where the limit is restricted to $z$ in the open unit
disk.  This is an immediate consequence of (\ref{P(z, w) = ..., 2}),
which also shows that we have uniform convergence for $w \in {\bf T}$
that satisfy $|w - z_0| \ge \delta$ for some $\delta > 0$.

        Note that
\begin{equation}
        \phi(z) - f(z_0) = \int_{\bf T} P(z, w) \, (f(w) - f(z_0)) \, |dw|
\end{equation}
for every $z_0 \in {\bf T}$ and $z$ in the open unit disk, because of
(\ref{int_{bf T} P(z, w) |dw| = 1}), and hence
\begin{equation}
        |\phi(z) - f(z_0)| \le \int_{\bf T} P(z, w) \, |f(w) - f(z_0)| \, |dw|.
\end{equation}
Using this and the continuity of $f$, one can check that $\phi(z) \to
f(z_0)$ as $z \to z_0$ in the open unit disk.  More precisely, $f(w) -
f(z_0)$ is small when $w$ is close to $z_0$, while $P(z, w)$ is small
when $w$ is not too close to $z_0$ and $z$ is very close to $z_0$.

        It follows that the function defined on the closed unit disk
by taking $\phi$ on the open unit disk and $f$ on the unit circle is
continuous.  In particular, if $\widehat{f}(j) = 0$ when $j < 0$, then
$\phi = \phi_1$ is holomorphic, as mentioned at the end of Section
\ref{fourier series}.

\section{Cauchy products}
\label{cauchy products}
\setcounter{equation}{0}

        If $\sum_{j = 0}^\infty a_j \, z^j$, $\sum_{j = 0}^\infty b_l
\, z^l$ are power series with complex coefficients, then
\begin{equation}
        \Big(\sum_{j = 0}^\infty a_j \, z^j\Big) \,
 \Big(\sum_{l = 0}^\infty b_l \, z^l\Big) = \sum_{n = 0}^\infty c_n \, z^n
\end{equation}
formally, where
\begin{equation}
        c_n = \sum_{j = 0}^n a_j \, b_{n - j}.
\end{equation}
In particular,
\begin{equation}
\label{(sum_{j = 0}^infty a_j) (sum_{l = 0}^infty b_l) = sum_{n = 0}^infty c_n}
 \Big(\sum_{j = 0}^\infty a_j\Big) \, \Big(\sum_{l = 0}^\infty b_l\Big)
        = \sum_{n = 0}^\infty c_n
\end{equation}
formally.  These identities clearly hold when $a_j = b_l = 0$ for all
but finitely many $j$, $l$, for instance.

        If $a_j$, $b_l$ are nonnegative real numbers, then it is easy
to see that
\begin{equation}
        \sum_{n = 0}^N c_n \le \Big(\sum_{j = 0}^N a_j\Big) \,
                                        \Big(\sum_{l = 0}^N b_l\Big)
\end{equation}
for every nonnegative integer $N$.  Similarly,
\begin{equation}
 \Big(\sum_{j = 0}^N a_j\Big) \, \Big(\sum_{l = 0}^N b_l\Big)
        \le \sum_{n = 0}^{2 N} c_n.
\end{equation}
Hence $\sum_{n = 0}^\infty c_n$ converges and satisfies (\ref{(sum_{j
= 0}^infty a_j) (sum_{l = 0}^infty b_l) = sum_{n = 0}^infty c_n}) when
$\sum_{j = 0}^\infty a_j$, $\sum_{l = 0}^\infty b_l$ converge.

        If $a_j$, $b_l$ are arbitrary real or complex numbers, then
\begin{equation}
        |c_n| \le \sum_{j = 0}^n |a_j| \, |b_{n - j}|
\end{equation}
for each $n$.  If $\sum_{j = 0}^\infty a_j$, $\sum_{l = 0}^\infty b_l$
converge absolutely, then it follows that $\sum_{n = 0}^\infty c_n$
converges absolutely too, by the remarks in the previous paragraph.
In this case, one can check that (\ref{(sum_{j = 0}^infty a_j) (sum_{l
= 0}^infty b_l) = sum_{n = 0}^infty c_n}) holds, by expressing these
series as linear combinations of convergent series of nonnegative real
numbers, and using the remarks in the previous paragraph.
Alternatively, one can approximate these series by ones with only
finitely many nonzero terms, and estimate the remainders using
absolute convergence.

        Suppose now that $\sum_{j = 0}^\infty a_j \, z^j$, $\sum_{l =
0}^\infty b_l \, z^l$ are power series that converge when $|z| < 1$,
and hence converge absolutely when $|z| < 1$, by standard results.
Thus $\sum_{n = 0}^\infty c_n \, z^n$ converges absolutely when $|z| <
1$, and is equal to the product of the other two series.  The partial
sums of these series also converge uniformly for $|z| \le r$ when $r <
1$, by standard results.

        Put $f(z) = \sum_{j = 0}^\infty a_j \, z^j$, $g(z) = \sum_{l =
0}^\infty b_l \, z^l$, and $h(z) = \sum_{n = 0}^\infty c_n \, z^n$
when $|z| < 1$, so that
\begin{equation}
        f(z) \, g(z) = h(z),
\end{equation}
as in the preceding paragraph.  If $f(z)$, $g(z)$ have continuous
extensions to the closed unit disk, then it follows that $h(z)$ does
as well.

        Note that
\begin{equation}
        a_j \, r^j = \frac{1}{2 \pi} \int_{\bf T} f(r \, z) \, z^{-j} \, |dz|
\end{equation}
for each $j \ge 0$ and $0 < r < 1$, and similarly for $g$, $h$.  This
is because $f(r \, z)$ is defined by an absolutely convergent Fourier
series, so that we can reduce to the usual identities for the integral
of a power of $z$ on the unit circle by interchaning the order of
integration and summation.  If $f$ extends continuously to the closed
unit disk, then this formula also holds with $r = 1$.

        If $\sum_{j = -\infty}^\infty a_j$, $\sum_{l = -\infty}^\infty
b_l$ are doubly-infinite series of complex numbers, then we have again that
\begin{equation}
\label{sum_{n = -infty}^infty c_n = ...}
        \Big(\sum_{j = -\infty}^\infty a_j\Big) \,
 \Big(\sum_{l = -\infty}^\infty b_l\Big) = \sum_{n = -\infty}^\infty c_n
\end{equation}
with $c_n = \sum_{j = -\infty}^\infty a_j \, b_{n - j}$, and similarly
\begin{equation}
\label{sum_{n = -infty}^infty c_n z^n = ...}
        \Big(\sum_{j = -\infty}^\infty a_j \, z^j\Big) \,
         \Big(\sum_{l = -\infty}^\infty b_l \, z^l\Big)
          = \sum_{n = -\infty}^\infty c_n \, z^n,
\end{equation}
at least formally.  As before, there is no problem with these
identities when $a_j = b_l = 0$ for all but finitely many $j$, $l$.
Otherwise, even the definition of $c_n$ requires some convergence
conditions.  If the $a_j$'s are absolutely summable and the $b_l$'s
are bounded, or vice-versa, then the series defining $c_n$ converges
absolutely, and
\begin{equation}
        |c_n| \le \sum_{j = -\infty}^\infty |a_j| \, |b_{n - j}|
\end{equation}
for each $n$.  If both the $a_j$'s and $b_l$'s are absolutely
summable, then it is easy to see that $c_n$'s are absolutely summable
too, with
\begin{equation}
 \sum_{n = -\infty}^\infty |c_n| \le \Big(\sum_{j = -\infty}^\infty |a_j|\Big)
                                  \, \Big(\sum_{l = -\infty}^\infty |b_l|\Big).
\end{equation}
This follows from the previous estimate for $|c_n|$ by interchanging
the order of summation.  One can also check that (\ref{sum_{n =
-infty}^infty c_n = ...}) holds under these conditions, in the same
way as in the earlier situation for sums over nonnegative integers.
Of course, this implies that (\ref{sum_{n = -infty}^infty c_n z^n =
...})  holds as well when $|z| = 1$, which is basically the same as
(\ref{widehat{(a * b)}(z) = widehat{a}(z) widehat{b}(z)}).

\section{Inner product spaces}
\label{inner products}
\setcounter{equation}{0}

        Let $V$ be a vector space over the real or complex numbers.
An \emph{inner product} on $V$ is a function $\langle v, w \rangle$
defined for $v, w \in V$ with values in ${\bf R}$ or ${\bf C}$, as
appropriate, that satisfies the following three conditions.  First,
\begin{equation}
        \lambda_w(v) = \langle v, w \rangle
\end{equation}
is linear as a function of $v$ for each $w \in V$.  Second,
\begin{equation}
        \langle w, v \rangle = \langle v, w \rangle
\end{equation}
for every $v, w \in V$ in the real case, and
\begin{equation}
        \langle w, v \rangle = \overline{\langle v, w \rangle}
\end{equation}
for every $v, w \in V$ in the complex case.  This implies that
$\langle v, w \rangle$ is linear in $w$ in the real case, and
conjugate-linear in $w$ in the complex case.  It also implies that
\begin{equation}
        \langle v, w \rangle = \overline{\langle v, v \rangle} \in {\bf R}
\end{equation}
for every $v \in V$ in the complex case.  The third condition is that
$\langle v, v \rangle \ge 0$ for every $v \in V$ in both the real and
complex cases, with equality only when $v = 0$.

        Put
\begin{equation}
        \|v\| = \langle v, v \rangle^{1/2}
\end{equation}
for every $v \in V$.  This satisfies the positivity and homogeneity
requirements of a norm, and we would like to show that it also
satisfies the triangle inequality.  Observe that
\begin{eqnarray}
 0 \le \|v + t \, w\|^2 & = & \langle v, v \rangle + t \, \langle v, w \rangle
                      + t \, \langle w, v \rangle + t^2 \langle w, w \rangle \\
      & = & \|v\|^2 + 2 \, t \, \langle v, w \rangle + t^2 \, \|w\|^2 \nonumber
\end{eqnarray}
for every $v, w \in V$ and $t \in {\bf R}$ in the real case, and similarly
\begin{eqnarray}
 0 \le \|v + t \, w\|^2 & = & \langle v, v \rangle + t \, \langle v, w \rangle
      + \overline{t} \, \langle w, v \rangle + |t|^2 \, \langle w, w \rangle \\
 & = & \|v\|^2 + t \, \langle v, w \rangle
+ \overline{t} \, \overline{\langle v, w\rangle} + |t|^2 \, \|w\|^2 \nonumber\\
 & = & \|v\|^2 + 2 \re t \, \langle v, w \rangle + |t|^2 \, \|w\|^2 \nonumber
\end{eqnarray}
for every $v, w \in V$ and $t \in {\bf C}$ in the complex case.  In both
cases, we get that
\begin{equation}
        0 \le \|v\|^2 - 2 \, r \, |\langle v, w \rangle| + r^2 \|w\|^2
\end{equation}
for every $v, w \in V$ and $r \ge 0$, by taking $t = - r \, \alpha$,
where $|\alpha| = 1$ and
\begin{equation}
        \alpha \, \langle v, w \rangle = |\langle v, w \rangle|.
\end{equation}
Equivalently,
\begin{equation}
        2 \, r \, |\langle v, w \rangle| \le \|v\|^2 + r^2 \, \|w\|^2
\end{equation}
for every $v, w \in V$ and $r \ge 0$, and hence
\begin{equation}
 |\langle v, w \rangle| \le \frac{1}{2} \, (r^{-1} \, \|v\|^2 + r \, \|w\|^2)
\end{equation}
when $r > 0$.  If $v, w \ne 0$, then we can take $r = \|v\| / \|w\|$
to get that
\begin{equation}
        |\langle v, w \rangle| \le \|v\| \, \|w\|.
\end{equation}
This is the \emph{Cauchy--Schwarz inequality}, which also holds
trivially when $v = 0$ or when $w = 0$.

        As before,
\begin{equation}
        \|v + w\|^2 = \|v\|^2 + 2 \, \langle v, w \rangle + \|w\|^2
\end{equation}
for every $v, w \in V$ in the real case, and
\begin{equation}
        \|v + w\|^2 = \|v\|^2 + 2 \, \re \langle v, w \rangle + \|w\|^2
\end{equation}
for every $v, w \in V$ in the complex case.  In both case,
\begin{eqnarray}
 \|v + w\|^2 & \le & \|v\|^2 + 2 \, |\langle v, w \rangle| + \|w\|^2 \\
             & \le & \|v\|^2 + 2 \, \|v\| \, \|w\| + \|w\|^2
                                         = (\|v\| + \|w\|)^2, \nonumber
\end{eqnarray}
using the Cauchy--Schwarz inequality in the second step.  This implies that
\begin{equation}
        \|v + w\| \le \|v\| + \|w\|
\end{equation}
for every $v, w \in V$, so that $\|v\|$ defines a norm on $V$, as desired.

        The standard inner products on ${\bf R}^n$ and ${\bf C}^n$ are given by
\begin{equation}
        \langle v, w \rangle = \sum_{j = 1}^n v_j \, w_j
\end{equation}
and
\begin{equation}
        \langle v, w \rangle = \sum_{j = 1}^n v_j \, \overline{w_j},
\end{equation}
respectively.  In both cases, the corresponding norm is given by
\begin{equation}
        \|v\| = \Big(\sum_{j = 1}^n |v_j|^2\Big)^{1/2}.
\end{equation}
This is the standard Euclidean norm on ${\bf R}^n$, ${\bf C}^n$, for
which the corresponding topology is the standard topology.

\section{$\ell^2(E)$}
\label{ell^2(E)}
\setcounter{equation}{0}

        Let $E$ be a nonempty set, and let $\ell^2(E)$ be the space of
real or complex-valued functions $f(x)$ on $E$ such that $|f(x)|^2$ is
a summable function on $E$, as in Section \ref{summable functions}.
As usual, this may also be denoted $\ell^2(E, {\bf R})$ or $\ell^2(E,
{\bf C})$, to indicate whether real or complex-valued functions are
being used.  Remember that
\begin{equation}
        a \, b \le \frac{a^2 + b^2}{2}
\end{equation}
for every $a, b \ge 0$, since
\begin{equation}
        0 \le (a - b)^2 = a^2 - 2 \, a \, b + b^2.
\end{equation}
If $f, g \in \ell^2(E)$, then it follows that
\begin{eqnarray}
 |f(x) + g(x)|^2 & \le & (|f(x)| + |g(x)|)^2 \\
                  & = & |f(x)|^2 + 2 \, |f(x)| \, |g(x)| + |g(x)|^2 \nonumber\\
                 & \le & 2 \, |f(x)|^2 + 2 \, |g(x)|^2 \nonumber
\end{eqnarray}
for every $x \in E$.  Hence $f + g \in \ell^2(E)$, because $|f(x)|^2$,
$|g(x)|^2$ are summable on $E$ by hypothesis.

        Similarly,
\begin{equation}
        |f(x)| \, |g(x)| \le \frac{1}{2} \, |f(x)|^2 + \frac{1}{2} \, |g(x)|^2
\end{equation}
is a summable function on $E$ when $f, g \in \ell^2(E)$.  Put
\begin{equation}
        \langle f, g \rangle = \sum_{x \in E} f(x) \, g(x)
\end{equation}
in the real case, and
\begin{equation}
        \langle f, g \rangle = \sum_{x \in E} f(x) \, \overline{g(x)}
\end{equation}
in the complex case.  Thus
\begin{equation}
        \langle f, f \rangle = \sum_{x \in E} |f(x)|^2
\end{equation}
in both cases.  It is easy to see that $\ell^2(E)$ is a vector space
with respect to pointwise addition and scalar multiplication, and that
$\langle f, g \rangle$ defines an inner product on $\ell^2(E)$.  The
norm associated to this inner product is denoted $\|f\|_2$.

        If $f \in \ell^1(E)$, then $f$ is bounded, and $\|f\|_\infty
\le \|f\|_1$.  This implies that
\begin{equation}
        \sum_{x \in E} |f(x)|^2 \le \|f\|_\infty \, \sum_{x \in E} |f(x)|
                                 = \|f\|_\infty \, \|f\|_1 \le \|f\|_1^2,
\end{equation}
so that $f \in \ell^2(E)$ and
\begin{equation}
        \|f\|_2 \le \|f\|_1.
\end{equation}
Similarly, if $f \in \ell^2(E)$, then $f$ is bounded on $E$, and
\begin{equation}
        \|f\|_\infty \le \|f\|_2.
\end{equation}
One can also check that $f \in c_0(E)$, for the same reasons as for
summable functions, and hence
\begin{equation}
        \ell^1(E) \subseteq \ell^2(E) \subseteq c_0(E).
\end{equation}
As in the case of $\ell^1(E)$, one can show that functions with finite
support on $E$ are dense in $\ell^2(E)$.

        If $(V, \langle v, w \rangle)$ is a real or complex inner
product space, then $\lambda_w(v) = \langle v, w \rangle$ defines a
continuous linear functional on $V$ for every $w \in V$.  This uses
the Cauchy--Schwarz inequality, which implies that the dual norm of
$\lambda_w$ is less than or equal to the norm of $w$.  The dual norm
of $\lambda_w$ is actually equal to the norm of $w$, as one can check
by taking $v = w$.  If $V = \ell^2(E)$ with the inner product defined
before, then one can show that every continuous linear functional is
of this form, using arguments like those in Sections \ref{c_0(E)} and
\ref{dual of ell^1}.  An inner product space $(V, \langle v, w
\rangle)$ is said to be a \emph{Hilbert space} if $V$ is complete as a
metric space with respect to the metric determined by the norm
associated to the inner product.  It is well known that $\ell^2(E)$ is
complete with respect to the $\ell^2$ norm, and hence is a Hilbert
space.  Conversely, it can be shown that every Hilbert space is
isometrically equivalent to $\ell^2(E)$ for some set $E$.  This is
simpler when $V$ is separable, in the sense that it has a countable
dense set, in which case $E$ has only finitely or countably many
elements.  One can also show more directly that every continuous
linear functional on a Hilbert space can be expressed as
$\lambda_w(v)$ for some $w \in V$.

\section{Orthogonality}
\label{orthogonality}
\setcounter{equation}{0}

        Let $(V, \langle v, w \rangle)$ be a real or complex inner
product space.  We say that $v, w \in V$ are \emph{orthogonal} if
\begin{equation}
        \langle v, w \rangle = 0,
\end{equation}
which implies that
\begin{equation}
        \|v + w\|^2 = \|v\|^2 + \|w\|^2.
\end{equation}
A collection of vectors $v_1, \ldots, v_n \in V$ is said to be
\emph{orthonormal} if $v_j$ is orthogonal to $v_l$ when $j \ne l$, and
$\|v_j\| = 1$ for each $j$.  This implies that
\begin{equation}
 \bigg\langle \sum_{j = 1}^n a_j \, v_j, \sum_{l = 1}^n b_l \, v_l \bigg\rangle
        = \sum_{j = 1}^n a_j \, b_j
\end{equation}
for every $a_1, \ldots, a_n, b_1, \ldots, b_n \in {\bf R}$ in the real case,
and
\begin{equation}
 \bigg\langle \sum_{j = 1}^n a_j \, v_j, \sum_{l = 1}^n b_l \, v_l \bigg\rangle
        = \sum_{j = 1}^n a_j \, \overline{b_l}
\end{equation}
for every $a_1, \ldots, a_n, b_1, \ldots, b_n \in {\bf C}$ in the complex case.

        Suppose that $v_1, \ldots, v_n \in V$ are orthonormal, and put
\begin{equation}
        P(v) = \sum_{j = 1}^n \langle v, v_j \rangle v_j
\end{equation}
for each $v \in V$.  Thus $P(v)$ is an element of the linear span of
$v_1, \ldots, v_n$ for each $v \in V$, and $P(v) = v$ when $v$ is in
the linear span of $v_1, \ldots, v_n$.  Moreover,
\begin{equation}
        \langle P(v), v_l \rangle = \langle v, v_l \rangle
\end{equation}
for every $v \in V$ and $l = 1, \ldots, n$, which implies that
\begin{equation}
        \langle v - P(v), v_l \rangle = 0
\end{equation}
for $l = 1, \ldots, n$.  Hence $v - P(v)$ is orthogonal to every
element of the linear span of $v_1, \ldots, v_n$.  In particular,
$v - P(v)$ is orthogonal to $P(v)$, which implies that
\begin{equation}
\label{||v||^2 = ... = ||v - P(v)||^2 + sum_{j = 1}^n |langle v, v_j rangle|^2}
        \|v\|^2 = \|v - P(v)\|^2 + \|P(v)\|^2 
                = \|v - P(v)\|^2 + \sum_{j = 1}^n |\langle v, v_j \rangle|^2.
\end{equation}

        Let $w$ be any element of the linear span of $v_1, \ldots,
v_n$.  Thus $v - P(v)$ is orthogonal to $w$, and hence $v - P(v)$ is
orthogonal to $P(v) - w$.  This implies that
\begin{equation}
\label{||v - w||^2 = ||v - P(v)||^2 + ||P(v) - w||^2 ge ||v - P(v)||^2}
        \|v - w\|^2 = \|v - P(v)\|^2 + \|P(v) - w\|^2 \ge \|v - P(v)\|^2,
\end{equation}
so that $P(v)$ is the element of the linear span of $v_1, \ldots, v_n$
closest to $v$.

        Let $A$ be a nonempty set, and suppose that for each $\alpha
\in A$ we have a vector $v_\alpha \in V$ such that $\|v_\alpha\| = 1$
and $v_\alpha$ is orthogonal to $v_\beta$ when $\beta \in A$ and
$\alpha \ne \beta$.  Thus $v_\alpha$, $\alpha \in A$, is an
orthonormal family of vectors in $V$.  If $v \in V$ and $\alpha_1,
\ldots, \alpha_n$ are distinct elements of $A$, then (\ref{||v||^2 =
... = ||v - P(v)||^2 + sum_{j = 1}^n |langle v, v_j rangle|^2})
implies that
\begin{equation}
        \sum_{j = 1}^n |\langle v, v_{\alpha_j}\rangle|^2 \le \|v\|^2.
\end{equation}
It follows that $\langle v, v_\alpha \rangle$ is an element of
$\ell^2(A)$ as a function of $\alpha$, with
\begin{equation}
        \sum_{\alpha \in A} |\langle v, v_\alpha\rangle|^2 \le \|v\|^2.
\end{equation}
If $v$ is in the closure of the linear span of the $v_\alpha$'s,
$\alpha \in A$, with respect to the norm associated to the inner
product on $V$, then one can check that
\begin{equation}
\label{sum_{alpha in A} |langle v, v_alpha rangle|^2 = ||v||^2}
        \sum_{\alpha \in A} |\langle v, v_\alpha \rangle|^2 = \|v\|^2.
\end{equation}

\section{Parseval's formula}
\label{parseval}
\setcounter{equation}{0}

        Let $C({\bf T})$ be the space of continuous complex-valued
functions on the unit circle.  It is easy to see that
\begin{equation}
\label{langle f, g rangle = frac{1}{2 pi} int_{bf T} f(z) overline{g(z)} |dz|}
        \langle f, g \rangle
              = \frac{1}{2 \pi} \int_{\bf T} f(z) \, \overline{g(z)} \, |dz|
\end{equation}
defines an inner product on $C({\bf T})$, for which the corresponding
norm is given by
\begin{equation}
        \|f\| = \Big(\frac{1}{2 \pi} \int_{\bf T} |f(z)|^2 \, |dz|\Big)^{1/2}.
\end{equation}
As in Section \ref{fourier series}, the functions on ${\bf T}$ of the
form $z^j$, $j \in {\bf Z}$, are orthonormal with respect to this
inner product.  The Fourier coefficients of a continuous function $f$
on ${\bf T}$ can also be expressed as
\begin{equation}
        \widehat{f}(j) = \langle f, z^j \rangle.
\end{equation}
\emph{Parseval's formula} states that
\begin{equation}
        \sum_{j = -\infty}^\infty |\widehat{f}(j)|^2
               = \frac{1}{2 \pi} \int_{\bf T} |f(z)|^2 \, |dz|.
\end{equation}
That the sum on the left is less than or equal to the integral on the
right follows immediately from the orthonormality of $z^j$, $j \in
{\bf Z}$, as in the previous section.  In order to show that equality
holds, it suffices to check that $f$ can be approximated by finite
linear combinations of the $z^j$'s with respect to the norm associated
to the inner product.  In fact, a continuous function $f$ on the unit
circle can be approximated uniformly by a finite linear combinations
of the $z^j$'s, $j \in {\bf Z}$.  To see this, one can use the
function $\phi(z)$ on the open unit disk discussed in Section
\ref{poisson kernel}.  Remember that $\phi$ extends to a continuous
function on the closed unit disk, which is equal to $f$ on the unit
circle.  It follows that $\phi(r \, z)$ converges uniformly to $f(z)$
for $z \in {\bf T}$ as $r \to 1$, because continuous functions on
compact sets are uniformly continuous.  It is easy to see that $\phi(r
\, z)$ can be approximated uniformly on ${\bf T}$ by a finite linear
combination of the $z^j$'s for each $r < 1$, because of the absolute
convergence of the series defining $\phi(r \, z)$ when $r < 1$.  This
implies that $f$ can be approximated uniformly by finite linear
combinations of the $z^j$'s on ${\bf T}$, as desired.

\section{$\ell^p(E)$}
\label{ell^p(E)}
\setcounter{equation}{0}

        Let $E$ be a nonempty set, and let $p$ be a positive real
number.  A real or complex-valued function $f(x)$ on $E$ is said to be
\emph{$p$-summable} if $|f(x)|^p$ is a summable function on $E$.  The
space of $p$-summable functions on $E$ is denoted $\ell^p(E)$, or
$\ell^p(E, {\bf R})$, $\ell^p(E, {\bf C})$ to indicate whether real or
complex-valued functions are being used.  This is consistent with
previous definitions when $p = 1, 2$.

        Observe that
\begin{equation}
 (a + b)^p \le (2 \max (a, b))^p = 2^p \max(a^p, b^p) \le 2^p \, (a^p + b^p)
\end{equation}
for any pair of nonnegative real numbers $a$, $b$.  If $f$, $g$ are
$p$-summable functions on $E$, then it follows that $f + g$ is also
$p$-summable, with
\begin{eqnarray}
 \sum_{x \in E} |f(x) + g(x)|^p & \le & \sum_{x \in E} (|f(x)| + |g(x)|)^p \\
   & \le & 2^p \sum_{x \in E} |f(x)|^2 + 2^p \sum_{x \in E} |g(x)|^p. \nonumber
\end{eqnarray}
This implies that $\ell^p(E)$ is a vector space with respect to
pointwise addition and scalar multiplication over the real or complex
numbers, as appropriate.

        If $f$ is a $p$-summable function on $E$, then we put
\begin{equation}
        \|f\|_p = \Big(\sum_{x \in E} |f(x)|^p\Big)^{1/p}.
\end{equation}
It is easy to see that $f$ vanishes at infinity on $E$, as in the $p =
1$ case.  In particular, $f$ is bounded, and we have that
\begin{equation}
        \|f\|_\infty \le \|f\|_p.
\end{equation}
This implies that $f$ is $q$-summable when $p \le q < \infty$, since
\begin{equation}
 \sum_{x \in E} |f(x)|^q \le \|f\|_\infty^{q - p} \sum_{x \in E} |f(x)|^p.
\end{equation}
More precisely, we get that
\begin{equation}
        \|f\|_q^q \le \|f\|_\infty^{q - p} \, \|f\|_p^p \le \|f\|_p^q,
\end{equation}
and hence
\begin{equation}
\label{||f||_q le ||f||_p}
        \|f\|_q \le \|f\|_p.
\end{equation}

        If $0 < p \le 1$, then
\begin{equation}
        a + b \le (a^p + b^p)^{1/p}
\end{equation}
for every $a, b \ge 0$.  This follows from (\ref{||f||_q le ||f||_p})
with $q = 1$, using a set $E$ with two elements.  Equivalently,
\begin{equation}
        (a + b)^p \le a^p + b^p.
\end{equation}
If $f$, $g$ are $p$-summable functions on $E$, then we get that
\begin{eqnarray}
 \sum_{x \in E} |f(x) + g(x)|^p & \le & \sum_{x \in E} (|f(x)| + |g(x)|)^p \\
  & \le & \sum_{x \in E} |f(x)|^p + \sum_{x \in E} |g(x)|^p. \nonumber
\end{eqnarray}
Thus
\begin{equation}
\label{||f + g||_p^p le ||f||_p^p + ||g||_p^p}
        \|f + g\|_p^p \le \|f\|_p^p + \|g\|_p^p.
\end{equation}
This is a bit better than what we had before, since there is no longer
an extra factor of $2^p$.  Note that $\|f\|_p$ does not satisfy the
ordinary triangle inequality when $0 < p < 1$ and $E$ has at least two
elements, and hence is not a norm on $\ell^p(E)$.  However, $\|f -
g\|_p^p$ defines a metric on $\ell^p(E)$ when $0 < p \le 1$, by
(\ref{||f + g||_p^p le ||f||_p^p + ||g||_p^p}).

\section{Convexity}
\label{comvexity}
\setcounter{equation}{0}

        It is well known that $\phi_p(r) = r^p$ defines a convex
function of $r \ge 0$ when $p \ge 1$.  Therefore
\begin{equation}
        (t \, a + (1 - t) \, b)^p \le t \, a^p + (1 - t) \, b^p
\end{equation}
for every $a, b \ge 0$ and $0 \le t \le 1$ when $p \ge 1$.  In
particular, if we take $t = 1/2$, then we get that
\begin{equation}
        (a + b)^p \le 2^{p - 1} \, (a^p + b^p).
\end{equation}
This improves an inequality in the previous section by a factor of $2$.

        If $f$, $g$ are $p$-summable functions on a set $E$, $0 \le t
\le 1$, and $p \ge 1$, then it follows that
\begin{eqnarray}
 \sum_{x \in E} |t \, f(x) + (1 - t) \, g(x)|^p
          & \le & \sum_{x \in E} (t \, |f(x)| + (1 - t) \, |g(x)|)^p \\
 & \le & t \sum_{x \in E} |f(x)|^p + (1 - t) \sum_{x \in E} |g(x)|^p. \nonumber
\end{eqnarray}
Equivalently,
\begin{equation}
\label{||t f + (1 - t) g||_p^p le t ||f||_p^p + (1 - t) ||g||_p^p}
 \|t \, f + (1 - t) \, g\|_p^p \le t \, \|f\|_p^p + (1 - t) \, \|g\|_p^p.
\end{equation}
\emph{Minkowski's inequality} states that
\begin{equation}
        \|f + g\|_p \le \|f\|_p + \|g\|_p
\end{equation}
for every $f, g \in \ell^p(E)$ when $p \ge 1$.  This implies that
$\|f\|_p$ is a norm on $\ell^p(E)$ when $p \ge 1$, because $\|f\|_p$
satisfies the positivity and homogeneity conditions of a norm for
every $p > 0$.

        To prove Minkowski's inequality, we may as well suppose that
neither $f$ nor $g$ is identically $0$ on $E$, since it is trivial
otherwise.  Put $f' = f/\|f\|_p$, $g' = g/\|g\|_p$, so that $\|f'\|_p
= \|g'\|_p = 1$.  Thus
\begin{equation}
\label{||t f' + (1 - t) g'||_p le 1}
        \|t \, f' + (1 - t) \, g'\|_p \le 1
\end{equation}
when $0 \le t \le 1$, by (\ref{||t f + (1 - t) g||_p^p le t ||f||_p^p
+ (1 - t) ||g||_p^p}).  If
\begin{equation}
\label{t = frac{||f||_p}{(||f||_p + ||g||_p)}}
        t = \frac{\|f\|_p}{(\|f\|_p + \|g\|_p)},
\end{equation}
then $1 - t = \|g\|_p / (\|f\|_p + \|g\|_p)$, and Minkowski's
inequality follows from (\ref{||t f' + (1 - t) g'||_p le 1}).

        Remember that a subset $A$ of a vector space $V$ is said to be
\emph{convex} if
\begin{equation}
        t \, v + (1 - t) \, w \in A
\end{equation}
for every $v, w \in A$ and $0 \le t \le 1$.  If $N(v)$ is a seminorm
on $V$, then it is easy to see that the corresponding closed unit ball
\begin{equation}
\label{B = {v in V : N(v) le 1}}
        B = \{v \in V : N(v) \le 1\}
\end{equation}
is a convex set in $V$.  Conversely, if a nonnegative real-valued
function $N(v)$ on $V$ satisfies the homogeneity condition of a
seminorm and $B$ is convex, then one can check $N(v)$ is a seminorm on
$V$.  This is basically the same as the argument in the previous
paragraph for $\|f\|_p$, at least when $N(v)$ satisfies the positivity
condition of a norm.  Otherwise, some minor adjustments are needed to
deal with $v \in V$ such that $N(v) = 0$ but $v \ne 0$.

\section{H\"older's inequality}
\label{holder's inequality}
\setcounter{equation}{0}

        Let $1 < p, q < \infty$ be conjugate exponents, in the sense that
\begin{equation}
        \frac{1}{p} + \frac{1}{q} = 1.
\end{equation}
If $E$ is a nonempty set, $f \in \ell^p(E)$, and $g \in \ell^q(E)$,
then \emph{H\"older's inequality} states that $f \, g \in \ell^1(E)$, and
\begin{equation}
        \|f \, g\|_1 \le \|f\|_p \, \|g\|_q.
\end{equation}
This also works when $p = 1$ and $q = \infty$, or the other way
around, and is much simpler.  The $p = q = 2$ case can be reduced to
the Cauchy--Schwarz inequality.

        Using the convexity of the exponential function, one can check that
\begin{equation}
        a \, b \le \frac{a^p}{p} + \frac{b^q}{q}
\end{equation}
for every $a, b \ge 0$.  Applying this to $a = |f(x)|$, $b = |g(x)|$,
and summing over $x \in E$, we get that
\begin{equation}
\label{sum_{x in E} |f(x)| |g(x)| le ...}
 \sum_{x \in E} |f(x)| \, |g(x)| \le p^{-1} \sum_{x \in E} |f(x)|^p
                                      + q^{-1} \sum_{x \in E} |g(x)|^q.
\end{equation}
In particular, $f \, g \in \ell^1(E)$, and
\begin{equation}
\label{||f g||_1 le p^{-1} ||f||_p^p + q^{-1} ||g||_q^q}
        \|f \, g\|_1 \le p^{-1} \, \|f\|_p^p + q^{-1} \, \|g\|_q^q,
\end{equation}
which implies H\"older's inequality in the special case where $\|f\|_p
= \|g\|_q = 1$.  If $f$ and $g$ are not identically $0$ on $E$, then
one can reduce to this case, by considering $f' = f / \|f\|_p$, $g' =
g/\|g\|_q$.  Otherwise, if $f$ or $g$ is identically $0$ on $E$, then
the result is trivial.

        If $f \in \ell^p(E)$, $g \in \ell^q(E)$, then put
\begin{equation}
        \lambda_g(f) = \sum_{x \in E} f(x) \, g(x).
\end{equation}
H\"older's inequality implies that
\begin{equation}
        |\lambda_g(f)| \le \|f\|_p \, \|g\|_q,
\end{equation}
so that $\lambda_g(f)$ defines a continuous linear functional on
$\ell^p(E)$ for each $g \in \ell^q(E)$, with dual norm less than or
equal to $\|g\|_q$.  One can check that the dual norm of $\lambda$ on
$\ell^p(E)$ is actually equal to $\|g\|_q$, by choosing $g$ such that
\begin{equation}
\label{f(x) g(x) = |f(x)|^p = |g(x)|^q}
        f(x) \, g(x) = |f(x)|^p = |g(x)|^q
\end{equation}
for every $x \in E$.  These conditions on $g$ are consistent with each
other, because $p$ and $q$ are conjugate exponents.

        Conversely, if $\lambda$ is a continuous linear functional on
$\ell^p(E)$, then one can show that $\lambda = \lambda_g$ for some $g
\in \ell^q(E)$.  As usual, one can start by putting $g(x) =
\lambda(\delta_x)$, where $\delta_x$ is the function on $E$ equal to
$1$ at $x$ and to $0$ elsewhere.  This permits $\lambda_g(f)$ to be
defined as in the previous paragraph when $f$ has finite support on
$E$, in which cas it agrees with $\lambda(f)$, by linearity.  The next
step is to show that
\begin{equation}
        \Big(\sum_{x \in A} |g(x)|^q\Big)^{1/q}
\end{equation}
is bounded by the dual norm of $\lambda$ on $\ell^p(E)$ when $A$ is a
finite subset of $E$.  This can be done by choosing $f$ such that
(\ref{f(x) g(x) = |f(x)|^p = |g(x)|^q}) holds when $x \in A$, and
$f(x) = 0$ when $x \in E \backslash A$.  This implies that $g \in
\ell^q(E)$, and that $\|g\|_q$ is less than or equal to the dual norm
of $\lambda$ on $\ell^p(E)$.  The remaining point is that $\lambda(f)
= \lambda_g(f)$ for every $f \in \ell^p(E)$.  We already know that
this holds when $f$ has finite support on $E$, which implies that it
holds for every $f \in \ell^p(E)$, because functions with finite
support are dense in $\ell^p(E)$, and because $\lambda$ and
$\lambda_g$ are continuous on $\ell^p(E)$.

\section{$p < 1$}
\label{p < 1}
\setcounter{equation}{0}

        Let $E$ be a nonempty set, and let $p$ be a positive real
number strictly less than $1$.  As in Section \ref{ell^p(E)},
\begin{equation}
        d_p(f, g) = \|f - g\|_p^p
\end{equation}
defines a metric on $\ell^p(E)$.  It is easy to see that addition and
scalar multiplication are continuous with respect to the topology
associated to this metric, so that $\ell^p(E)$ becomes a topological
vector space.  If $E$ has only finitely many elements, then
$\ell^p(E)$ can be identified with ${\bf R}^n$ or ${\bf C}^n$, as
appropriate, where $n$ is the number of elements of $E$, and the
topology on $\ell^p(E)$ determined by this metric corresponds exactly to
the standard topology on ${\bf R}^n$ or ${\bf C}^n$.

        If $f \in \ell^p(E)$ and $g \in \ell^\infty(E)$, then $f \, g
\in \ell^p(E) \subseteq \ell^1(E)$, and we can put
\begin{equation}
        \lambda_g(f) = \sum_{x \in E} f(x) \, g(x).
\end{equation}
Moreover,
\begin{equation}
        |\lambda_g(f)| \le \|f\|_1 \, \|g\|_\infty \le \|f\|_p \, \|g\|_\infty.
\end{equation}
Using this estimate, it is easy to see that $\lambda_g$ is a
continuous linear functional on $\ell^p(E)$ with respect to the
topology associated to the metric defined in the previous paragraph.

        Conversely, suppose that $\lambda$ is a continuous linear functional
on $\ell^p(E)$.  This implies that there is a $\delta > 0$ such that
\begin{equation}
        |\lambda(f)| \le 1
\end{equation}
for all $f \in \ell^p(E)$ such that $d_p(f, 0) = \|f\|_p^p < \delta$.
Equivalently, there is a $C \ge 0$ such that
\begin{equation}
        |\lambda(f)| \le C \, \|f\|_p
\end{equation}
for every $f \in \ell^p(E)$, because of linearity.  Put $g(x) =
\lambda(\delta_x)$ for each $x \in E$, where $\delta_x$ is the
function on $E$ equal to $1$ at $x$ and to $0$ elsewhere.  Thus
$|g(x)| \le C$ for every $x \in E$, because $\|\delta_x\|_p = 1$.
This permits us to define $\lambda_g$ as in the preceding paragraph.
By construction, $\lambda(f) = \lambda_g(f)$ when $f$ has finite
support on $E$.  It is easy to see that functions with finite support
on $E$ are dense in $\ell^p(E)$, for basically the same reasons as
when $1 \le p < \infty$.  Hence $\lambda(f) = \lambda_g(f)$ for every
$f \in \ell^p(E)$, since $\lambda$, $\lambda_g$ are both continuous on
$\ell^p(E)$.

        If $E$ has at least two elements, then the unit ball in
$\ell^p(E)$ is not convex, unlike the situation when $p \ge 1$.  If
$E$ has infinitely many elements, then the convex hull of the unit
ball in $\ell^p(E)$ is not even bounded with respect to $\|f\|_p$,
since it contains all functions $f$ on $E$ with finite support such
that $\|f\|_1 \le 1$, for instance.  However, if $f, g \in \ell^p(E)$,
$0 \le t \le 1$, and $h$ is another function on $E$ that satisfies
\begin{equation}
        |h(x)| \le |f(x)|^t \, |g(x)|^{1 - t}
\end{equation}
for every $x \in E$, then $h \in \ell^p(E)$, and
\begin{equation}
        \|h\|_p \le \|f\|_p^t \, \|g\|_p^{1 - t}.
\end{equation}
This follows from H\"older's inequality, and works for all $p > 0$.
In particular, $\|h\|_p \le 1$ when $\|f\|_p, \|g\|_p \le 1$, which is
a multiplicative convexity property of the unit ball in $\ell^p(E)$.

\section{Bounded linear mappings, revisited}
\label{bounded linear mappings, revisited}
\setcounter{equation}{0}

        Let $V$ be a real or complex vector space with a norm
$\|v\|_V$, and consider the space $\mathcal{BL}(V) = \mathcal{BL}(V,
V)$ of bounded linear mappings from $V$ into itself.  This is an
associative algebra, with composition of linear operators as
multiplication, and the identity operator $I$ on $V$ as the
multiplicative identity element.  Note that $\|I\|_{op} = 1$, except
in the trivial case where $V$ consists of only the zero element.  If
$V$ is complete, then $\mathcal{BL}(V)$ is also complete with respect
to the operator norm, as in Section \ref{bounded linear mappings}.
Thus $\mathcal{BL}(V)$ is a Banach algebra when $V$ is a Banach space
and $V \ne \{0\}$.  If $V$ is finite-dimensional, then
$\mathcal{BL}(V)$ is the same as the algebra of all linear
transformations on $V$.  In particular, $\mathcal{BL}(V)$ is not
commutative when the dimension of $V$ is greater than or equal to $2$.
This includes the case where $V$ is infinite-dimensional, since the
Hahn--Banach theorem may be used to get plenty of bounded linear
operators on $V$ with finite rank.

        As an example, let $V$ be the space of real or complex-valued
continuous functions on $[0, 1]$, equipped with the supremum norm.  If
$f$ is a continuous function on $[0, 1]$, then let $T(f)$ be the
function defined on $[0, 1]$ by
\begin{equation}
\label{T(f)(x) = int_0^x f(y) dy}
        T(f)(x) = \int_0^x f(y) \, dy.
\end{equation}
Note that $T(f)$ is continuously-differentiable on $[0, 1]$, with
derivative equal to $f$.  In particular, $T(f)$ is continuous on $[0,
1]$.  Moreover,
\begin{equation}
 |T(f)(x)| \le \int_0^x |f(y)| \, dy \le \int_0^1 |f(y)| \, dy \le \|f\|_{sup}
\end{equation}
for every $f \in C([0, 1])$ and $x \in [0, 1]$, which implies that
\begin{equation}
        \|T(f)\|_{sup} \le \int_0^1 |f(y)| \, dy \le \|f\|_{sup}.
\end{equation}
It follows that $T$ is a bounded linear mapping from $C([0, 1])$ into
itself, with operator norm less than or equal to $1$.  It is easy to
see that $\|T\|_{op} = 1$, by considering the case where $f$ is the
constant function equal to $1$ on $[0, 1]$.

        Let $n$ be a positive integer, and let $T^n = T \circ \cdots
\circ T$ be the $n$-fold composition of $T$.  This can be expressed by
the $n$-fold integral
\begin{equation}
\label{T^n(f)(x) = ...}
 T^n(f)(x) = \int_0^x \int_0^{y_n} \cdots \int_0^{y_2} f(y_1) \,
                                           dy_1 \cdots dy_{n - 1} \, dy_n.
\end{equation}
Thus
\begin{eqnarray}
 |T^n(f)(x)| & \le & \int_0^x \int_0^{y_n} \cdots \int_0^{y_2} |f(y_1)| \,
                                           dy_1 \cdots dy_{n - 1} \, dy_n \\
             & \le & \int_0^1 \int_0^{y_n} \cdots \int_0^{y_2} |f(y_1)| \,
                                  dy_1 \cdots dy_{n - 1} \, dy_n. \nonumber
\end{eqnarray}
If
\begin{equation}
        \sigma(n) = \int_0^1 \int_0^{y_n} \cdots \int_0^{y_2} dy_1
                                            \cdots dy_{n - 1} \, dy_n,
\end{equation}
then we get that
\begin{equation}
        \|T^n(f)\|_{sup} \le \sigma(n) \, \|f\|_{\sup}.
\end{equation}
This shows that the operator norm of $T^n$ on $C([0, 1])$ is less than
or equal to $\sigma(n)$, and it is again easy to see that
$\|T^n\|_{op} = \sigma(n)$, by considering the case where $f$ is the
constant function equal to $1$ on $[0, 1]$.

        In fact, if ${\bf 1}_{[0, 1]}$ denotes the constant function
equal to $1$ on $[0, 1]$, then it is easy to check that
\begin{equation}
        T^n({\bf 1}_{[0, 1]})(x) = \frac{x^n}{n!},
\end{equation}
using induction on $n$.  In particular,
\begin{equation}
\label{sigma(n) = T^n({bf 1}_{[0, 1]})(1) = frac{1}{n!}}
        \sigma(n) = T^n({\bf 1}_{[0, 1]})(1) = \frac{1}{n!}.
\end{equation}
Alternatively, $\sigma(n)$ is the same as the $n$-dimensional volume
of the $n$-dimensional simplex
\begin{equation}
\label{def of Sigma(n)}
        \Sigma(n) = \{y \in {\bf R}^n :
                 0 \le y_1 \le y_2 \le \cdots \le y_{n - 1} \le y_n \le 1\}.
\end{equation}
That the volume of $\Sigma(n)$ is equal to $1/n!$ can also be seen
geometrically, by decomposing the unit cube in ${\bf R}^n$ into $n!$
copies of $\Sigma(n)$ with disjoint interiors.  These copies of
$\Sigma(n)$ are obtained by permuting the standard coordinates of
${\bf R}^n$, using the $n!$ permutations on the set $\{1, \ldots,
n\}$.  Each copy of $\Sigma(n)$ has the same $n$-dimensional volume as
$\Sigma(n)$, and the intersection of any two distinct copies has
measure $0$.  Thus the sum of the volumes of all of these copies of
$\Sigma(n)$ is equal to $n!$ times the volume of $\Sigma(n)$, and is
also equal to the volume of the unit cube, which is equal to $1$.

        Observe that $n! \ge k^{n - k + 1}$ for each positive integer
$k$ when $n \ge k$, so that
\begin{equation}
\label{(n!)^{-1/n} le k^{(k - 1)/n - 1}}
        (n!)^{-1/n} \le k^{(k - 1)/n - 1}
\end{equation}
when $n \ge k$.  In particular,
\begin{equation}
        (n!)^{-1/n} \le k^{-1/2}
\end{equation}
when $n \ge 2 k$, which implies that
\begin{equation}
        \lim_{n \to \infty} (n!)^{-1/n} = 0,
\end{equation}
since the previous statement works for every positive integer $k$.
It follows that
\begin{equation}
        \lim_{n \to \infty} \|T^n\|_{op}^{1/n} = 0,
\end{equation}
because $\|T^n\|_{op} = \sigma(n) = 1/n!$.

        Equivalently, this shows that $r(T) = 0$, in the notation of
Section \ref{spectral radius}.  This would be trivial if $T^n = 0$ for
some positive integer $n$, which is clearly not the case in this
example.

        Let $\mathcal{A}$ be an associative algebra over the real or
complex numbers with a multiplicative identity element, such as the
algebra of bounded linear operators on a vector space with a norm.  If
$x \in \mathcal{A}$, then let $\mathcal{A}(x)$ be the subalgebra of
$\mathcal{A}(x)$ generated by $x$, consisting of linear combinations
of the multiplicative identity element and positive powers of $x$.  It
is easy to see that this is a commutative subalgebra of $\mathcal{A}$,
even if $\mathcal{A}$ is not commutative.  If $\mathcal{A}$ is a
topological algebra, then the closure of a commutative subalgebra of
$\mathcal{A}$ is also commutative.  If $\mathcal{A}$ is a Banach
algebra, then closed subalgebras of $\mathcal{A}$ are Banach algebras
too.

\section{Involutions}
\label{involutions}
\setcounter{equation}{0}

        Let $\mathcal{A}$ be an associative algebra over the real or
complex numbers.  A mapping
\begin{equation}
\label{x mapsto x^*}
        x \mapsto x^*
\end{equation}
on $\mathcal{A}$ is said to be an \emph{involution} if it satisfies
the following three conditions.  First, (\ref{x mapsto x^*}) should be
linear in the real case, and conjugate-linear in the complex case.
This means that
\begin{equation}
        (x + y)^* = x^* + y^*
\end{equation}
for every $x, y \in \mathcal{A}$ in both cases,
\begin{equation}
        (t \, x)^* = t \, x^*
\end{equation}
for every $x \in \mathcal{A}$ and $t \in {\bf R}$ in the real case, and
\begin{equation}
        (t \, x)^* = \overline{t} \, x^*
\end{equation}
in the complex case.  Second, (\ref{x mapsto x^*}) should be
compatible with multiplication in $\mathcal{A}$, in the sense that
\begin{equation}
\label{(x y)^* = y^* x^*}
        (x \, y)^* = y^* \, x^*
\end{equation}
for every $x, y \in \mathcal{A}$.  Of course, (\ref{(x y)^* = y^*
x^*}) is the same as
\begin{equation}
\label{(x y)^* = x^* y^*}
        (x \, y)^* = x^* \, y^*
\end{equation}
when $\mathcal{A}$ is commutative.  The third condition is that
\begin{equation}
        (x^*)^* = x
\end{equation}
for every $x \in \mathcal{A}$.  In particular, this implies that
(\ref{x mapsto x^*}) is a one-to-one mapping of $\mathcal{A}$ onto
itself.  If $\mathcal{A}$ has a multiplicative identity element $e$,
then it follows from the multiplicativity condition (\ref{(x y)^* =
y^* x^*}) that
\begin{equation}
\label{e^* = e}
        e^* = e.
\end{equation}
If $\mathcal{A}$ is equipped with a norm, then one normally asks also
that the involution be isometric, so that
\begin{equation}
\label{||x^*|| = ||x||}
        \|x^*\| = \|x\|
\end{equation}
for every $x \in \mathcal{A}$.

        If $\mathcal{A}$ is the algebra of continuous complex-valued
functions on a topological space, then
\begin{equation}
\label{f(p) mapsto overline{f(p)}}
        f(p) \mapsto \overline{f(p)}
\end{equation}
defines an involution on $\mathcal{A}$.  This would not work for
holomorphic functions, because the complex-conjugate of a holomorphic
function $f$ is also holomorphic only when $f$ is constant.  If
$\mathcal{A}$ is the algebra of $n \times n$ matrices of real numbers
with respect to matrix multiplication, then the transpose of a matrix
defines an involution on $\mathcal{A}$.  If instead $\mathcal{A}$ is
the algebra of $n \times n$ matrices of complex numbers with respect
to matrix multiplication, then one can get an involution on
$\mathcal{A}$ by taking the complex conjugates of the entries of the
transpose of a matrix.  

        If $(V, \langle v, w \rangle)$ is a real or complex Hilbert
space and $T$ is a bounded linear operator on $V$, then it is well
known that there is a unique bounded linear operator $T^*$ on $V$ such that
\begin{equation}
        \langle T(v), w \rangle = \langle v, T^*(w) \rangle
\end{equation}
for every $v, w \in V$, known as the \emph{adjoint} of $T$.  It is
easy to see that this defines an involution on the algebra
$\mathcal{BL}(V)$ of bounded linear operators on $V$.  The adjoint of
$T$ corresponds exactly to the transpose of a real matrix or the
complex conjugate of the transpose of a complex matrix when $T$ is
represented by a matrix with respect to an orthonormal basis for $V$.

        Using the definition of the norm associated to an inner
product and the Cauchy--Schwarz inequality, one can check that
\begin{equation}
 \|T\|_{op} = \sup \{|\langle T(v), w \rangle| : v, w \in V,
                                                 \, \|v\|, \|w\| \le 1\}
\end{equation}
for every bounded linear operator $T$ on $V$.  This implies that
\begin{equation}
        \|T^*\|_{op} = \|T\|_{op}
\end{equation}
for every $T \in \mathcal{BL}(V)$, using the symmetry properties of
the inner product and interchanging the roles of $v$ and $w$ in the
previous expression for the operator norm of $T^*$.  Moreover,
\begin{equation}
\label{||T^* circ T||_{op} = ||T||_{op}^2}
        \|T^* \circ T\|_{op} = \|T\|_{op}^2.
\end{equation}
Of course,
\begin{equation}
        \|T^* \circ T\|_{op} \le \|T^*\|_{op} \, \|T\|_{op} = \|T\|_{op}^2,
\end{equation}
and so it suffices to show the opposite inequality.  Observe that
\begin{equation}
 \langle (T^*(T(v)), v \rangle = \langle T(v), T(v) \rangle = \|T(v)\|^2,
\end{equation}
by the definition of the adjoint operator $T^*$.  This implies that
\begin{equation}
 \|T(v)\|^2 \le \|(T^*(T(v))\| \, \|v\| \le \|T^* \circ T\|_{op} \, \|v\|^2,
\end{equation}
by the Cauchy--Schwarz inequality and the definition of the operator norm.
Thus
\begin{equation}
        \|T\|_{op}^2 \le \|T^* \circ T\|_{op},
\end{equation}
as desired.

        A Banach algebra $(\mathcal{A}, \|x\|)$ equipped with an
isometric involution $x \mapsto x^*$ is said to be a \emph{$C^*$ algebra} if
\begin{equation}
\label{||x^* x|| = ||x||^2}
        \|x^* \, x\| = \|x\|^2
\end{equation}
for every $x \in \mathcal{A}$.  This includes the algebras of bounded
linear operators on real or complex Hilbert spaces, as in the previous
paragraphs.  This also includes the algebra of real or complex-valued
bounded continuous functions on a topological space $X$ with respect
to the supremum norm, where the involution is given by complex
conjugation as in (\ref{f(p) mapsto overline{f(p)}}) in the complex
case, and by the identity operator in the real case.  The same
involutions are defined and isometric on the algebras of real and
complex-valued continuously-differentiable functions on the unit
interval, as in Section \ref{C^1 functions}, but the $C^1$ norm
does not satisfy the $C^*$ condition (\ref{||x^* x|| = ||x||^2}).

        Suppose that $\tau$ is a continuous involution on a
topological space $X$, which is to say a continuous mapping from $X$
into itself such that
\begin{equation}
        \tau(\tau(p)) = p
\end{equation}
for every $p \in X$.  Equivalently, $\tau$ is its own inverse, and
hence a homeomorphism from $X$ onto itself.  Under these conditions,
\begin{equation}
        f(p) \mapsto f(\tau(p))
\end{equation}
is an involution on the algebra of real-valued continuous functions on $X$, and
\begin{equation}
        f(p) \mapsto \overline{f(\tau(p))}
\end{equation}
is an involution on the algebra of complex-valued continuous functions
on $X$.  These involutions also preserve the supremum norms of bounded
continuous functions on $X$.  However, the $C^*$ condition (\ref{||x^*
x|| = ||x||^2}) does not work when $\tau$ is not the identity mapping
on $X$, at least when $X$ is sufficiently regular to have enough
continuous functions.

        As a variant of this, let $U$ be the open disk in the complex
plane.  If $f(z)$ is a holomorphic function on $U$, then it is well
known that
\begin{equation}
\label{overline{f(overline{z})}}
        \overline{f(\overline{z})}
\end{equation}
is also holomorphic on $U$.  It is easy to see that this defines an
involution on the algebra of holomorphic functions on $U$, which
preserves the supremum norm of bounded holomorphic functions on $U$.
However, if
\begin{equation}
        f(z) = z + i,
\end{equation}
then the supremum norm of $f$ on $U$ is equal to $2$, and the supremum
norm of
\begin{equation}
\label{overline{f(overline{z})} f(z) = (z - i) (z + i) = z^2 + 1}
        \overline{f(\overline{z})} \, f(z) = (z - i) (z + i) = z^2 + 1
\end{equation}
is equal to $2$ as well.  Thus the $C^*$ condition (\ref{||x^* x|| =
||x||^2}) does not work in this case either, when we restrict our
attention to bounded holomorphic functions on $U$, since the supremum
norm of (\ref{overline{f(overline{z})} f(z) = (z - i) (z + i) = z^2 +
1}) on $U$ is strictly less than the square of the supremum norm of $f$.

        Let $(\mathcal{A}, \|x\|, x^*)$ be a real or complex $C^*$
algebra, and suppose that $x \in \mathcal{A}$ satisfies
\begin{equation}
        x^* = x.
\end{equation}
In this case, the $C^*$ condition (\ref{||x^* x|| = ||x||^2}) reduces to
\begin{equation}
        \|x^2\| = \|x\|^2.
\end{equation}
If $l$ is a positive integer, then
\begin{equation}
        (x^l)^* = (x^*)^l = x^l,
\end{equation}
and so we can apply the previous statement to $x^l$ to get that
\begin{equation}
        \|x^{2 l}\| = \|x^l\|^2.
\end{equation}
Applying this repeatedly, we get that
\begin{equation}
        \|x^{2^n}\| = \|x\|^{2^n}
\end{equation}
for each positive integer $n$.  Of course,
\begin{equation}
        \|x^l\| \le \|x\|^l
\end{equation}
for any positive integer $n$, by the submultiplicative property of the norm.
If we choose a positive integer $n$ such that $l \le 2^n$, then we get that
\begin{equation}
        \|x\|^{2^n} = \|x^{2^n}\| \le \|x^l\| \, \|x\|^{2^n - l},
\end{equation}
using the submultiplicative property of the norm again.  This implies that
\begin{equation}
        \|x\|^l \le \|x^l\|,
\end{equation}
and hence that
\begin{equation}
        \|x^l\| = \|x\|^l
\end{equation}
for each positive integer $l$.

        If $y$ is any element of $\mathcal{A}$, then $x = y^* \, y$
satisfies $x^* = x$.  Thus we get
\begin{equation}
        \|(y^* \, y)^l\| = \|y^* \, y\|^l = \|y\|^{2 \, l}
\end{equation}
for each positive integer $l$.  Suppose that $y^*$ commutes with $y$,
so that
\begin{equation}
        (y^* \, y)^l = (y^*)^l \, y^l = (y^l)^* \, y^l
\end{equation}
for each $l$, and hence
\begin{equation}
        \|(y^* \, y)^l\| = \|(y^l)^* \, y^l\| = \|y^l\|^2.
\end{equation}
This implies that
\begin{equation}
        \|y^l\| = \|y\|^l
\end{equation}
for each positive integer $l$, as before.

\part{Several variables}

\section{Power series}
\label{power series}
\setcounter{equation}{0}

        Let $n$ be a positive integer, and let
\begin{equation}
        \sum_\alpha a_\alpha \, z^\alpha
\end{equation}
be a power series in $n$ complex variables.  More precisely, the sum
is taken over all multi-indices $\alpha = (\alpha_1, \ldots,
\alpha_n)$, $z^\alpha = z_1^{\alpha_1} \cdots z_n^{\alpha_n}$ is the
corresponding monomial, and the coefficients $a_\alpha$ are complex
numbers.  Let $A$ be the set of $z = (z_1, \ldots, z_n) \in {\bf C}^n$
for which this series converges absolutely, in the sense that
$a_\alpha \, z^\alpha$ is a summable function of $\alpha$ on the set
of multi-indices.  Thus $0 \in A$ trivially, and $w \in A$ whenever
there is a $z \in A$ such that $|w_j| \le |z_j|$ for $j = 1, \ldots,
n$, by the comparison test.

        Let $\sum_\alpha b_\alpha \, z^\alpha$ be another power
series, and let $B$ be the set of $z \in {\bf C}^n$ on which this
series converges absolutely, as before.  Note that
\begin{equation}
        \sum_\alpha (a_\alpha + b_\alpha) \, z^\alpha
\end{equation}
converges absolutely for every $z \in A \cap B$.  The product of these
two power series can be expressed formally as
\begin{equation}
\label{(sum a_alpha z^alpha) (sum b_beta z^beta) = sum c_gamma z^gamma}
        \Big(\sum_\alpha a_\alpha \, z^\alpha\Big) \,
 \Big(\sum_\beta b_\beta \, z^\beta\Big) = \sum_\gamma c_\gamma \, z^\gamma,
\end{equation}
where
\begin{equation}
\label{c_gamma = sum_{alpha + beta = gamma} a_alpha b_beta}
        c_\gamma = \sum_{\alpha + \beta = \gamma} a_\alpha \, b_\beta.
\end{equation}
More precisely, the sum on the right is taken over all multi-indices
$\alpha$, $\beta$ such that $\alpha + \beta = \gamma$, of which there
are only finitely many.  If $z \in A \cap B$, then one can check that
$\sum_\gamma c_\gamma \, z^\gamma$ converges absolutely, and that the
sum satisfies (\ref{(sum a_alpha z^alpha) (sum b_beta z^beta) = sum
c_gamma z^gamma}).  As a first step, one can verify that
\begin{equation}
\label{sum |c_gamma| |z^gamma| le ...}
        \sum_\gamma |c_\gamma| \, |z^\gamma|
           \le \Big(\sum_\alpha |a_\alpha| \, |z^\alpha|\Big) \,
              \Big(\sum_\beta |b_\beta| \, |z^\beta|\Big),
\end{equation}
by estimating the sum over finitely many $\gamma$'s in terms of the
product of sums over finitely many $\alpha$'s and $\beta$'s.  This
implies that $\sum_\gamma c_\gamma \, z^\gamma$ converges absolutely
when $z \in A \cap B$, and one can show that (\ref{(sum a_alpha
z^alpha) (sum b_beta z^beta) = sum c_gamma z^gamma}) holds by
approximating infinite sums by sums with only finitely many nonzero terms.
It suffices to consider the case where $z = (1, \ldots, 1)$, since otherwise
the monomials in $z$ can be absorbed into the coefficients.  One can also
use linearity to reduce to the case where the coefficients are nonnegative
real numbers, and estimate products of sums of finitely many $a_\alpha$'s
and $b_\beta$'s in terms of sums of finitely many $c_\gamma$'s.

        Let us return to a single power series $\sum_\alpha a_\alpha
\, z^\alpha$, and suppose that $w, z \in A$ and $u \in {\bf C}^n$ satisfy
\begin{equation}
\label{|u_j| le |w_j|^t |z_j|^{1 - t}}
        |u_j| \le |w_j|^t \, |z_j|^{1 - t}
\end{equation}
for some $t \in {\bf R}$, $0 < t < 1$, and each $j = 1, \ldots, n$.  Hence
\begin{equation}
        |u^\alpha| \le |w^\alpha|^t \, |z^\alpha|^{1 - t}
\end{equation}
for each multi-index $\alpha$.  The convexity of the exponential
function on the real line implies that
\begin{equation}
\label{k^t l^{1 - t} le t k + (1 - t) l}
        k^t \, l^{1 - t} \le t \, k + (1 - t) \, l
\end{equation}
for every $k, l \ge 0$.  Applying this to $k = |w^\alpha|$, $l =
|z^\alpha|$ and summing over $\alpha$, we get that $u \in A$, because
\begin{equation}
        \sum_\alpha |a_\alpha| \, |u^\alpha|
         \le t \sum_\alpha |a_\alpha| \, |w^\alpha|^t
            + (1 - t) \sum_\alpha |a_\alpha| \, |z^\alpha|^{1 - t}.
\end{equation}

\section{Power series, continued}
\label{power series, continued}
\setcounter{equation}{0}

        Let $n$ be a positive integer, and let $\sum_\alpha a_\alpha
\, z^\alpha$ be a power series with complex coefficients in $z = (z_1,
\ldots, z_n)$.  If $l$ is a nonnegative integer, then
\begin{equation}
        p_l(z) = \sum_{|\alpha| = l} a_\alpha \, z^\alpha
\end{equation}
is a homogeneous polynomial of degree $l$ in $z$, where more precisely
the sum is taken over the finitely many multi-indices $\alpha$ such
that $|\alpha| = l$.  Of course,
\begin{equation}
\label{sum_{l = 0}^infty p_l(z) = sum_alpha a_alpha z^alpha}
        \sum_{l = 0}^\infty p_l(z) = \sum_\alpha a_\alpha \, z^\alpha
\end{equation}
formally, which gives another way to look at the convergence of
$\sum_\alpha a_\alpha \, z^\alpha$.  In particular, if $\sum_\alpha
a_\alpha \, z^\alpha$ converges absolutely for some $z \in {\bf C}^n$,
then $\sum_{l = 0}^\infty p_l(z)$ converges absolutely, and the
two sums are the same.  This uses the fact that
\begin{equation}
        |p_l(z)| \le \sum_{|\alpha| = l} |a_\alpha| \, |z^\alpha|
\end{equation}
for each $l$.

        Let $\sum_\alpha b_\alpha \, z^\alpha$ be another power series,
with the corresponding polynomials
\begin{equation}
        q_l(z) = \sum_{|\alpha| = l} b_\alpha \, z^\alpha.
\end{equation}
Thus $p_l(z) + q_l(z)$ are the polynomials associated to $\sum_\alpha
(a_\alpha + b_\alpha) \, z^\alpha$.  Suppose that $\sum_\gamma
c_\gamma \, z^\gamma$ is the power series obtained by formally
multiplying $\sum_\alpha a_\alpha \, z^\alpha$ and $\sum_\beta b_\beta
\, z^\beta$, so that
\begin{equation}
        c_\gamma = \sum_{\alpha + \beta = \gamma} a_\alpha \, b_\beta.
\end{equation} 
It is easy to check that the corresponding polynomials
\begin{equation}
        r_l = \sum_{|\gamma| = l} c_\gamma \, z^\gamma
\end{equation}
are also given by
\begin{equation}
        r_l = \sum_{j = 0}^l p_j(z) \, q_{l - j}(z).
\end{equation}
This shows that $r_l$ is the Cauchy product of the $p_j$'s and $q_k$'s.

        Note that
\begin{equation}
\label{sum_{l = 0}^infty p_l(t z) = sum_{l = 0}^infty t^l p_l(z)}
        \sum_{l = 0}^\infty p_l(t \, z) = \sum_{l = 0}^\infty t^l \, p_l(z)
\end{equation}
may be considered as an ordinary power series in $t \in {\bf C}$ for
each $z \in {\bf C}^n$.  This gives another way to look at the Cauchy
product in the preceding paragraph, as the coefficients of the product
of two power series in $t$.  If $\sum_{l = 0}^\infty p_l(z)$ converges
for some $z \in {\bf C}^n$, then $\{p_l(z)\}_{l = 1}^\infty$ converges
to $0$, and hence $\{p_l(z)\}_{l = 1}^\infty$ is bounded.  This
implies that (\ref{sum_{l = 0}^infty p_l(t z) = sum_{l = 0}^infty t^l
p_l(z)}) converges absolutely when $|t| < 1$, by the comparison test.

        Consider
\begin{equation}
        p^*(z) = \limsup_{l \to \infty} |p_l(z)|^{1/l},
\end{equation}
which takes values in $[0, \infty]$.  Observe that
\begin{equation}
\label{p^*(t z) = |t| rho(z)}
        p^*(t \, z) = |t| \, p^*(z)
\end{equation}
for each $t \in {\bf C}$ and $z \in {\bf C}^n$, because $p_l(z)$ is
homogeneous of degree $l$.  The right side of (\ref{p^*(t z) = |t|
rho(z)}) should be interpreted as being $0$ when $t = 0$, even when
$p^*(z) = +\infty$, because $p^*(0) = 0$.  The \emph{root test} states
that $\sum_{l = 0}^\infty p_l(z)$ converges absolutely when $p^*(z) <
1$, and diverges when $p^*(z) > 1$.  It follows that the radius of
convergence of (\ref{sum_{l = 0}^infty p_l(t z) = sum_{l = 0}^infty
t^l p_l(z)}) as a power series in $t$ is equal to $1/p^*(z)$.

\section{Linear transformations}
\label{linear transformations}
\setcounter{equation}{0}

        Let $n$ be a positive integer, and let $T$ be a one-to-one
linear transformation from ${\bf C}^n$ onto itself.  Consider the mapping
$\rho_T$ acting on complex-valued functions on ${\bf C}^n$ defined by
\begin{equation}
        \rho_T(f)(z) = f(T^{-1}(z)).
\end{equation}
Thus
\begin{equation}
        \rho_T(f + g) = \rho_T(f) + \rho_T(g)
\end{equation}
and
\begin{equation}
        \rho_T(f \, g) = \rho_T(f) \, \rho_T(g).
\end{equation}
for any pair of functions $f$, $g$ on ${\bf C}^n$.
If $f$ is a polynomial on ${\bf C}^n$, then it is easy to see that
$\rho_T(f)$ is a polynomial too.  If $f$ is a homogeneous polynomial,
then $\rho_T(f)$ is a homogeneous polynomial as well, of the same
degree.

        Of course, $\rho_T(f) = f$ for every function $f$ on ${\bf
C}^n$ when $T$ is the identity transformation on ${\bf C}^n$.  If $R$,
$T$ are arbitrary invertible linear transformation on ${\bf C}^n$,
then
\begin{eqnarray}
 \rho_R(\rho_T(f))(z) & = & \rho_T(f)(R^{-1}(z)) = f(T^{-1}(R^{-1}(z))) \\
        & = & f((R \circ T)^{-1}(z)) = \rho_{R \circ T}(f)(z). \nonumber
\end{eqnarray}
In particular, $\rho_{T^{-1}} = (\rho_T)^{-1}$.  Let $GL({\bf C}^n)$
be the group of invertible linear transformations on ${\bf C}^n$, with
composition of mappings as the group operation.  It follows that $T
\mapsto \rho_T$ is a homomorphism from $GL({\bf C}^n)$ into the group
of invertible linear transformations on the space of functions on
${\bf C}^n$, which is to say a representation of $GL({\bf C}^n)$ on
the space of functions on ${\bf C}^n$.

        Let $f(z) = \sum_\alpha a_\alpha \, z^\alpha$ be a formal
power series with complex coefficients.  This can also be expressed as
$\sum_{l = 0}^\infty p_l(z)$, where $p_l(z)$ is a homogeneous
polynomial of degree $l$ for each $l \ge 0$.  If $T$ is an invertible
linear transformation on ${\bf C}^n$, then we can take $\rho_T(f)$ to
be the formal power series that corresponds to $\sum_{l = 0}^\infty
\rho_T(p_l)$.  It is easy to see that this preserves sums and products
of power series, just as for ordinary functions.  In particular, this
defines a representation of $GL({\bf C}^n)$ on the space of formal
power series.

        If $\sum_{l = 0}^\infty p_l(z)$ converges for some $z \in {\bf
C}^n$, then $\sum_{l = 0}^\infty \rho_T(p_l)(T(z))$ converges and has
the same sum, because it is the same series of complex numbers.  If
$\sum_{l = 0}^\infty \rho_l(z)$ converges for every $z \in {\bf C}^n$,
then $\sum_{l = 0}^\infty \rho_T(p_l)(T(z))$ converges for every $z
\in {\bf C}^n$, and has the same sum.  Hence the formal and pointwise
definitions of $\rho(f)$ are consistent with each other in this case.

\section{Abel summability}
\label{abel summability}
\setcounter{equation}{0}

        Let $\sum_{j = 0}^\infty a_j$ be an infinite series of complex
numbers, and put
\begin{equation}
\label{A(r) = sum_{j = 0}^infty a_j r^j}
        A(r) = \sum_{j = 0}^\infty a_j \, r^j
\end{equation}
when $0 \le r < 1$.  More precisely, we suppose that the sum on the
right converges for each $r < 1$, which implies that $\{a_j \,
r^j\}_{j = 0}^\infty$ converges to $0$ for each $r < 1$, and hence
that $\{a_j \, r^j\}_{j = 0}^\infty$ is bounded for each $r < 1$.
Conversely, if $\{a_j \, r^j\}_{j = 0}^\infty$ is bounded for each $r
< 1$, then $\sum_{j = 0}^\infty a_j \, t^j$ converges absolutely for
each $t < 1$, as one can see by taking $t < r < 1$ and using the
comparison test, since $\sum_{j = 0}^\infty (t/r)^j$ is a convergent
geometric series under these conditions.  The expressions $A(r)$ are
known as the \emph{Abel sums} associated to $\sum_{j = 0}^\infty a_j$,
and we say that $\sum_{j = 0}^\infty a_j$ is \emph{Abel summable} if
\begin{equation}
        \lim_{r \to 1-} A(r)
\end{equation}
exists.

        If $\sum_{j = 0}^\infty a_j$  converges in the usual sense,
then it is Abel summable.  To see this, let
\begin{equation}
\label{s_n = sum_{j = 0}^n a_j}
        s_n = \sum_{j = 0}^n a_j
\end{equation}
be the $n$th partial sum of $\sum_{j = 0}^\infty a_j$ when $n \ge 0$,
and put $s_{-1} = 0$.  Thus $a_j = s_j - s_{j - 1}$ for each $j \ge
0$, and hence
\begin{equation}
 A(r) = \sum_{j = 0}^\infty (s_j - s_{j - 1}) \, r^j
     = \sum_{j = 0}^\infty s_j \, r^j - \sum_{j = 0}^\infty s_{j - 1} \, r^j
\end{equation}
when $0 \le r < 1$.  There is no problem with the convergence of the
series on the right, because the convergence of $\sum_{j = 0}^\infty
a_j$ implies that $\{a_j\}_{j = 0}^\infty$ converges to $0$ and is
therefore bounded, which implies that $s_n = O(n)$.  Of course,
\begin{equation}
 \sum_{j = 0}^\infty s_{j - 1} \, r^j = \sum_{j = 1}^\infty s_{j - 1} \, r^j
        = \sum_{j = 0}^\infty s_j \, r^{j + 1},
\end{equation}
because $s_{-1} = 0$, which implies that
\begin{equation}
        A(r) = \sum_{j = 0}^\infty s_j (r^j - r^{j + 1})
             = (1 - r) \sum_{j = 0}^\infty s_j \, r^j
\end{equation}
when $r < 1$.

        We would like to show that
\begin{equation}
        \lim_{r \to 1-} A(r) = \lim_{j \to \infty} s_j
\end{equation}
when the limit on the right side exists.  Put $s = \lim_{j \to \infty}
s_j$, let $\epsilon > 0$ be given, and choose $L \ge 0$ such that
\begin{equation}
\label{|s_j - s| < frac{epsilon}{2}}
        |s_j - s| < \frac{\epsilon}{2}
\end{equation}
for every $j \ge L$.  Observe that
\begin{equation}
        A(r) - s = (1 - r) \sum_{j = 0}^\infty (s_j - s) \, r^j
\end{equation}
when $r < 1$, because $(1 - r) \sum_{j = 0}^\infty r^j = 1$.  It
follows that
\begin{eqnarray}
 |A(r) - s| & \le & (1 - r) \sum_{j = 0}^\infty |s_j - s| \, r^j   \\
            & < & (1 - r) \sum_{j = 0}^{L - 1} |s_j - s| \, r^j
                 + (1 - r) \sum_{j = L}^\infty (\epsilon/2) \, r^j \nonumber \\
 & \le & (1 - r) \sum_{j = 0}^{L - 1} |s_j - s| \, r^j + \frac{\epsilon}{2}
                                                                   \nonumber
\end{eqnarray}
for each $r < 1$.  If $r$ is sufficiently close to $1$, then
\begin{equation}
        (1 - r) \sum_{j = 0}^{L - 1} |s_j - s| \, r^j
               \le (1 - r) \sum_{j = 0}^{L - 1} |s_j - s| < \frac{\epsilon}{2}
\end{equation}
so that $|A(r) - s| < \epsilon / 2 + \epsilon / 2 = \epsilon$, as
desired.

        If $a \in {\bf C}$ satisfies $|a| = 1$, then
\begin{equation}
        \sum_{j = 0}^\infty a^j \, r^j = \frac{1}{1 - a \, r}
\end{equation}
when $0 \le r < 1$.  Hence $\sum_{j = 0}^\infty a^j$ is Abel summable
when $a \ne 1$, with the sum equal to $(1 - a)^{-1}$.  Let $\sum_{j =
0}^\infty a_j$, $\sum_{j = 0}^\infty b_j$ be infinite series of
complex numbers with Abel sums $A(r)$, $B(r)$, respectively, and note
that $\sum_{j = 0}^\infty (a_j + b_j)$ has Abel sums given by $A(r) +
B(r)$.  If $\sum_{j = 0}^\infty a_j$, $\sum_{j = 0}^\infty b_j$ are
Abel summable, then it follows that $\sum_{j = 0}^\infty (a_j + b_j)$
is Abel summable, with the Abel sum of the latter equal to the sum of
the Abel sums of the first two series.  Suppose now that $c_n =
\sum_{j = 0}^n a_j \, b_{n - j}$ is the Cauchy product of the $a_j$'s
and $b_j$'s, and let $C(r)$ be the corresponding Abel sums.  As in
Section \ref{cauchy products},
\begin{equation}
\label{C(r) = A(r) B(r)}
        C(r) = A(r) \, B(r)
\end{equation}
when $0 \le r < 1$.  More precisely, if the series defining $A(r)$,
$B(r)$ converge absolutely, then the series defining $C(r)$ also
converges absolutely, and satisfies (\ref{C(r) = A(r) B(r)}).  The
existence of the Abel sums for $\sum_{j = 0}^\infty a_j$, $\sum_{j =
0}^\infty b_j$ for each $r < 1$ implies that this condition holds for
every $r < 1$, as discussed at the beginning of this section.  If
$\sum_{j = 0}^\infty a_j$, $\sum_{j = 0}^\infty b_j$ are Abel
summable, then it follows that $\sum_{n = 0}^\infty c_n$ is Abel
summable, and that the Abel sum of the latter equal to the product of
the Abel sums of the former.

\section{Multiple Fourier series}
\label{multiple fourier series}
\setcounter{equation}{0}

        Let $n$ be a positive integer, and let ${\bf T}^n$ be the
$n$-dimensional torus, consisting of $z = (z_1, \ldots, z_n) \in {\bf
C}^n$ such that $|z_j| = 1$ for $j = 1, \ldots, n$.  If $\alpha =
(\alpha_1, \ldots, \alpha_n)$ is an $n$-tuple of integers, then put
\begin{equation}
        z^\alpha = z_1^{\alpha_1} \cdots z_n^{\alpha_n},
\end{equation}
with the usual convention that $z_j^{\alpha_j} = 1$ when $\alpha_j =
0$.  Thus
\begin{equation}
 \frac{1}{(2 \pi)^n} \int_{{\bf T}^n} z^\alpha \, |dz|
 = \prod_{j = 1}^n \frac{1}{2 \, \pi} \int_{\bf T} z_j^{\alpha_j} \, |dz_j|
\end{equation}
is equal to $0$ when $\alpha \ne 0$, and is equal to $1$ when $\alpha
= 0$.  Here $|dz|$ is the $n$-dimensional element of integration on
${\bf T}^n$ corresponding to the element $|dz_j|$ of arc length in
each variable.

        If $f$ is a continuous complex-valued function on ${\bf T}^n$
and $\alpha \in {\bf Z}^n$, then we put
\begin{equation}
\label{widehat{f}(alpha) = frac{1}{(2 pi)^n} int_{T^n} f(z) z^{-alpha} |dz|}
        \widehat{f}(\alpha)
          = \frac{1}{(2 \pi)^n} \int_{{\bf T}^n} f(z) \, z^{-\alpha} \, |dz|.
\end{equation}
The corresponding Fourier series is given by
\begin{equation}
\label{sum_{alpha in {bf Z}^n} widehat{f}(alpha) z^alpha}
        \sum_{\alpha \in {\bf Z}^n} \widehat{f}(\alpha) \, z^\alpha.
\end{equation}
For example, if $f(z) = z^\beta$ for some $\beta \in {\bf Z}^n$, then
$\widehat{f}(\alpha) = 1$ when $\alpha = \beta$ and is equal to $0$
otherwise.  Thus (\ref{sum_{alpha in {bf Z}^n} widehat{f}(alpha)
z^alpha}) reduces to $f$ in this case, or when $f$ is a finite linear
combination of $z^\beta$'s.  Note that
\begin{equation}
\label{|widehat{f}(alpha)| le frac{1}{(2 pi)^n} int_{{bf T}^n} |f(z)| |dz|}
 |\widehat{f}(\alpha)| \le \frac{1}{(2 \pi)^n} \int_{{\bf T}^n} |f(z)| \, |dz|
\end{equation}
for any continuous function $f$ on ${\bf T}^n$ and $\alpha \in {\bf
Z}^n$.

        Let $U^n$ be the $n$-dimensional open unit polydisk,
consisting of $z \in {\bf C}^n$ with $|z_j| < 1$ for $j = 1, \ldots,
n$.  The $n$-dimensional Poisson kernel $P_n(z, w)$ can be defined for
$z \in U^n$ and $w \in {\bf T}^n$ by
\begin{equation}
        P_n(z, w) = \prod_{j = 1}^n P(z_j, w_j),
\end{equation}
where $P(z_j, w_j)$ is the ordinary Poisson kernel evaluated at $z_j$,
$w_j$, as in Section \ref{poisson kernel}.  If $f$ is a continuous
function on ${\bf T}^n$, then its Poisson integral is defined on $U^n$ by
\begin{equation}
        \phi(z) = \int_{{\bf T}^n} P_n(z, w) \, f(w) \, |dw|.
\end{equation}
As before, one can show that $\phi(z) \to f(z_0)$ as $z \in U^n$ tends
to $z_0 \in {\bf T}^n$, but one can also do more than this.

        Let $\overline{U}^n$ be the $n$-dimensional closed unit
polydisk, consisting of $z \in {\bf C}^n$ such that $|z_j| \le 1$ for
each $j$.  Of course, this is the same as the closure of $U^n$ in
${\bf C}^n$.  The boundary $\partial U^n$ of $U^n$ in ${\bf C}^n$
consists of $z \in {\bf C}^n$ such that $|z_j| \le 1$ for each $j$ and
$|z_j| = 1$ for at least one $j$.  In particular, ${\bf T}^n \subseteq
\partial U^n$, but ${\bf T}^n$ is a relatively small subset of
$\partial U^n$ when $n > 1$.  More precisely, $U^n$ has complex
dimension $n$ and hence real dimension $2 \, n$, $\partial U^n$ has
real dimension $2 n - 1$, and ${\bf T}^n$ has real dimension $n$.

        One can extend $\phi(z)$ to $z \in \overline{U}^n$ in the
following way.  If $z \in {\bf T}^n$, then we simply put $\phi(z) =
f(z)$.  If $z \in \partial U^n \backslash {\bf T}^n$, then $|z_j| < 1$
for at least one $j$, and we define $\phi(z)$ by taking the Poisson
integral of $f$ in the $j$th variable when $|z_j| < 1$, and simply
evaluating $f$ at $z_j$ in the $j$th variable when $|z_j| = 1$.  It is
not too difficult to show that this defines a continuous function on
$\overline{U}^n$.

        If $z \in {\bf C}^n$ and $\alpha \in {\bf Z}^n$, then put
\begin{equation}
\label{widetilde{z}^alpha = ...}
        \widetilde{z}^\alpha
           = \widetilde{z}_1^{\alpha_1} \cdots \widetilde{z}_n^{\alpha_n},
\end{equation}
where $\widetilde{z}_j^{\alpha_j} = z_j^{\alpha_j}$ when $\alpha_j \ge
0$ and $\widetilde{z}_j^{\alpha_j} = \overline{z_j}^{-\alpha_j}$ when
$\alpha_j < 0$.  Thus $\widetilde{z}^\alpha = z^\alpha$ when $z \in
{\bf T}^n$, and
\begin{equation}
\label{|widetilde{z}^alpha| = |z_1|^{|alpha_1|} cdots |z_n|^{|alpha_n|}}
        |\widetilde{z}^\alpha| = |z_1|^{|\alpha_1|} \cdots |z_n|^{|\alpha_n|}
\end{equation}
for every $z \in {\bf C}^n$.  If $z \in U^n$, then it is easy to see that
\begin{equation}
\label{phi(z) = sum_{alpha in {bf Z}^n} widehat{f}(alpha) widetilde{z}^alpha}
        \phi(z) = \sum_{\alpha \in {\bf Z}^n} \widehat{f}(\alpha) \,
                                                       \widetilde{z}^\alpha,
\end{equation}
using the analogous expansion for the Poisson kernel in one variable.
Note that this series converges absolutely for every $z \in U^n$,
since the Fourier coefficients $\widehat{f}(\alpha)$ are bounded, as
in (\ref{|widehat{f}(alpha)| le frac{1}{(2 pi)^n} int_{{bf T}^n}
|f(z)| |dz|}).

        If $z \in {\bf T}^n$ and $0 \le r < 1$, then put
\begin{equation}
\label{f_r(z) = phi(r z) = ...}
        f_r(z) = \phi(r \, z)
 = \sum_{\alpha \in {\bf Z}^n} \widehat{f}(\alpha) \, r^{|\alpha|} \, z^\alpha,
\end{equation}
where $|\alpha| = |\alpha_1| + \cdots + |\alpha_n|$.  One can check
that $f_r \to f$ as $r \to 1$ uniformly on ${\bf T}^n$, using the fact
that continuous functions on compact sets are uniformly continuous.
The sum on the right side of (\ref{f_r(z) = phi(r z) = ...}) can be
approximated by finite subsums uniformly on ${\bf T}^n$ for each $r <
1$, as in Weierstrass' M-test.  It follows that every continuous
function $f$ on ${\bf T}^n$ can be approximated uniformly by finite
linear combinations of $z^\alpha$'s, $\alpha \in {\bf Z}^n$.

        Observe that $\phi(z)$ is ``polyharmonic'', in the sense that
it is harmonic as a function of $z_j$ on the set where $|z_j| < 1$ for
each $j$.  This follows from the remarks about harmonic functions of
one complex variable in Section \ref{poisson kernel}.  In addition,
\begin{equation}
\label{sup_{z in overline{U}^n} |phi(z)| = sup_{z in {bf T}^n} |f(z)|}
        \sup_{z \in \overline{U}^n} |\phi(z)| = \sup_{z \in {\bf T}^n} |f(z)|.
\end{equation}
More precisely, the right side of (\ref{sup_{z in overline{U}^n}
|phi(z)| = sup_{z in {bf T}^n} |f(z)|}) is less than or equal to the
left side because ${\bf T}^n \subseteq \overline{U}^n$ and $\phi = f$
on ${\bf T}^n$.  To get the opposite inequality, one can use the fact
that the Poisson kernel is positive and has integral equal to $1$.

        If $f$, $g$ are continuous complex-valued functions on ${\bf
T}^n$, then put
\begin{equation}
\label{langle f, g rangle = (2 pi)^{-n} int_{T^n} f(z) overline{g(z)} |dz|}
 \langle f, g \rangle
   = \frac{1}{(2 \pi)^n} \int_{{\bf T}^n} f(z) \, \overline{g(z)} \, |dz|.
\end{equation}
This defines an inner product on the vector space $C({\bf T}^n)$ of
continuous complex-valued functions on ${\bf T}^n$, for which the
corresponding norm is given by
\begin{equation}
\label{||f|| = (frac{1}{(2 pi)^n} int_{{bf T}^n} |f(z)|^2 |dz|)^{1/2}}
 \|f\| = \Big(\frac{1}{(2 \pi)^n} \int_{{\bf T}^n} |f(z)|^2 \, |dz|\Big)^{1/2}.
\end{equation}
It is easy to see that the functions $z^\alpha$, $\alpha \in {\bf Z}^n$ are
orthonormal with respect to this inner product, and that the Fourier
coefficients of a continuous function $f$ on ${\bf T}^n$ can be expressed by
\begin{equation}
\label{widehat{f}(alpha) = langle f, z^alpha rangle}
        \widehat{f}(\alpha) = \langle f, z^\alpha \rangle.
\end{equation}
The $n$-dimensional version of Parseval's formula states that
\begin{equation}
\label{sum_{alpha in {bf Z}^n} |widehat{f}(alpha)|^2 = ...}
        \sum_{\alpha \in {\bf Z}^n} |\widehat{f}(\alpha)|^2
          = \frac{1}{(2 \pi)^n} \int_{{\bf T}^n} |f(z)|^2 \, |dz|,
\end{equation}
where the summability of the sum on the left is part of the
conclusion.  This follows from the orthonormality of the $z^\alpha$'s
and the fact that their finite linear combinations are dense in
$C({\bf T})$, as in the one-dimensional case.

        Suppose that $f$, $g$ are continuous functions on ${\bf T}^n$,
and let us check that
\begin{equation}
\label{widehat{(f g)}(alpha) = ...}
        \widehat{(f \, g)}(\alpha) =
 \sum_{\beta \in {\bf Z}^n} \widehat{f}(\alpha - \beta) \, \widehat{g}(\beta).
\end{equation}
This can be derived formally by multiplying the Fourier series for
$f$, $g$ and collecting terms.  To make this rigorous, observe first
that $\widehat{f}(\alpha - \beta) \, \widehat{g}(\beta)$ is summable
in $\beta$, because $\widehat{f}, \widehat{g} \in \ell^2({\bf Z}^n)$,
as in the previous paragraph.  If $g(z) = z^\gamma$ for some $\gamma
\in {\bf Z}^n$, then it is easy to see that both sides of
(\ref{widehat{(f g)}(alpha) = ...})  are equal to $\widehat{f}(\alpha
- \gamma)$.  It follows that (\ref{widehat{(f g)}(alpha) = ...}) holds
when $g$ is a finite linear combination of $z^\gamma$'s, and the same
conclusion for an arbitrary continuous function $g$ on ${\bf T}^n$ can
be obtained by approximation by linear combinations of $z^\gamma$'s.

        If $a(\alpha)$, $b(\alpha)$ are summable functions on ${\bf
Z}^n$, then their convolution can be defined by
\begin{equation}
\label{(a * b)(alpha) = sum_{beta in {bf Z}^n} a(alpha - beta) b(beta)}
 (a * b)(\alpha) = \sum_{\beta \in {\bf Z}^n} a(\alpha - \beta) \, b(\beta),
\end{equation}
as in the one-dimensional case.  More precisely, $a * b$ is also
summable on ${\bf Z}^n$, and satisfies
\begin{equation}
\label{||a * b||_1 le ||a||_1 ||b||_1}
        \|a * b\|_1 \le \|a\|_1 \, \|b\|_1,
\end{equation}
where $\|a\|_1$ is the $\ell^1$ norm of $a$ on ${\bf Z}^n$.  This
follows by interchanging the order of summation, as before, and one
can also check that $\ell^1({\bf Z}^n)$ is a commutative Banach
algebra with respect to convolution.  The Fourier transform of $a$ in
$\ell^1({\bf Z}^n)$ is defined by
\begin{equation}
\label{widehat{a}(z) = sum_{alpha in {bf Z}^n} a(alpha) z^alpha}
        \widehat{a}(z) = \sum_{\alpha \in {\bf Z}^n} a(\alpha) \, z^\alpha
\end{equation}
for $z \in {\bf T}^n$.  The sum on the right is absolutely summable
for each $z \in {\bf T}^n$, because $a(\alpha)$ is summable, and can
be approximated by finite subsums uniformly on ${\bf T}^n$, as in
Weierstrass' M-test.  This implies that $\widehat{a}(z)$ is continuous
on ${\bf T}^n$, and it is easy to see that
\begin{equation}
\label{widehat{(a * b)}(z) = widehat{a}(z) widehat{b}(z), n dimensions}
        \widehat{(a * b)}(z) = \widehat{a}(z) \, \widehat{b}(z)
\end{equation}
for every $a, b \in \ell^1({\bf Z}^n)$ and $z \in {\bf T}^n$, as before.
Conversely, every nonzero multiplicative homomorphism on $\ell^1({\bf
Z}^n)$ with respect to convolution can be represented as $a \mapsto
\widehat{a}(z)$ for some $z \in {\bf T}^n$, as in the one-dimensional
situation.  Note that the Fourier coefficients of $\widehat{a}$ are
given by $a(\alpha)$ for every $a \in \ell^1({\bf Z}^n)$, because of
the orthogonality properties of the $z^\alpha$'s.  Observe also that
\begin{equation}
\label{sum_{alpha in {bf Z}^n} a(alpha) widetilde{z}^alpha}
        \sum_{\alpha \in {\bf Z}^n} a(\alpha) \, \widetilde{z}^\alpha
\end{equation}
is absolutely summable for every $z \in \overline{U}^n$, which is the
analogue of the function $\phi$ discussed earlier.  As usual,
(\ref{sum_{alpha in {bf Z}^n} a(alpha) widetilde{z}^alpha}) can be
approximated by finite subsums uniformly on $\overline{U}^n$ under
these conditions, which implies more directly that it defines a
continuous function on $\overline{U}^n$ than in the earlier discussion.

\section{Functions of analytic type}
\label{functions of analytic type}
\setcounter{equation}{0}

        Let $A({\bf T}^n)$ be the collection of continuous functions
$f$ on ${\bf T}^n$ such that
\begin{equation}
\label{widehat{f}(alpha) = 0}
        \widehat{f}(\alpha) = 0
\end{equation}
when $\alpha \in {\bf Z}^n$ satisfies $\alpha_j < 0$ for some $j$.  If
$f, g \in A({\bf T}^n)$, then it is easy to see from (\ref{widehat{(f
g)}(alpha) = ...}) that their product $f \, g$ is in $A({\bf T}^n)$
too.  Note that the sum on the right side of (\ref{widehat{(f
g)}(alpha) = ...}) has only finitely many nonzero terms in this
situation.  It follows that $A({\bf T}^n)$ is a subalgebra of $C({\bf
T}^n)$, since the former is clearly a linear subspace of the latter.

        If $f \in A({\bf T}^n)$, then (\ref{phi(z) = sum_{alpha in {bf
Z}^n} widehat{f}(alpha) widetilde{z}^alpha}) reduces to
\begin{equation}
        \phi(z) = \sum_\alpha \widehat{f}(\alpha) \, z^\alpha,
\end{equation}
where now the sum is taken over all multi-indices $\alpha$, which is
to say $\alpha \in {\bf Z}^n$ such that $\alpha_j \ge 0$ for each $j$.
Thus we get an ordinary power series in this case, in the sense that
the $z^\alpha$'s are the usual monomials, instead of the modified
monomials $\widetilde{z}^\alpha$ that may include complex conjugation.
In particular, this implies that $f$ can be approximated uniformly on
${\bf T}^n$ by a finite linear combinations of $z^\alpha$'s, where the
$\alpha$'s are multi-indices, by the same type of argument as in the
previous section.  Of course, $z^\alpha \in A({\bf T}^n)$ for every
multi-index $\alpha$, and $A({\bf T}^n)$ is a closed set in $C({\bf
T}^n)$ with respect to the supremum norm.  It follows that $A({\bf
T}^n)$ is the same as the closure in $C({\bf T}^n)$ of the linear span
of the $z^\alpha$'s, where $\alpha$ is a multi-index.

        Let $\phi_f$ be the continuous function $\phi$ on the closed
unit polydisk $\overline{U}^n$ associated to $f \in C({\bf T}^n)$ as
in the previous section.  If $f, g \in A({\bf T}^n)$, then
\begin{equation}
\label{phi_{f g} = phi_f phi_g}
        \phi_{f g} = \phi_f \, \phi_g.
\end{equation}
This follows by multiplying the series expansions for $\phi_f$,
$\phi_g$ in the previous paragraph and collecting terms, as in
(\ref{(sum a_alpha z^alpha) (sum b_beta z^beta) = sum c_gamma
z^gamma}) and (\ref{c_gamma = sum_{alpha + beta = gamma} a_alpha
b_beta}).  This also uses the formula (\ref{widehat{(f g)}(alpha) =
...}) for the Fourier coefficients of the product $f \, g$.  The main
point is that
\begin{equation}
        z^\alpha \, z^\beta = z^{\alpha + \beta}
\end{equation}
while $\widetilde{z}^\alpha \, \widetilde{z}^\beta$ is not necessarily
the same as $\widetilde{z}^{\alpha + \beta}$.  More precisely, this
argument works on the open unit polydisk $U^n$, where the series
expansions for $\phi_f$, $\phi_g$ are absolutely summable.  This
implies that (\ref{phi_{f g} = phi_f phi_g}) holds on
$\overline{U}^n$, by continuity.

        Observe that $A({\bf T}^n)$ is a commutative Banach algebra
with respect to the supremum norm, since it is a closed subalgebra of
$C({\bf T}^n)$ that contains the constant functions, and hence the
multiplicative identity element.  If $p \in \overline{U}^n$, then $f
\mapsto \phi_f(p)$ defines a nonzero homomorphism from $A({\bf T}^n)$
into the complex numbers.  Conversely, suppose that $h$ is a nonzero
homomorphism on $A({\bf T}^n)$, and let us show that there is a $p \in
\overline{U}^n$ such that $h(f) = \phi_f(p)$ for every $f \in A({\bf
T}^n)$.  As usual, $h({\bf 1}_{{\bf T}^n}) = 1$, where ${\bf 1}_{{\bf
T}^n}$ is the constant function equal to $1$ on ${\bf T}^n$, and
\begin{equation}
\label{|h(f)| le sup_{z in {bf T}^n} |f(z)|}
        |h(f)| \le \sup_{z \in {\bf T}^n} |f(z)|
\end{equation}
for every $f \in A({\bf T}^n)$.  Consider $f_j(z) = z_j$, $j = 1,
\ldots, n$, as an element of $A({\bf T}^n)$.  If $p_j = h(f_j)$, then
$|p_j| \le 1$ for each $j$, by (\ref{|h(f)| le sup_{z in {bf T}^n}
|f(z)|}).  Hence $p = (p_1, \ldots, p_n) \in \overline{U}^n$.  By
construction, $h(f) = \phi_f(p)$ when $f = f_j$ for some $j$, and it
follows that this also holds when $f$ is a polynomial, because $h$ is
a homomorphism.  Using (\ref{|h(f)| le sup_{z in {bf T}^n} |f(z)|})
again, we get that $h(f) = \phi_f(p)$ for every $f \in A({\bf T}^n)$,
because polynomials are dense in $A({\bf T}^n)$.

        Let $\ell^1_A({\bf Z}^n)$ be the set of $a \in \ell^1({\bf
Z}^n)$ such that $a(\alpha) = 0$ whenever $\alpha \in {\bf Z}^n$
satisfies $\alpha_j < 0$ for some $j$.  It is easy to see that this is
a closed subalgebra of $\ell^1({\bf Z}^n)$ with respect to
convolution.  If $a \in \ell^1_A({\bf Z}^n)$, then (\ref{sum_{alpha in
{bf Z}^n} a(alpha) widetilde{z}^alpha}) reduces to an ordinary power series
\begin{equation}
\label{sum_alpha a(alpha) z^alpha}
        \sum_\alpha a(\alpha) \, z^\alpha,
\end{equation}
where the sum is taken over all multi-indices $\alpha$.  If $b \in
\ell^1_A({\bf Z}^n)$ too, then
\begin{equation}
        \Big(\sum_\alpha a(\alpha) \, z^\alpha\Big) \,
         \Big(\sum_\beta b(\beta) \, z^\beta\Big)
           = \sum_\gamma (a * b)(\gamma) \, z^\gamma
\end{equation}
for every $z \in \overline{U}^n$, which is basically the same as
(\ref{(sum a_alpha z^alpha) (sum b_beta z^beta) = sum c_gamma
z^gamma}) again.  Thus the mapping from $a$ to (\ref{sum_alpha
a(alpha) z^alpha}) defines a homomorphism from $\ell^1({\bf Z}^n)$
into the complex numbers for each $z \in \overline{U}^n$, using
convolution as multiplication on $\ell^1_A({\bf Z}^n)$.  

        Conversely, let us check that any nonzero homomorphism $h$
from $\ell^1_A({\bf Z}^n)$ into the complex numbers is of this form.
If $\alpha \in {\bf Z}^n$, then let $\delta_\alpha$ be the function on
${\bf Z}^n$ defined by $\delta_\alpha(\beta) = 1$ when $\alpha =
\beta$ and $\delta_\alpha(\beta) = 0$ otherwise.  Thus $\delta_\alpha
\in \ell^1_A({\bf Z}^n)$ when $\alpha_j \ge 0$ for each $j$.  In particular,
$\delta_0 \in \ell^1_A({\bf Z}^n)$, which is the multiplicative identity
element for $\ell^1({\bf Z}^n)$, and hence for $\ell^1_A({\bf Z}^n)$.
It follows that $h(\delta_0) = 1$, and we also have that
\begin{equation}
\label{|h(a)| le ||a||_1}
        |h(a)| \le \|a\|_1
\end{equation}
for every $a \in \ell^1_A({\bf Z}^n)$, since $\ell^1_A({\bf Z}^n)$ is
a Banach algebra.  Let $\alpha(l)$ be the element of ${\bf Z}^n$ with
$l$th component equal to $1$ and other components equal to $0$, for $l
= 1, \ldots, n$.  Put $z_l = h(\delta_{\alpha(l)})$, so that $|z_l|
\le 1$, since $\|\delta_{\alpha(l)}\|_1 = 1$.  Thus $z = (z_1, \ldots,
z_n) \in \overline{U}^n$, and $h(\delta_\alpha) = z^\alpha$ for every
$\alpha \in {\bf Z}^n$ with $\alpha_j \ge 0$ for each $j$, because $h$
is a homomorphism with respect to convolution on $\ell^1_A({\bf
Z}^n)$.  More precisely, this uses the fact that
\begin{equation}
        \delta_\alpha * \delta_\beta = \delta_{\alpha + \beta}
\end{equation}
for every $\alpha, \beta \in {\bf Z}^n$.  This implies that $h(a)$ is
equal to (\ref{sum_alpha a(alpha) z^alpha}) for every $a$ in
$\ell^1_A({\bf Z}^n)$ and this choice of $z$, by the linearity and
continuity of $h$.

        If $h$ were a homomorphism on all of $\ell^1({\bf Z}^n)$, then
we would have (\ref{|h(a)| le ||a||_1}) for every $a \in \ell^1({\bf
Z}^n)$, which would imply that $|z_l| = 1$ for each $l$.  This is
because $\delta_{\alpha(l)} * \delta_{-\alpha(l)} = \delta_0$, so that
\begin{equation}
        z_l \, h(\delta_{-\alpha(l)})
         = h(\delta_{\alpha(l)}) \, h(\delta_{-\alpha(l)}) = 1,
\end{equation}
while $|h(\delta_{-\alpha(l)})| \le 1$ by (\ref{|h(a)| le ||a||_1}).
In this case, we would get that $h(a)$ is equal to $\widehat{a}(z)$ as
in (\ref{widehat{a}(z) = sum_{alpha in {bf Z}^n} a(alpha) z^alpha})
for every $a \in \ell^1({\bf Z}^n)$, by essentially the same argument
as before.  Of course, $a \mapsto \widehat{a}(z)$ defines a
homomorphism on $\ell^1({\bf Z}^n)$ for every $z \in {\bf T}^n$, as in
the previous section.

        Similarly, if $h$ is a nonzero homomorphism on all of $C({\bf
T}^n)$ and $f_j(z) = z_j$, then $|h(f_j)| = 1$ for each $j$, because
$z_j^{-1}$ is also a continuous function on ${\bf T}^n$ with supremum
norm equal to $1$.  Using this, one can show that $h(f) = f(p)$ for
every $f \in C({\bf T}^n)$, where $p = (p_1, \ldots, p_n) \in {\bf
T}^n$ is defined by $p_j = h(f_j)$, in the same way as before.
Although this is a special case of the results discussed in Section
\ref{compact spaces}, the present approach has the advantage of making
the relationship with $A({\bf T}^n)$ more clear.

\section{The maximum principle}
\label{maximum principle}
\setcounter{equation}{0}

        Let $D$ be a nonempty bounded connected open set in the
complex plane ${\bf C}$.  If $f$ is a continuous complex-valued
function on the closure $\overline{D}$ of $D$ in ${\bf C}$, then the
extreme value theorem implies that the $|f(z)|$ attains its maximum on
$\overline{D}$.  If $f$ is also holomorphic on $D$, then the
\emph{maximum modulus principle} implies that the maximum of $|f(z)|$
on $\overline{D}$ is attained on the boundary $\partial D$ of $D$.
More precisely, if $|f(z)|$ has a local maximum on $D$, then $f$ is
constant.  This follows from the fact that a nonconstant holomorphic
function on a connected open set in ${\bf C}$ is an open mapping, in
the sense that it maps open sets to open sets.

        Alternatively, suppose that $z \in D$, and that the closed
disk centered at $z$ with radius $r > 0$ is contained in $D$.  If $f$
is holomorphic on $D$, then
\begin{equation}
\label{f(z) = frac{1}{2 pi r} int_{|w - z| = r} f(w) |dw|}
        f(z) = \frac{1}{2 \pi r} \int_{|w - z| = r} f(w) \, |dw|,
\end{equation}
by the Cauchy integral formula.  If $|f|$ has a local maximum at $z$,
then one can use this to show that $f(w) = f(z)$ when $|w - z|$ is
sufficiently small, and hence that $f$ is constant on $D$ when $D$ is
connected.  This identity is known as the ``mean value property'',
since it says that the value of $f$ at $z$ is given by the average of
$f$ on the circle $|w - z| = r$.  This also works for harmonic
functions on $D$, and the analogous statement for harmonic functions
on open subsets of ${\bf R}^n$ holds for every $n$.  In particular, if
$f$ is a harmonic function on a connected open set $D \subseteq {\bf
R}^n$, and if $|f|$ has a local maximum on $D$, then one can show that
$f$ is constant on $D$. Similarly, if $f$ is a real-valued harmonic
function on $D$ with a local maximum on $D$, then $f$ is constant on
$D$.

        A holomorphic function of several complex variables is
holomorphic in each variable separately.  In particular, such a
function is harmonic, but one can get stronger versions of the maximum
principle by considering restrictions of the function to complex
lines, or even ``analytic disks'' that do not have to be flat.

        Suppose for instance that $D$ is the unit polydisk $U^n$.  If
$f$ is a continuous complex-valued function on $\overline{U}^n$ that
is holomorphic on $U^n$, then one can show that the maximum of
$|f(z)|$ on $\overline{U}^n$ is actually attained on ${\bf T}^n$.
This is the same as the boundary of $U^n$ when $n = 1$, but otherwise
is significantly smaller, as mentioned previously.  This version of
the maximum principle was implicitly given already in (\ref{sup_{z in
overline{U}^n} |phi(z)| = sup_{z in {bf T}^n} |f(z)|}), using Poisson
integrals.  This also works for functions that are polyharmonic
instead of holomorphic, which is to say harmonic in $z_j$ for $j = 1,
\ldots, n$.  This can also be derived from the maximum principle for
the unit disk, by looking at restrictions of the function to disks in
which all but one variable is constant.

\section{Convex hulls}
\label{convex hulls}
\setcounter{equation}{0}

        Let $A$ be a nonempty subset of ${\bf R}^n$ for some positive
integer $n$.  The \emph{convex hull} of $A$ is denoted $\con(A)$ and is
defined to be the set of $x \in {\bf R}^n$ for which there are
finitely many elements $y_1, \ldots, y_l$ of $A$ and nonnegative real
numbers $t_1, \ldots, t_l$ such that $\sum_{j = 1}^l t_j = 1$ and
\begin{equation}
        x = \sum_{j = 1}^l t_j \, y_j.
\end{equation}
It is easy to see that $\con(A)$ is a convex set in ${\bf R}^n$, and
that $\con(A) \subseteq B$ whenever $A \subseteq B$ and $B \subseteq
{\bf R}^n$ is convex.  Thus $\con(A)$ is the smallest convex set in
${\bf R}^n$ that contains $A$.

        It is well known that every element of $\con(A)$ can be
expressed as a convex combination of less than or equal to $n + 1$
elements of $A$.  This uses the fact that ${\bf R}^n$ is an
$n$-dimensional real vector space, while the definition of the convex
hull and the other remarks in the previous paragraph would work just
as well in any real vector space.  Using this, one can show that $\con
(A)$ is compact when $A \subseteq {\bf R}^n$ is compact.  Otherwise, the
\emph{closed convex hull} of $A$ is defined to be the closure of the
convex hull of $A$, and is automatically convex, because the closure
of any convex set in ${\bf R}^n$ is also convex.  This is the smallest
closed convex set that contains $A$, because any closed convex set
that contains $A$ also contains $\con(A)$ and hence $\overline{\con(A)}$.

        If $E$ is a nonempty closed convex set in ${\bf R}^n$ and $x
\in {\bf R}^n \backslash E$, then a well-known separation theorem
states that there is a linear function $\lambda(y)$ on ${\bf R}^n$
such that
\begin{equation}
        \sup_{y \in E} \lambda(y) < \lambda(x).
\end{equation}
To see this, observe first that there is an element $u$ of $E$ that
minimizes the distance to $x$ with respect to the standard Euclidean
metric, so that
\begin{equation}
        \sum_{j = 1}^n (y_j - x_j)^2 \ge \sum_{j = 1}^n (u_j - x_j)^2
\end{equation}
for every $y \in E$.  This follows immediately from the extreme value
theorem when $E$ is compact, and otherwise one can reduce to that case
by considering the intersection of $E$ with a closed ball centered at
$x$ with sufficiently large radius.  Without loss of generality, we
may suppose that $u = 0$, since otherwise we can translate everything
by $-u$ to reduce to this case.  Thus the previous inequality becomes
\begin{equation}
        \sum_{j = 1}^n (y_j - x_j)^2 \ge \sum_{j = 1}^n x_j^2,
\end{equation}
which holds for every $y \in E$.  Equivalently,
\begin{equation}
        \sum_{j = 1}^n y_j^2 \ge \sum_{j = 1}^n y_j \, x_j
\end{equation}
for every $y \in E$.  Because $u = 0 \in E$ and $E$ is convex, $t \, y
\in E$ for every $y \in E$ and $t \in [0, 1]$.  Hence
\begin{equation}
        \sum_{j = 1}^n (t \, y_j)^2 \ge \sum_{j = 1}^n (t \, y_j) \, x_j
\end{equation}
for every $y \in E$ and $0 \le t \le 1$.  This implies that
\begin{equation}
        t \sum_{j = 1}^n y_j^2 \ge \sum_{j = 1}^n y_j \, x_j
\end{equation}
when $y \in E$ and $0 < t \le 1$.  Taking the limit as $t \to 0$, we get that
\begin{equation}
        \sum_{j = 1}^n y_j \, x_j \le 0
\end{equation}
for every $y \in E$.  Put $\lambda(y) = \sum_{j = 1}^n y_j \, x_j$, so
that $\lambda(y) \le 0$ for every $y \in E$, by the preceding
inequality.  Note that $x \ne u = 0$, because $x \not\in E$ and $u \in
E$.  Thus we also have that $\lambda(x) = \sum_{j = 1}^n x_j^2 > 0$,
as desired.

        Let $A$ be a nonempty set in ${\bf R}^n$, and let $\lambda$ be
a linear function on ${\bf R}^n$.  Suppose that $x \in \con(A)$, so that
there are $y_1, \ldots y_l \in A$ and $t_1, \ldots, t_l \ge 0$ such
that $\sum_{j = 1}^j t_j = 1$ and $x = \sum_{j = 1}^l t_j \, y_j$.
In particular,
\begin{equation}
        \lambda(x) = \sum_{j = 1}^l t_j \, \lambda(y_j) 
         \le \max_{1 \le j \le l} \lambda(y_j).
\end{equation}
This implies that
\begin{equation}
\label{lambda(x) le sup_{y in A} lambda(y)}
        \lambda(x) \le \sup_{y \in A} \lambda(y),
\end{equation}
where the supremum on the right side may be $+\infty$, in which case
the inequality is trivial.  If $x \in \overline{\con(A)}$, then it is
easy to see that (\ref{lambda(x) le sup_{y in A} lambda(y)}) also
holds, by continuity.  However, if $x \in {\bf R}^n \backslash
\overline{\con(A)}$, then there is a linear function $\lambda$ on ${\bf
R}^n$ for which (\ref{lambda(x) le sup_{y in A} lambda(y)}) does not
hold, as in the previous paragraph.  Thus the closed convex hull of
$A$ is the same as the set of $x \in {\bf R}^n$ such that
(\ref{lambda(x) le sup_{y in A} lambda(y)}) holds for every linear
function $\lambda$ on ${\bf R}^n$.

\section{Polynomial hulls}
\label{polynomial hulls}
\setcounter{equation}{0}

        Let $E$ be a nonempty subset of ${\bf C}^n$ for some positive
integer $n$.  The \emph{polynomial hull} of $E$ in ${\bf C}^n$ is
denoted $\pol(E)$ and defined to be the set of $z \in {\bf C}^n$ such that
\begin{equation}
\label{|p(z)| le sup_{w in E} |p(w)|}
        |p(z)| \le \sup_{w \in E} |p(w)|
\end{equation}
for every polynomial $p$ on ${\bf C}^n$.  More precisely, to say that
$p$ is a polynomial on ${\bf C}^n$ means that $p$ can be expressed as
\begin{equation}
        p(w) = \sum_{|\alpha| \le N} a_\alpha \, w^\alpha
\end{equation}
for some nonnegative integer $N$, where the sum is taken over all multi-indices
$\alpha$ with $|\alpha| \le N$, and $a_\alpha \in {\bf C}$ for each $\alpha$.
If $E$ is unbounded, then $p$ may be unbounded on $E$, so that the supremum
in (\ref{|p(z)| le sup_{w in E} |p(w)|}) is $+\infty$, and the inequality
is trivial.  

        Of course,
\begin{equation}
        E \subseteq \pol(E)
\end{equation}
by definition.  If $E_1 \subseteq E_2 \subseteq {\bf C}^n$, then
\begin{equation}
        \pol(E_1) \subseteq \pol(E_2).
\end{equation}
It is easy to see that $\pol(E)$ is always a closed set in ${\bf
C}^n$, because polynomials are continuous.  Similarly,
\begin{equation}
\label{pol(E) = pol(overline{E})}
        \pol(E) = \pol(\overline{E}),
\end{equation}
and so we may as well restrict our attention to closed sets $E
\subseteq {\bf C}^n$.

        If $E$ is any nonempty subset of ${\bf C}^n$ and $p$ is a
polynomial on ${\bf C}^n$, then
\begin{equation}
        \sup_{z \in \pol(E)} |p(z)| = \sup_{w \in E} |p(w)|.
\end{equation}
More precisely, the right side is less than or equal to the left side
because $E \subseteq \pol(E)$, while the opposite inequality follows
from the definition of $\pol(E)$.  If $\zeta \in \pol(\pol(E))$, then
we get that
\begin{equation}
 |p(\zeta)| \le \sup_{z \in \pol(E)} |p(z)| = \sup_{w \in E} |p(w)|
\end{equation}
for every polynomial $p$ on ${\bf C}^n$, which implies that $\zeta \in
\pol(E)$.  Thus $\pol(\pol(E))$ is contained in $\pol(E)$, and hence
\begin{equation}
        \pol(\pol(E)) = \pol(E),
\end{equation}
because $\pol(E) \subseteq \pol(\pol(E))$ automatically.

        As an example, let us check that
\begin{equation}
        \pol({\bf T}^n) = \overline{U}^n.
\end{equation}
If $z \in \overline{U}^n$, then $z \in \pol({\bf T}^n)$, as in Section
\ref{maximum principle}, and so $\overline{U}^n \subseteq \pol({\bf
T}^n)$.  However, if $z \in {\bf C}^n \backslash \overline{U}^n$, then
$|z_j| > 1$, and one can check that $z \not\in \pol({\bf T}^n)$, by
taking $p(w) = w_j$.  Thus $\pol({\bf T}^n) \subseteq \overline{U}^n$,
as desired.

        If $E$ is any nonempty bounded subset of ${\bf C}^n$, then
\begin{equation}
\label{pol(E) subseteq overline{con(E)}}
        \pol(E) \subseteq \overline{\con(E)}.
\end{equation}
To see this, we identify ${\bf C}^n$ with ${\bf R}^{2n}$ as a real
vector space, so that the results in the previous section are
applicable.  If $z \in {\bf C}^n \backslash \overline{\con(E)}$, then
there is a real-valued real-linear function $\lambda$ on ${\bf C}^n
\cong {\bf R}^{2 n}$ such that
\begin{equation}
        \sup_{w \in E} \lambda(w) < \lambda(z),
\end{equation}
as in the previous section.  Equivalently, $\lambda$ can be expressed
as the real part of a complex-linear function $\mu$ on ${\bf C}^n$,
and
\begin{equation}
\label{sup_{w in E} re mu(w) < re mu(z)}
        \sup_{w \in E} \re \mu(w) < \re \mu(z).
\end{equation}
We would like to show that
\begin{equation}
\label{sup_{w in E} |1 + t mu(w)| < |1 + t mu(z)|}
        \sup_{w \in E} |1 + t \, \mu(w)| < |1 + t \, \mu(z)|
\end{equation}
when $t$ is a sufficiently small positive real number, so that $z
\not\in \pol(E)$.

        Note that
\begin{eqnarray}
\, |1 + t \, \mu(w)|^2 & = & (1 + t \, \re \mu(w))^2 + t^2 \, (\im \mu(w))^2 \\
 & = & 1 + 2 \, t \, \re \mu(w) + t^2 \, (\re \mu(w))^2 + t^2 \, (\im \mu(w))^2
                                                                   \nonumber
\end{eqnarray}
for every $t \in {\bf R}$ and $w \in {\bf C}^n$.  Because $E$ is bounded,
\begin{equation}
        |\mu(w)|^2 = (\re \mu(w))^2 + (\im \mu(w))^2 \le C
\end{equation}
for some $C \ge 0$ and every $w \in E$.  Thus
\begin{equation}
\label{sup_{w in E} |1 + t mu(w)|^2 le 1 + 2 t sup_{w in E} re mu(w) + C t^2}
        \sup_{w \in E} |1 + t \, \mu(w)|^2
         \le 1 + 2 \, t \, \sup_{w \in E} \re \mu(w) + C \, t^2
\end{equation}
for every $t > 0$.  Using this and (\ref{sup_{w in E} re mu(w) < re
mu(z)}), it is easy to see that
\begin{equation}
\label{sup_{w in E} |1 + t mu(w)|^2 < |1 + t mu(z)|^2}
        \sup_{w \in E} |1 + t \, \mu(w)|^2 < |1 + t \, \mu(z)|^2
\end{equation}
when $t > 0$ is sufficiently small, as desired.

        If $n = 1$ and $E \subseteq {\bf C}$ is unbounded, then every
nonconstant polynomial $p$ on ${\bf C}$ is unbounded on $E$.  Thus
$\pol(E) = {\bf C}$ in this case.  In particular, $\pol(E)$ may not be
contained in $\overline{\con(E)}$ when $E$ is unbounded.

        Suppose that $E$ is a nonempty set in ${\bf C}^n$ with only
finitely many elements.  If $z \in {\bf C}^n \backslash E$, then it is
easy to see that there is a polynomial $p$ on ${\bf C}^n$ such that
$p(w) = 0$ for each $w \in E$ and $p(z) \ne 0$, by taking a product of
affine functions that vanish at the elements of $E$, one at a time,
and are nonzero at $z$.  This implies that $z \not\in \pol(E)$, so
that $\pol(E) \subseteq E$.  Hence $\pol(E) = E$ when $E$ has only
finitely many elements, since $E \subseteq \pol(E)$ automatically.  By
contrast, the convex hull of a finite set may be much larger.

\section{Algebras and homomorphisms}
\label{algebras, homomorphisms}
\setcounter{equation}{0}

        Let $E$ be a nonempty compact set in ${\bf C}^n$, and let
$C(E)$ be the algebra of continuous complex-valued functions on $E$.
Let $PC(E)$ be the subalgebra of $C(E)$ consisting of the restrictions
to $E$ of polynomials on ${\bf C}^n$, and let $AC(E)$ be the closure
of $PC(E)$ in $C(E)$ with respect to the supremum norm.  Thus $AC(E)$
is a closed subalgebra of $C(E)$, and hence a commutative Banach
algebra with respect to the supremum norm, since $C(E)$ is.  Of
course, the constant function equal to $1$ on $E$ is the
multiplicative identity element in $C(E)$, and is contained in $PC(E)
\subseteq AC(E)$.

        Suppose that $h$ is a nonzero homomorphism from $AC(E)$ into
the complex numbers.  Let $f_j$ be the function on $E$ defined by
$f_j(w) = w_j$ for $j = 1, \ldots, n$, so that $f_j \in PC(E)
\subseteq AC(E)$ for each $j$.  Put
\begin{equation}
        z_j = h(f_j)
\end{equation}
for each $j$, and consider $z = (z_1, \ldots, z_n) \in {\bf C}^n$.  If
$p$ is any polynomial on ${\bf C}^n$, and $\widetilde{p}$ is the
restriction of $p$ to $E$, then $\widetilde{p} \in PC(E) \subseteq
AC(E)$, and
\begin{equation}
        h(\widetilde{p}) = p(z).
\end{equation}
As in Section \ref{banach algebras},
\begin{equation}
        |h(f)| \le \sup_{w \in E} |f(w)|
\end{equation}
for every $f \in AC(E)$.  It follows that
\begin{equation}
\label{|p(z)| = |h(widetilde{p})| le ... = sup_{w in E} |p(w)|}
 |p(z)| = |h(\widetilde{p})| \le \sup_{w \in E} |\widetilde{p}(w)|
                               = \sup_{w \in E} |p(w)|
\end{equation}
for every polynomial $p$ on ${\bf C}^n$.  Thus $z \in \pol(E)$.

        Conversely, suppose that $z \in \pol(E)$, so that
\begin{equation}
\label{|p(z)| le sup_{w in E} |p(w)| = sup_{w in E} |widetilde{p}(w)|}
        |p(z)| \le \sup_{w \in E} |p(w)|
                   = \sup_{w \in E} |\widetilde{p}(w)|
\end{equation}
for every polynomial $p$ on ${\bf C}^n$.  In particular, if $p(w) = 0$
for every $w \in E$, then $p(z) = 0$.  This implies that
\begin{equation}
\label{h_z(widetilde{p}) = p(z)}
        h_z(\widetilde{p}) = p(z)
\end{equation}
is well-defined on $PC(E)$, and in fact it is a homomorphism from
$PC(E)$ into the complex numbers.  Moreover, (\ref{|p(z)| le sup_{w in
E} |p(w)| = sup_{w in E} |widetilde{p}(w)|}) implies that $h_z$ is a
continuous linear functional on $PC(E)$ with respect to the supremum
norm, so that $h_z$ has a unique extension to a continuous linear
functional on $AC(E)$.  It is easy to see that this extension is also
a homomorphism with respect to multiplication.

        The argument in the preceding paragraph would work just as well if
\begin{equation}
        |p(z)| \le C \, \sup_{w \in E} |p(w)|
\end{equation}
for some $C \ge 0$ and every polynomial $p$ on ${\bf C}^n$.  Note that
$p^l$ is also a polynomial on ${\bf C}^n$ for every polynomial $p$ and
positive integer $l$.  Applying the previous condition to $p^l$, we
get that
\begin{equation}
        |p(z)|^l \le C \, \sup_{w \in E} |p(w)|^l.
\end{equation}
Equivalently,
\begin{equation}
        |p(z)| \le C^{1/l} \, \sup_{w \in E} |p(w)|
\end{equation}
for each $l \ge 1$ and polynomial $p$ on ${\bf C}^n$.  Taking the
limit as $l \to \infty$, it follows that the initial inequality holds
with $C = 1$.  Hence this apprently weaker condition implies that $z
\in \pol(E)$.  This could also be derived from the earlier discussion,
but this approach is more direct.

\section{The exponential function}
\label{exponential function}
\setcounter{equation}{0}

        Put
\begin{equation}
\label{E(z) = sum_{j = 0}^infty frac{z^j}{j!}}
        E(z) = \sum_{j = 0}^\infty \frac{z^j}{j!}
\end{equation}
for each $z \in {\bf C}$, where $j!$ is ``$j$ factorial'', the product
of $1, \ldots, j$.  As usual, this is interpreted as being equal to
$1$ when $j = 0$.  It is easy to see that this series converges
absolutely for every $z \in {\bf C}$, by the ratio test, for instance.

        If $z, w \in {\bf C}$, then
\begin{equation}
\label{E(z) E(w) = ...}
 E(z) \, E(w) = \Big(\sum_{j = 0}^\infty \frac{z^j}{j!}\Big) \,
                 \Big(\sum_{l = 0}^\infty \frac{w^l}{l!}\Big)
              = \sum_{n = 0}^\infty \Big(\sum_{j = 0}^n
                             \frac{z^j \, w^{n - j}}{j! \, (n - j)!}\Big),
\end{equation}
as in Section \ref{cauchy products}.  This uses the absolute
convergence of the series defining $E(z)$ and $E(w)$.  The binomial
theorem states that
\begin{equation}
\label{sum_{j = 0}^n frac{n!}{j! (n - j)!} z^j w^{n - j} = (z + w)^n}
 \sum_{j = 0}^n \frac{n!}{j! \, (n - j)!} \, z^j \, w^{n - j} = (z + w)^n,
\end{equation}
so that
\begin{equation}
\label{E(z) E(w) = sum_{n = 0}^infty frac{(z + w)^n}{n!} = E(z + w)}
        E(z) \, E(w) = \sum_{n = 0}^\infty \frac{(z + w)^n}{n!} = E(z + w).
\end{equation}
In particular,
\begin{equation}
        E(z) \, E(-z) = E(0) = 1
\end{equation}
for every $z \in {\bf C}$.  Equivalently, $E(z) \ne 0$ for every $z
\in {\bf C}$, and $1/E(z) = E(-z)$.  If $x$ is a nonnegative real
number, then it is clear from the definition of $E(x)$ that $E(x) \in
{\bf R}$ and $E(x) \ge 1$.  It follows that $E(x)$ is a positive real
number for every $x \in {\bf R}$, and that $E(x) \le 1$ when $x \le
0$.  Similarly, it is easy to see from the definition that $E(x)$ is
strictly increasing when $x \ge 0$, and one can extend this to the
whole real line using the fact that $E(-x) = 1/E(x)$.

        Observe that
\begin{equation}
        \overline{E(z)} = E(\overline{z})
\end{equation}
for every $z \in {\bf C}$, by the definition of $E(z)$.  This implies that
\begin{equation}
        |E(z)|^2 = E(z) \, \overline{E(z)} = E(z) \, E(\overline{z})
                 = E(z + \overline{z}) = E(2 \re z),
\end{equation}
and hence
\begin{equation}
\label{|E(z)| = E(re z)}
        |E(z)| = E(\re z)
\end{equation}
for every $z \in {\bf C}$.

        If $z = i \, y$ for some $y \in {\bf R}$, then (\ref{|E(z)| =
E(re z)}) implies that
\begin{equation}
        |E(i \, y)| = 1.
\end{equation}
It is well known that
\begin{equation}
        E(i \, y) = \cos y + i \sin y
\end{equation}
for every $y \in {\bf R}$.  One way to see this is to use the standard
power series expansions for the sine and cosine.  Alternatively,
\begin{equation}
        \frac{d}{dy} E(i \, y) = i \, E(i \, y),
\end{equation}
as one can check using the series expansion for $E(i \, y)$.  We
already know that $E(i \, y)$ maps the real line into the unit circle
${\bf T}$ and sends $y = 0$ to $1$.  This formula for the derivative
of $E(i \, y)$ shows that it goes around the circle at unit speed in
the positive orientation.  It follows that the real and imaginary
parts of $E(i \, y)$ are given by the cosine and sine, respectively,
by the geometric definitions of the cosine and sine.

\section{Entire functions}
\label{entire functions}
\setcounter{equation}{0}

        Suppose that
\begin{equation}
\label{sum_alpha a_alpha z^alpha}
        \sum_\alpha a_\alpha \, z^\alpha
\end{equation}
is a power series with complex coefficients that is absolutely summable
for every $z \in {\bf C}^n$, and let $f(z)$ be the sum of this series.
If $E$ is a nonempty bounded set in ${\bf C}^n$, then
\begin{equation}
        |f(z)| \le \sup_{w \in E} |f(w)|
\end{equation}
for every $z \in \pol(E)$.  This follows from the definition of the
polynomial hull, and the fact that $f$ can be approximated uniformly
by polynomials corresponding to finite subsums of (\ref{sum_alpha
a_alpha z^alpha}) on bounded subsets of ${\bf C}^n$.

        Let $u = (u_1, \ldots, u_n) \in {\bf C}^n$ be given, and put
\begin{equation}
        \mu(z) = \sum_{j = 1}^n u_j \, z_j.
\end{equation}
Observe that $f_\mu(z) = E(\mu(z))$ can be expressed as
\begin{equation}
\label{prod_{j = 1}^n E(u_j z_j) = sum_alpha frac{u^alpha}{alpha!} z^alpha}
        \prod_{j = 1}^n E(u_j \, z_j)
          = \sum_\alpha \frac{u^\alpha}{\alpha!} \, z^\alpha,
\end{equation}
where $\alpha! = \alpha_1! \cdots \alpha_n!$.  In particular, this
power series is absolutely summable for every $z \in {\bf C}^n$.  Moreover,
\begin{equation}
        |f_\mu(z)| = E(\re \mu(z)),
\end{equation}
as in the previous section.

        If $E$ is a nonempty subset of ${\bf C}^n$ and $z \in {\bf
C}^n \backslash \overline{\con(E)}$, then there is a complex-linear
function $\mu$ on ${\bf C}^n$ such that
\begin{equation}
        \sup_{w \in E} \re \mu(w) < \re \mu(z),
\end{equation}
as in Section \ref{polynomial hulls}.  Hence
\begin{equation}
        \sup_{w \in E} |f_\mu(w)| < |f_\mu(z)|.
\end{equation}
This has the advantage of working for both bounded and unbounded sets
$E$, in exchange for allowing a larger class of functions than
polynomials, as before.

\section{The three lines theorem}
\label{three lines theorem}
\setcounter{equation}{0}

        Let $D$ be the open unit strip in the complex plane,
\begin{equation}
        D = \{z \in {\bf C} : 0 < \re z < 1\},
\end{equation}
so that the closure of $D$ is the closed unit strip,
\begin{equation}
        \overline{D} = \{z \in {\bf C} : 0 \le \re z \le 1\}.
\end{equation}
Also let $f$ be a continuous complex-valued function on $\overline{D}$
which is holomorphic on $D$, and suppose that $A_0$, $A_1$ are
positive real numbers such that
\begin{equation}
        |f(x + i \, y)| \le A_x
\end{equation}
for $x = 0, 1$ and every $y \in {\bf R}$.  If $f$ is also bounded on
$D$, then the \emph{three lines theorem} states that
\begin{equation}
\label{|f(x + i y)| le A_0^{1 - x} A_1^x}
        |f(x + i \, y)| \le A_0^{1 - x} \, A_1^x
\end{equation}
when $0 < x < 1$ and $y \in {\bf R}$.

        To do this, we would first like to show that $f$ satisfies the
maximum principle, so that
\begin{equation}
\label{|f(x + i y)| le max (A_0, A_1)}
        |f(x + i \, y)| \le \max (A_0, A_1)
\end{equation}
when $0 < x < 1$ and $y \in {\bf R}$.  However, $D$ is not bounded,
$\overline{D}$ is not compact, and so we cannot use the ordinary
maximum principle in quite the usual way.  Let us begin with the
case where $f$ satisfies the additional condition that $f(z)
\to 0$ uniformly on $\overline{D}$ as $|\im z| \to \infty$.  Put
\begin{equation}
        D_R = \{z \in D : |\im z| < R\}
\end{equation}
for each $R > 0$.  Thus $D_R$ is bounded, and we can apply the maximum
principle to $f$ on $D_R$.  If
\begin{equation}
        B_R = \sup \{|f(x + i \, y)| : 0 \le x \le 1, y = \pm R\},
\end{equation}
then we get that
\begin{equation}
        |f(x + i \, y)| \le \max (A_0, A_1, B_R)
\end{equation}
when $0 < x < 1$ and $|y| < R$.  By hypothesis, $B_R \to 0$ as $R \to
\infty$, and so (\ref{|f(x + i y)| le max (A_0, A_1)}) follows easily
in this case.

        If $f$ is bounded but does not necessarily tend to $0$ at
infinity, then we can approximate it by functions that do.  Consider
\begin{equation}
        f_\epsilon(z) = f(z) \, E(\epsilon \, z^2)
\end{equation}
for each $\epsilon > 0$.  Observe that
\begin{equation}
\label{|f_epsilon(z)| = ... = |f(z)| E(epsilon (x^2 - y^2))}
        |f_\epsilon(z)| = |f(z)| \, |E(\epsilon \, z^2)|
                        = |f(z)| \, E(\epsilon \, (x^2 - y^2)),
\end{equation}
where $z = x + i \, y$, and hence $\re z^2 = x^2 - y^2$.  Thus
$f_\epsilon(z)$ is continuous on $\overline{D}$, holomorphic on $D$,
and tends to $0$ uniformly as $|y| \to \infty$ for each $\epsilon >
0$, because $f$ is bounded on $\overline{D}$ by hypothesis and
$E(\epsilon \, y^2) \to +\infty$ as $y \to + \infty$ for each $\epsilon > 0$.
We also have that
\begin{equation}
        |f_\epsilon(i \, y)| \le A_0, \quad
         |f_\epsilon(i \, y)| \le A_1 \, E(\epsilon)
\end{equation}
for each $y \in {\bf R}$, since $E(-\epsilon \, y^2) \le 1$ for every
$y \in {\bf R}$.  This permits us to use the version of the maximum
principle in the previous paragraph, to get that
\begin{equation}
\label{|f_epsilon(x + i y)| le max(A_0, A_1 E(epsilon))}
        |f_\epsilon(x + i \, y)| \le \max(A_0, A_1 \, E(\epsilon))
\end{equation}
when $0 < x < 1$ and $y \in {\bf R}$.  Of course, $f_\epsilon(z) \to
f(z)$ for each $z \in \overline{D}$ as $\epsilon \to 0$, and it
follows that $f$ satisfies (\ref{|f(x + i y)| le max (A_0, A_1)}), by
taking the limit as $\epsilon \to 0$ in (\ref{|f_epsilon(x + i y)| le
max(A_0, A_1 E(epsilon))}).

        Note that $E(t) \ge 1 + t$ for every nonnegative real number
$t$, by the definition of $E(t)$.  Thus $E(t) \to +\infty$ as $t \to
+\infty$, as in the previous paragraph, and hence $E(-t) = 1/E(t) \to
0$ as $t \to +\infty$.  If $A$ is a positive real number, then it
follows that there is a unique real number $\log A$ such that $E(\log
A) = A$, because $E(t)$ is a strictly increasing continuous function
on the real line.  Put $A^z = E(z \, \log A)$, and observe that $|A^z|
= A^{\re z}$, by the properties of the exponential function.  In order
to get (\ref{|f(x + i y)| le A_0^{1 - x} A_1^x}), consider
\begin{equation}
        g(z) = f(z) \, A_0^{z - 1} \, A_1^{-z}.
\end{equation}
This is a bounded continuous function on $\overline{D}$ which is
holomorphic on $D$ and satisfies
\begin{equation}
\label{|g(x + i y)| le 1}
        |g(x + i \, y)| \le 1
\end{equation}
when $x = 0, 1$ and $y \in {\bf R}$, by the corresponding properties of $f$.
The analogue of (\ref{|f(x + i y)| le max (A_0, A_1)}) for $g$ implies
that (\ref{|g(x + i y)| le 1}) holds for every $0 < x < 1$ and $y \in {\bf
R}$, which is the same as (\ref{|f(x + i y)| le A_0^{1 - x} A_1^x}).

\section{Completely circular sets}
\label{completely circular sets}
\setcounter{equation}{0}

        A set $E \subseteq {\bf C}^n$ is said to be \emph{completely
circular} if
\begin{equation}
        (u_1 \, z_1, \ldots, u_n \, z_n) \in E
\end{equation}
for every $z = (z_1, \ldots, z_n) \in E$ and $u = (u_1, \ldots, u_n)
\in {\bf C}^n$ such that $|u_j| \le 1$ for each $j$.  Equivalently, $w
\in E$ whenever $w \in {\bf C}^n$ satisfies $|w_j| \le |z_j|$ for some
$z \in E$ and each $j$.  In particular, $0 \in E$ when $E \ne
\emptyset$.

        Suppose that $w, z \in E$, $0 < t < 1$, and $v \in {\bf C}^n$ satisfy
\begin{equation}
        |v_j| \le |z_j|^t \, |w_j|^{1 - t}
\end{equation}
for each $j$.  We would like to show that $v \in \pol(E)$ when $E$ is
completely circular.  Thus we would like to show that
\begin{equation}
        |p(v)| \le \sup_{\zeta \in E} |p(\zeta)|
\end{equation}
for every polynomial $p$ on ${\bf C}^n$.  To do this, we shall use the
version of the maximum principle discussed in the previous section.

        By hypothesis, we can express $v$ as
\begin{equation}
        v_j = u_j \, |z_j|^t \, |w_j|^{1 - t},
\end{equation}
where $|u_j| \le 1$ for each $j$.  Put
\begin{equation}
        g_j(\tau) = u_j \, |z_j|^\tau \, |w_j|^{1 - \tau}
\end{equation}
for each $\tau \in {\bf C}$ and $1 \le j \le n$.  This uses the
definition of $A^\tau$ for any positive real number $A$ and complex
number $\tau$ as $E(\tau \, \log A)$, as in the previous section, and
we put $A^\tau = 0$ for every $\tau \in {\bf C}$ when $A = 0$.  Note that
\begin{equation}
        |g_j(\tau)| \le |z_j|^{\re \tau} \, |w_j|^{1 - \re \tau}
\end{equation}
for each $\tau$ and $j$, and hence
\begin{equation}
        g(\tau) = (g_1(\tau), \ldots, g_n(\tau)) \in E
\end{equation}
when $\re \tau = 0, 1$.  We also have that $g(t) = v$, by
construction.

        Let $p$ be a polynomial on ${\bf C}^n$, and consider
\begin{equation}
        f(\tau) = p(g(\tau)).
\end{equation}
This is a holomorphic function on the complex plane ${\bf C}$, and in
particular it is a holomorphic function on the open unit strip $D$
that extends continuously to the closure $\overline{D}$.  Moreover,
$f$ is bounded on $\overline{D}$, because $g$ is bounded on
$\overline{D}$, and $p$ is bounded on bounded subsets of ${\bf C}^n$.
It follows that
\begin{equation}
        |f(t)| \le \sup \{|f(\tau)| : \tau \in {\bf C}, \ \re \tau = 0, 1\},
\end{equation}
as in the previous section.  This is the same as saying that
\begin{equation}
        |p(v)| \le \sup \{|p(g(\tau))| : \tau \in {\bf C}, \ \re \tau = 0, 1\},
\end{equation}
which is exactly what we wanted, since $g(\tau) \in E$ when $\re \tau = 0, 1$.

\section{Completely circular sets, continued}
\label{completely circular sets, continued}
\setcounter{equation}{0}

        Let $E$ be a nonempty bounded completely circular set in ${\bf
C}^n$, and let $z \in {\bf C}^n$ be given.  Suppose that $z \ne 0$,
and let $I$ be the set of $j = 1, \ldots, n$ such that $z_j \ne 0$.
Let $\alpha(I) = (\alpha_1(I), \ldots, \alpha_n(I))$ be the
multi-index defined by $\alpha_j(I) = 1$ when $j \in I$, and
$\alpha_j(I) = 0$ otherwise.  Thus $z^{\alpha(I)} \ne 0$, and if
$w^{\alpha(I)} = 0$ for every $w \in E$, then $z \not\in \pol(E)$.

        Let $E_I$ be the set of $w \in E$ such that $w_j \ne 0$ when
$j \in I$, and suppose from now on in this section that $E_I \ne
\emptyset$.  Also let ${\bf R}^I$ be the set of real-valued functions
on $I$, which is basically the same as ${\bf R}^l$, where $l$ is the
number of elements of $I$.  Thus $\log |w_j|$, $j \in I$, determines
an element of ${\bf R}^I$ for each $w \in E_I$, and we let $A_I$ be
the subset of ${\bf R}^I$ corresponding to elements of $E_I$ in this
way.  If $r \in {\bf R}^I$, $t \in A_I$, and
\begin{equation}
\label{r_j le t_j}
        r_j \le t_j
\end{equation}
for each $j \in I$, then $r \in A_I$, because $E$ is completely circular.

        Let $\zeta$ be the element of ${\bf R}^I$ given by $\zeta_j =
\log |z_j|$ for $j \in I$, and suppose that $\zeta$ is not an element
of the closure of the convex hull of $A_I$ in ${\bf R}^I$.  This implies
that there is a linear function $\lambda$ on ${\bf R}^I$ such that
\begin{equation}
        \sup_{r \in A_I} \lambda(r) < \lambda(\zeta).
\end{equation}
More precisely, $\lambda(r)$ can be given as
\begin{equation}
        \lambda(r) = \sum_{j \in I} \lambda_j \, r_j
\end{equation}
for some $\lambda_j \in {\bf R}$, $j \in I$, and every $r \in {\bf R}^I$.
If $\lambda_j < 0$ for some $j \in I$, then there are $r \in A_I$ for which
$\lambda(r)$ is arbitrarily large, because $A_I$ is nonempty and satisfies
the condition  mentioned at the end of the previous paragraph.  Thus
\begin{equation}
        \lambda_j \ge 0
\end{equation}
for each $j \in I$, since $\lambda(r)$ is bounded from above for $r \in A_I$.

        If $z \in \pol(E)$, then
\begin{equation}
        |z^\alpha| \le \sup_{w \in E} |w^\alpha|
\end{equation}
for every multi-index $\alpha$.  Let us restrict our attention to
multi-indices $\alpha$ such that $\alpha_j \ge 1$ when $j \in I$ and
$\alpha_j = 0$ otherwise, so that $w^\alpha = 0$ when $w \in E
\backslash E_I$.  In this case, the previous inequality reduces to
\begin{equation}
        \sum_{j \in I} \alpha_j \log |z_j|
         \le \sup_{w \in E_I} \sum_{j \in I} \alpha_j \, \log |w_j|.
\end{equation}
Equivalently,
\begin{equation}
        \sum_{j \in I} \alpha_j \, \zeta_j
         \le \sup_{r \in A_I} \sum_{j \in I} \alpha_j \, r_j.
\end{equation}
This inequality holds for arbitrary positive integers $\alpha_j$, $j
\in I$, and hence for arbitrary positive rational numbers $\alpha_j$,
by dividing both sides by a positive integer.  It follows that this
inequality also holds for arbitrary nonnegative real numbers, by
approximation.  This uses the hypothesis that $E$ be bounded, so that
$r_j$ has an upper bound for each $j \in I$ and $r \in A_I$, and more
precisely one should approximate nonnegative real numbers $\alpha_j$
by positive rational numbers $\alpha_j'$ such that $\alpha_j \le
\alpha_j'$ for each $j$.  Combining this with the discussion in the
previous paragraph, we get that $\zeta$ is in the closure of the
convex hull of $A_I$ in ${\bf R}^I$ when $z \in \pol(E)$.

\section{The torus action}
\label{torus action}
\setcounter{equation}{0}

        Let ${\bf T}^n$ be the set of $t = (t_1, \ldots, t_n) \in {\bf
C}^n$ such that $|t_j| = 1$ for each $j$, as usual.  If $t \in {\bf
T}^n$ and $z \in {\bf C}^n$, then put
\begin{equation}
        T_t(z) = (t_1 \, z_1, \ldots, t_n \, z_n),
\end{equation}
so that $T_t$ is an invertible linear transformation on ${\bf C}^n$
for each $t \in {\bf T}^n$.  Note that ${\bf T}^n$ is a commutative
group with respect to coordinatewise multiplication, and that $t
\mapsto T_t$ is a homomorphism from ${\bf T}^n$ into the group of
invertible linear transformations on ${\bf C}^n$.  Suppose that $E$ is
a nonempty subset of ${\bf C}^n$ such that
\begin{equation}
\label{T_t(E) subseteq E}
        T_t(E) \subseteq E
\end{equation}
for every $t \in {\bf T}^n$.  This implies that
\begin{equation}
\label{T_t(E) = E}
        T_t(E) = E
\end{equation}
for each $t \in {\bf T}^n$, because $T_t^{-1}(E) = T_{t^{-1}}(E)
\subseteq E$, where $t^{-1} = (t_1^{-1}. \ldots, t_n^{-1})$.

        Suppose that $z \in \pol(E)$, so that
\begin{equation}
        |p(z)| \le \sup_{w \in E} |p(w)|
\end{equation}
for every polynomial $p$ on ${\bf C}^n$.  If $p_t(w) = p(T_t(w))$,
then $p_t$ is also a polynomial on ${\bf C}^n$ for each $t \in {\bf
T}^n$, and hence
\begin{equation}
\label{|p_t(z)| le sup_{w in E} |p_t(w)|}
        |p_t(z)| \le \sup_{w \in E} |p_t(w)|.
\end{equation}
We also have that
\begin{equation}
        \sup_{w \in E} |p_t(w)| = \sup_{w \in E} |p(w)|
\end{equation}
for every $t \in {\bf T}^n$, because of (\ref{T_t(E) = E}).  Thus
\begin{equation}
\label{|p(T_t(z))| = |p_t(z)| le sup_{w in E} |p_t(w)| = sup_{w in E} |p(w)|}
 |p(T_t(z))| = |p_t(z)| \le \sup_{w \in E} |p_t(w)| = \sup_{w \in E} |p(w)|
\end{equation}
for every polynomial $p$ on ${\bf C}^n$ and $t \in {\bf T}^n$, which
implies that $T_t(z) \in \pol(E)$ for every $t \in {\bf T}^n$.  This
shows that $T_t(\pol(E)) \subseteq \pol(E)$ for every $t \in {\bf
T}^n$, and hence $T_t(\pol(E)) = \pol(E)$ for every $t \in {\bf T}^n$,
as before.

        Let $U^n$ be the open unit polydisk in ${\bf C}^n$, consisting
of $u \in {\bf C}^n$ such that $|u_j| < 1$ for each $j$.  If $p$ is a
polynomial on ${\bf C}^n$, $u \in U^n$, and $w \in {\bf C}^n$, then
\begin{equation}
        |p(u_1 \, w_1, \ldots, u_n \, w_n)|
         \le \sup_{t \in {\bf T}^n} |p(t_1 \, w_1, \ldots, t_n \, w_n)|
\end{equation}
More precisely, we can think of $p(u_1 \, w_1, \ldots, u_n \, w_n)$ as
a polynomial in $u$ for each $w \in {\bf C}^n$, and apply the maximum
principle as in Section \ref{maximum principle}.  If $z \in \pol(E)$
and $u \in U^n$, then we get that
\begin{equation}
\label{|p(u_1 z_1, ldots, u_n z_n)| le ... le sup_{w in E} |p(w)|}
 |p(u_1 \, z_1, \ldots, u_n \, z_n)| \le \sup_{t \in {\bf T}^n} |p_t(z)|
                                      \le \sup_{w \in E} |p(w)|
\end{equation}
for every polynomial $p$ on ${\bf C}^n$, where the second step is as
in the previous paragraph.  Thus
\begin{equation}
\label{(u_1 z_1, ldots, u_n z_n) in pol(E)}
        (u_1 \, z_1, \ldots, u_n \, z_n) \in \pol(E)
\end{equation}
for every $z \in \pol(E)$ and $u \in U^n$, and it follows that $\pol(E)$
is completely circular in this case.

\section{Another condition}
\label{another condition}
\setcounter{equation}{0}

        Let $E$ be a nonempty completely circular closed set in ${\bf
C}^n$ such that
\begin{equation}
\label{E^* = {w in E : w_j ne 0 hbox{ for each } j}}
        E^* = \{w \in E : w_j \ne 0 \hbox{ for each } j\}
\end{equation}
is dense in $E$.  This happens when $E$ is the closure of a nonempty
completely circular open set in ${\bf C}^n$, for instance.  Put
\begin{equation}
        A = \{y \in {\bf R}^n : y_j = \log |w_j| \hbox{ for some } w \in E^*
                                                  \hbox{ and each } j\}.
\end{equation}
so that
\begin{equation}
 \quad E^* = \{w \in {\bf C}^n : w_j \ne 0 \hbox{ for each } j, \hbox{ and }
                                (\log |w_1|, \ldots, \log |w_n|) \in A\},
\end{equation}
because $E$ is completely circular.  Thus $E = \overline{E^*}$ is
uniquely determined by $A$ under these conditions.  If $x \in {\bf
R}^n$, $y \in A$, and
\begin{equation}
        x_j \le y_j
\end{equation}
for each $j$, then we also have that $x \in A$, since $E$ is
completely circular.

        Let $I$ be a nonempty subset of $\{1, \ldots, n\}$, and let
${\bf R}^I$ be the set of real-valued functions on $I$, as before.
There is a natural projection from ${\bf R}^n$ onto ${\bf R}^I$, in
which one keeps the coordinates corresoponding to $j \in I$ and drops
the others.  Let $E_I$ be the set of $w \in E$ such that $w_j \ne 0$
when $j \in I$, and let $A_I$ be the subset of ${\bf R}^I$ whose
elements correspond to $\log |w_j|$, $j \in I$, with $w \in E_I$.
Observe that
\begin{equation}
        \pi_I(A) \subseteq A_I \subseteq \overline{\pi_I(A)},
\end{equation}
because $E^* \subseteq E_I$ and $E^*$ is dense in $E$.

        If $z \in \pol(E)$, then
\begin{equation}
        |z^\alpha| \le \sup_{w \in E} |w^\alpha|
\end{equation}
for every multi-index $\alpha$.  Moreover,
\begin{equation}
        \sup_{w \in E^*} |w^\alpha| = \sup_{w \in E} |w^\alpha|,
\end{equation}
since $E^*$ is dense in $E$.  Hence
\begin{equation}
\label{sup_{w in E_I} |w^alpha| = sup_{w in E} |w^alpha|}
        \sup_{w \in E_I} |w^\alpha| = \sup_{w \in E} |w^\alpha|
\end{equation}
for any $I \subseteq \{1, \ldots, n\}$, because $E_I \subseteq E^*
\subseteq E$.  Of course, (\ref{sup_{w in E_I} |w^alpha| = sup_{w in
E} |w^alpha|}) is trivial when $\alpha_j \ge 1$ for each $j \in I$, so
that $w^\alpha = 0$ when $w \in E \backslash E_I$.

        Let us consider some examples in ${\bf C}^2$ where $E^*$ is
not dense in $E$.  If
\begin{equation}
        E = ({\bf C} \times \{0\}) \cup (\{0\} \times {\bf C}),
\end{equation}
then $E$ is closed and completely circular, and $E^* = \emptyset$.
Equivalently,
\begin{equation}
        E = \{z = (z_1, z_2) \in {\bf C}^2 : z_1 \, z_2 = 0\},
\end{equation}
and it is easy to see that $\pol(E) = E$ in this case.

        Put
\begin{equation}
        D(r) = \{\zeta \in {\bf C} : |\zeta| \le r\}
\end{equation}
for each $r > 0$, and consider
\begin{equation}
        E = (D(r_1) \times \{0\}) \cup (\{0\} \times D(r_2))
\end{equation}
for some $r_1, r_2 > 0$.  Thus $E$ is closed and completely circular
again, $E^* = \emptyset$, and $E$ is also bounded in this case.  One
can check that $\pol(E) = E$ as well, using the polynomials $p_1(w) =
w_1$, $p_2(w) = w_2$, and $p(w) = w_1 \, w_2$.

        If $0 < r < R$ and
\begin{equation}
        E(r, R) = (D(r) \times {\bf C}) \cup (D(R) \times \{0\}),
\end{equation}
then $E(r, R)$ is closed and completely circular, and
\begin{equation}
        \overline{E(r, R)^*} = D(r) \times {\bf C}.
\end{equation}
If $p$ is a polynomial on ${\bf C}$ that is bounded on $E(r, R)$, then
$p(w_1, w_2)$ is bounded as a polynomial in $w_2$ for each $w_1 \in
D(r)$.  This implies that $p(w_1, w_2)$ is constant in $w_2$ for each
$w_1 \in D(r)$, and hence for every $w_1 \in {\bf C}$.  Thus $p(w_1,
w_2)$ reduces to a polynomial in $w_1$, and one can use this to show
that the polynomial hull of $E(r, R)$ is equal to $D(R) \times {\bf
C}$.

        Put
\begin{equation}
        E(r) = (D(r) \times {\bf C}) \cup ({\bf C} \times \{0\})
\end{equation}
for $r > 0$, which is the analogue of $E(r, R)$ with $R = +\infty$.
As before, $E(r)$ is closed and completely circular, and
\begin{equation}
        \overline{E(r)^*} = D(r) \times \{0\}.
\end{equation}
If $p$ is a polynomial on ${\bf C}^2$ that is bounded on $E(r)$, then $p$ is
constant, as in the previous paragraph, so that $\pol(E(r)) = {\bf C}^2$.

        Of course,
\begin{equation}
        E = D(r) \times {\bf C}
\end{equation}
is closed and completely circular for each $r > 0$, and satisfies
$\overline{E^*} = E$.  It is also easy to see that $\pol(E) = E$ in
this case, using the polynomial $p(w) = w_1$.

        If $r_1, r_2, R > 0$ and $r_1 < R$, then put
\begin{equation}
        E(r_1, r_2, R) = (D(r_1) \times D(r_2)) \cup (D(R) \times \{0\}).
\end{equation}
Thus $E(r_1, r_2, R)$ is closed, bounded, and completely circular, and
\begin{equation}
        \overline{E(r_1, r_2, R)^*} = D(r_1) \times D(r_2).
\end{equation}
If $z = (z_1, z_2) \in \pol(E(r_1, r_2, R))$, then it is easy to see
that $|z_1| \le R$ and $|z_2| \le r_2$, using the polynomials $p_1(w)
= w_1$ and $p_2(w) = w_2$.  If $z_2 \ne 0$, then one can show that
$|z_1| \le r_1$, using the polynomials $q_n(w) = w_1^n \, w_2$ for
each positive integer $n$.  More precisely,
\begin{equation}
        |q_n(z)| \le \sup \{|q_n(w)| : w \in E(r_1, r_2, R)\}
\end{equation}
implies that $|z_1|^n \, |z_2| \le r_1^n \, r_2$ for each $n$, and
hence that
\begin{equation}
        |z_1| \, |z_2|^{1/n} \le r_1 \, r_2^{1/n}.
\end{equation}
If $z_2 \ne 0$, then we can take the limit as $n \to \infty$ to get
that $|z_1| \le r_1$, as desired.  Thus $z \in E(r_1, r_2, R)$, which
implies that the polynomial hull of $E(r_1, r_2, R)$ is itself.

        If $E$ is a closed bi-disk
\begin{equation}
        D(r_1) \times D(r_2)
\end{equation}
for some $r_1, r_2 > 0$, then $E$ is completely circular,
$\overline{E^*} = E$, and $\pol(E) = E$.  If $E$ is the union of two
closed bi-disks
\begin{equation}
        (D(r_1) \times D(r_2)) \cup (D(t_1) \times D(t_2))
\end{equation}
for some $r_1, r_2, t_1, t_2 > 0$, then $E$ is completely circular and
$\overline{E^*} = E$ again.  Of course, this reduces to the single
bi-disk $D(t_1) \times D(t_2)$ when $r_1 \le t_1$ and $r_2 \le t_2$,
and to $D(r_1) \times D(r_2)$ when $t_1 \le r_1$ and $t_2 \le r_2$.
Otherwise, $E \ne \pol(E)$, because $E$ is not multiplicatively
convex.  Note that all of the other examples mentioned in this section
are multiplicatively convex.

\section{Multiplicative convexity}
\label{multiplicative convexity}
\setcounter{equation}{0}

        Suppose that $E \subseteq {\bf C}^n$ has the property that
\begin{equation}
\label{(t_1 z_1, ldots, t_n z_n) in E}
        (t_1 \, z_1, \ldots, t_n \, z_n) \in E
\end{equation}
when $z = (z_1, \ldots, z_n) \in E$, $t = (t_1, \ldots, t_n) \in {\bf
C}^n$, and $|t_j| = 1$ for each $j$, so that $t \in {\bf T}^n$.  Let
us say that $E$ is \emph{multiplicatively convex} if for each $v, w \in E$
and $a \in (0, 1)$, we have that $u \in E$ whenever $u \in {\bf C}^n$ and
\begin{equation}
        |u_j| = |v_j|^a \, |w_j|^{1 - a}
\end{equation}
for each $j$.  Similarly, if $E$ is completely circular and
multiplicatively convex, then $u \in E$ whenever
\begin{equation}
        |u_j| \le |v_j|^a \, |w_j|^{1 - a}
\end{equation}
for some $v, w \in E$, $0 < a < 1$, and each $j$.  If $E$ is completely
circular and convex, then $E$ is multiplicatively convex, because
\begin{equation}
\label{|v_j|^a |w_j|^{1 - a} le a |v_j| + (1 - a) |w_j|}
        |v_j|^a \, |w_j|^{1 - a} \le a \, |v_j| + (1 - a) \, |w_j|
\end{equation}
when $0 < a < 1$, by the convexity of the exponential function.  More
precisely, if $E$ is invariant under the usual action of ${\bf T}^n$,
as in (\ref{(t_1 z_1, ldots, t_n z_n) in E}), and $E$ is also nonempty
and convex, then it is easy to see that $0 \in E$.  This implies that
$r \, z \in E$ when $z \in E$ and $0 < r < 1$, and hence that $E$ is
completely circular.  We have also seen examples of sets that are
completely circular and multiplicatively convex, but not convex.

        If $E \subseteq {\bf C}^n$ is completely circular and $\pol(E)
= E$, then $E$ is multiplicatively convex, as in Section
\ref{completely circular sets}.  Conversely, if $E$ is closed,
bounded, completely circular, and multiplicatively convex, then
$\pol(E) = E$.  This basically follows from the discussion in Section
\ref{completely circular sets, continued}, with a few extra details.
The main point is that the sets $A_I$ considered there are closed and
convex in this case.  The convexity of the $A_I$'s corresponds exactly
to the multiplicative convexity of $E$.  To see that each $A_I$ is
closed, one can use the fact that $E$ is closed, and that for each $y
\in A_I$ there is a $w = w(y) \in E$ such that $w_j > 0$ and $\log
|w_j|$ when $j \in I$, and $w_j = 0$ when $j \not\in I$.  This also
uses the complete circularity of $E$, and otherwise there is a
standard argument based on the compactness of $E$.  Although the
boundedness of $E$ is not necessary for this step, it is important for
the approximation argument in Section \ref{completely circular sets,
continued}.

        Note that the polynomial hull of a bounded set $E \subseteq
{\bf C}^n$ is also bounded.  More precisely, if $|w_j| \le r_j$ for
some $r_j \ge 0$ and each $w \in E$, then $|z_j| \le r_j$ for each $z
\in \pol(E)$, as one can see by considering the polynomial $p_j(w) =
w_j$.

\section{Coefficients}
\label{coefficients}
\setcounter{equation}{0}

        Let
\begin{equation}
        p(z) = \sum_{|\alpha| \le N} a_\alpha \, z^\alpha
\end{equation}
be a polynomial with complex coefficients on ${\bf C}^n$, where the
sum is taken over all multi-indices $\alpha$ with $|\alpha| \le N$ for
some $N$.  Thus
\begin{equation}
        p(t_1 \, z_1, \ldots, t_n \, z_n)
         = \sum_{|\alpha| \le N} a_\alpha \, t^\alpha \, z^\alpha
\end{equation}
for each $t \in {\bf T}^n$, and so
\begin{equation}
        a_\beta \, z^\beta = \frac{1}{(2 \pi)^n} \int_{{\bf T}^n}
                  p(t_1 \, z_1, \ldots, t_n \, z_n) \, t^{-\beta} \, |dt|
\end{equation}
for every multi-index $\beta$, as in Section \ref{multiple fourier series}.
In particular,
\begin{eqnarray}
\label{|a_beta| |z^beta| le ... le sup_{t in T^n} |p(t_1 z_1, ldots, t_n z_n)|}
 |a_\beta| \, |z^\beta| & \le & \frac{1}{(2 \pi)^n} \int_{{\bf T}^n}
                                |p(t_1 \, z_1, \ldots, t_n \, z_n)| \, |dt| \\
 & \le & \sup_{t \in {\bf T}^n} |p(t_1 \, z_1, \ldots, t_n \, z_n)|. \nonumber
\end{eqnarray}

        Let $E$ be a nonempty subset of ${\bf C}^n$ which is
completely circular, or at least invariant under the usual action of
${\bf T}^n$.  If $p(z)$ is bounded on $E$, then it follows from the
discussion in the previous pargraph that each term $a_\beta \,
z^\beta$ in $p(z)$ is bounded on $E$.  Equivalently, the monomial
$z^\beta$ is bounded on $E$ whenever its coefficient $a_\beta$ in
$p(z)$ is not equal to $0$.  If $E$ is unbounded, then it may be that
$z^\beta$ is not bounded on $E$ for any nonzero multi-index $\beta$.
This implies that the only polynomials on ${\bf C}^n$ that are bounded
on $E$ are constant, and hence that $\pol(E) = {\bf C}^n$.

        As a nice family of examples in ${\bf C}^n$, consider
\begin{equation}
        E(b) = \{(z_1, z_2) \in {\bf C}^2 : |z_1|^b \, |z_2| \le 1\},
\end{equation}
where $b$ is a positive real number.  Thus $E(b)$ is closed,
completely circular, and multiplicatively convex for each $b > 0$, and
\begin{equation}
        E(b)^* = \{(z_1, z_2) \in {\bf C}^2 : 0 < |z_1|^b \, |z_2| \le 1\}
\end{equation}
is dense in $E(b)$ for each $b$ as well, as in Section \ref{another
condition}.  If $b$ is rational, so that $b = \beta_1 / \beta_2$ for
some positive integers $\beta_1$, $\beta_2$, then
\begin{equation}
 E(b) = \{(z_1, z_2) \in {\bf C}^2 : |z_1|^{\beta_1} \, |z_2|^{\beta_2} \le 1\}
           = \{z \in {\bf C}^2 : |z^\beta| \le 1\},
\end{equation}
where $\beta = (\beta_1, \beta_2)$.  In this case, it is easy to see
that $\pol(E(b)) = E(b)$, using the polynomial $p(z) = z^\beta$.
Otherwise, if $b$ is irrational, then one can check that $z^\beta$ is
unbounded on $E(b)$ for every nonzero multi-index $\beta$, which
implies that every nonconstant polynomial on ${\bf C}^n$ is unbounded
on $E(b)$, as before, and hence that $\pol(E(b)) = {\bf C}^2$.

\section{Polynomial convexity}
\label{polynomial convexity}
\setcounter{equation}{0}

        A set $E \subseteq {\bf C}^n$ is said to be \emph{polynomially
convex} if $\pol(E) = E$.  Thus $E$ has to be closed in this case,
since the polynomial hull of any set is closed.  Of course, $E
\subseteq \pol(E)$ automatically, and so $E$ is polynomially convex
when $\pol(E)$ is contained in $E$.  We have seen before that finite
subsets of ${\bf C}^n$ are polynomially convex, as are compact convex
sets.  A closed, bounded, and completely circular set is polynomially
convex if and only if it is multiplicatively convex, as in Section
\ref{multiplicative convexity}.  The polynomial hull of any set $E
\subseteq {\bf C}^n$ is polynomially convex, because $\pol(\pol(E)) =
\pol(E)$.  If $p$ is a polynomial on ${\bf C}^n$ and $k$ is a
nonnegative real number, then it is easy to see that
\begin{equation}
\label{E(p, k) = {z in {bf C}^n : |p(z)| le k}}
        E(p, k) = \{z \in {\bf C}^n : |p(z)| \le k\}
\end{equation}
is polynomially convex.  In particular, one can take $k = 0$, so that
the zero set of any polynomial is polynomially convex.

        If $E_\alpha$, $\alpha \in A$, is any collection of subsets of
${\bf C}^n$, then
\begin{equation}
        \pol\Big(\bigcap_{\alpha \in A} E_\alpha\Big)
             \subseteq \bigcap_{\alpha \in A} \pol(E_\alpha),
\end{equation}
because $\bigcap_{\alpha \in A} E_\alpha \subseteq E_\beta$ for each
$\beta \in A$, so that $\pol\Big(\bigcap_{\alpha \in A} E_\alpha\Big)
\subseteq \pol(E_\beta)$ for each $\beta \in A$.  If $E_\alpha$ is
polynomially convex for each $\alpha \in A$, then we get that
\begin{equation}
        \pol\Big(\bigcap_{\alpha \in A} E_\alpha\Big)
            \subseteq \bigcap_{\alpha \in A} \pol(E_\alpha)
              = \bigcap_{\alpha \in A} E_\alpha.
\end{equation}
This implies that $\bigcap_{\alpha \in A} E_\alpha$ is also
polynomially convex, since it is automatically contained in its
polynomial hull, as in the previous paragraph.  The polynomial hull of
any set $E \subseteq {\bf C}^n$ may be described as the intersection
of all sets $E(p, k)$ such that
\begin{equation}
        E \subseteq E(p, k),
\end{equation}
where $p$ is a polynomial on ${\bf C}^n$ and $k$ is a nonnegative real
number, as before.  It follows that $E$ is polynomially convex if and
only if it can be expressed as the intersection of some collection of
sets of the form $E(p, k)$, since these sets are all polynomially
convex, and the intersection of any collection of polynomially convex
sets is also polynomially convex.

        Alternatively, to avoid technical problems with unbounded
sets, one can expand the definition to say that a closed set $E
\subseteq {\bf C}^n$ is polynomially convex if for every compact set
$K \subseteq E$ we have that $\pol(K) \subseteq E$.  Of course, this
still implies that $\pol(E) = E$ when $E$ is compact.  With this
expanded definition, it is easy to see that every closed convex set in
${\bf C}^n$ is polynomially convex, for essentially the same reasons
as before.  Similarly, a closed completely circular set $E \subseteq
{\bf C}^n$ is polynomially convex in this expanded sense if and only
if it is multiplicatively convex.

\section{Entire functions, revisited}
\label{entire functions, revisited}
\setcounter{equation}{0}

        Let $E$ be a nonempty subset of ${\bf C}^n$, and let $\hol(E)$
be the set of $z \in {\bf C}^n$ such that
\begin{equation}
\label{|f(z)| le sup_{w in E} |f(w)|}
        |f(z)| \le \sup_{w \in E} |f(w)|
\end{equation}
for every complex-valued function $f$ on ${\bf C}^n$ that can be expressed as
\begin{equation}
\label{f(w) = sum_alpha a_alpha w^alpha}
        f(w) = \sum_\alpha a_\alpha \, w^\alpha.
\end{equation}
More precisely, the $a_\alpha$'s are supposed to be complex numbers,
and the sum is taken over all multi-indices $\alpha$ and is supposed
to be absolutely convergent for every $w \in {\bf C}^n$.  This
includes the case of polynomials, for which $a_\alpha = 0$ for all but
finitely many $\alpha$, and so
\begin{equation}
        \hol(E) \subseteq \pol(E).
\end{equation}
If $E$ is bounded, then $f$ can be approximated uniformly on $E$ by
finite subsums of (\ref{f(w) = sum_alpha a_alpha w^alpha}), which are
polynomials, and hence
\begin{equation}
        \hol(E) = \pol(E).
\end{equation}
If $E$ is not bounded, then $f$ may be unbounded on $E$, so that the
supremum in (\ref{|f(z)| le sup_{w in E} |f(w)|}) is $+\infty$, and
(\ref{|f(z)| le sup_{w in E} |f(w)|}) holds vacuously.

        Of course,
\begin{equation}
        E \subseteq \hol(E)
\end{equation}
automatically.  If $E_1 \subseteq E_2$, then
\begin{equation}
        \hol(E_1) \subseteq \hol(E_2).
\end{equation}
Note that functions on ${\bf C}^n$ expressed by absolutely summable
power series are continuous, because of uniform convergence on compact
sets, and continuity of polynomials.  This implies that $\hol(E)$ is
always a closed set in ${\bf C}^n$, and that
\begin{equation}
        \hol(\overline{E}) = \hol(E).
\end{equation}
As in the case of polynomial hulls, one can check that
\begin{equation}
        \hol(\hol(E)) = \hol(E).
\end{equation}
Using exponential functions as in Section \ref{entire functions}, one
also gets that
\begin{equation}
        \hol(E) \subseteq \overline{\con(E)}.
\end{equation}
More precisely, this works for both bounded and unbounded sets $E$.

        If $E$ is invariant under the torus action, as in Section
\ref{torus action}, then it is easy to see that $\hol(E)$ is too, as
before.  One can also use the maximum principle to show that $\hol(E)$
is completely circular in this case, as before.  One can use the three
lines theorem to show that $\hol(E)$ is multiplicatively convex in
this situation as well.  However, there are many examples where $E$ is
closed, completely circular, and multiplicatively convex, but $\hol(E)
\ne E$.  This uses the same type of arguments as in Sections
\ref{another condition} and \ref{coefficients}, and of course it is
important that $E$ be unbounded in these examples.

        If $E_\alpha$, $\alpha \in A$, is any collection of subsets of
${\bf C}^n$, then
\begin{equation}
        \hol\Big(\bigcap_{\alpha \in A} E_\alpha\Big)
          \subseteq \bigcap_{\alpha \in A} \hol(E_\alpha),
\end{equation}
as in the previous section.  If $\hol(E_\alpha) = E_\alpha$ for each
$\alpha \in A$, then it follows that
\begin{equation}
        \hol\Big(\bigcap_{\alpha \in A} E_\alpha\Big)
          = \bigcap_{\alpha \in A} E_\alpha,
\end{equation}
as before.  Let $g(w)$ be a complex-valued function on ${\bf C}^n$
that can be expressed by a power series that is absolutely summable
for each $w \in {\bf C}^n$, and put
\begin{equation}
        E(g, k) = \{w \in {\bf C}^n : |g(w)| \le k\}
\end{equation}
for each nonnegative real number $k$.  As in the previous section, it
is easy to see that
\begin{equation}
        \hol(E(g, k)) = E(g, k).
\end{equation}
One can also check that $\hol(E)$ is the same as the intersection of
all sets of the form $E(g, k)$ such that $E \subseteq E(g, k)$ for any
$E \subseteq {\bf C}^n$, as before.

\section{Power series expansions}
\label{power series expansions}
\setcounter{equation}{0}

        Let $R$ be a positive real number, and put
\begin{equation}
        D(R) = \{w \in {\bf C} : |w| < R\},
\end{equation}
as before.  Suppose that $f(w)$ is a holomorphic function on $D(R)$,
which one can take to mean that $f(w)$ is continuously-differentiable
and satisfies the Cauchy--Riemann equations.  Of course, it is well known
that one can also start with significantly weaker regularity conditions
on $f$.  If $|z| < r < R$, then Cauchy's integral formula implies that
\begin{equation}
\label{f(z) = frac{1}{2 pi i} oint_{partial D(r)} frac{f(w)}{w - z} dw}
 f(z) = \frac{1}{2 \pi i} \oint_{\partial D(r)} \frac{f(w)}{w - z} \, dw.
\end{equation}
More precisely, this uses an oriented contour integral over the circle
centered at $0$ with radius $r$, which is the boundary of the
corresponding disk $D(r)$.

        Let us briefly review the standard argument for obtaining a
power series expansion for $f(z)$ from (\ref{f(z) = frac{1}{2 pi i}
oint_{partial D(r)} frac{f(w)}{w - z} dw}).  If $|z| < r = |w|$, then
\begin{equation}
\label{frac{1}{w - z} = ... = w^{-1} sum_{j = 0}^infty w^{-j} z^j}
        \frac{1}{w - z} = \frac{1}{w \, (1 - w^{-1} \, z)}
                        = w^{-1} \sum_{j = 0}^\infty w^{-j} \, z^j,
\end{equation}
where the series on the right is an absolutely convergent geometric
series under these conditions.  The partial sums of this series also
converge uniformly as a function of $w$ on $\partial D(r)$ for each $z
\in D(r)$, by Weierstrass' M-test.  This permits us to interchange the
order of summation and integration in (\ref{f(z) = frac{1}{2 pi i}
oint_{partial D(r)} frac{f(w)}{w - z} dw}), to get that
\begin{equation}
\label{f(z) = sum_{j = 0}^infty a_j z^j}
        f(z) = \sum_{j = 0}^\infty a_j \, z^j,
\end{equation}
where $|z| < r < R$ and
\begin{equation}
\label{a_j = frac{1}{2 pi i} oint_{partial D(r)} f(w) w^{- j - 1} dw}
 a_j = \frac{1}{2 \pi i} \oint_{\partial D(r)} f(w) \, w^{- j - 1} \, dw
\end{equation}
for each $j \ge 0$.

        Although this expression for $a_j$ implicitly depends on $r$,
different choices of $r < R$ lead to the same value of $a_j$.  This is
an immediate consequence of Cauchy's theorem, and one can also observe
that $a_j$ is equal to $1/j!$ times the $j$th derivative of $f$ at
$0$, which obviously does not depend on $r$.  Alternatively, once one
has this power series expansion for $f$ on $D(r)$, one can use it to
evaluate integrals of $f$ over circles of radius less than $r$.  In
particular, the coefficients of the power series are given by the
corresponding integrals over circles of radius less than $r$, because
of the usual orthogonality properties of the $w^j$'s with respect to
integration over the unit circle.  This also uses the fact that the
partial sums of the power series converge uniformly on compact subsets
of $D(r)$, to interchange the order of integration and summation.

        Note that
\begin{equation}
\label{|a_j| le frac{1}{2 pi r^{j + 1}} int_{partial D(r)} |f(w)| |dw|}
 |a_j| \le \frac{1}{2 \pi r^{j + 1}} \int_{\partial D(r)} |f(w)| \, |dw|
\end{equation}
for each $j$, where the integral is now taken with respect to the
element of arc length $|dw|$.  In particular,
\begin{equation}
\label{|a_j| le r^{-j} (sup_{|w| = r} |f(w)|)}
        |a_j| \le r^{-j} \, \Big(\sup_{|w| = r} |f(w)|\Big).
\end{equation}
This works for each $r < R$, since $a_j$ does not depend on $r$, as in
the previous paragraph.

\section{Power series expansions, continued}
\label{power series expansions, continued}
\setcounter{equation}{0}

        Let $n$ be a positive integer, and let $R = (R_1, \ldots,
R_n)$ be an $n$-tuple of positive real numbers.  Also let
\begin{equation}
        D_n(R) = D(R_1) \times \cdots \times D(R_n)
\end{equation}
be the corresponding polydisk in ${\bf C}^n$.  To say that a
complex-valued function $f(w)$ on $D(R)$ is holomorphic, we mean that
$f(w)$ is continuously-differentiable on $D_n(R)$ and holomorphic as a
function of $w_j$ for each $j$, which is to say that $f(w)$ satisfies
the Cauchy--Riemann equations as a function of $w_j$ for each $j$.  As
in the one-variable case, one can start with weaker regularity
conditions on $f$, but we shall not pursue this here.  One might at
least note that it would be sufficient in this section to ask that $f$
be continuous on $D_n(R)$ and holomorphic in each variable separately.

        If $z \in D(R)$ and $|z_1| < r_1 < R_1$, then we can apply Cauchy's
integral formula to $f(w)$ as a holomorphic function of $w_1$ to get that
\begin{equation}
        f(z) = \frac{1}{2 \pi i} \oint_{\partial D(r_1)}
                       \frac{f(w_1, z_2, \ldots, z_n)}{w_1 - z_1} \, dw_1,
\end{equation}
as in the previous section.  Repeating the process, if $|z_j| < r_j < R_j$
for each $j$, then we get that
\begin{equation}
\label{f(z) = frac{1}{(2 pi i)^n} oint ... dw_1 cdots dw_n}
 \quad  f(z) = \frac{1}{(2 \pi i)^n} \oint_{\partial D(r_1)} \cdots
                                                        \oint_{\partial D(r_n)}
       f(w) \, \Big(\prod_{j = 1}^n (w_j - z_j)^{-1}\Big) \, dw_1 \cdots dw_n,
\end{equation}
which is an $n$-dimensional version of Cauchy's integral formula.

        Let us pause for a moment to consider ``multiple geometric
series''.  If $\zeta \in {\bf C}^n$ and $|\zeta_j| < 1$ for each $j$, then
\begin{equation}
\label{prod_{j = 1}^n (1 - zeta_j)^{-1} = ...  = sum_alpha zeta^alpha}
        \prod_{j = 1}^n (1 - \zeta_j)^{-1}
 = \prod_{j = 1}^n \Big(\sum_{\ell_j = 0}^\infty \zeta_j^{\ell_j}\Big)
   = \sum_\alpha \zeta^\alpha,
\end{equation}
where the last sum is taken over all multi-indices $\alpha$, and
$\zeta^\alpha = \zeta_1^{\alpha_1} \cdots \zeta_n^{\alpha_n}$ is the
usual monomial.  All of these sums converge absolutely under these conditions.

        If $|z_j| < r_j = |w_j|$ for each $j$, then
\begin{equation}
\label{prod_{j = 1}^n (w_j - z_j)^{-1} = sum_alpha w^{-alpha - 1} z^alpha}
        \prod_{j = 1}^n (w_j - z_j)^{-1}
         = \prod_{j = 1}^n w_j^{-1} \, (1 - w_j^{-1} \, z_j)^{-1}
         = \sum_\alpha w^{-\alpha - 1} \, z^\alpha,
\end{equation}
where $w^{-\alpha - 1} = w_1^{- \alpha_1 - 1} \cdots w_n^{- \alpha_n -
1}$.  As usual, this sum is absolutely convergent under these
conditions, and is uniformly approximated by finite subsums as a
function of $w$ on $\partial D(r_1) \times \cdots \times \partial
D(r_n)$ for each $z \in D_n(r)$, $r = (r_1, \ldots, r_n)$.

        If $r_j < R_j$ for each $j$, then put
\begin{equation}
\label{a_alpha = frac{1}{(2 pi i)^n} oint cdots dw_1 cdots dw_n}
        a_\alpha = \frac{1}{(2 \pi i)^n} \oint_{\partial D(r_1)} \cdots
          \oint_{\partial D(r_n)} f(w) \, w^{- \alpha - 1} \, dw_1 \cdots dw_n
\end{equation}
for each multi-index $\alpha$.  Thus
\begin{equation}
\label{|a_alpha| le frac{r^{-alpha - 1}}{(2 pi)^n} int ... |dw_1| cdots |dw_n|}
 |a_\alpha| \le \frac{r^{-\alpha - 1}}{(2 \pi)^n} \int_{\partial D(r_1)}
                 \cdots \int_{\partial D(r_n)} |f(w)| \, |dw_1| \cdots |dw_n|,
\end{equation}
where $r^{-\alpha - 1}$ is as in the previous paragraph, and hence
\begin{equation}
\label{|a_alpha| le r^{-alpha} sup {|f(w)| : |w_j| = r_j for j = 1, ldots, n}}
 |a_\alpha| \le r^{-\alpha} \sup \{|f(w)| : |w_j| = r_j
                                                \hbox{ for } j = 1, \ldots, n\}
\end{equation}
for each $\alpha$.

        If $|z_j| < r_j < R_j$ for each $j$, then we get that
\begin{equation}
\label{f(z) = sum_alpha a_alpha z^alpha}
        f(z) = \sum_\alpha a_\alpha \, z^\alpha.
\end{equation}
More precisely, it is easy to see that the sum on the right converges
absolutely under these conditions, by comparison with a convergent
multiple geometric series.  To get (\ref{f(z) = sum_alpha a_alpha
z^alpha}), one can plug (\ref{prod_{j = 1}^n (w_j - z_j)^{-1} =
sum_alpha w^{-alpha - 1} z^alpha}) into (\ref{f(z) = frac{1}{(2 pi
i)^n} oint ... dw_1 cdots dw_n}), and interchange the order of
summation and integration.  This uses the fact that the sum in
(\ref{prod_{j = 1}^n (w_j - z_j)^{-1} = sum_alpha w^{-alpha - 1}
z^alpha}) can be approximated uniformly by finite subsums for $w \in
\partial D(r_1) \times \cdots \times \partial D(r_n)$.

        As in the previous section, the coefficients $a_\alpha$ do not
depend on the choice of $r = (r_1, \ldots, r_n)$, as long as $0 < r_j
< R_j$ for each $j$.  Thus (\ref{f(z) = sum_alpha a_alpha z^alpha})
holds on all of $D_n(R)$, with absolute convergence of the sum for
every $z \in D_n(R)$.

\section{Holomorphic functions, revisited}
\label{holomorphic functions, revisited}
\setcounter{equation}{0}

        Let us say that a complex-valued function $f(z)$ on a nonempty
open set $U$ in ${\bf C}^n$ is holomorphic if it is
continuously-differentiable in the real-variable sense and holomorphic
in each variable separately.  As in the previous section, this implies
that $f$ can be represented by an absolutely convergent power series
on a neighborhood of any point in $U$.  In particular, $f$ is
automatically continuously-differentiable of all orders on $U$.  This
would also work if we only asked that $f$ be continuous on $U$ and
holomorphic in each variable separately, but we shall not try to deal
with weaker regularity conditions here.

        Let $C(U)$ be the algebra of continuous complex-valued
functions on $U$, and let $\mathcal{H}(U)$ be the subspace of $C(U)$
consisting of holomorphic functions.  More precisely, $\mathcal{H}(U)$
is a subalgebra of $C(U)$, because the sum and product of two
holomorphic functions on $U$ are also holomorphic.  Remember that
there is also a natural topology on $C(U)$, defined by the supremum
seminorms associated to nonempty compact subsets of $U$.  As in the
one-variable case, one can check that $\mathcal{H}(U)$ is closed in
$C(U)$ with respect to this topology, using the $n$-dimensional
version of the Cauchy integral formula.

        Let $f$ be a holomorhic function on $U$, and let $U_0$ be the
set of $p \in U$ such that $f = 0$ at every point in a neighborhood of
$p$, so that $U_0$ is an open set in $U$, by construction.  If $Z$ is
the set of $p \in U$ such that $f$ and all of its derivatives are
equal to $0$ at $p$, then $Z$ is relatively closed in $U$, because $f$
and its derivatives are continuous on $U$.  Clearly $U_0 \subseteq Z$,
and $Z \subseteq U_0$ because of the local power series representation
of $f$ at each point in $U$.  Thus $U_0 = Z$ is both open and
relatively closed in $U$.  It follows that $U_0 = U$ when $U_0 \ne
\emptyset$ and $U$ is connected.

        Suppose that $h$ is a continuous complex-valued function on a
closed disk in the complex plane which is holomorphic in the interior
and not equal to $0$ at any point on the boundary.  Let $a$ be the
number of points in the interior at which $h$ is equal to $0$, counted
with their appropriate multiplicity.  The argument principle implies
that $a$ is the same as the winding number of the boundary values of
$h$ around $0$ in the range.  This winding number is not changed by
small perturbations of $h$ on the boundary with respect to the
supremum norm, and hence $a$ is not changed by small perturbations of
$h$ as a continuous function on the closed disk which is holomorphic
in the interior with respect to the supremum norm.  This implies that
a holomorphic function $f$ in $n \ge 2$ complex variables cannot have
isolated zeros, by considering $f$ as a continuous family of
holomorphic functions in one variable parameterized by the other $n -
1$ variables.

\section{Laurent expansions}
\label{laurent expansions}
\setcounter{equation}{0}

        Let $R$, $T$ be nonnegative real numbers with $R < T$, and let
\begin{equation}
\label{A(R, T) = {z in {bf C} : R < |w| < T}}
        A(R, T) = \{z \in {\bf C} : R < |w| < T\}
\end{equation}
be the open annulus in the complex plane with inner radius $R$ and
outer radius $T$.  If $f(w)$ is a holomorphic function on $A(R, T)$
and $R < r < |z| < t < T$, then Cauchy's integral formula implies that
\begin{equation}
 f(z) = \frac{1}{2 \pi i} \oint_{\partial A(r, t)} \frac{f(w)}{w - z} \, dw.
\end{equation}
The boundary of $A(r, t)$ consists of the circles centered at $0$ with
radii $r$, $t$ and opposite orientations, and the integral over
$\partial A(r, t)$ may be re-expressed as
\begin{equation}
        \oint_{|w| = t} \frac{f(w)}{w - z} \, dw
         - \oint_{|w| = r} \frac{f(w)}{w - z} \, dw,
\end{equation}
where these circles have their usual positive orientations in both
integrals.

        As in Section \ref{power series expansions},
\begin{equation}
        \frac{1}{2 \pi i} \oint_{|w| = t} \frac{f(w)}{w - z} \, dw
         = \sum_{j = 0}^\infty a_j \, z^j,
\end{equation}
where
\begin{equation}
        a_j = \frac{1}{2 \pi i} \oint_{|w| = t} f(w) \, w^{- j - 1} \, dw.
\end{equation}
Note that
\begin{equation}
        |a_j| \le \frac{1}{2 \pi t^{j + 1}} \oint_{|w| = t} |f(w)| \, |dw|
               \le t^{-j} \Big(\sup_{|w| = t} |f(w)|\Big)
\end{equation}
for each $j \ge 0$, so that $\sum_{j = 0}^\infty a_j \, z^j$ converges
absolutely when $|z| < t$.

        The other term is a bit different, because $|z| > |w| = r$.
This time we use
\begin{equation}
        \frac{-1}{w - z} = \frac{1}{z \, (1 - z^{-1} \, w)}
                         = z^{-1} \sum_{j = 0}^\infty z^{-j} \, w^j
\end{equation}
to get that
\begin{equation}
        - \frac{1}{2 \pi i} \oint_{|w| = r} \frac{f(w)}{w - z} \, dw
           = \sum_{j = -1}^{-\infty} a_j \, z^j,
\end{equation}
where
\begin{equation}
        a_j = \frac{1}{2 \pi i} \oint_{|w| = r} f(w) \, w^{- j - 1} \, dw
\end{equation}
for $j \le -1$.  Thus
\begin{equation}
        |a_j| \le \frac{1}{2 \pi r^{j + 1}} \int_{|w| = r} |f(w)| \, |dw|
               \le r^{-j} \Big(\sup_{|w| = r} |f(w)|\Big)
\end{equation}
for each $j \le -1$, so that $\sum_{j = -1}^{-\infty} a_j \, z^j$
converges absolutely when $|z| > r$.

        Combining these two series, we get that
\begin{equation}
\label{f(z) = sum_{j = -infty}^infty a_j z^j}
        f(z) = \sum_{j = -\infty}^\infty a_j \, z^j
\end{equation}
when $r < |z| < t$, where the coefficients $a_j$ are given as above
for $j \ge 0$ and $j \le -1$, respectively.  As in Section \ref{power
series expansions}, these coefficients do not actually depend on the
choices of radii $r, t \in (R, T)$.

\section{Laurent expansions, continued}
\label{laurent expansions, continued}
\setcounter{equation}{0}

        Let $R$, $T$ be nonnegative real numbers with $R < T$, and let
$V$ be a nonempty open set in ${\bf C}^{n - 1}$ for some $n \ge 2$.
If $z = (z_1, z_2, \ldots, z_n) \in {\bf C}^n$, then we put $z' =
(z_2, \ldots, z_n) \in {\bf C}^{n - 1}$, and identify $z$ with $(z_1,
z') \in {\bf C} \times {\bf C}^{n - 1}$, so that
\begin{equation}
        U = A(R, T) \times V
\end{equation}
is identified with an open set in ${\bf C}^n$.

        Let $f$ be a holomorphic function on $U$, and let $z$ be an
element of $U$, with $r < |z_1| < t$ for some $r, t \in (R, T)$.
Applying the discussion in the previous section to $f(z_1, z')$ as a
function of $z_1$, we get that
\begin{equation}
        f(z_1, z') = \sum_{j = -\infty}^\infty a_j(z') \, z_1^j,
\end{equation}
where
\begin{equation}
 a_j(z') = \frac{1}{2 \pi i} \oint_{|w| = t} f(w, z') \, w^{- j - 1} \, dw
\end{equation}
when $j \ge 0$, and
\begin{equation}
 a_j(z') = \frac{1}{2 \pi i} \oint_{|w| = r} f(w, z') \, w^{- j - 1} \, dw
\end{equation}
when $j \le -1$.  It follows from these expressions that $a_j(z')$ is
holomorphic as a function of $z'$ on $V$ for each $j$, because $f$ is
holomorphic.

        Suppose that $V_1$ is a nonempty open subset of $V$, and that
$f$ is actually a holomorphic function on the open set
\begin{equation}
\label{(A(R, T) times V) cup (D(T) times V_1)}
        (A(R, T) \times V) \cup (D(T) \times V_1)
\end{equation}
in ${\bf C}^n$.  Thus $f(w, z')$ is holomorphic as a function of $w$
on the open disk $D(T)$ for each $z' \in V_1$.  This implies that
\begin{equation}
        a_j(z') = 0
\end{equation}
when $z' \in V_1$ and $j \le -1$.  If $V$ is connected, then it
follows that the same conclusion holds for every $z' \in V$ and $j \le
-1$, because $a_j(z')$ is holomorphic as a function of $z'$ on $V$ for
each $j$.

        Under these conditions, we get that
\begin{equation}
        f(z_1, z') = \sum_{j = 0}^\infty a_j(z') \, z_1^j
\end{equation}
for every $z = (z_1, z')$ in (\ref{(A(R, T) times V) cup (D(T) times V_1)}).
This series actually converges absolutely when $|z_1| < T$ and $z' \in V$,
as one can see by choosing $t$ such that $|z_1| < t < T$, and applying
the estimate for $|a_j|$ in the previous section.  Similarly, the partial
sums of this series converge uniformly on compact subsets of $D(T) \times V$.
The partial sums are also holomorphic in $z_1$ and $z'$, and it follows
that the series defines a holomorphic function on $D(T) \times V$.  Thus
$f$ extends to a holomorphic function on $D(T) \times V$ in this case.

\section[\ Completely circular domains]{Completely circular domains}
\label{completely circular domains}
\setcounter{equation}{0}

        Let $U$ be a nonempty complete circular open subset of ${\bf C}^n$.
If $z \in U$, then there is an $n$-tuple $R = (R_1, \ldots, R_n)$ of
positive real numbers such that
\begin{equation}
        z \in D_n(R) \subseteq U,
\end{equation}
where $D_n(R) = D(R_1) \times \cdots \times D(R_n)$ is the polydisk
associated to $R$, as before.  Thus $U$ can be expressed as a union of
open polydisks.

        As in Section \ref{another condition}, let $U^*$ be the set of
$w \in U$ such that $w_j \ne 0$ for each $j$, and let $A$ be the set
of $y \in {\bf R}^n$ for which there is a $w \in U^*$ such that $y_j =
\log |w_j|$ for each $j$.  Note that $A$ is an open set in ${\bf
R}^n$, and that for each $z \in U$ there is a $w \in U^*$ such that
$|z_j| < |w_j|$, because $U$ is an open set in ${\bf C}^n$.  As
before, if $x \in {\bf R}^n$, $y \in A$, and $x_j \le y_j$ for each
$j$, then $x \in A$, because $U$ is completely circular.  Similarly,
if $\zeta \in {\bf C}^n$, $x \in A$, and $|z_j| \le \exp x_j$ for each
$j$, then $\zeta \in U$.  Conversely, for each $\zeta \in U$ there is
an $x \in A$ with this property, so that $U$ is completely determined
by $A$ under these conditions.

        Let $f$ be a holomorphic function on $U$.  If $R$ is an
$n$-tuple of positive real numbers such that $D_n(R) \subseteq U$,
then $f$ can be represented by a power series on $D_n(R)$, as in
Section \ref{power series expansions, continued}.  More precisely,
there are complex numbers $a_\alpha$ for each multi-index $\alpha$
such that
\begin{equation}
        f(z) = \sum_\alpha a_\alpha \, z^\alpha
\end{equation}
for each $z \in D_n(R)$, where the sum converges absolutely.  The
coefficients $a_\alpha$ can be given by the derivatives of $f$ at $0$
in the usual way, since
\begin{equation}
 \frac{\partial^{|\alpha|} f}{\partial z^\alpha}(0) = \alpha ! \cdot a_\alpha,
\end{equation}
where $\alpha ! = \alpha_1 ! \cdots \alpha_n !$.  In particular, the
coefficients $a_\alpha$ do not depend on $R$, and so this power series
representation for $f(z)$ holds for every $z \in U$.

        Remember that $\con(A)$ denotes the convex hull of $A$ in
${\bf R}^n$, which is an open set in ${\bf R}^n$ in this case, because
$A$ is open.  Similarly, if $x \in {\bf R}^n$, $y \in \con(A)$, and
$x_j \le y_j$ for each $j$, then $x \in \con(A)$, because of the
corresponding property of $A$.  Consider
\begin{eqnarray}
 V & = & \{\zeta \in {\bf C}^n : \hbox{ there is an } x \in \con(A) 
         \hbox{ such that } \\
 & & \qquad\qquad\quad |\zeta_j| \le \exp x_j \hbox{ for } j = 1, \ldots, n\}.
                                                                 \nonumber
\end{eqnarray}
It is easy to see that $V$ is open, completely circular, and
multiplicatively convex under these conditions.  We also have that $U
\subseteq V$, with $U = V$ exactly when $U$ is multiplicatively
convex.  As in Section \ref{power series}, the set of $z \in {\bf
C}^n$ for which $\sum_\alpha a_\alpha \, z^\alpha$ is absolutely
summable is completely circular and multiplicatively convex.  It is
not difficult to check that this happens for each $z \in V$, so that
$f$ extends to a holomorphic function on $V$.

\section[\ Convex sets]{Convex sets}
\label{convex sets}
\setcounter{equation}{0}

        Let $A$ be a nonempty convex set in ${\bf R}^n$.  As in
Section \ref{convex hulls}, if $x \in {\bf R}^n \backslash \overline{A}$,
then there is a linear function $\lambda$ on ${\bf R}^n$ such that
\begin{equation}
\label{sup_{y in A} lambda(y) < lambda(x)}
        \sup_{y \in A} \lambda(y) < \lambda(x).
\end{equation}
More precisely, we can express $\lambda$ as
\begin{equation}
        \lambda(y) = \sum_{j = 1}^n a_j \, y_j
\end{equation}
for some $a \in {\bf R}^n$.  Of course, $a \ne 0$, and we can
normalize $a$ so that
\begin{equation}
\label{max_{1 le j le n} |a_j| = 1}
        \max_{1 \le j \le n} |a_j| = 1,
\end{equation}
by multiplying $a$ by a positive real number.

        Suppose now that $x \in \partial A$, and let us show that
there is a nonzero linear functional $\lambda$ on ${\bf R}^n$ such that
\begin{equation}
\label{lambda(y) le lambda(x)}
        \lambda(y) \le \lambda(x)
\end{equation}
for every $y \in A$.  By hypothesis, there is a sequence $\{x(l)\}_{l
= 1}^\infty$ of elements of ${\bf R}^n \backslash \overline{A}$ that
converges to $x$.  As in the previous paragraph, for each $l$ there is
an $a(l) \in {\bf R}^n$ such that
\begin{equation}
        \max_{1 \le j \le n} |a_j(l)| = 1
\end{equation}
and $\lambda_l(y) = \sum_{j = 1}^n a_j(l) \, y_j$ satisfies
\begin{equation}
        \sup_{y \in A} \lambda_l(y) < \lambda_l(x(l)).
\end{equation}
Passing to a subsequence if necessary, we may suppose that
$\{a(l)\}_{l = 1}^\infty$ converges to some $a \in {\bf R}^n$, which
also satisfies (\ref{max_{1 le j le n} |a_j| = 1}).  If $\lambda$ is
the linear functional on ${\bf R}^n$ corresponding to $a$ as before,
then it is easy to see that $\lambda$ satisfies (\ref{lambda(y) le
lambda(x)}), as desired.

        If in addition $A$ is an open set in ${\bf R}^n$, then we get that
\begin{equation}
        \lambda(y) < \lambda(x)
\end{equation}
for every $y \in A$.  Otherwise, if $\lambda(y) = \lambda(x)$ for some
$y \in A$, then one can use the facts that $A$ is open and $\lambda
\ne 0$ to get that $\lambda(z) > \lambda(x)$ for some $z \in A$.

        As another special case, suppose that $A$ has the property
that for each $u \in {\bf R}^n$ and $y \in A$ with $u_j \le y_j$ for
each $j$ we have that $u \in A$ too.  If $\lambda(y) = \sum_{j = 1}^n
a_j \, y_j$ satisfies (\ref{lambda(y) le lambda(x)}), then $a_j \ge 0$
for each $j$.

\section[\ Completely circular domains, continued]{Completely circular domains, continued}
\label{completely circular domains, continued}
\setcounter{equation}{0}

        Let $U$ be a nonempty open subset of ${\bf C}^n$ that is also
completely circular and multiplicatively convex, and let $w$ be an
element of the boundary of $U$.  Note that $w \ne 0$, because $0 \in
U$.  Let $I$ be the set of $j = 1, \ldots, n$ such that $w_j \ne 0$,
and let $U_I$ be the set of $z \in U$ such that $z_j \ne 0$ when $j
\in I$.  Also let ${\bf R}^I$ be the set of real-valued functions on
$I$, and let $A_I$ be the set of elements of ${\bf R}^I$ of the form
$\log |z_j|$, $j \in I$, with $z \in U_I$.  If $u \in {\bf R}^I$, $v
\in A_I$, and $u_j \le v_j$ for each $j \in I$, then $u \in A_I$ too,
because $U$ is completely circular.  It is easy to see that $A_I$ is
open and convex in ${\bf R}^I$, because $U$ is open and
multiplicatively convex.  One can also check that $\log |w_j|$, $j \in
I$, corresponds to an element of the boundary of $A_I$ in ${\bf R}^I$
under these conditions.

        As in the previous section, there is an $a \in {\bf R}^I$ such
that $a_j \ge 0$ for each $j \in I$, $\max_{j \in I} a_j = 1$, and
\begin{equation}
\label{sum_{j in I} a_j v_j < sum_{j in I} a_j log |w_j|}
        \sum_{j \in I} a_j \, v_j < \sum_{j \in I} a_j \, \log |w_j|
\end{equation}
for each $v \in A_I$.  If $j \in I$ and $l$ is a positive integer, then
let $\alpha_j(l)$ be the smallest positive integer such that
\begin{equation}
        a_j \, l \le \alpha_j(l).
\end{equation}
Put $\alpha_j(l) = 0$ when $j \not\in I$, so that $\alpha(l) =
(\alpha_1(l), \ldots, \alpha_n(l))$ is a multi-index for each positive
integer.  By construction, $a_{j_0} = 1$ for some $j_0 \in I$, which
implies that $\alpha_{j_0}(l) = l$ for each $l$.  In particular, the
multi-indices $\alpha(l)$ are all distinct.

        Consider
\begin{equation}
\label{f_w(z) = sum_{l = 1}^infty w^{-alpha(l)} z^{alpha(l)}}
        f_w(z) = \sum_{l = 1}^\infty w^{-\alpha(l)} \, z^{\alpha(l)}.
\end{equation}
This is a power series in $z$, with coefficients $w^{-\alpha(l)} =
\prod_{j \in I} w_j^{-\alpha_j(l)}$, and we would like to show that it
converges absolutely when $z \in U$.  If $z \in U \backslash U_I$, so
that $z_j = 0$ for some $j \in I$, then $z^{\alpha(l)} = 0$ for each
$l$, because $\alpha_j(l) \ge 1$ for every $j \in I$ and $l \ge 1$ by
construction.  Thus we may as well suppose that $z \in U_I$, so that
$\log |z_j|$, $j \in I$, determines an element of $A_I$, and hence
\begin{equation}
        \sum_{j \in I} a_j \, \log |z_j| < \sum_{j \in I} a_j \, \log |w_j|.
\end{equation}
Equivalently,
\begin{equation}
\label{prod_{j in I} |z_j|^{a_j} < prod_{j in I} |w_j|^{a_j}}
        \prod_{j \in I} |z_j|^{a_j} < \prod_{j \in I} |w_j|^{a_j}.
\end{equation}

        Observe that
\begin{equation}
        0 \le \alpha_j(l) - a_j \, l \le 1
\end{equation}
for each $j \in I$ and $l \ge 1$.  Remember that $\alpha_j(l)$ is the
smallest positive integer greater than or equal to $a_j \, l$, so that
$\alpha_j(l) - a_j \, l \ge 0$ in particular.  If $a_j > 0$, then
$\alpha_j(l) - a_j < 1$ for each $l$.  Otherwise, if $a_j = 0$, then
$\alpha_j(l) = 1$ for each $l$.

        Of course,
\begin{equation}
        |w^{-\alpha(l)}| \, |z^{\alpha(l)}|
         = \prod_{j \in I} \Big(\frac{|z_j|}{|w_j|}\Big)^{\alpha_j(l)}.
\end{equation}
Using the observation in the previous paragraph, we get that
\begin{equation}
 \prod_{j \in I} \Big(\frac{|z_j|}{|w_j|}\Big)^{\alpha_j(l) - a_j \, l} \le C
\end{equation}
for some $C \ge 0$, where $C$ depends on $w$ and $z$ but not $l$.  Hence
\begin{equation}
\label{|w^{-alpha(l)}| |z^{alpha(l)}| le C prod_{j in I} ...}
        |w^{-\alpha(l)}| \, |z^{\alpha(l)}|
         \le C \,\prod_{j \in I} \Big(\frac{|z_j|}{|w_j|}\Big)^{a_j \, l}
\end{equation}
for each $l$.

        Equivalently,
\begin{equation}
\label{|w^{-alpha(l)}| |z^{alpha(l)}| le C (prod_{j in I} ....)^l}
        |w^{-\alpha(l)}| \, |z^{\alpha(l)}|
         \le C \, \Big(\prod_{j \in I} \frac{|z_j|^{a_j}}{|w_j|^{a_j}}\Big)^l
\end{equation}
for each $l$.  Note that the quantity in parentheses on the right side
is strictly less than $1$, by (\ref{prod_{j in I} |z_j|^{a_j} <
prod_{j in I} |w_j|^{a_j}}).  It follows that the series in
(\ref{f_w(z) = sum_{l = 1}^infty w^{-alpha(l)} z^{alpha(l)}})
converges absolutely when $z \in U_I$, by comparison with a convergent
geometric series, as desired.

        Thus $f_w(z)$ defines a holomorphic function of $z$ on $U$.
If $z = w$, then the series in (\ref{f_w(z) = sum_{l = 1}^infty
w^{-alpha(l)} z^{alpha(l)}}) does not converge, because every term in
the series is equal to $1$.  It is easy to see that $t \, w \in U$
when $t$ is a nonnegative real number strictly less than $1$, because
$w \in \partial U$ and $U$ is completely circular.  In this case,
\begin{equation}
        f_w(t \, w) = \sum_{l = 1}^\infty t^{|\alpha(l)|},
\end{equation}
which tends to $+\infty$ as $t \to 1$.  It follows that $f_w(z)$ does
not have a holomorphic extension to a neighborhood of $w$, since it is
not even bounded on $U$ near $w$.

\section[\ Convex domains]{Convex domains}
\label{convex domains}
\setcounter{equation}{0}

        Let $U$ be a nonempty convex open set in ${\bf C}^n$, and let
$w$ be an element of the boundary of $U$.  As in Section \ref{convex
sets}, there is a complex-linear function $\mu$ on ${\bf C}^n$ such
that
\begin{equation}
        \re \mu(z) < \re \mu(w)
\end{equation}
for every $z \in U$.  In particular,
\begin{equation}
        \mu(z) \ne \mu(w)
\end{equation}
for every $z \in U$.  It follows that
\begin{equation}
        g_w(z) = \frac{1}{\mu(z) - \mu(w)}
\end{equation}
is a holomorphic function on $U$ that is unbounded on the intersection
of $U$ with any neighborhood of $w$, and hence does not have a
holomorphic extension to the union of $U$ with any neighborhood of $w$.

\section[\ Planar domains]{Planar domains}
\label{Planar domains}
\setcounter{equation}{0}

        Let $U$ be a nonempty open set in the complex plane, and let
$w$ be an element of the boundary of $U$.  Observe that
\begin{equation}
        h_w(z) = \frac{1}{z - w}
\end{equation}
is a holomorphic function on $U$ that is unbounded on the intersection
of $U$ with any neighborhood of $w$, and hence cannot be extended to a
holomorphic function on the union of $U$ with any neighborhood of $w$.
In particular, holomorphic functions in one complex variable can have
isolated zeros, and thus isolated singularities.  We have seen before
that holomorphic functions in two or more complex variables cannot
have isolated zeros, and they also cannot have isolated singularities,
by the earlier discussion about Laurent expansions.

\part{Convolution}

\section[\ Convolution on ${\bf T}^n$]{Convolution on ${\bf T}^n$}
\label{convolution on T^n}
\setcounter{equation}{0}

        Let $f$, $g$ be continuous complex-valued functions on the
$n$-dimensional torus ${\bf T}^n$.  The \emph{convolution} $f * g$ is
the function defined on ${\bf T}^n$ by
\begin{equation}
\label{(f * g)(z) = frac{1}{(2 pi)^n} int_{T^n} f(z diamond w^{-1}) g(w) |dw|}
        (f * g)(z) = \frac{1}{(2 \pi)^n} \int_{{\bf T}^n} f(z \diamond w^{-1})
                                                        \, g(w) \, |dw|.
\end{equation}
As before, $|dw|$ is the $n$-dimensional element of integration on
${\bf T}^n$ corresponding to the element $|dw_j|$ of arc length in
each variable.  Alternatively, $|dw|$ represents the appropriate
version of Lebesgue measure on ${\bf T}^n$.  If $z = (z_1, \ldots,
z_n)$ and $w = (w_1, \ldots, w_n)$ are elements of ${\bf T}^n$, then
we put
\begin{equation}
        w^{-1} = (w_1^{-1}, \ldots, w_n^{-1})
\end{equation}
and
\begin{equation}
        z \diamond w = (z_1 \, w_1, \ldots, z_n \, w_n),
\end{equation}
so that $z \diamond w^{-1}$ is also defined.

        It is easy to see that $f * g$ is also a continuous function
on ${\bf T}^n$ when $f$, $g$ are continuous, using the fact that
continuous functions on ${\bf T}^n$ are uniformly continuous, since
${\bf T}^n$ is compact.  Observe that
\begin{equation}
\label{f * g = g * f}
        f * g = g * f,
\end{equation}
as one can see using the change of variables $w \mapsto w^{-1} \diamond z$ in
(\ref{(f * g)(z) = frac{1}{(2 pi)^n} int_{T^n} f(z diamond w^{-1}) g(w) |dw|}).
More precisely, this also uses the fact that the measure on ${\bf T}^n$
is invariant under the mappings $w \mapsto w^{-1}$ and $w \mapsto u \diamond w$
for each $u \in {\bf T}^n$.  Similarly, one can check that
\begin{equation}
        (f * g) * h = f * (g * h)
\end{equation}
for all continuous functions $f$, $g$, and $h$ on ${\bf T}^n$.

        If $\alpha = (\alpha_1, \ldots, \alpha_n)$ is an $n$-tuple of
integers, then the corresponding Fourier coefficient of a continuous
function $f$ on ${\bf T}^n$ is defined as usual by
\begin{equation}
 \widehat{f}(\alpha) = \frac{1}{(2 \pi)^n} \int_{{\bf T}^n} f(z) \, z^{-\alpha}
                                                                       \, |dz|.
\end{equation}
It is easy to check that
\begin{equation}
\label{widehat{(f * g)}(alpha) = widehat{f}(alpha) widehat{g}(alpha)}
        \widehat{(f * g)}(\alpha) = \widehat{f}(\alpha) \, \widehat{g}(\alpha)
\end{equation}
for all continuous functions $f$, $g$ on ${\bf T}^n$ and $\alpha \in
{\bf Z}^n$.  More precisely, if we substitute the definition of $f *
g$ into the definition of the Fourier coefficient, then we get a
double integral in $z$ and $w$.  This double integral can be evaluated
by integrating in $z$ first, using the change of variables $z \mapsto
z \diamond w$ and the fact that
\begin{equation}
        (z \diamond w)^{-\alpha} = z^{-\alpha} \, w^{-\alpha}
\end{equation}
for all $z, w \in {\bf T}^n$ and $\alpha \in {\bf Z}^n$.  The double
integral then splits into a product of integrals over $z$ and $w$
separately, which leads to (\ref{widehat{(f * g)}(alpha) =
widehat{f}(alpha) widehat{g}(alpha)}).

        Note that the convolution $f * g$ can be defined as before when $f$
is continuous on ${\bf T}^n$ and $g$ is Lebesgue integrable, and satisfies
\begin{equation}
\label{sup_{z in {bf T}^n} |(f * g)(z)| le ...}
\sup_{z \in {\bf T}^n} |(f * g)(z)| \le \Big(\sup_{z \in {\bf T}^n} |f(z)|\Big)
              \, \Big(\frac{1}{(2 \pi)^n} \int_{{\bf T}^n} |g(w)| \, |dw|\Big).
\end{equation}
In this case, it is easy to see that $f * g$ is still continuous,
because $f$ is uniformly continuous on ${\bf T}^n$.  Of course, the
analogous statements also hold when the roles of $f$ and $g$ are
reversed, because convolution is commutative.  If $f$ is bounded and
measurable on ${\bf T}^n$ and $g$ is integrable, then the convolution
$(f * g)(z)$ can be defined in the same way for each $z \in {\bf
T}^n$, and satisfies (\ref{sup_{z in {bf T}^n} |(f * g)(z)| le ...}).
The convolution $f * g$ is actually continuous in this case as well,
as one can show by approximating $g$ by continuous functions with
respect to the $L^1$ norm on ${\bf T}^n$, so that $f * g$ is
approximated uniformly by continuous functions on ${\bf T}^n$ by
(\ref{sup_{z in {bf T}^n} |(f * g)(z)| le ...}) and the previous remarks.

        Suppose that $f$, $g$ are nonnegative real-valued integrable
functions on ${\bf T}^n$.  In this case,
\begin{eqnarray}
\label{frac{1}{(2 pi)^n} int_{{bf T}^n} (f * g)(z) |dz| = ...}
\lefteqn{\frac{1}{(2 \pi)^n} \int_{{\bf T}^n} (f * g)(z) \, |dz| =} \\
 & & \Big(\frac{1}{(2 \pi)^n} \int_{{\bf T}^n} f(z) \, |dz|\Big)
    \Big(\frac{1}{(2 \pi)^n} \int_{{\bf T}^n} g(z) \, |dw|\Big). \nonumber
\end{eqnarray}
To see this, one can substitute the definition of $(f * g)(z)$ into
the integral on the left, which leads to a double integral in $w$ and
$z$.  One can then interchange the order of integration and use the
change of variable $z \mapsto z \diamond w$ to split the double integral
into a product of integrals in $z$ and $w$, as before.  In particular,
it follows that $(f * g)(z)$ is finite for almost every $z \in {\bf
T}^n$.

        Now let $f$, $g$ be integrable complex-valued functions on
${\bf T}^n$.  Observe that
\begin{equation}
        \int_{{\bf T}^n} |f(z \diamond w^{-1})| \, |g(w)| \, |dw| < \infty
\end{equation}
for almost every $z \in {\bf T}^n$, by the argument in the previous
paragraph applied to $|f|$, $|g|$.  Thus $(f * g)(z)$ is defined for
almost every $z \in {\bf T}^n$, and satisfies
\begin{equation}
\label{|(f * g)(z)| le ...}
 |(f * g)(z)| \le \frac{1}{(2 \pi)^n} \int_{{\bf T}^n} |f(z \diamond w^{-1})| \,
                                                                 |g(w)| \, dw.
\end{equation}
Using Fubini's theorem, one may conclude that $f * g$ is an integrable
function on ${\bf T}^n$, and that
\begin{eqnarray}
\lefteqn{\frac{1}{(2 \pi)^n} \int_{{\bf T}^n} |(f * g)(z)| \, |dz| \le} \\
 & & \Big(\frac{1}{(2 \pi)^n} \int_{{\bf T}^n} |f(z)| \, |dz|\Big)
      \Big(\frac{1}{(2 \pi)^n} \int_{{\bf T}^n} |g(w)| \, |dw|\Big). \nonumber
\end{eqnarray}
One can also check that convolution is commutative and associative on
$L^1({\bf T}^n)$, as before.

        If $f$ is an integrable function on ${\bf T}^n$, then the
Fourier coefficients $\widehat{f}(\alpha)$ can be defined in the usual
way, and satisfy
\begin{equation}
 |\widehat{f}(\alpha)| \le \frac{1}{(2 \pi)^n} \int_{{\bf T}^n} |f(z)| \, |dz|
\end{equation}
for each $\alpha \in {\bf Z}^n$.  If $f$ and $g$ are integrable functions
on ${\bf T}^n$, so that their convolution $f * g$ is also integrable, as
in the preceding paragraph, then the Fourier coefficients of $f * g$ are
equal to the product of the Fourier coefficients of $f$ and $g$, as in
(\ref{widehat{(f * g)}(alpha) = widehat{f}(alpha) widehat{g}(alpha)}).
This follows from Fubini's theorem, as before.

\section[\ Convolution on ${\bf R}^n$]{Convolution on ${\bf R}^n$}
\label{convolution on R^n}
\setcounter{equation}{0}

        Let $f$ and $g$ be nonnegative real-valued integrable
functions on ${\bf R}^n$, and put
\begin{equation}
\label{(f * g)(x) = int_{{bf R}^n} f(x - y) g(y) dy}
        (f * g)(x) = \int_{{\bf R}^n} f(x - y) \, g(y) \, dy,
\end{equation}
where $dy$ denotes Lebesgue measure on ${\bf R}^n$, as usual.  It is
easy to see that
\begin{equation}
\label{int_{{bf R}^n} (f * g)(x) dx = ...}
 \int_{{\bf R}^n} (f * g)(x) \, dx = \Big(\int_{{\bf R}^n} f(x) \, dx\Big)
                                      \Big(\int_{{\bf R}^n} g(y) \, dy\Big),
\end{equation}
by interchanging the order of integration and using the change of
variables $x \mapsto x + y$, as in the previous section.  Thus $f * g$
is integrable on ${\bf R}^n$ under these conditions, and finite almost
everywhere on ${\bf R}^n$ in particular.

        If $f$ and $g$ are arbitrary real or complex-valued integrable
functions on ${\bf R}^n$, then it follows that
\begin{equation}
        \int_{{\bf R}^n} |f(x - y)| \, |g(y)| \, dy < \infty
\end{equation}
for almost every $x \in {\bf R}^n$, by applying the preceding argument
to $|f|$ and $|g|$.  This shows that the definition (\ref{(f * g)(x) =
int_{{bf R}^n} f(x - y) g(y) dy}) of $(f * g)(x)$ also makes sense in
this case for almost every $x \in {\bf R}^n$, and satisfies
\begin{equation}
\label{|(f * g)(x)| le int_{{bf R}^n} |f(x - y)| |g(y)| dy}
        |(f * g)(x)| \le \int_{{\bf R}^n} |f(x - y)| \, |g(y)| \, dy.
\end{equation}
One can also check that $f * g$ is measurable, using Fubini's theorem.
Integrating in $x$ as before, we get that
\begin{equation}
\label{int_{{bf R}^n} |(f * g)(x)| dx le ...}
\int_{{\bf R}^n} |(f * g)(x)| \, dx \le \Big(\int_{{\bf R}^n} |f(x)| \, dx\Big)
                                       \Big(\int_{{\bf R}^n} |g(y)| \, dy\Big),
\end{equation}
and that $f * g$ is integrable in particular.

        As in the previous section, it is easy to see that
\begin{equation}
        f * g = g * f,
\end{equation}
using the change of variables $y \mapsto x - y$.  Similarly, one can
verify that
\begin{equation}
\label{(f * g) * h = f * (g * h)}
        (f * g) * h = f * (g * h)
\end{equation}
for any integrable functions $f$, $g$, and $h$ on ${\bf R}^n$.

        If $f$ is an integrable function on ${\bf R}^n$ and $g$ is
bounded and measurable, then the convolution $f * g$ can be defined
using (\ref{(f * g)(x) = int_{{bf R}^n} f(x - y) g(y) dy}) as before,
and satisfies
\begin{equation}
\label{sup_{x in {bf R}^n} |(f * g)(x)| le ...}
\sup_{x \in {\bf R}^n} |(f * g)(x)| \le \Big(\int_{{\bf R}^n} |f(x)| \, dx\Big)
                                       \Big(\sup_{y \in {\bf R}^n} |g(y)|\Big).
\end{equation}
One can also check that $f * g$ is uniformly continuous under these
conditions, as follows.  If $f$ is a continuous function on ${\bf
R}^n$ with compact support, then $f$ is uniformly continuous, and it
is easy to see that $f * g$ is uniformly continuous directly from the
definitions.  Otherwise, if $f$ is any integrable function on ${\bf
R}^n$, then it is well known that $f$ can be approximated by
continuous functions on ${\bf R}^n$ with compact support in the $L^1$
norm.  This implies that $f * g$ can be approximated by uniformly
continuous functions on ${\bf R}^n$ with respect to the supremum norm,
and hence that $f * g$ is uniformly continuous as well.

\section[\ The Fourier transform]{The Fourier transform}
\label{fourier transform}
\setcounter{equation}{0}

        If $f$ is an integrable complex-valued function on ${\bf
R}^n$, then the \emph{Fourier transform} $\widehat{f}$ of $f$ is
defined by
\begin{equation}
\label{widehat{f}(xi) = int_{{bf R}^n} f(x) exp (- i xi cdot x) dx}
\widehat{f}(\xi) = \int_{{\bf R}^n} f(x) \, \exp (- i \xi \cdot x) \, dx.
\end{equation}
Here $\xi \in {\bf R}^n$, and $\xi \cdot x$ is the usual dot product,
given by
\begin{equation}
        \xi \cdot x = \sum_{j = 1}^n \xi_j \, x_j.
\end{equation}
Also, $\exp (- i \xi \cdot x)$ refers to the complex exponential
function, which satisifies $|\exp (i t)| = 1$ for every $t \in {\bf
R}$.  Thus the integrand in (\ref{widehat{f}(xi) = int_{{bf R}^n} f(x)
exp (- i xi cdot x) dx}) is an integrable function, and
\begin{equation}
\label{|widehat{f}(xi)| le int_{{bf R}^n} |f(x)| dx}
        |\widehat{f}(\xi)| \le \int_{{\bf R}^n} |f(x)| \, dx
\end{equation}
for every $\xi \in {\bf R}^n$.

        Let $R$ be a positive real number, and put $f_R(x) = f(x)$
when $|x| \le R$ and $f_R(x) = 0$ when $|x| > R$.  Thus
\begin{equation}
        \widehat{f}_R(\xi)
             = \int_{|x| \le R} f(x) \, \exp (- i \xi \cdot x) \, dx,
\end{equation}
and
\begin{equation}
\label{|widehat{f}(xi) - widehat{f}_R(xi)| le int_{|x| > R} |f(x)| dx}
        |\widehat{f}(\xi) - \widehat{f}_R(\xi)| \le \int_{|x| > R} |f(x)| \, dx
\end{equation}
for every $\xi \in {\bf R}^n$ and $R > 0$.  In particular,
$\widehat{f}_R \to \widehat{f}$ uniformly on ${\bf R}^n$ as $R \to
\infty$.  It is easy to see that $\widehat{f}_R(\xi)$ is uniformly
continuous on ${\bf R}^n$ for each $R > 0$, using the fact that $\exp
(i t)$ is uniformly continuous on the real line.  It follows that
$\widehat{f}(\xi)$ is also uniformly continuous on ${\bf R}^n$, since
it is the uniform limit of uniformly continuous functions on ${\bf
R}^n$.

        Now let $f$, $g$ be integrable functions on the real line, so
that their convolution $f * g$ is also integrable, as in the previous
section.  The Fourier transform of $f * g$ is given by
\begin{eqnarray}
 \widehat{(f * g)}(\xi) & = & \int_{{\bf R}^n} (f * g)(x) \,
                                  \exp (- i \xi \cdot x) \, dx \\
 & = & \int_{{\bf R}^n} \int_{{\bf R}^n} f(x - y) \, g(y) \,
                              \exp (- i \xi \cdot x) \, dy \, dx. \nonumber
\end{eqnarray}
This is the same as
\begin{equation}
        \int_{{\bf R}^n} \int_{{\bf R}^n} f(x) \, g(y) \,
                          \exp (- i \xi \cdot (x + y)) \, dx \, dy,
\end{equation}
by interchanging the order of integration and using the change of
variables $x \mapsto x + y$.  Because $\exp (i (r + t)) = \exp (i r)
\, \exp (i t)$ for every $r, t \in {\bf R}$, this double integral
reduces to
\begin{equation}
        \Big(\int_{{\bf R}^n} f(x) \, \exp (- i \xi \cdot x) \, dx\Big)
         \Big(\int_{{\bf R}^n} g(y) \, \exp (- i \xi \cdot y) \, dy\Big).
\end{equation}
Thus
\begin{equation}
        \widehat{(f * g)}(\xi) = \widehat{f}(\xi) \, \widehat{g}(\xi)
\end{equation}
for every $\xi \in {\bf R}^n$.

\section[\ Holomorphic extensions]{Holomorphic extensions}
\label{holomorphic extensions}
\setcounter{equation}{0}

        Let $L^1({\bf R}^n)$ be the space of Lebesgue integrable
functions on ${\bf R}^n$ equipped with the $L^1$ norm
\begin{equation}
\label{||f||_1 = int_{{bf R}^n} |f(x)| dx}
        \|f\|_1 = \int_{{\bf R}^n} |f(x)| \, dx,
\end{equation}
as usual.  Let us say that $f \in L^1({\bf R}^n)$ has support
contained in a closed set $E \subseteq {\bf R}^n$ if $f(x) = 0$ almost
everywhere on ${\bf R}^n \backslash E$.  The space $L^1_{com}({\bf
R}^n)$ of $f \in L^1({\bf R}^n)$ with compact support is a dense
linear subspace of $L^1({\bf R}^n)$ which is closed under convolution,
in the sense that $f * g \in L^1_{com}({\bf R}^n)$ for every $f$, $g$ in
$L^1_{com}({\bf R}^n)$.  If $f \in L^1_{com}({\bf R}^n)$ is supported
in a compact set $K$, then the Fourier transform $\widehat{f}(\xi)$ extends
to a holomorphic function $\widehat{f}(\zeta)$ on ${\bf C}^n$, given by
\begin{equation}
\label{widehat{f}(zeta) = int_K f(x) exp (- i zeta cdot x) dx}
        \widehat{f}(\zeta) = \int_K f(x) \, \exp (- i \zeta \cdot x) \, dx.
\end{equation}
Here $\zeta \in {\bf C}^n$ may be expressed as $\xi + i \eta$, with
$\xi, \eta \in {\bf R}^n$, and
\begin{equation}
        \zeta \cdot x = \sum_{j = 1}^n \zeta_j \, x_j,
\end{equation}
as before.  Thus (\ref{widehat{f}(zeta) = int_K f(x) exp (- i zeta
cdot x) dx}) reduces to (\ref{widehat{f}(xi) = int_{{bf R}^n} f(x) exp
(- i xi cdot x) dx}) when $\zeta = \xi \in {\bf R}^n$, and otherwise
it is easy to check that $\widehat{f}(\zeta)$ is a holomorphic
function on ${\bf C}^n$, since the exponential function is holomorphic.
In addition,
\begin{equation}
\label{widehat{(f * g)}(zeta) = widehat{f}(zeta) widehat{g}(zeta)}
        \widehat{(f * g)}(\zeta) = \widehat{f}(\zeta) \, \widehat{g}(\zeta)
\end{equation}
for every $f, g \in L^1_{com}({\bf R}^n)$ and $\zeta \in {\bf C}^n$,
for the same reasons as in the previous section.

        Let $L^1_+({\bf R})$ be the space of $f \in L^1({\bf R})$ that
are supported in $[0, \infty)$, and let $L^1_-({\bf R})$ be the space
of $f \in L^1({\bf R})$ that are supported in $(-\infty, 0]$.  These
are closed linear subspaces of $L^1({\bf R})$ that are closed under
convolution, in the sense that $f * g \in L^1_+({\bf R})$ when $f, g
\in L^1_+({\bf R})$, and similarly for $L^1_-({\bf R})$.  Let $H_+$,
$H_-$ be the upper and lower open half-planes in the complex plane,
consisting of complex numbers with positive and negative imaginary
parts, respectively.  Thus their closures $\overline{H}_+$,
$\overline{H}_-$ are the upper and lower closed half-planes in ${\bf
C}$, consisting of complex numbers with imaginary part greater than or
equal to $0$ and less than or equal to $0$, respectively.
If $f \in L^1_+({\bf R})$, then
\begin{equation}
\label{widehat{f}(zeta) = int_0^infty f(x) exp (- i zeta x) dx = ...}
 \widehat{f}(\zeta) = \int_0^\infty f(x) \, \exp (- i \zeta \, x) \, dx
     = \int_0^\infty f(x) \, \exp (- i \xi \, x + \eta \, x) \, dx
\end{equation}
is defined for all $\zeta = \xi + i \, \eta \in \overline{H}_-$.
In this case, $\eta \le 0$, so that
\begin{equation}
        |\exp (- i \xi \, x + \eta \, x)| = \exp (\eta \, x) \le 1
\end{equation}
for every $x \ge 0$, and hence
\begin{equation}
        |\widehat{f}(\zeta)| \le \|f\|_1
\end{equation}
for every $\zeta \in \overline{H}_-$.  As in the previous section, one
can check that $\widehat{f}(\zeta)$ is uniformly continuous on
$\overline{H}_-$.  This uses the fact that $\exp (- i \zeta \, x)$ is
uniformly continuous as a function of $\zeta$ on $\overline{H}_-$ for
each $x \ge 0$, and it is easier to first consider the case where $f$
has compact support in $[0, \infty)$, and then get the same conclusion
for any $f \in L_+({\bf R})$ by approximation.  One can also check
that $\widehat{f}(\zeta)$ is holomorphic on $H_-$, using the
holomorphicity of the exponential function and the integrability of
the expressions in (\ref{widehat{f}(zeta) = int_0^infty f(x) exp (- i
zeta x) dx = ...}).  If $f, g \in L^1_+({\bf R})$, then
\begin{equation}
        \widehat{(f * g)}(\zeta) = \widehat{f}(\zeta) \, \widehat{g}(\zeta)
\end{equation}
for every $\zeta \in \overline{H}_-$, for the same reasons as before.
In the same way, the Fourier transform of a function in $L^1_-({\bf
R})$ has a natural extension to a bounded uniformly continuous
function on $\overline{H}_+$ that is holomorphic on $H_+$, and with
the analogous property for convolutions.

        Let $\epsilon = (\epsilon_1, \ldots, \epsilon_n)$ be an
$n$-tuple with $\epsilon_j \in \{1, -1\}$ for each $j$, which is
to say an element of $\{1, -1\}^n$.  Put
\begin{equation}
\label{Q_{n, epsilon} = ...}
        Q_{n,\epsilon} = \{x \in {\bf R}^n : \epsilon_j \, x_j \ge 0
                                          \hbox{ for } j = 1, \ldots, n\},
\end{equation}
which is the closed ``quadrant'' in ${\bf R}^n$ associated to $\epsilon$.
Let $L^1_\epsilon({\bf R}^n)$ be the set of $f \in L^1({\bf R}^n)$
which are supported in $Q_{n, \epsilon}$.  It is easy to see that
$L^1_\epsilon({\bf R}^n)$ is a closed linear subspace of $L^1({\bf
R}^n)$ that is closed with respect to convolution, as before.
Consider
\begin{equation}
\label{H_{n, epsilon} = ...}
 H_{n, \epsilon} = \{\zeta = \xi + i \, \eta \in {\bf C}^n :
                      \epsilon_j \, \eta_j > 0 \hbox{ for } j = 1, \ldots, n\},
\end{equation}
so that the closure $\overline{H}_{n, \epsilon}$ of $H_{n, \epsilon}$
consists of the $\zeta = \xi + \eta \in {\bf C}^n$ with $\eta \in
Q_{n, \epsilon}$.  If $f \in L^1_\epsilon({\bf R}^n)$, then
\begin{eqnarray}
\label{widehat{f}(zeta), f in L^1_epsilon({bf R}^n)}
        \widehat{f}(\zeta)
 & = & \int_{Q_{n, \epsilon}} f(x) \, \exp (- i \zeta \cdot x) \, dx \\
 & = & \int_{Q_{n, \epsilon}} f(x) \, \exp (- i \xi \cdot x + \eta \cdot x)
                                                               \, dx \nonumber
\end{eqnarray}
is defined for every $\zeta = \xi + \eta \in H_{n, -\epsilon}$, where
$-\epsilon = (-\epsilon_1, \ldots, -\epsilon_n)$.  In this case, $\eta
\cdot x \le 0$ for every $x \in Q_{n, \epsilon}$, so that
\begin{equation}
 |\exp (- i \xi \cdot x + \eta \cdot x)| = \exp (\eta \cdot x) \le 1,
\end{equation}
and hence $|\widehat{f}(\zeta)| \le \|f\|_1$ for every $\zeta \in
\overline{H}_{n, - \epsilon}$.  As before, one can check that
$\widehat{f}(\zeta)$ is uniformly continuous on $\overline{H}_{n, -
\epsilon}$, and holomorphic on $H_{n, -\epsilon}$.  If $f, g \in
L^1_\epsilon({\bf R}^n)$, then the extension of the Fourier transform
of $f * g$ to $H_{n, - \epsilon}$ is equal to the product of the
extensions of the Fourier transforms of $f$ and $g$ to $H_{n, -
\epsilon}$, as usual.

\section[\ The Riemann--Lebesgue lemma]{The Riemann--Lebesgue lemma}
\label{riemann--lebesgue lemma}
\setcounter{equation}{0}

        If $a$, $b$ are real numbers with $a < b$, then the Fourier
transform of the indicator function ${\bf 1}_{[a, b]}$ of the interval
$[a, b]$ in the real line is equal to
\begin{equation}
\label{int_a^b exp (- i xi x) dx = i xi^{-1} (exp (- i xi b) - exp (- i xi a))}
        \int_a^b \exp (- i \xi \, x) \, dx 
          = i \xi^{-1} \, (\exp (- i \, \xi \, b) - \exp (- i \xi \, a))
\end{equation}
when $\xi \ne 0$, and to $b - a$ when $\xi = 0$.  In particular,
this tends to $0$ as $|\xi| \to \infty$.

        If $f \in L^1({\bf R})$, then the Riemann--Lebesgue lemma states that
\begin{equation}
\label{lim_{|xi| to infty} widehat{f}(xi) = 0}
        \lim_{|\xi| \to \infty} \widehat{f}(\xi) = 0.
\end{equation}
This follows immediately from the remarks in the previous paragraph
when $f$ is a step function, which is to say a finite linear
combination of indicator functions of intervals in the real line.
Otherwise, any integrable function $f$ on the real line can be
approximated by step functions in the $L^1$ norm, which leads to an
approximation of the Fourier transform $\widehat{f}$ of $f$ by Fourier
transforms of step functions in the supremum norm, by
(\ref{|widehat{f}(xi)| le int_{{bf R}^n} |f(x)| dx}).  This permits
one to derive (\ref{lim_{|xi| to infty} widehat{f}(xi) = 0}) for $f$
from the corresponding statement for step functions.

       This also works for integrable functions on ${\bf R}^n$.  In
this case, we can start with a rectangular box $B$ in ${\bf R}^n$,
which is to say the Cartesian product of $n$ intervals in the real
line.  The indicator function of $B$ on ${\bf R}^n$ is the same as the
product of the $n$ indicator functions of the corresponding intervals
in ${\bf R}$, as functions of $x_1, \ldots, x_n$.  Thus the Fourier
transform of the indicator function of $B$ is the same as the product
of the $n$ one-dimensional Fourier transforms of these indicator
functions of intervals in ${\bf R}$, as functions of $\xi_1, \ldots,
\xi_n$.  This implies that the Fourier transform of the indicator
function of $B$ tends to $0$ at infinity, as before.  Hence the
Fourier transform of any finite linear combination of indicator
functions of rectangular boxes in ${\bf R}^n$ also tends to $0$ at
infinity.  Any integrable function $f$ on ${\bf R}^n$ can be
approximated by a finite linear combination of indicator functions of
rectangular boxes in the $L^1$ norm, which implies (\ref{lim_{|xi| to
infty} widehat{f}(xi) = 0}) as in the one-dimensional case.

        As in the previous section, the Fourier transform of the
indicator function ${\bf 1}_{[a, b]}$ of an interval $[a, b]$ in the
real line extends to a holomorphic function on the complex plane,
given by
\begin{equation}
        \int_a^b \exp (- i \zeta \, x) \, dx
         = i \, \zeta^{-1} \, (\exp (- i \zeta \, b) - \exp (- i \zeta \, a))
\end{equation}
when $\zeta \ne 0$, and equal to $b - a$ when $\zeta = 0$.  If $a, b
\ge 0$, then it is easy to see that this tends to $0$ as $|\zeta| \to
\infty$ when $\zeta$ is in the closed lower half-plane
$\overline{H}_-$.  If $f \in L^1_+({\bf R})$, so that the Fourier
transform of $f$ has a natural extension $\widehat{f}(\zeta)$ to
$\zeta \in \overline{H}_-$, as in the preceding section, then one can
use this to show that $\widehat{f}(\zeta) \to 0$ as $|\zeta| \to \infty$
in $\overline{H}_-$, by approximating $f$ by step functions as before.
Of course, there is an analogous statement for the extension to the
closed upper half-plane $\overline{H}_+$ of the Fourier transform of a
function in $L^1_-({\bf R}^n)$.  There is also an analogous statement
for the extension to $\overline{H}_{n, -\epsilon}$ of the Fourier transform
of a function in $L^1_\epsilon({\bf R}^n)$, as in the previous section.

\section[\ Translation and multiplication]{Translation and multiplication}
\label{translation, multiplication}
\setcounter{equation}{0}

        If $f \in L^1({\bf R}^n)$ and $t \in {\bf R}^n$, then let
$T_t(f)$ be the function on ${\bf R}^n$ obtained by translating $f$ by
$t$, so that
\begin{equation}
\label{T_t(f)(x) = f(x - t)}
        T_t(f)(x) = f(x - t).
\end{equation}
Thus $T_t(f) \in L^1({\bf R}^n)$ too, and $\|T_t(f)\|_1 = \|f\|_1$.
It is easy to see that
\begin{equation}
\label{widehat{(T_t(f))}(xi) = exp (- i xi cdot t) widehat{f}(xi)}
        \widehat{(T_t(f))}(\xi) = \exp (- i \xi \cdot t) \, \widehat{f}(\xi),
\end{equation}
for each $\xi \in {\bf R}^n$, using the change of variable $x \mapsto
x + t$ in the definition of $\widehat{T_t(f)}$.  Similarly, if $f$ has
compact support in ${\bf R}^n$, then $T_t(f)$ does too, and the
natural extension of the Fourier transform of $T_t(f)$ to a
holomorphic function on ${\bf C}^n$ satisfies
\begin{equation}
\label{widehat{(T_t(f))}(zeta) = exp (- i zeta cdot t) widehat{f}(zeta)}
 \widehat{(T_t(f))}(\zeta) = \exp (- i \zeta \cdot t) \, \widehat{f}(\zeta)
\end{equation}
for each $\zeta \in {\bf C}^n$.

        Suppose now that $\epsilon \in \{1, -1\}^n$, and that $f \in
L^1_\epsilon({\bf R}^n)$, as in Section \ref{holomorphic extensions}.
Thus $f$ is supported in the ``quadrant'' $Q_{n, \epsilon}$ defined in
(\ref{Q_{n, epsilon} = ...}).  If $t \in Q_{n, \epsilon}$, then it is
easy to see that $T_t(f)$ is supported in $Q_{n, \epsilon}$ as well,
so that $T_t(f) \in L^1_\epsilon({\bf R}^n)$.  As in Section
\ref{holomorphic extensions}, the Fourier transform of $f$ and
$T_t(f)$ have natural extensions to $\overline{H}_{n, -\epsilon}$,
which are related by the same expression (\ref{widehat{(T_t(f))}(zeta)
= exp (- i zeta cdot t) widehat{f}(zeta)}) as in the previous paragraph.
Note that
\begin{equation}
\label{|exp (- i zeta cdot t)| le 1}
        |\exp (- i \zeta \cdot t)| \le 1
\end{equation}
for each $\zeta \in H_{n, -\epsilon}$ and $t \in Q_{n, \epsilon}$,
as in Section \ref{holomorphic extensions}.

        If $w \in {\bf R}^n$ and $f \in L^1({\bf R}^n)$, then let
$M_w(f)$ be the function on ${\bf R}^n$ defined by multiplying $f$ by
$\exp (i w \cdot x)$, so that
\begin{equation}
\label{(M_w(f))(x) = exp (i w cdot x) f(x)}
        (M_w(f))(x) = \exp (i w \cdot x) \, f(x).
\end{equation}
Thus $M_w(f) \in L^1({\bf R}^n)$ and $\|M_w(f)\|_1 = \|f\|_1$, since
$|\exp (i w \cdot x)| = 1$ for every $x, w \in {\bf R}^n$.
It is easy to see that
\begin{equation}
\label{widehat{(M_w(f))}(xi) = widehat{f}(xi - w)}
        \widehat{(M_w(f))}(\xi) = \widehat{f}(\xi - w)
\end{equation}
for every $\xi, w \in {\bf R}^n$, directly from the definition of the
Fourier transform.  If $w \in {\bf C}^n$, then we can still define
$M_w(f)$ for $f \in L^1({\bf R}^n)$ by (\ref{(M_w(f))(x) = exp (i w
cdot x) f(x)}), and $M_w(f)$ will be locally integrable on ${\bf
R}^n$, but it may not be integrable on ${\bf R}^n$.  However, if $f$
has compact support in ${\bf R}^n$, then $M_w(f)$ also has compact
support in ${\bf R}^n$ for every $w \in {\bf C}^n$, and $M_w(f)$ is
integrable on ${\bf R}^n$ for every $w \in {\bf C}^n$.  In this case,
the Fourier transform of $f$ extends to a holomorphic function on
${\bf C}^n$, as in Section \ref{holomorphic extensions}, and the
Fourier transform of $M_w(f)$ is defined and extends to a holomorphic
function on ${\bf C}^n$ for each $w \in {\bf C}^n$.  As before,
we have that
\begin{equation}
\label{widehat{(M_w(f))}(zeta) = widehat{f}(zeta - w)}
        \widehat{(M_w(f))}(\zeta) = \widehat{f}(\zeta - w)
\end{equation}
for every $\zeta, w \in {\bf C}^n$ when $f \in L^1_{com}({\bf R}^n)$.

        Let $\epsilon$ be an element of $\{1, -1\}^n$ again, and
suppose that $f \in L^1_\epsilon({\bf R}^n)$.  As before, $M_w(f)$ is
a locally integrable function on ${\bf R}^n$ with support contained in
$Q_{n, \epsilon}$ for every $w \in {\bf C}^n$.  If $w \in
\overline{H}_{n, \epsilon}$, then $|\exp (i w \cdot x)| \le 1$ for
every $x \in Q_{n, \epsilon}$, and hence $M_w(f) \in L^1_\epsilon({\bf
R}^n)$, with $\|M_w(f)\|_1 \le \|f\|_1$.  As in Section \ref{holomorphic 
extensions}, the Fourier transforms of $f$ and $M_w(f)$ have natural
extensions to $\overline{H}_{n, -\epsilon}$ under these conditions,
and one can check that they are related as in (\ref{widehat{(M_w(f))}(zeta)
= widehat{f}(zeta - w)}) for each $\zeta \in \overline{H}_{n, -\epsilon}$.
Note that $\widehat{f}(\zeta - w)$ is defined in this case, because
$-w$ and hence $\zeta - w$ is in $\overline{H}_{n, -\epsilon}$.

\section[\ Some examples]{Some examples}
\label{some examples}
\setcounter{equation}{0}

        Let $a$ be a positive real number, and put
\begin{equation}
        q_{a, +}(x) = \exp (- a \, x)
\end{equation}
when $x \ge 0$, and $q_{a, +}(x) = 0$ when $x < 0$.  Thus $q_{a, +}
\in L^1_+({\bf R})$, and so the Fourier transform of $q_{a, +}$ should
have a natural extension to the closed lower half-plane in ${\bf C}$,
as in Section \ref{holomorphic extensions}.  More precisely,
\begin{equation}
\label{widehat{q_{a, +}}(zeta) = ...}
 \widehat{q_{a, +}}(\zeta) = \int_0^\infty \exp (- a \, x - i \zeta \, x) \, dx
                           = \frac{-1}{- a - i \zeta} = \frac{1}{a + i \zeta}
\end{equation}
for every $\zeta \in \overline{H}_-$.  Note that $\re (a + i \zeta)
\ge a > 0$ when $\zeta \in \overline{H}_-$ and $a > 0$.

        Similarly, put
\begin{equation}
        q_{a, -}(x) = \exp (a \, x) = \exp (- a \, |x|)
\end{equation}
when $x \le 0$, and $q_{a, -}(x) = 0$ when $x > 0$.  In this case,
$q_{a, -} \in L^1_-({\bf R}^n)$, so that the Fourier transform of
$q_{a, -}$ should have a natural extension to the closed upper
half-plane in ${\bf C}$, as in Section \ref{holomorphic extensions}.
Indeed,
\begin{equation}
\label{widehat{q_{a, -}}(zeta) = ...}
\widehat{q_{a, -}}(\zeta) = \int_{-\infty}^0 \exp (a \, x - i \zeta \, x) \, dx
                          = \frac{1}{a - i \zeta}
\end{equation}
for every $\zeta \in \overline{H}_+$.  As before, $\re (a - i \zeta)
\ge a > 0$ when $\zeta \in \overline{H}_+$ and $a > 0$.

        Now let $a = (a_1, \ldots, a_n)$ be an $n$-tuple of positive
real numbers, and let $\epsilon$ be an element of $\{1, -1\}^n$.  Put
\begin{equation}
\label{q_{n, a, epsilon}(x) = ...}
        q_{n, a, \epsilon}(x)
             = \exp \Big(- \sum_{j = 1}^n a_j \, \epsilon_j \, x_j\Big)
             = \exp \Big(- \sum_{j = 1}^n a_j |x_j|\Big)
\end{equation}
when $x \in Q_{n, \epsilon}$, and $q_{n, a, \epsilon}(x) = 0$ when
$x \in {\bf R}^n \backslash Q_{n, \epsilon}$.  Equivalently,
\begin{equation}
\label{q_{n, a, epsilon}(x) = prod_{j = 1}^n q_{a_j, epsilon_j}(x_j)}
        q_{n, a, \epsilon}(x) = \prod_{j = 1}^n q_{a_j, \epsilon_j}(x_j),
\end{equation}
where the subscript $\epsilon_j$ on the right should be interpreted as
$+$ when $\epsilon_j = 1$ and as $-$ when $\epsilon_j = -1$.  Of
course, $q_{n, a, \epsilon} \in L^1_\epsilon({\bf R}^n)$, and so its
Fourier transform should have a natural extension to $\overline{H}_{n,
-\epsilon}$, as usual.  In fact, the Fourier transform of $q_{n, a,
\epsilon}$ can be given as the product of the one-dimensional Fourier
transforms of the factors $q_{a_j, \epsilon_j}$, so that
\begin{equation}
        \widehat{q_{n, a, \epsilon}}(\zeta)
          = \prod_{j = 1}^n \widehat{q_{a_j, \epsilon_j}}(\zeta_j)
          = \prod_{j = 1}^n \frac{1}{(a_j + i \epsilon_j \, \zeta_j)}
\end{equation}
for each $\zeta \in \overline{H}_{n, -\epsilon}$.

\section[\ Some examples, continued]{Some examples, continued}
\label{some examples, continued}
\setcounter{equation}{0}

        Let $a$ be a positive real number, and put
\begin{equation}
\label{p_a(x) = exp (- a |x|) = q_{a, +}(x) + q_{a, -}(x)}
        p_a(x) = \exp (- a \, |x|) = q_{a, +}(x) + q_{a, -}(x).
\end{equation}
This defines an integrable function on the real line, whose Fourier
transform is given by
\begin{equation}
\label{widehat{p_a}(xi) = ...}
 \widehat{p_a}(\xi) = \widehat{q_{a, +}}(\xi) + \widehat{q_{a, -}}(\xi)
            = \frac{1}{a + i \xi} + \frac{1}{a - i \xi}
            = 2 \re \Big(\frac{1}{a + i \xi}\Big)
\end{equation}
for each $\xi \in {\bf R}$.  Of course,
\begin{equation}
\label{frac{1}{a + i xi} = ...}
        \frac{1}{a + i \xi} = \frac{a - i \xi}{(a + i \xi) (a - i \xi)}
                            = \frac{a - i \xi}{a^2 + \xi^2},
\end{equation}
and so (\ref{widehat{p_a}(xi) = ...}) is the same as
\begin{equation}
\label{widehat{p_a}(xi) = frac{2 a}{a^2 + xi^2}}
        \widehat{p_a}(\xi) = \frac{2 a}{a^2 + \xi^2}.
\end{equation}

        It follows easily from (\ref{widehat{p_a}(xi) = frac{2 a}{a^2
+ xi^2}}) that $\widehat{p_a}(\xi)$ is an integrable function of $\xi$
on the real line.  In order to compute its integral, observe that
\begin{equation}
\label{int_{bf R} widehat{p_a}(xi) d xi = ...}
        \int_{\bf R} \widehat{p_a}(\xi) \, d\xi
 = \lim_{R \to \infty} \int_{-R}^R \widehat{p_a}(\xi) \, d\xi
 = \lim_{R \to \infty} 2 \re \int_{-R}^R \frac{1}{a + i \xi} \, d\xi.
\end{equation}
Using the change of variables $\xi \mapsto R \, \xi$, we get that
\begin{equation}
        \int_{-R}^R \frac{1}{a + i \xi} \, d\xi
          = \int_{-1}^1 \frac{1}{a + R \, \xi} \, R \, d\xi
          = \int_{-1}^1 \frac{1}{a \, R^{-1} + i \xi} \, d\xi
\end{equation}
for each $R > 0$.  Hence
\begin{equation}
        \int_{\bf R} \widehat{p_a}(\xi) \, d\xi
          = \lim_{r \to 0+} 2 \re \int_{-1}^1 \frac{1}{r + i \xi} \, d\xi.
\end{equation}

        Let $\log z$ be the principal branch of the logarithm.
Remember that this is a holomorphic function defined on the set of $z
\in {\bf C}$ such that $z$ is not a real number less than or equal to
$0$, which agrees with the ordinary natural logarithm of $z$ when $z$
is a positive real number, and whose derivative is equal to $1/z$.  Thus
\begin{equation}
        \int_{-1}^1 \frac{1}{r + i \xi} \, i d\xi = \log (r + i) - \log (r - i)
\end{equation}
for each $r > 0$, which implies that
\begin{eqnarray}
        2 \re \int_{-1}^1 \frac{1}{r + i \xi} \, d\xi
         & = & 2 \im \int_{-1}^1 \frac{1}{r + i \xi} \, i d\xi \\
         & = & 2 \im (\log (r + i) - \log (r - i)). \nonumber
\end{eqnarray}
Taking the limit as $r \to 0+$, we get that
\begin{equation}
\label{int_{bf R} widehat{p_a}(xi) d xi = 2 im (log i - log (-i)) = 2 pi}
        \int_{\bf R} \widehat{p_a}(\xi) \, d\xi
          = 2 \im (\log i - \log (-i)) = 2 \pi,
\end{equation}
since $\log i = (\pi/2) i$ and $\log (-i) = - (\pi/2) i$.

        Similarly, if $a = (a_1, \ldots, a_n)$ is an $n$-tuple of
positive real numbers, then
\begin{equation}
        p_{n, a}(x) = \prod_{j = 1}^n p_{a_j}(x_j)
                    = \exp \Big(- \sum_{j = 1}^n a_j \, |x_j|\Big)
\end{equation}
is an integrable function on ${\bf R}^n$.  The Fourier transform of
$p_{n, a}$ is the product of the one-dimensional Fourier transforms
of its factors, given by
\begin{equation}
        \widehat{p_{n, a}}(\xi) = \prod_{j = 1}^n \widehat{p_{a_j}}(\xi_j)
               = \prod_{j = 1}^n \frac{2 \, a_j}{(a_j^2 + \xi_j^2)}.
\end{equation}
The integral of this is equal to the product of the one-dimensional
integrals of its factors, so that
\begin{equation}
 \int_{{\bf R}^n} \widehat{p_{n, a}}(\xi) \, d\xi = (2 \pi)^n.
\end{equation}

\section[\ The multiplication formula]{The multiplication formula}
\label{multiplication formula}
\setcounter{equation}{0}

        Let $f$, $g$ be integrable functions on ${\bf R}^n$.  The
\emph{multiplication formula} states that
\begin{equation}
\label{int_{R^n} widehat{f}(xi) g(xi) d xi = int_{R^n} f(x) widehat{g}(x) dx}
        \int_{{\bf R}^n} \widehat{f}(\xi) \, g(\xi) \, d\xi
          = \int_{{\bf R}^n} f(x) \, \widehat{g}(x) \, dx.
\end{equation}
Note that both sides of this equation make sense, because the Fourier
transforms of $f$ and $g$ are bounded.  Equivalently, (\ref{int_{R^n}
widehat{f}(xi) g(xi) d xi = int_{R^n} f(x) widehat{g}(x) dx}) states that
\begin{eqnarray}
\lefteqn{\int_{{\bf R}^n} \Big(\int_{{\bf R}^n} f(x) \, \exp (- i \xi \cdot x)
                                          \, dx\Big) \, g(\xi) \, d\xi} \\
 & = & \int_{{\bf R}^n} \Big(\int_{{\bf R}^n} g(\xi) \, \exp (- i x \cdot \xi)
                                            \, d\xi\Big) f(x) \, dx, \nonumber
\end{eqnarray}
which follows from Fubini's theorem.

        Let $h$ be integrable function on ${\bf R}^n$, and let $w$ be
an element of ${\bf R}^n$.  If
\begin{equation}
        g(\xi) = \exp (i \xi \cdot w) \, h(\xi),
\end{equation}
then
\begin{equation}
        \widehat{g}(x) = \widehat{h}(x - w),
\end{equation}
as in Section \ref{translation, multiplication}.  If $f \in L^1({\bf
R}^n)$, then the multiplication formula implies that
\begin{equation}
\label{int_{R^n} widehat{f}(xi) exp (i xi cdot w) h(xi) = ...}
 \int_{{\bf R}^n} \widehat{f}(\xi) \, \exp (i \xi \cdot w) \, h(\xi) \, d\xi
                        = \int_{{\bf R}^n} f(x) \, \widehat{h}(x - w) \, dx.
\end{equation}
The right side is similar to $(f * \widehat{h})(w)$, but not quite the
same.

        If $h_1(\xi) = h(-\xi)$, then
\begin{eqnarray}
\label{widehat{h_1}(x) = ...}
 \widehat{h_1}(x) & = &
          \int_{{\bf R}^n} h(-\xi) \, \exp (- i \xi \cdot x) \, d\xi \\
 & = & \int_{{\bf R}^n} h(\xi) \, \exp (i \xi \cdot x) \, dx = \widehat{h}(-x),
                                                               \nonumber
\end{eqnarray}
using the change of variable $x \mapsto -x$.  Hence
\begin{equation}
\label{int_{R^n} f(x) widehat{h}(x - w) dx = ... = (f * widehat{h_1})(w)}
        \int_{{\bf R}^n} f(x) \, \widehat{h}(x - w) \, dx
          = \int_{{\bf R}^n} f(x) \, \widehat{h_1}(w - x) \, dx
          = (f * \widehat{h_1})(w).
\end{equation}

        Suppose now that $h$ is an even function on ${\bf R}^n$, so
that $h_1 = h$.  Thus $\widehat{h}$ is even too, by
(\ref{widehat{h_1}(x) = ...}).  In this case, (\ref{int_{R^n}
widehat{f}(xi) exp (i xi cdot w) h(xi) = ...}) reduces to
\begin{equation}
\label{int_{R^n} widehat{f}(xi) exp (i xi cdot w) h(xi) d xi = ..., 2}
 \int_{{\bf R}^n} \widehat{f}(\xi) \, \exp (i \xi \cdot w) \, h(\xi) \, d\xi
         = (f * \widehat{h})(w).
\end{equation}

\section[\ Convergence]{Convergence}
\label{convergence}
\setcounter{equation}{0}

        Let $a = (a_1, \ldots, a_n)$ be an $n$-tuple of positive real
numbers, and put
\begin{equation}
        P_{n, a}(x) = \pi^{-n} \prod_{j = 1}^n \frac{a_j}{(a_j^2 + x_j^2)}.
\end{equation}
Thus $P_{n, a}$ is a nonnegative integrable function on ${\bf R}^n$
that satisfies
\begin{equation}
        \int_{{\bf R}^n} P_{n, a}(x) \, dx = 1
\end{equation}
for each $a$, as in Section \ref{some examples, continued}.

        Let $f$ be a bounded continuous function on ${\bf R}^n$.  By
standard arguments,
\begin{equation}
\label{lim_{a to 0} (P_{n, a} * f)(x) = f(x)}
        \lim_{a \to 0} (P_{n, a} * f)(x) = f(x)
\end{equation}
for each $x \in {\bf R}^n$.  Because $f$ is uniformly continuous on
compact subsets of ${\bf R}^n$, one also gets uniform convergence on
compact subsets of ${\bf R}^n$ in (\ref{lim_{a to 0} (P_{n, a} * f)(x)
= f(x)}).  If $f$ is uniformly continuous on ${\bf R}^n$, then one
gets uniform convergence on ${\bf R}^n$.

        If $f$ is a continuous function on ${\bf R}^n$ with compact
support, then $f$ is bounded and uniformly continuous in particular,
so that $P_{n, a} * f \to f$ uniformly on ${\bf R}^n$ as $a \to 0$, as
in the previous paragraph.  In this case, it is easy to check that
$P_{n, a} * f \to f$ as $a \to 0$ in the $L^1$ norm on ${\bf R}^n$ too.

        If $f$ is any integrable function on ${\bf R}^n$, then
\begin{equation}
        \|P_{n, a} * f\|_1 \le \|P_{n, a}\|_1 \, \|f\|_1 = \|f\|_1
\end{equation}
for each $a$, as in Section \ref{convolution on R^n}.  One can also check that
$P_{n, a} * f \to f$ as $a \to 0$ in the $L^1$ norm on ${\bf R}^n$, since this
holds on a dense subset of $L^1({\bf R}^n)$, as in the preceding paragraph.

\section[\ Inversion]{Inversion}
\label{inversion}
\setcounter{equation}{0}

        If $f$ is an integrable function on ${\bf R}^n$, then
\begin{equation}
\label{inversion formula}
        \int_{{\bf R}^n} \widehat{f}(\xi) \, \exp (i \xi \cdot w)
                \, \exp \Big(- \sum_{j = 1}^n a_j \, |\xi_j|\Big) \, d\xi
          = (2 \pi)^n \, (P_{n, a} * f)(w)
\end{equation}
for every $n$-tuple $a = (a_1, \ldots, a_n)$ of positive real numbers
and $w \in {\bf R}^n$.  This follows from (\ref{int_{R^n}
widehat{f}(xi) exp (i xi cdot w) h(xi) d xi = ..., 2}), with $h$ equal
to $p_{n, a}$ from Section \ref{some examples, continued}.  This also
uses the fact that $p_{n, a}$ is even and satisfies
\begin{equation}
\label{widehat{p}_{n, a} = (2 pi)^n P_{n, a}}
        \widehat{p}_{n, a} = (2 \pi)^n \, P_{n, a},
\end{equation}
where $P_{n, a}$ is as in the previous section.

        If $\widehat{f}$ is also integrable on ${\bf R}^n$, then it is
easy to see that
\begin{eqnarray}
\lefteqn{\lim_{a \to 0} \int_{{\bf R}^n} \widehat{f}(\xi) \,
\exp (i \xi \cdot w) \, \exp \Big(- \sum_{j = 1}^n a_j \, |\xi_j|\Big) \, d\xi}
                                                                         \\
 & = & \int_{{\bf R}^n} \widehat{f}(\xi) \, \exp (i \xi \cdot w) \, d\xi
                                                                   \nonumber
\end{eqnarray}
for every $w \in {\bf R}^n$.  More precisely,
\begin{equation}
\label{widehat{f}(xi) exp (- sum_{j = 1}^n a_j |xi_j|) to widehat{f}(xi)}
        \widehat{f}(\xi) \, \exp \Big(- \sum_{j = 1}^n a_j \, |\xi_j|\Big)
                                                        \to \widehat{f}(\xi)
\end{equation}
as $a \to 0$ in the $L^1$ norm on ${\bf R}^n$, so that one has uniform
convergence in $w$ in the previous statement.  This can be derived
from the dominated convergence theorem, but one can also use the same
type of argument a bit more directly.  The main points are that
\begin{equation}
        \exp \Big(- \sum_{j = 1}^n a_j \, |\xi_j|\Big) \le 1
\end{equation}
for every $a$ and $\xi$, and that
\begin{equation}
        \exp \Big(- \sum_{j = 1}^n a_j \, |\xi_j|\Big) \to 1
\end{equation}
as $a \to 0$ uniformly on compact subsets of ${\bf R}^n$.

        It follows that
\begin{equation}
\label{int_{R^n} widehat{f}(xi) exp (i xi cdot w) d xi = (2 pi)^n f(w)}
        \int_{{\bf R}^n} \widehat{f}(\xi) \, \exp (i \xi \cdot w) \, d\xi 
                                                        = (2 \pi)^n \, f(w)
\end{equation}
for almost every $w \in {\bf R}^n$ when $f$ and $\widehat{f}$ are
integrable functions on ${\bf R}^n$, since $P_{n, a} * f \to f$ in
$L^1({\bf R}^n)$ as $a \to 0$, as in the preceding section.  In
particular, $f = 0$ almost everywhere on ${\bf R}^n$ when $\widehat{f}
= 0$.

\section[\ Measures on ${\bf T}^n$]{Measures on ${\bf T}^n$}
\label{measures on T^n}
\setcounter{equation}{0}

        There are two basic ways to think about Borel measures on
${\bf T}^n$.  The first is as countably-additive real or
complex-valued functions on the $\sigma$-algebra of Borel subsets of
${\bf T}^n$.  The second way is to look at countinuous linear
functionals on the space $C({\bf T}^n)$ of continuous real or
complex-valued functions on ${\bf T}^n$, with respect to the supremum
norm on $C({\bf T}^n)$.

        If $\mu$ is a countably-additive real or complex Borel measure
on ${\bf T}^n$, then there is a finite nonnegative Borel measure
$|\mu|$ on ${\bf T}^n$ associated to it, known as the total variation
measure corresponding to $\mu$.  This is characterized by the fact that
\begin{equation}
\label{|mu(E)| le |mu|(E)}
        |\mu(E)| \le |\mu|(E)
\end{equation}
for every Borel set $E \subseteq {\bf T}^n$, and that $|\mu|$ is the
smallest nonnegative Borel measure on ${\bf T}^n$ with this property.
More precisely, if $\nu$ is a nonnegative Borel measure on ${\bf T}^n$
such that $|\mu(E)| \le \nu(E)$ for every Borel set $E \subseteq {\bf
T}^n$, then $|\mu|(E) \le \nu(E)$ for every Borel set $E \subseteq
{\bf T}^n$.

        If $f$ is a real or complex-valued Borel measurable function
on ${\bf T}^n$ which is integrable with respect to $|\mu|$, then the
integral of $f$ with respect to $\mu$ can also be defined, and satisfies
\begin{equation}
\label{|int_{{bf T}^n} f d mu| le int_{{bf T}^n} |f| d |mu|}
\biggl|\int_{{\bf T}^n} f \, d\mu\biggr| \le \int_{{\bf T}^n} |f| \, d|\mu|.
\end{equation}
In particular, this applies to any bounded Borel measurable function $f$
on ${\bf T}^n$, in which case we get that
\begin{equation}
\label{|int_{{bf T}^n} f d mu| le (sup_{z in {bf T}^n} |f(z)|) |mu|({bf T}^n)}
        \biggl|\int_{{\bf T}^n} f \, d\mu\biggr|
          \le \Big(\sup_{z \in {\bf T}^n} |f(z)|\Big) \, |\mu|({\bf T}^n).
\end{equation}
Continuous functions on ${\bf T}^n$ are obviously Borel measurable, so
that
\begin{equation}
\label{lambda_mu(f) = int_{{bf T}^n} f d mu}
        \lambda_\mu(f) = \int_{{\bf T}^n} f \, d\mu
\end{equation}
defines a bounded linear functional on $C({\bf T}^n)$, with dual norm
less than or equal to $|\mu|({\bf T}^n)$ with respect to the supremum
norm on $C({\bf T}^n)$.

        Conversely, a version of the Riesz representation theorem
states that every continuous linear functional $\lambda$ on $C({\bf
T}^n)$ can be expressed as (\ref{lambda_mu(f) = int_{{bf T}^n} f d
mu}) for a unique Borel measure $\mu$ on ${\bf T}^n$.  The dual norm
of $\lambda$ with respect to the supremum norm on $C({\bf T}^n)$ is
also equal to $|\mu|({\bf T}^n)$.  Normally one asks that $\mu$ be
Borel regular, which means by definition that $|\mu|$ is Borel
regular, but this is automatic in this case, because open subsets of
${\bf T}^n$ are $\sigma$-compact.  An important advantage of looking
at measures on ${\bf T}^n$ in terms of continuous linear functionals
on $C({\bf T}^n)$ is that we can use the weak$^*$ topology on the dual
of $C({\bf T}^n)$, as in Section \ref{weak^* topology}.

\section[\ Convolution of measures]{Convolution of measures}
\label{convolution of measures}
\setcounter{equation}{0}

        Let $\mu$, $\nu$ be real or complex Borel measures on ${\bf
T}^n$.  Their convolution $\mu * \nu$ may be defined as the Borel
measure on ${\bf T}^n$ given by
\begin{equation}
\label{(mu * nu)(E) = ...}
 (\mu * \nu)(E) = (\mu \times \nu)(\{(z, w) \in {\bf T}^n \times {\bf T}^n :
                                                        z \diamond w \in E\}).
\end{equation}
Here $z \diamond w = (z_1 \, w_1, \ldots, z_n \, w_n)$, as in Section
\ref{convolution on T^n}, and $\mu \times \nu$ is the product measure
on ${\bf T}^n \times {\bf T}^n$ associated to $\mu$, $\nu$.  Note that
\begin{equation}
\label{{(z, w) in {bf T}^n times {bf T}^n : z diamond w in E}}
        \{(z, w) \in {\bf T}^n \times {\bf T}^n : z \diamond w \in E\}
\end{equation}
is a relatively open set in ${\bf T}^n \times {\bf T}^n$ when $E$ is a
relatively open set in ${\bf T}^n$, because $(z, w) \mapsto z \diamond
w$ is continuous as a mapping from ${\bf T}^n \times {\bf T}^n$ into
${\bf T}^n$.  This implies that (\ref{{(z, w) in {bf T}^n times {bf
T}^n : z diamond w in E}}) is a Borel set in ${\bf T}^n \times {\bf
T}^n$ when $E$ is a Borel set in ${\bf T}^n$.  Equivalently, if $f$ is
a bounded Borel measurable function on ${\bf T}^n$, then
\begin{equation}
\label{int_{{bf T}^n} f d(mu * nu) = ...}
        \int_{{\bf T}^n} f \, d(\mu * \nu) =
 \int_{{\bf T}^n \times {\bf T}^n} f(z \diamond w) \, d(\mu \times \nu)(z, w).
\end{equation}
It is easy to see that
\begin{equation}
        \mu * \nu = \nu * \mu,
\end{equation}
and that
\begin{equation}
        (\mu * \nu) * \rho = \mu * (\nu * \rho)
\end{equation}
for any three Borel measures $\mu$, $\nu$, and $\rho$ on ${\bf T}^n$.

        Observe that
\begin{equation}
        |(\mu * \nu)(E)| \le (|\mu| * |\nu|)(E)
\end{equation}
for every Borel set $E \subseteq {\bf T}^n$, and hence
\begin{equation}
        |\mu * \nu|(E) \le (|\mu| * |\nu|)(E).
\end{equation}
This implies that
\begin{equation}
        |\mu * \nu|({\bf T}^n) \le (|\mu| * |\nu|)({\bf T}^n)
                                  = |\mu|({\bf T}^n) \, |\nu|({\bf T}^n).
\end{equation}
Of course, $\|\mu\| = |\mu|({\bf T}^n)$ is a natural norm on the space
of Borel measures on ${\bf T}^n$, also known as the total variation of $\mu$.

        If one looks at measures on ${\bf T}^n$ in terms of continuous
linear functionals on $C({\bf T}^n)$, then convolution can be defined
more directly, basically using (\ref{int_{{bf T}^n} f d(mu * nu) =
...}).  To do this, the product $\lambda_1 \times \lambda_2$ of two
continuous linear functionals $\lambda_1$, $\lambda_2$ on $C({\bf
T}^n)$ should first be defined as a continuous linear functional on
$C({\bf T}^n \times {\bf T}^n)$.  This is not too difficult to do, but
there are some details to be checked.  If $f(z, w)$ is a continuous
function on ${\bf T}^n \times {\bf T}^n$, then one can apply
$\lambda_1$ to $f(z, w)$ as a function of $z$ for each $w \in {\bf
T}^n$, to get a function of $w$ on ${\bf T}^n$.  It is easy to see
that this is a continuous function of $w$, using the fact that $f(z,
w)$ is uniformly continuous on ${\bf T}^n \times {\bf T}^n$, because
${\bf T}^n$ and hence ${\bf T}^n \times {\bf T}^n$ is compact, and
using the continuity of $\lambda_1$ on $C({\bf T}^n)$.  Thus one can
apply $\lambda_2$ to the resulting function of $w$, to get a real or
complex number, as appropriate.  This defines $\lambda_1 \times
\lambda_2$ as a linear functional on $C({\bf T}^n \times {\bf T}^n)$.
By construction,
\begin{equation}
\label{|(lambda_1 times lambda_2)(f)| le ...}
 |(\lambda_1 \times \lambda_2)(f)| \le \|\lambda_1\|_* \, \|\lambda_2\|_* \,
                                \Big(\sup_{z, w \in {\bf T}^n} |f(z, w)|\Big),
\end{equation}
where $\|\lambda_1\|_*$, $\|\lambda_2\|_*$ are the dual norms of
$\lambda_1$, $\lambda_2$ with respect to the supremum norm on $C({\bf
T}^n)$.  This shows that $\lambda_1 \times \lambda_2$ is continuous
with respect to the supremum norm on $C({\bf T}^n \times {\bf T}^n)$,
with the dual norm less than or equal to $\|\lambda_1\|_* \,
\|\lambda_2\|_*$.  If $f(z, w) = f_1(z) \, f_2(w)$ for some continuous
functions $f_1$, $f_2$ on ${\bf T}^n$, then it follows directly from
the definition of $\lambda_1 \times \lambda_2$ that
\begin{equation}
\label{(lambda_1 times lambda_2)(f) = lambda_1(f_1) lambda_2(f_2)}
        (\lambda_1 \times \lambda_2)(f) = \lambda_1(f_1) \, \lambda_2(f_2).
\end{equation}
This implies that the dual norm of $\lambda_1 \times \lambda_2$ on
$C({\bf T}^n \times {\bf T}^n)$ is equal to $\|\lambda_1\|_* \,
\|\lambda_2\|_*$.  Every continuous function on ${\bf T}^n \times {\bf
T}^n$ can be approximated uniformly by a finite sum of products of
continuous functions of $z$ and $w$ on ${\bf T}^n$, and hence
$\lambda_1 \times \lambda_2$ may be characterized as the unique
continuous linear functional on $C({\bf T}^n \times {\bf T}^n)$ that
satisfies (\ref{(lambda_1 times lambda_2)(f) = lambda_1(f_1)
lambda_2(f_2)}) for all $f_1, f_2 \in C({\bf T}^n)$.  In particular,
suppose that $\lambda_1 \times \lambda_2$ was defined instead by first
applying $\lambda_2$ to a continuous function $f(z, w)$ on ${\bf T}^n
\times {\bf T}^n$ as a function of $w$ for each $z \in {\bf T}^n$, and
then applying $\lambda_1$ to the resulting function of $z$.  This
would also determine a continuous linear functional on $C({\bf T}^n
\times {\bf T}^n)$ that satisfies (\ref{(lambda_1 times lambda_2)(f) =
lambda_1(f_1) lambda_2(f_2)}), and which would therefore be equivalent
to the previous definition of $\lambda_1 \times \lambda_2$.

\section[\ Functions and measures]{Functions and measures}
\label{functions, measures}
\setcounter{equation}{0}

        If $g$ is an real or complex-valued function on ${\bf T}^n$
which is integrable with respect to Lebesgue measure, then
\begin{equation}
\label{mu_g(E) = frac{1}{(2 pi)^n} int_E g(z) |dz|}
        \mu_g(E) = \frac{1}{(2 \pi)^n} \int_E g(z) \, |dz|
\end{equation}
defines a Borel measure on ${\bf T}^n$.  As usual,
\begin{equation}
\label{int_{{bf T}^n} f d mu_g = frac{1}{(2 pi)^n} int_{T^n} f(z) g(z) |dz|}
        \int_{{\bf T}^n} f \, d\mu_g
          = \frac{1}{(2 \pi)^n} \int_{{\bf T}^n} f(z) \, g(z) \, |dz|
\end{equation}
for every bounded measurable function $f$ on ${\bf T}^n$.  It is also
well known that $|\mu_g| = \mu_{|g|}$, and hence
\begin{equation}
\label{||mu_g|| = |mu_g|({bf T}^n) = frac{1}{(2 pi)^n} int_{T^n} |g(z)| |dz|}
        \|\mu_g\| = |\mu_g|({\bf T}^n)
                  = \frac{1}{(2 \pi)^n} \int_{{\bf T}^n} |g(z)| \, |dz|.
\end{equation}
If $h$ is another Lebesgue integrable function on ${\bf T}^n$, then
the convolution $g * h$ is also defined as a Lebesgue integrable
function on ${\bf T}^n$, as in Section \ref{convolution on T^n}.  It
is not difficult to check that this is compatible with the definition
of convolution of measures in the previous section, in the sense that
\begin{equation}
\label{mu_g * mu_h = mu_{g * h}}
        \mu_g * \mu_h = \mu_{g * h}.
\end{equation}

       If $\nu$ is a real or complex Borel measure on ${\bf T}^n$,
then the convolution of $\mu_g$ and $\nu$ can be defined as a measure
on ${\bf T}^n$ as in the previous section.  Alternatively, $g * \nu$
can be defined as a Lebesgue integrable function on ${\bf T}^n$ by
\begin{equation}
\label{(g * nu)(z) = int_{{bf T}^n} g(z diamond w^{-1}) d nu(w)}
        (g * \nu)(z) = \int_{{\bf T}^n} g(z \diamond w^{-1}) \, d\nu(w),
\end{equation}
where $w^{-1} = (w_1^{-1}, \ldots, w_n^{-1})$, as before.  The
existence of this integral for almost every $z \in {\bf T}^n$ with
respect to Lebesgue measure uses Fubini's theorem, as in Section
\ref{convolution on T^n}.  More precisely, if $g$ and $\nu$ are
nonnegative and real-valued, then Fubini's theorem implies that
\begin{equation}
\label{frac{1}{(2 pi)^n} int_{{bf T}^n} (g * nu)(z) |dz| = ...}
 \frac{1}{(2 \pi)^n} \int_{{\bf T}^n} (g * \nu)(z) \, |dz| = 
 \Big(\frac{1}{(2 \pi)^n} \int_{{\bf T}^n} g(z) \, |dz|\Big) \, \nu({\bf T}^n).
\end{equation}
In particular, $(g * \nu) (z) < \infty$ for almost every $z \in {\bf
T}^n$ with respect to Lebesgue measure.  Otherwise, if $g$ and $\nu$
are real or complex-valued, then one can apply this to $|g|$ and
$|\nu|$.  This implies that the integral in (\ref{(g * nu)(z) =
int_{{bf T}^n} g(z diamond w^{-1}) d nu(w)}) makes sense for almost
every $z \in {\bf T}^n$ with respect to Lebesgue measure, and that
\begin{equation}
\label{frac{1}{(2 pi)^n} int_{{bf T}^n} |(g * nu)(z)| |dz| le ...}
        \frac{1}{(2 \pi)^n} \int_{{\bf T}^n} |(g * \nu)(z)| \, |dz|
           \le \Big(\frac{1}{(2 \pi)^n} \int_{{\bf T}^n} |g(z)| \, |dz|\Big)
                                                          \, |\nu|({\bf T}^n).
\end{equation}
Of course, if $\nu = \mu_h$ for some Lebesgue integrable function $h$
on ${\bf T}^n$, then this definition of $g * \nu$ reduces to the
earlier definition of $g * h$.  Similarly, if $\nu$ is any Borel
measure on ${\bf T}^n$, then this definition of $g * \nu$ is
compatible with the definition of convolution of measures in the
previous section, in the sense that
\begin{equation}
        \mu_g * \nu = \mu_{g * \nu}.
\end{equation}

        If $g$ is continuous on ${\bf T}^n$, and hence uniformly
continuous, then it is easy to see that $g * \nu$ also defines a
continuous function on ${\bf T}^n$.  In this case, we also have that
\begin{equation}
        \sup_{z \in {\bf T}^n} |(g * \nu)(z)|
         \le \Big(\sup_{z \in {\bf T}^n} |g(z)|\Big) \, |\nu|({\bf T}^n).
\end{equation}

\section[\ Fourier coefficients]{Fourier coefficients}
\label{fourier coefficients}
\setcounter{equation}{0}

        Let $\mu$ be a complex Borel measure on ${\bf T}^n$.  If
$\alpha \in {\bf Z}^n$, then the corresponding Fourier coefficient of
$\mu$ is defined by
\begin{equation}
\label{widehat{mu}(alpha) = int_{{bf T}^n} z^{-alpha} d mu(z)}
        \widehat{\mu}(\alpha) = \int_{{\bf T}^n} z^{-\alpha} \, d\mu(z).
\end{equation}
This reduces to the earlier definition of the Fourier coefficients of
a Lebesgue integrable function $g$ on ${\bf T}^n$ when $\mu = \mu_g$,
as in the previous section.  The Fourier coefficients of any complex
Borel measure $\mu$ on ${\bf T}^n$ are bounded, with
\begin{equation}
\label{|widehat{mu}(alpha)| le |mu|({bf T}^n)}
        |\widehat{\mu}(\alpha)| \le |\mu|({\bf T}^n)
\end{equation}
for each $\alpha \in {\bf Z}^n$.  If $\nu$ is another complex Borel
measure on ${\bf T}^n$, then it is easy to see that
\begin{equation}
\label{widehat{(mu * nu)}(alpha) = widehat{mu}(alpha) widehat{nu}(alpha)}
 \widehat{(\mu * \nu)}(\alpha) = \widehat{\mu}(\alpha) \, \widehat{\nu}(\alpha)
\end{equation}
for every $\alpha \in {\bf Z}^n$.

        Let $U^n$ be the open unit polydisk in ${\bf C}^n$, and let
$\widetilde{z}^\alpha$ be defined for $\alpha \in {\bf Z}^n$ and $z
\in {\bf C}^n$ as in Section \ref{multiple fourier series}.  If $\mu$
is a complex Borel measure on ${\bf T}^n$ and $z \in U^n$, then put
\begin{equation}
\label{phi_mu(z) = sum_{alpha in Z^n} widehat{mu}(alpha) widetilde{z}^alpha}
        \phi_\mu(z) = \sum_{\alpha \in {\bf Z}^n}
                        \widehat{\mu}(\alpha) \, \widetilde{z}^\alpha.
\end{equation}
As in Section \ref{multiple fourier series}, the sum converges
absolutely for every $z \in U^n$, because of the boundedness of the
Fourier coefficients of $\mu$.  This can also be expressed as
\begin{equation}
\label{phi_mu(z) = (2 pi)^n int_{{bf T}^n} P_n(z, w) d mu(w)}
        \phi_\mu(z) = (2 \pi)^n \int_{{\bf T}^n} P_n(z, w) \, d\mu(w),
\end{equation}
where $P_n(z, w)$ is the $n$-dimensional Poisson kernel, discussed
in Section \ref{multiple fourier series}.  More precisely,
\begin{equation}
\label{(2 pi)^n P_n(z, w) = sum_{alpha in Z^n} widetilde{z}^alpha w^{-alpha}}
 (2 \pi)^n \, P_n(z, w) = \sum_{\alpha \in {\bf Z}^n} \widetilde{z}^\alpha
                                                             \, w^{-\alpha}
\end{equation}
for each $z \in U^n$ and $w \in {\bf T}^n$.  This sum can be
approximated uniformly by finite subsums as a function of $w \in {\bf
T}^n$ for each $z \in U^n$, which permits one to interchange the order
of summation and integration in (\ref{phi_mu(z) = sum_{alpha in Z^n}
widehat{mu}(alpha) widetilde{z}^alpha}) to get (\ref{phi_mu(z) = (2
pi)^n int_{{bf T}^n} P_n(z, w) d mu(w)}).  Of course, the extra factor
of $(2 \pi)^n$ here simply comes from slightly different
normalizations being used.

        If $r \in [0, 1)^n$ and $z, w \in {\bf T}^n$, then $r \diamond
z \in U^n$, and
\begin{equation}
 (2 \pi)^n \, P_n(r \diamond z, w) = (2 \pi)^n \, P_n(r, w \diamond z^{-1}).
\end{equation}
Put
\begin{equation}
        \rho_{n, r}(w) = (2 \pi)^n \, P_n(r, w)
\end{equation}
for each $r \in [0, 1)^n$ and $w \in {\bf T}^n$.  It is easy to see that
\begin{equation}
\label{rho_{n, r}(w^{-1}) = rho_{n, r}(w)}
        \rho_{n, r}(w^{-1}) = \rho_{n, r}(w),
\end{equation}
using the change or variables $\alpha \mapsto - \alpha$ in (\ref{(2
pi)^n P_n(z, w) = sum_{alpha in Z^n} widetilde{z}^alpha w^{-alpha}}).
It follows that
\begin{equation}
\label{phi_mu(r diamond z) = (rho_{n, r} * mu)(z)}
        \phi_\mu(r \diamond z) = (\rho_{n, r} * \mu)(z)
\end{equation}
for every $r \in [0, 1)^n$ and $z \in {\bf T}^n$, by (\ref{phi_mu(z) =
(2 pi)^n int_{{bf T}^n} P_n(z, w) d mu(w)}).

        Note that $\rho_{n, r}(w) \ge 0$ and
\begin{equation}
\label{frac{1}{(2 pi)^n} int_{{bf T}^n} rho_{n, r}(w) |dw| = 1}
        \frac{1}{(2 \pi)^n} \int_{{\bf T}^n} \rho_{n, r}(w) \, |dw| = 1
\end{equation}
for each $r \in [0, 1)$, by the corresponding properties of the Poisson
kernel.  If $\mu = \mu_g$ for some continuous function $g$ on ${\bf T}^n$, then
\begin{equation}
\label{phi_mu(r diamond z) to g(z)}
        \phi_\mu(r \diamond z) \to g(z)
\end{equation}
as $r \to (1, \ldots, 1)$ for each $z \in {\bf T}^n$, as in previous
discussions of Poisson integrals.  As usual, the convergence is also
uniform over $z \in {\bf T}^n$, because $g$ is uniformly continuous on
${\bf T}^n$.  If $\mu = \mu_g$ for a Lebesgue integrable function $g$
on ${\bf T}^n$, then one can show that there is convergence in the
$L^1$ norm on ${\bf T}^n$.  More precisely, this follows by
approximating $g$ by continuous functions on ${\bf T}^n$ in the $L^1$
norm, and using uniform bounds for the $L^1$ norm of $\phi_\mu(r
\diamond z)$ as a function of $z \in {\bf T}^n$ over $r \in [0, 1)^n$.
If $\mu$ is any complex Borel measure on ${\bf T}^n$ and $f$ is a
continuous function on ${\bf T}^n$, then
\begin{equation}
 \frac{1}{(2 \pi)^n} \int_{{\bf T}^n} \phi_\mu(r \diamond z) \, f(z) \, |dz|
          = \int_{{\bf T}^n} \rho_{n, r} * f \, d\mu,
\end{equation}
by Fubini's theorem and (\ref{rho_{n, r}(w^{-1}) = rho_{n, r}(w)}).
Hence
\begin{equation}
 \frac{1}{(2 \pi)^n} \int_{{\bf T}^n} \phi_\mu(r \diamond z) \, f(z) \, |dz| 
          \to \int_{{\bf T}^n} f \, d\mu
\end{equation}
as $r \to (1, \ldots, 1)$, because $\rho_{n, r} * f \to f$ uniformly
on ${\bf T}^n$ as $r \to (1, \ldots, 1)$, as before.  This says that the
measure on ${\bf T}^n$ associated to $\phi_\mu(r \diamond z)$ as in the
preceding section converges to $\mu$ in the weak$^*$ topology on the dual of
$C({\bf T}^n)$ as $r \to (1, \ldots, 1)$, when we identify Borel measures
on ${\bf T}^n$ with continuous linear functionals on $C({\bf T}^n)$.

\section[\ Measures on ${\bf R}^n$]{Measures on ${\bf R}^n$}
\label{measures on R^n}
\setcounter{equation}{0}

        Let $\mu$ be a real or complex Borel measure on ${\bf R}^n$,
which is to say a countably-additive real or complex valued function
on the $\sigma$-algebra of Borel sets in ${\bf R}^n$.  As before,
there is a finite nonnegative Borel measure $|\mu|$ on ${\bf R}^n$
associated to $\mu$ such that
\begin{equation}
        |\mu(E)| \le |\mu|(E)
\end{equation}
for every Borel set $E \subseteq {\bf R}^n$, and which is less than or
equal to every other nonnegative Borel measure on ${\bf R}^n$ with
this property.  If $f$ is a Borel measurable function on ${\bf R}^n$
which is integrable with respect to $|\mu|$, then the integral of $f$
with respect to $\mu$ can also be defined, and satisfies
\begin{equation}
\label{|int_{{bf R}^n} f d mu| le int_{{bf R}^n} |f| d |mu|}
 \biggl|\int_{{\bf R}^n} f \, d\mu\biggr| \le \int_{{\bf R}^n} |f| \, d|\mu|.
\end{equation}
In particular, this works when $f$ is a bounded Borel measurable
function on ${\bf R}^n$, for which we have that
\begin{equation}
\label{|int_{{bf R}^n} f d mu| le (sup_{x in {bf R}^n} |f(x)|) |mu|({bf R}^n)}
 \biggl|\int_{{\bf R}^n} f \, d\mu\biggr|
        \le \Big(\sup_{x \in {\bf R}^n} |f(x)|\Big) \, |\mu|({\bf R}^n).
\end{equation}

        Of course, continuous functions on ${\bf R}^n$ are Borel
measurable, and so
\begin{equation}
\label{lambda_mu(f) = int_{{bf R}^n} f d mu}
        \lambda_\mu(f) = \int_{{\bf R}^n} f \, d\mu
\end{equation}
defines a bounded linear functional on the space $C_b({\bf R}^n)$ of
bounded continuous functions on ${\bf R}^n$ with respect to the
supremum norm, with dual norm less than or equal to $|\mu|({\bf R}^n)$.
The restriction of $\lambda_\mu$ to the space $C_0({\bf R}^n)$ of continuous
functions on ${\bf R}^n$ that vanish at infinity is also bounded with
respect to the supremum norm, with dual norm less than or equal to
$|\mu|({\bf R}^n)$.  Conversely, a version of the Riesz representation theorem
states that every bounded linear functional $\lambda$ on $C_0({\bf
R}^n)$ corresponds to a unique Borel measure $\mu$ in this way, where
the dual norm of $\lambda$ with respect to the supremum norm on
$C_0({\bf R}^n)$ is equal to $\|\mu\| = |\mu|({\bf R}^n)$.  Normally
one also asks $\mu$ to satisfy some additional regularity conditions,
but these hold automatically on ${\bf R}^n$, since open sets in ${\bf R}^n$
are $\sigma$-compact.

        A small part of this theorem implies that a bounded linear
functional $\lambda$ on $C_0({\bf R}^n)$ has a natural extension to
$C_b({\bf R}^n)$.  This extension is characterized by the following
additional continuity condition, which is a mild version of the
dominated convergence theorem.  Namely, if $\{f_j\}_{j = 1}^\infty$ is
a sequence of bounded continuous functions on ${\bf R}^n$ that are
uniformly bounded on ${\bf R}^n$ and converge uniformly on compact
subsets to a function $f$ on ${\bf R}^n$, then $\{\lambda(f_j)\}_{j =
1}^\infty$ converges to $\lambda(f)$.  Note that $f$ is bounded and
continuous under these conditions, and that any bounded continuous
function on ${\bf R}^n$ is the limit of a uniformly bounded sequence
of continuous functions with compact support on ${\bf R}^n$ that
converges uniformly on compact subsets of ${\bf R}^n$.  Hence the
extension of $\lambda$ to $C_b({\bf R}^n)$ is uniquely determined by
$\lambda$ on $C_0({\bf R}^n)$ when the extension satisfies this
additional continuity condition.

        If $\lambda$ is a bounded linear functional on $C_0({\bf
R}^n)$ with compact support, so that $\lambda(f)$ only depends on the
restriction of $f$ to a compact set in ${\bf R}^n$, then this
extension of $\lambda$ to $C_b({\bf R}^n)$ is basically trivial.
Otherwise, it is not too difficult to show that a bounded linear
functional $\lambda$ on $C_0({\bf R}^n)$ can be approximated by
bounded linear functionals on $C_0({\bf R}^n)$ with compact support
with respect to the dual norm.  One can then use this approximation to
show more directly that $\lambda$ can be extended to a bounded linear
functional on $C_b({\bf R}^n)$ that satisfies the additional
continuity condition mentioned in the previous paragraph.  Note that
the dual norm of the extension of $\lambda$ to $C_b({\bf R}^n)$ with
respect to the supremum norm on $C_b({\bf R}^n)$ is equal to the dual
norm of $\lambda$ on $C_0({\bf R}^n)$.

\section[\ Convolution of measures, continued]{Convolution of measures, continued}
\label{convolution of measures, continued}
\setcounter{equation}{0}

        If $\mu$, $\nu$ are real or complex Borel measures on ${\bf
R}^n$, then their convolution $\mu * \nu$ may be defined as a Borel
measure on ${\bf R}^n$ by
\begin{equation}
\label{(mu * nu)(E) = (mu times nu)({(x, y) in R^n times R^n : x + y in E}}
 (\mu * \nu)(E) = (\mu \times \nu)(\{(x, y) \in {\bf R}^n \times {\bf R}^n :
                                                                 x + y \in E\},
\end{equation}
where $\mu \times \nu$ is the product measure on ${\bf R}^n \times
{\bf R}^n$ corresponding to $\mu$, $\nu$.  Note that
\begin{equation}
\label{{(x, y) in {bf R}^n times {bf R}^n : x + y in E}}
        \{(x, y) \in {\bf R}^n \times {\bf R}^n : x + y \in E\}
\end{equation}
is an open set in ${\bf R}^n \times {\bf R}^n$ for every open set $E
\subseteq {\bf R}^n$, by continuity of addition, which implies that
(\ref{{(x, y) in {bf R}^n times {bf R}^n : x + y in E}}) is a Borel set
in ${\bf R}^n \times {\bf R}^n$ when $E$ is a Borel set in ${\bf R}^n$.
If $f$ is a bounded Borel measurable function on ${\bf R}^n$, then
we get that
\label{int_{{bf R}^n} f d(mu * nu) = ...}
\begin{equation}
 \int_{{\bf R}^n} f \, d(\mu * \nu)
   = \int_{{\bf R}^n \times {\bf R}^n} f(x + y) \, d (\mu \times \nu)(x, y).
\end{equation}
As usual,
\begin{equation}
\label{nu * mu = mu * nu and (mu * nu) * rho = mu * (nu * rho)}
        \nu * \mu = \mu * \nu \quad\hbox{and}\quad
           (\mu * \nu) * \rho = \mu * (\nu * \rho)
\end{equation}
for any Borel measures $\mu$, $\nu$ and $\rho$ on ${\bf R}^n$.

        As before,
\begin{equation}
        |(\mu * \nu)(E)| \le (|\mu| * |\nu|)(E)
\end{equation}
for any Borel set $E$ in ${\bf R}^n$.  This implies that
\begin{equation}
        |\mu * \nu|(E) \le (|\mu| * |\nu|)(E)
\end{equation}
for every Borel set $E \subseteq {\bf R}^n$.  In particular,
\begin{equation}
        \|\mu * \nu\| \le \|\mu\| \, \|\nu\|,
\end{equation}
where $\|\mu\| = |\mu|({\bf R}^n)$, as in the previous section.

        One can also look at convolution in terms of bounded linear
functionals on spaces of continuous functions, as in Section
\ref{convolution of measures}.  If $\lambda_1$, $\lambda_2$ are
bounded linear functionals on $C_0({\bf R}^n)$, then the product
linear functional $\lambda_1 \times \lambda_2$ can be defined as a
bounded linear functional on $C_0({\bf R}^n \times {\bf R}^n)$, in
basically the same way as before.  In order to define the convolution
$\lambda_1 * \lambda_2$ as a bounded linear functional on $C_0({\bf
R}^n)$, one would like to let $\lambda_1 * \lambda_2$ act on $f(x +
y)$ as a continuous function on ${\bf R}^n \times {\bf R}^n$, where
$f$ is a continuous function on ${\bf R}^n$ tht vanishes at infinity.
However, $f(x + y)$ does not vanish at infinity on ${\bf R}^n \times
{\bf R}^n$ unless $f \equiv 0$, and so it is better to use the natural
extension of $\lambda_1 \times \lambda_2$ to bounded continuous
functions on ${\bf R}^n \times {\bf R}^n$, as in the preceding
section.

\section[\ Functions and measures, continued]{Functions and measures, continued}
\label{functions, measures, continued}
\setcounter{equation}{0}

        If $g$ is a real or complex-valued function on ${\bf R}^n$
that is integrable with respect to Lebesgue measure, then
\begin{equation}
\label{mu_g(E) = int_E g(x) dx}
        \mu_g(E) = \int_E g(x) \, dx
\end{equation}
defines a Borel measure on ${\bf R}^n$.  As before,
\begin{equation}
        \int_{{\bf R}^n} f \, d\mu_g = \int_{{\bf R}^n} f(x) \, g(x) \, dx
\end{equation}
for every bounded measurable function $f$ on ${\bf R}^n$.  Also,
$|\mu_g| = \mu_{|g|}$, so that $\|\mu_g\|$ is the same as the $L^1$
norm of $g$ on ${\bf R}^n$.  If $h$ is another integrable function on
${\bf R}^n$, then one can check that
\begin{equation}
        \mu_g * \mu_h = \mu_{g * h},
\end{equation}
where $g * h$ is the integrable function on ${\bf R}^n$ defined as in
Section \ref{convolution on R^n}.

        If $\nu$ is a real or complex Borel measure on ${\bf R}^n$,
then the convolution of $g$ and $\nu$ can be defined as a Lebesgue
integrable function on ${\bf R}^n$ by
\begin{equation}
\label{(g * nu)(x) = int_{{bf R}^n} g(x - y) d nu(y)}
        (g * \nu)(x) = \int_{{\bf R}^n} g(x - y) \, d\nu(y).
\end{equation}
As usual, the existence of this integral almost everywhere on ${\bf
R}^n$ uses Fubini's theorem.  If $g$ and $\nu$ are nonnegative and
real-valued, then
\begin{equation}
\label{int_{{bf R}^n} (g * nu)(x) dx = (int_{{bf R}^n} g(x) dx) nu({bf R}^n)}
 \int_{{\bf R}^n} (g * \nu)(x) \, dx = \Big(\int_{{\bf R}^n} g(x) \, dx\Big)
                                        \, \nu({\bf R}^n),
\end{equation}
and in particular $(g * \nu)(x) < \infty$ for almost every $x \in {\bf
R}^n$ with respect to Lebesgue measure.  Otherwise, if $g$ and $\nu$
are real or complex-valued, then one can apply this to $|g|$ and
$|\nu|$, to get that the integral in (\ref{(g * nu)(x) = int_{{bf
R}^n} g(x - y) d nu(y)}) makes sense for almost every $x \in {\bf
R}^n$ with respect to Lebesgue measure, and that
\begin{equation}
\label{int_{R^n} |(g * nu)(x)| dx le (int_{R^n} |g(x)| dx) |nu|(R^n)}
        \int_{{\bf R}^n} |(g * \nu)(x)| \, dx
            \le \Big(\int_{{\bf R}^n} |g(x)| \, dx\Big) \, |\nu|({\bf R}^n).
\end{equation}
If $\nu = \mu_h$ for some integrable function $h$ on ${\bf R}^n$, then
$g * \nu$ reduces to the usual definition of $g * h$, while if $\nu$
is any Borel measure on ${\bf R}^n$, then $\mu_g * \nu = \mu_{g * \nu}$.

        If $\nu$ is a real or complex Borel measure on ${\bf R}^n$ and
$g$ is a bounded continuous function on ${\bf R}^n$, then $(g *
\nu)(x)$ is defined for every $x \in {\bf R}^n$, and satisfies
\begin{equation}
\label{sup_{x in R^n} |(g * nu)(x)| le (sup_{x in R^n} |g(x)|) |nu|(R^n)}
        \sup_{x \in {\bf R}^n} |(g * \nu)(x)|
         \le \Big(\sup_{x \in {\bf R}^n} |g(x)|\Big) \, |\nu|({\bf R}^n),
\end{equation}
as before.  One can also check that $g * \nu$ is continuous on ${\bf
R}^n$, using the dominated convergence theorem.  If $g$ is bounded and
uniformly continuous, then it is easy to see that $g * \nu$ is
uniformly continuous too.  Alternatively, to show that $g * \nu$ is
continuous when $g$ is bounded and continuous, one can use the fact
that $g$ is uniformly continuous on compact sets, and approximate
$\nu$ by measures with compact support.

        If $g$ and $\nu$ have compact support in ${\bf R}^n$, then it
is easy to see that $g * \nu$ has compact support as well.  If $g$ is
a continuous function on ${\bf R}^n$ that vanishes at infinity and
$\nu$ has compact support, then it is easy to check that $g * \nu$
vanishes at infinity on ${\bf R}^n$ too.  This also works when $\nu$
does not have compact support, by approximating $\nu$ by measures with
compact support on ${\bf R}^n$.

\section[\ The Fourier transform, continued]{The Fourier transform, continued}
\label{fourier transform, continued}
\setcounter{equation}{0}

        The Fourier transform of a complex Borel measure $\mu$ on
${\bf R}^n$ can be defined by
\begin{equation}
\label{widehat{mu}(xi) = int_{{bf R}^n} exp (- i xi cdot x) d mu(x)}
        \widehat{\mu}(\xi) = \int_{{\bf R}^n} \exp (- i \xi \cdot x) \, d\mu(x)
\end{equation}
for each $\xi \in {\bf R}^n$.  This coincides with the earlier
definition for an integrable function $f$ on ${\bf R}^n$ when $\mu =
\mu_f$.  As before, 
\begin{equation}
\label{|widehat{mu}(xi)| le |mu|({bf R}^n)}
        |\widehat{\mu}(\xi)| \le |\mu|({\bf R}^n)
\end{equation}
for every $\xi \in {\bf R}^n$, and one can also check that
$\widehat{\mu}(\xi)$ is uniformly continuous on ${\bf R}^n$.  This is
easier to do when $\mu$ has compact support in ${\bf R}^n$, and
otherwise one can approximate $\mu$ by measures with compact support.
If $\nu$ is another complex Borel measure on ${\bf R}^n$, then it is
easy to see that
\begin{equation}
\label{widehat{(mu * nu)}(xi) = widehat{mu}(xi) widehat{nu}(xi)}
        \widehat{(\mu * \nu)}(\xi) = \widehat{\mu}(\xi) \, \widehat{\nu}(\xi)
\end{equation}
for every $\xi \in {\bf R}^n$.

        The analogue of the multiplication formula in this context states that
\begin{equation}
\label{int_{R^n} widehat{mu}(xi) d nu(xi) = int_{R^n} widehat{nu}(x) d mu(x)}
        \int_{{\bf R}^n} \widehat{\mu}(\xi) \, d\nu(\xi)
          = \int_{{\bf R}^n} \widehat{\nu}(x) \, d\mu(x)
\end{equation}
for any pair of complex Borel measures $\mu$, $\nu$ on ${\bf R}^n$.
This follows from Fubini's theorem, as before.  In particular,
\begin{equation}
\label{int_{R^n} widehat{mu}(xi) g(xi) d xi = int_{R^n} widehat{g}(x) d mu(x)}
        \int_{{\bf R}^n} \widehat{\mu}(\xi) \, g(\xi) \, d\xi
          = \int_{{\bf R}^n} \widehat{g}(x) \, d\mu(x)
\end{equation}
for every Lebesgue integrable function $g$ on ${\bf R}^n$.  As in
Section \ref{multiplication formula}, this implies that
\begin{equation}
\label{int_{{bf R}^n} widehat{mu}(xi) exp (i xi cdot w) h(xi) d xi = ...}
 \int_{{\bf R}^n} \widehat{\mu}(\xi) \, \exp (i \xi \cdot w) \, h(\xi) \, d\xi
                         = \int_{{\bf R}^n} \widehat{h}(x - w) \, d\mu(w)
\end{equation}
for every Lebesgue integrable function $h$ on ${\bf R}^n$ and every $w
\in {\bf R}^n$.  If $h$ is a even function on ${\bf R}^n$, then this
reduces to
\begin{equation}
\label{int_{{bf R}^n} widehat{mu}(xi) exp (i xi cdot w) h(xi) d xi = ..., 2}
 \int_{{\bf R}^n} \widehat{\mu}(\xi) \, \exp (i \xi \cdot w) \, h(\xi) \, d\xi
                         = (\widehat{h} * \mu)(w).
\end{equation}

        Let $a = (a_1, \ldots, a_n)$ be an $n$-tuple of positive real
numbers, and let $P_{n, a}(x)$ be the function on ${\bf R}^n$
discussed in Section \ref{convergence}.  Also let $f$ be a continuous
function on ${\bf R}^n$ that vanishes at infinity, and observe that
\begin{equation}
\label{int_{R^n} (P_{n, a} * mu)(w) f(w) dw = int_{R^n} P_{n, a} * f d mu}
        \int_{{\bf R}^n} (P_{n, a} * \mu)(w) \, f(w) \, dw
          = \int_{{\bf R}^n} P_{n, a} * f \, d\mu
\end{equation}
by Fubini's theorem, using also the fact that $P_{n, a}$ is an
even function on ${\bf R}^n$.  It follows that
\begin{equation}
        \lim_{a \to 0} \int_{{\bf R}^n} (P_{n, a} * \mu)(w) \, f(w) \, dw
                                              = \int_{{\bf R}^n} f \, d\mu,
\end{equation}
because $P_{n, a} * f \to f$ uniformly on ${\bf R}^n$ as $a \to 0$, as
in Section \ref{convergence}.  This says that the measure on ${\bf
R}^n$ associated to $P_{n, a} * \mu$ converges to $\mu$ as $a \to 0$
with respect to the weak$^*$ topology on the dual of $C_0({\bf R}^n)$
when we identify complex Borel measures on ${\bf R}^n$ with bounded
linear functionals on $C_0({\bf R}^n)$.

        As in Section \ref{inversion}, we have that
\begin{equation}
\label{inversion formula, 2}
        \quad  \int_{{\bf R}^n} \widehat{\mu}(\xi) \, \exp (i \xi \cdot w) \, 
                        \exp \Big(- \sum_{j = 1}^n a_j \, |\xi_j|\Big) \, d\xi 
                                         = (2 \pi)^n \, (P_{n, a} * \mu)(w)
\end{equation}
for each $w \in {\bf R}^n$, by applying (\ref{int_{{bf R}^n}
widehat{mu}(xi) exp (i xi cdot w) h(xi) d xi = ..., 2}) with $h =
p_{n, a}$ as in Section \ref{some examples, continued}.  This
converges to $(2 \pi)^n \, \mu$ as $a \to 0$ with respect to the
weak$^*$ topology on the dual of $C_0({\bf R}^n)$, as in the previous
paragraph.  In particular, $\mu = 0$ when $\widehat{\mu} = 0$.

\section[\ Holomorphic extensions, continued]{Holomorphic extensions, continued}
\label{holomorphic extensions, continued}
\setcounter{equation}{0}

        Let us say that a complex Borel measure $\mu$ on ${\bf R}^n$
has support contained in a closed set $E \subseteq {\bf R}^n$ if
\begin{equation}
        \mu({\bf R}^n \backslash E) = 0.
\end{equation}
If $\mu$ has support contained in a compact set $K$ in ${\bf R}^n$,
then the Fourier transform $\widehat{\mu}(\xi)$ extends to a
holomorphic function $\widehat{\mu}(\zeta)$ on ${\bf C}^n$, given by
\begin{equation}
        \widehat{\mu}(\zeta) = \int_K \exp (- i \zeta \cdot x) \, d\mu(x),
\end{equation}
as in Section \ref{holomorphic extensions}.  If $\mu$, $\nu$ are
compactly supported complex Borel measures on ${\bf R}^n$, then one
can check that $\mu * \nu$ also has compact support, and that
\begin{equation}
 \widehat{(\mu * \nu)}(\zeta) = \widehat{\mu}(\zeta) \, \widehat{\nu}(\zeta)
\end{equation}
for each $\zeta \in {\bf C}^n$.  Similarly, let $\epsilon \in \{-1,
1\}^n$ be given, let $Q_{n, \epsilon}$ be the closed ``quadrant'' in
${\bf R}^n$ associated to $\epsilon$ as before, and let $H_{n,
\epsilon}$ be the corresponding region in ${\bf C}^n$.  If $n = 1$,
then $Q_{n, \epsilon}$ is a closed half-line in ${\bf R}$, and $H_{n,
\epsilon}$ is the open upper or lower half-plane in ${\bf C}$, as
appropriate.  If $\mu$ is a complex Borel measure on ${\bf R}^n$ with
support contained in $Q_{n, \epsilon}$, then the Fourier transform of
$\mu$ extends naturally to a bounded uniformly continuous function on
$\overline{H}_{n, - \epsilon}$ that is holomorphic on $H_{n, - \epsilon}$,
for basically the same reasons as for integrable functions.
If $\mu$, $\nu$ are complex Borel measures on ${\bf R}^n$ supported on
$Q_{n, \epsilon}$, then one can check that $\mu * \nu$ is also supported on
$Q_{n, \epsilon}$, and that the natural extension of the Fourier transform
of $\mu * \nu$ to $\overline{H}_{n, - \epsilon}$ is equal to the product of
the corresponding extensions of the Fourier transforms of $\mu$, $\nu$.

\section[\ Approximation and support]{Approximation and support}
\label{approximation, support}
\setcounter{equation}{0}

        Let $X$ be a locally compact Hausdorff topological space, and
let $\lambda$ be a bounded linear functional on the space $C_0(X)$ of
continuous functions on $X$ that vanish at infinity, with respect to
the supremum norm.  If $\phi$ is a bounded continuous function on $X$,
then
\begin{equation}
\label{lambda_phi(f) = lambda(phi f)}
        \lambda_\phi(f) = \lambda(\phi \, f)
\end{equation}
is also a bounded linear functional on $C_0(X)$, and
\begin{equation}
\label{||lambda_phi||_* le ||phi||_{sup} ||lambda||_*}
        \|\lambda_\phi\|_* \le \|\phi\|_{sup} \, \|\lambda\|_*.
\end{equation}
Here $\|\phi\|_{sup}$ denotes the supremum norm of $\phi$ on $X$, and
$\|\lambda\|_*$ is the dual norm of $\lambda$ with respect to the
supremum norm on $C_0(X)$.  Note that $\phi \, f \in C_0(X)$ when $f
\in C_0(X)$ and $\phi \in C_b(X)$.

        Let $\psi$ be another bounded continuous function on $X$, and
let us check that
\begin{equation}
\label{||lambda_phi||_* + ||lambda_psi||_* le ...}
        \|\lambda_\phi\|_* + \|\lambda_\psi\|_*
         \le \sup_{x \in X} (|\phi(x)| + |\psi(x)|) \, \|\lambda\|_*.
\end{equation}
Let $a$, $b$ be real or complex numbers, as appropriate, and let $f$,
$g$ be continuous functions on $X$ that vanish at infinity, with
\begin{equation}
        |a|, |b|, \|f\|_{sup}, \|g\|_{sup} \le 1.
\end{equation}
Observe that
\begin{equation}
        a \, \lambda_\phi(f) + b \, \lambda_\psi(g)
            = \lambda(a \, \phi \, f + b \, \psi \, g),
\end{equation}
and hence
\begin{equation}
        |a \, \lambda_\phi(f) + b \, \lambda_\psi(g)|
         \le \|a \, \phi \, f + b \, \psi \, g\|_{sup} \, \|\lambda\|_*.
\end{equation}
Our hypotheses on $a$, $b$, $f$, and $g$ imply that
\begin{equation}
        \|a \, \phi \, f + b \, \psi \, g\|_{sup}
               \le \sup_{x \in X} (|\phi(x)| + |\psi(x)|),
\end{equation}
so that
\begin{equation}
        |a \, \lambda_\phi(f) + b \, \lambda_\psi(g)|
              \le \sup_{x \in X} (|\phi(x)| + |\psi(x)|) \, \|\lambda\|_*.
\end{equation}
Using suitable choices of $a$ and $b$, we get that
\begin{equation}
        |\lambda_\phi(f)| + |\lambda_\psi(g)|
            \le \sup_{x \in X} (|\phi(x)| + |\psi(x)|) \, \|\lambda\|_*,
\end{equation}
which implies (\ref{||lambda_phi||_* + ||lambda_psi||_* le ...}), by
taking the supremum over $f$ and $g$.

        Suppose now that $\phi$ is a bounded real-valued continuous
function on $X$ such that $0 \le \phi(x) \le 1$ for each $x \in X$.
If we take $\psi = 1 - \phi$ in (\ref{||lambda_phi||_* +
||lambda_psi||_* le ...}), then we get that
\begin{equation}
        \|\lambda_\phi\|_* + \|\lambda_{1 - \phi}\|_* \le \|\lambda\|_*.
\end{equation}
Of course,
\begin{equation}
        \|\lambda\|_* \le \|\lambda_\phi\|_* + \|\lambda_{1 - \phi}\|_*,
\end{equation}
because $\lambda_\phi + \lambda_{1 - \phi} = \lambda$, and so
\begin{equation}
\label{||lambda_phi||_* + ||lambda_{1 - phi}||_* = ||lambda||_*}
        \|\lambda_\phi\|_* + \|\lambda_{1 - \phi}\|_* = \|\lambda\|_*.
\end{equation}

        Let $\epsilon > 0$ be given, and let $f$ be a continuous function
on $X$ that vanishes at infinity such that $\|f\|_{sup} \le 1$ and
\begin{equation}
        |\lambda(f)| > \|\lambda\|_* - \epsilon.
\end{equation}
We may also ask $f$ to have compact support in $X$, since continuous
functions with compact support are dense in $C_0(X)$.  Let $\phi$ be a
continuous real-valued function on $X$ with compact support such that
$\phi(x) = 1$ for every $x$ in the support of $f$ and $0 \le \phi(x)
\le 1$ for every $x \in X$, which exists by Urysohn's lemma.  Thus
$\lambda_\phi(f) = \lambda(f)$, so that
\begin{equation}
        \|\lambda_\phi\|_* > \|\lambda\|_* - \epsilon.
\end{equation}
This implies that
\begin{equation}
        \|\lambda - \lambda_\phi\|_* = \|\lambda_{1 - \phi}\|_* < \epsilon,
\end{equation}
by (\ref{||lambda_phi||_* + ||lambda_{1 - phi}||_* = ||lambda||_*}).

\section[\ Extensions to $C_b(X)$]{Extensions to $C_b(X)$}
\label{extensions to C_b(X)}
\setcounter{equation}{0}

        Let $X$ be a locally compact Hausdorff topological space, and
let $\lambda$ be a bounded linear functional on $C_0(X)$.  As in
Section \ref{measures on R^n}, there is a natural extension of
$\lambda$ to a bounded linear functional on $C_b(X)$ with some
additional continuity properties.  Of course, this is trivial when $X$
is compact, and so we may as well suppose that $X$ is not compact.
Remember that there is a natural topology on the space $C(X)$ of all
continuous real or complex-valued functions on $X$, which is
determined by the collection of supremum seminorms associated to
nonempty compact subsets of $X$.  If $X$ is $\sigma$-compact, as in
the case of $X = {\bf R}^n$, then we have seen that it suffices to
consider the supremum seminorms corresponding to a sequence of compact
subsets of $X$, which implies that this topology on $C(X)$ is
metrizable.

        If $L$ is a nonnegative real number, then let $C_{b, L}(X)$ be
the space of continuous functions $f$ on $X$ such that $|f(x)| \le L$
for every $x \in X$.  Similarly, let $C_{0, L}(X)$ be the intersection
of $C_0(X)$ and $C_{b, L}(X)$, consisting of all continuous functions
$f$ that vanish at infinity and satisfy $\|f\|_{sup} \le L$.  It is
easy to see that $C_{0, L}(X)$ is dense in $C_{b, L}(X)$ with respect
to the topology induced on $C_{b, L}(X)$ by the one on $C(X)$
described in the preceding paragraph.  More precisely, for each
bounded continuous function $f$ on $X$ with $\|f\|_{sup} \le L$ and
every nonempty compact set $K \subseteq X$ there is a continuous
function $g$ with compact support on $X$ such that $\|g\|_{sup} \le L$
and $g(x) = f(x)$ for every $x \in K$.  To see this, one can take $g =
\theta \, f$, where $\theta$ is a continuous real-valued function on
$X$ with compact support such that $\theta(x) = 1$ for every $x \in K$
and $0 \le \theta(x) \le 1$ for every $x \in X$.

        Thus we are actually interested in extending $\lambda$ to a
bounded linear functional on $C_b(X)$ with the additional property
that the restriction of $\lambda$ to $C_{b, L}(X)$ is continuous with
respect to the topology induced by the one on $C(X)$ described before
for each $L \ge 0$.  This extension would be unique, because $C_{0,
L}(X)$ is dense in $C_{b, L}(X)$ with respect to the topology induced
by the one on $C(X)$.  If $X$ is $\sigma$-compact, then this
additional continuity condition is equivalent to asking that
$\{\lambda(f_j)\}_{j = 1}^\infty$ converges to $\lambda(f)$ for each uniformly
bounded sequence $\{f_j\}_{j = 1}^\infty$ of continuous functions on $X$ that
converges uniformly on compact subsets of $X$ to a function $f$ on $X$.
Of course, a necessary condition for the existence of an extension of
$\lambda$ to $C_b(X)$ with this additional continuity property is
that the restriction of $\lambda$ to $C_{0, L}(X)$ be continuous with
respect to the topology induced by the one on $C(X)$ for each $L \ge 0$.
It is easy to see that $\lambda$ satisfies this condition, using the
approximation of $\lambda$ by bounded linear functionals on $C_0(X)$
with compact support, as in the previous section.

        The existence of the extension of $\lambda$ to $C_b(X)$ with
this additional continuity property can be obtained by approximating a
bounded continuous function $f$ on $X$ by uniformly bounded continuous
functions $g$ on $X$ with compact support, as before, and choosing
$\lambda(f)$ so that it is approximated by the $\lambda(g)$'s.  This
is analogous to the fact that a uniformly continuous real or
complex-valued function on a dense subset of a metric space $M$ has a
unique extension to a uniformly continuous function $M$.
Alternatively, let $\{\phi_j\}_{j = 1}^\infty$ be a sequence of
uniformly bounded continuous functions on $X$ with compact support
such that the corresponding linear functionals $\lambda_{\phi_j}$
converge to $\lambda$ with respect to the dual norm associated to the
supremum norm on $C_0(X)$, as in the previous section.  Each
$\lambda_{\phi_j}$ has an obvious extension to $C_b(X)$, and one can
check that these extensions converge as $j \to \infty$ to a bounded
linear functional on $C_b(X)$.  One can then take the desired
extension of $\lambda$ to $C_b(X)$ to be the limit of this sequence,
which amounts to approximating $\lambda(f)$ for $f \in C_b(X)$ by
$\lambda(g)$ with uniformly bounded continuous functions $g$ on $X$
with compact support, as before.

\section[\ Delta functions]{Delta functions}
\label{delta functions}
\setcounter{equation}{0}

        A Dirac delta function is not really a function in the usual
sense, but can easily be interpreted as a measure on ${\bf R}^n$.
Thus if $u \in {\bf R}^n$, then the corresponding measure $\delta_u$
is defined on ${\bf R}^n$ by
\begin{eqnarray}
 \delta_u(E) & = & 1 \quad\hbox{when } u \in E \\
             & = & 0 \quad\hbox{when } u \in {\bf R}^n \backslash E. \nonumber
\end{eqnarray}
Equivalently,
\begin{equation}
        \int_{{\bf R}^n} f \, d\delta_u = f(u)
\end{equation}
for any function $f$ on ${\bf R}^n$.

        The Fourier transform of $\delta_u$ is given by
\begin{equation}
        \widehat{\delta_u}(\xi) = \exp (- i \xi \cdot u)
\end{equation}
for every $\xi \in {\bf R}^n$.  In particular,
$\widehat{\delta_0}(\xi) = 1$ for every $\xi \in {\bf R}^n$, and
$|\widehat{\delta_u}(\xi)| = 1$ for every $u, \xi \in {\bf R}^n$.
This shows that the analogue of the Riemann--Lebesgue lemma for
measures instead of integrable functions does not work.  As in Section
\ref{holomorphic extensions, continued}, there is a natural extension
of $\widehat{\delta_u}$ to a holomorphic function on ${\bf C}^n$,
given by $\widehat{\delta_u}(\zeta) = \exp (- i \zeta \cdot u)$.  If
$\epsilon \in \{-1, 1\}^n$ and $u \in Q_{n, \epsilon}$, then
$|\widehat{\delta_u}(\zeta)| \le 1$ for every $\zeta \in
\overline{H}_{n, -\epsilon}$, as before.

        If $\mu$ is a real or complex Borel measure on ${\bf R}^n$, then
\begin{equation}
\label{(mu * delta_u)(E) = mu(E - u)}
        (\mu * \delta_u)(E) = \mu(E - u)
\end{equation}
for every Borel set $E \subseteq {\bf R}^n$, where $E - u$ is the set
of points in ${\bf R}^n$ of the form $x - u$ with $x \in E$.  In
particular,
\begin{equation}
\label{mu * delta_0 = mu}
        \mu * \delta_0 = \mu
\end{equation}
for every Borel measure $\mu$, and
\begin{equation}
\label{delta_u * delta_v = delta_{u + v}}
        \delta_u * \delta_v = \delta_{u + v}
\end{equation}
for every $u, v \in {\bf R}^n$.  If $f$ is a suitable function on
${\bf R}^n$, then
\begin{equation}
        (f * \delta_u)(x) = f(x - u).
\end{equation}


\begin{thebibliography}{199}
%\begin{thebibliography}{ABC}

\addcontentsline{toc}{section}{References}


\bibitem {a} K.~Adachi, {\it Several Complex Variables and Integral
Formulas}, World Scientific, 2007.

\bibitem {ah} L.~Ahlfors, {\it Complex Analysis: An Introduction to the
Theory of Analytic Functions of One Complex Variable}, 3rd edition,
McGraw-Hill, 1978.

\bibitem {a-w} H.~Alexander and J.~Wermer, {\it Several Complex
Variables and Banach Algebras}, 3rd edition, Springer-Verlag, 1998.

\bibitem {an} M.~Andersson, {\it Topics in Complex Analysis},
Springer-Verlag, 1997.

\bibitem {a1} W.~Arveson, {\it An Invitation to $C^*$ Algebras},
Springer-Verlag, 1976.

\bibitem {a2} W.~Arveson, {\it A Short Course on Spectral Theory},
Springer-Verlag, 2002.

\bibitem {a-b-r} S.~Axler, P.~Bourdon, and W.~Ramey, {\it Harmonic
Function Theory}, 2nd edition, Springer-Verlag, 2001.

\bibitem {bac-n} G.~Bachman and L.~Narici, {\it Functional Analysis},
Dover, 2000.

\bibitem {bak-n} J.~Bak and D.~Newman, {\it Complex Analysis}, 2nd
edition, Springer-Verlag, 1997.

\bibitem {b} V.~Balachandran, {\it Topological Algebras},
North-Holland, 2000.

\bibitem {b1} R.~Beals, {\it Advanced Mathematical Analysis},
Springer-Verlag, 1973.

\bibitem {b2} R.~Beals, {\it Analysis: An Introduction}, Cambridge
University Press, 2004.

\bibitem {b-n-s} E.~Beckenstein, L.~Narici, and C.~Suffel, {\it
Topological Algebras}, North-Holland, 1977.

\bibitem {sb} S.~Berberian, {\it Lectures in Functional Analysis and
Operator Theory}, Springer-Verlag, 1974.

\bibitem {b-g-1} C.~Berenstein and R.~Gay, {\it Complex Variables: An
Introduction}, Springer-Verlag, 1991.

\bibitem {b-g-2} C.~Berenstein and R.~Gay, {\it Complex Analysis and
Special Topics in Harmonic Analysis}, Springer-Verlag, 1995.

\bibitem {bh} R.~Bhatia, {\it Notes on Functional Analysis}, Hindustan
Book Agency, 2009.

\bibitem {boc} S.~Bochner, {\it Lectures on Fourier Integrals},
translated by M.~Tenenbaum and H.~Pollard, Princeton University Press,
1959.

\bibitem {boc-c} S.~Bochner and K.~Chandrasekharan, {\it Fourier
Transforms}, Princeton University Press, 1949.

\bibitem {boc-m} S.~Bochner and W.~Martin, {\it Several Complex
Variables}, Princeton University Press, 1948.

\bibitem {br1} A.~Browder, {\it Introduction to Function Algebras},
Benjamin, 1969.

\bibitem {br2} A.~Browder, {\it Mathematical Analysis: An
Introduction}, Springer-Verlag, 1996.

\bibitem {hc} H.~Cartan, {\it Elementary Theory of Analytic Functions
of One or Several Variables}, translated from the French, Dover, 1995.

\bibitem {c-n} W.~Comfort and S.~Negrepontis, {\it The Theory of
Ultrafilters}, Springer-Verlag, 1974.

\bibitem {c1} J.~Conway, {\it Functions of One Complex Variable}, 2nd
edition, Springer-Verlag, 1978.

\bibitem {c2} J.~Conway, {\it A Course in Functional Analysis}, 2nd
edition, Springer-Verlag, 1990.

\bibitem {c3} J.~Conway, {\it Functions of One Complex Variable, II},
Springer-Verlag, 1995.

\bibitem {dlb} P.~Dolbeault, {\it Analyse Complexe}, Masson, 1990.

\bibitem {dgl} R.~Douglas, {\it Banach Algebra Techniques in Operator
Theory}, 2nd edition, Springer-Verlag, 1998.

\bibitem {duo} J.~Duoandikoetxea, {\it Fourier Analysis}, translated
and revised from the 1995 Spanish original by D.~Cruz-Uribe, SFO,
American Mathematical Society, 2001.

\bibitem {drn} P.~Duren, {\it Theory of $H^p$ Spaces}, Academic Press,
1970.

\bibitem {ehr} L.~Ehrenpreis, {\it Fourier Analysis in Several Complex
Variables}, Wiley, 1970.

\bibitem {e-m-t} Y.~Eidelman, V.~Milman, and A.~Tsolomitis, {\it
Functional Analysis: An Introduction}, American Mathematical Society,
2004.

\bibitem {fld} M.~Field, {\it Several Complex Variables and Complex
Manifolds}, I, II, Cambridge University Press, 1982.

\bibitem {sf} S.~Fisher, {\it Function Theory on Planar Domains: A
Second Course in Complex Analysis}, Wiley, 1983.

\bibitem {gf1} G.~Folland, {\it Real Analysis}, 2nd edition, Wiley,
1999.

\bibitem {gf2} G.~Folland, {\it A Guide to Advanced Real Analysis},
Mathematical Association of America, 2009.

\bibitem {f-s} J.~Forn{\ae}ss and B.~Stens{\o}nes, {\it Lectures on
Counterexamples in Several Complex Variables}, AMS Chelsea, 2007.

\bibitem {fr} M.~Fragoulopoulou, {\it Topological Algebras with
Involution}, Elsevier, 2005.

\bibitem {f-b} E.~Freitag and R.~Busam, {\it Complex Analysis}, 2nd
edition, translated from the German by D.~Fulea, Springer-Verlag,
2009.

\bibitem {f-g} K.~Fritzsche and H.~Grauert, {\it From Holomorphic
Functions to Complex Manifolds}, Springer-Verlag, 2002.

\bibitem {g1} T.~Gamelin, {\it Uniform Algebras}, Prentice-Hall, 1969.

\bibitem {g2} T.~Gamelin, {\it Uniform Algebras and Jensen Measures},
Cambridge Univerisity Press, 1978.

\bibitem {g3} T.~Gamelin, {\it Complex Analysis}, Springer-Verlag,
2001.

\bibitem {g-g} T.~Gamelin and R.~Greene, {\it Introduction to
Topology}, 2nd edition, Dover, 1999.

\bibitem {jg} J.~Garnett, {\it Bounded Analytic Functions}, revised
first edition, Springer-Verlag, 2007.

\bibitem {bg} B.~Gelbaum, {\it Modern Real and Complex Analysis},
Wiley, 1995.

\bibitem {g-j} L.~Gillman and M.~Jerison, {\it Rings of Continuous
Functions}, Springer-Verlag, 1976.

\bibitem {g-k-r} J.~Gilman, I.~Kra, and R.~Rodr\'{\i}guez, {\it
Complex Analysis}, Springer-Verlag, 2007.

\bibitem {g-p} C.~Goffman and G.~Pedrick, {\it First Course in
Functional Analysis}, Prentice-Hall, 1965.

\bibitem {rg} R.~Goldberg, {\it Methods of Real Analysis}, 2nd
edition, Wiley, 1976.

\bibitem {g-f} H.~Grauert and K.~Fritzsche, {\it Several Complex
Variables}, translated from the German, Springer-Verlag, 1976.

\bibitem {g-k} R.~Greene and S.~Krantz, {\it Function Theory of One
Complex Variable}, 3rd edition, American Mathematical Society, 2006.

\bibitem {gn} R.~Gunning, {\it Introduction to Holomorphic Functions
of Several Variables}, Volume I, {\it Function Theory}, Volume II,
{\it Local Theory}, Volume III, {\it Homological Theory}, Wadsworth \&
Brooks / Cole, 1990.

\bibitem {g-r} R.~Gunning and H.~Rossi, {\it Analytic Functions of
Several Complex Variables}, AMS Chelsea, 2009.

\bibitem {vh} V.~Hansen, {\it Functional Analysis: Entering Hilbert
Space}, World Scientific, 2006.

\bibitem {h-l-1} G.~Henkin and J.~Leiterer, {\it Theory of Functions
on Complex Manifolds}, Birkh\"auser, 1984.

\bibitem {h-l-2} G.~Henkin and J.~Leiterer, {\it Andreotti--Grauert
Theory by Integral Formulas}, Birkh\"auser, 1988.

\bibitem {h-s} E.~Hewitt and K.~Stromberg, {\it Real and Abstract
Analysis}, Springer-Verlag, 1975.

\bibitem {h-y} J.~Hocking and G.~Young, {\it Topology}, 2nd edition,
Dover, 1988.

\bibitem {kh} K.~Hoffman, {\it Banach Spaces of Analytic Functions},
Dover, 1988.

\bibitem {h1} L.~H\"ormander, {\it An Introduction to Complex Analysis
in Several Variables}, 3rd edition, North-Holland, 1990.

\bibitem {h2} L.~H\"ormander, {\it Notions of Convexity},
Birkh\"auser, 2007.

\bibitem {jh} J.~Horv\'ath, {\it Topological Vector Spaces and
Distributions}, Addison-Wesley, 1966.

\bibitem {j1} I.~James, {\it Introduction to Uniform Spaces},
Cambridge University Press, 1990.

\bibitem {j2} I.~James, {\it Topologies and Uniformities},
Springer-Verlag, 1999.

\bibitem {j-p} M.~Jarnicki and P.~Pflug, {\it First Steps in Several
Complex Variables: Reinhardt Domains}, European Mathematical Society,
2008.

\bibitem {fj} F.~Jones, {\it Lebesgue Integration on Euclidean
Spaces}, Jones and Bartlett, 1993.

\bibitem {khn} J.-P.~Kahane, {\it S\'eries de Fourier Absolument
Convergentes}, Springer-Verlag, 1970.

\bibitem {kan} E.~Kaniuth, {\it A Course in Commutative Banach
Algebras}, Springer-Verlag, 2009.

\bibitem {kz} S.~Kantorovitz, {\it Introduction to Modern Analysis},
Oxford University Press, 2003.

\bibitem {kap} I.~Kaplansky, {\it Set Theory and Metric Spaces}, 2nd
edition, Chelsea, 1977.

\bibitem {kat} Y.~Katznelson, {\it An Introduction to Harmonic
Analysis}, 3rd edition, Cambridge University Press, 2004.

\bibitem {jk} J.~Kelley, {\it General Topology}, Springer-Verlag,
1975.

\bibitem {k-n} J.~Kelley, I.~Namioka, et al., {\it Linear Topological
Spaces}, Springer-Verlag, 1976.

\bibitem {k-s} J.~Kelley and T.~Srinivasan, {\it Measure and
Integral}, Springer-Verlag, 1988.

\bibitem {ke} S.~Kesevan, {\it Functional Analysis}, Hindustan Book
Agency, 2009.

\bibitem {ak1} A.~Knapp, {\it Basic Real Analysis}, Birkh\"auser,
2005.

\bibitem {ak2} A.~Knapp, {\it Advanced Real Analysis}, Birkh\"auser,
2005.

\bibitem {kod} K.~Kodaira, {\it Complex Analysis}, translated from the
1977 Japanese original by A.~Sevenster, edited by A.~Beardon and
T.~Carne, Cambridge University Press, 2007.

\bibitem {pk} P.~Koosis, {\it Introduction to $H_p$ Spaces}, 2nd
edition, with two appendices by V.~Havin, Cambridge University Press,
1998.

\bibitem {k1} G.~K\"othe, {\it Topological Vector Spaces, I},
translated from the German by D.~Garling, Springer-Verlag, 1969.

\bibitem {k2} G.~K\"othe, {\it Topological Vector Spaces, II},
Springer-Verlag, 1979.

\bibitem {sk1} S.~Krantz, {\it A Panorama of Harmonic Analysis},
Mathematical Association of America, 1999.

\bibitem {sk2} S.~Krantz, {\it Function Theory of Several Complex
Variables}, AMS Chelsea, 2001.

\bibitem {sk3} S.~Krantz, {\it Complex Analysis: The Geometric
Viewpoint}, 2nd edition, Mathematical Association of America, 2004.

\bibitem {sk4} S.~Krantz, {\it Real Analysis and Foundations}, 2nd
edition, Chapman \& Hall / CRC, 2005.

\bibitem {sk5} S.~Krantz, {\it Geometric Function Theory: Explorations
in Complex Analysis}, Birkh\"auser, 2006.

\bibitem {sk6} S.~Krantz, {\it A Guide to Complex Variables},
Mathematical Association of America, 2008.

\bibitem {sk7} S.~Krantz, {\it A Guide to Real Variables}, Mathematical
Association of America, 2009.

\bibitem {sk8} S.~Krantz, {\it A Guide to Topology}, Mathemtical
Association of America, 2009.

\bibitem {sk9} S.~Krantz, {\it Essentials of Topology with
Applications}, CRC Press, 2010.

\bibitem {k-p-1} S.~Krantz and H.~Parks, {\it A Primer of Real
Analytic Functions}, 2nd edition, Birkh\"auser, 2002.

\bibitem {k-p-2} S.~Krantz and H.~Parks, {\it The Implicit Function
Theorem: History, Theory, and Applications}, Birkh\"auser, 2002.

\bibitem {bl} B.~Lahiri, {\it Elements of Functional Analysis}, World
Press, 2005.

\bibitem {sl1} S.~Lang, {\it Real and Functional Analysis}, 3rd
edition, Springer-Verlag, 1993.

\bibitem {sl2} S.~Lang, {\it Undergraduate Analysis}, 2nd edition,
Springer-Verlag, 1997.

\bibitem {sl3} S.~Lang, {\it Complex Analysis}, 4th edition,
Springer-Verlag, 1999.

\bibitem {rl} R.~Larsen, {\it Functional Analysis: An Introduction},
Dekker, 1973.

\bibitem {lax} P.~Lax, {\it Functional Analysis}, Wiley, 2002.

\bibitem {l-g} P.~Lelong and L.~Gruman, {\it Entire Functions of
Several Complex Variables}, Springer-Verlag, 1986.

\bibitem {l-r} D.~Luecking and L.~Rubel, {\it Complex Analysis: A
Functional Analysis Approach}, Springer-Verlag, 1984.

\bibitem {mac} B.~MacCluer, {\it Elementary Functional Analysis},
Springer-Verlag, 2009.

\bibitem {mdx} I.~Maddox, {\it Elements of Functional Analysis}, 2nd
edition, Cambridge University Press, 1988.

\bibitem {am} A.~Mallios, {\it Topological Algebras: Selected Topics},
North-Holland, 1986.

\bibitem {m-v} R.~Meise and D.~Vogt, {\it Introduction to Functional
Analysis}, translated from the German by M.~Ramanujan and revised by
the authors, Oxford University Press, 1997.

\bibitem {mn} B.~Mendelson, {\it Introduction to Topology}, 3rd
edition, Dover, 1990.

\bibitem {tm} T.~Morrison, {\it Functional Analysis: An Introduction
to Banach Space Theory}, Wiley, 2001.

\bibitem {ln1} L.~Nachbin, {\it Holomorphic Functions, Domains of
Holomorphy, and Local Properties}, notes prepared by R.~Aron,
North-Holland, 1970.

\bibitem {ln2} L.~Nachbin, {\it Introduction to Functional Analysis:
Banach Spaces and Differential Calculus}, translated from the
Portuguese by R.~Aron, Dekker, 1981.

\bibitem {rn} R.~Narasimhan, {\it Several Complex Variables},
University of Chicago Press, 1995.

\bibitem {n-n} R.~Narasimhan and Y.~Nievergelt, {\it Complex Analysis
in One Variable}, 2nd edition, Birkh\"auser, 2001.

\bibitem {n} N.~Nikolski, {\it Operators, Functions, and Systems: An
Easy Reading}, Volumes 1, 2, translated from the French by A.~Hartmann
and revised by the author, American Mathematical Society, 2002.

\bibitem {tn} T.~Nishino, {\it Function Theory in Several Complex
Variables}, translated from the 1996 Japanese original by N.~Levenberg
and H.~Yamagchi, American Mathematical Society, 2001.

\bibitem {jn} J.~Noguchi, {\it Introduction to Complex Analysis},
translated from the 1993 Japanese original by the author, American
Mathematical Society, 1998.

\bibitem {oh} T.~Ohsawa, {\it Analysis of Several Complex Variables},
translated from the Japanese by S.~Nakamura, American Mathematical
Society, 2002.

\bibitem {p} E.~Packel, {\it Functional Analysis: A Short Course},
Krieger, 1980.

\bibitem {wp} W.~Page, {\it Topological Uniform Structures}, Dover,
1988.

\bibitem {p-w} R.~Paley and N.~Wiener, {\it Fourier Transforms in the
Complex Domain}, American Mathematical Society, 1987.

\bibitem {vp} V.~Peller, {\it Hankel Operators and their
Applications}, Springer-Verlag, 2003.

\bibitem {pr} S.~Promislow, {\it A First Course in Functional
Analysis}, Wiley, 2008.

\bibitem {r} M.~Range, {\it Holomorphic Functions and Integral
Representations in Several Complex Variables}, Springer-Verlag, 1986.

\bibitem {cr1} C.~Rickart, {\it General Theory of Banach Algebras},
van Nostrand, 1960.

\bibitem {cr2} C.~Rickart, {\it Natural Function Algebras},
Springer-Verlag, 1979.

\bibitem {r-r} A.~Robertson and W.~Robertson, {\it Topological Vector
Spaces}, 2nd edition, Cambridge University Press, 1980.

\bibitem {hr} H.~Royden, {\it Real Analysis}, 3rd edition, Macmillan,
1988.

\bibitem {mer} M.~Rudin, {\it Lectures on Set Theoretic Topology},
American Mathematical Society, 1975.

\bibitem {r1} W.~Rudin, {\it Function Theory in Polydisks}, Benjamin,
1969.

\bibitem {r2} W.~Rudin, {\it Lectures on the Edge-of-the-Wedge
Theorem}, American Mathematical Society, 1971.

\bibitem {r3} W.~Rudin, {\it Principles of Mathematical Analysis},
3rd edition, McGraw-Hill, 1976.

\bibitem {r4} W.~Rudin, {\it Real and Complex Analysis}, 3rd edition,
McGraw-Hill, 1987.

\bibitem {r5} W.~Rudin, {\it Fourier Analysis on Groups}, Wiley, 1990.

\bibitem {r6} W.~Rudin, {\it Functional Analysis}, 2nd edition,
McGraw-Hill, 1991.

\bibitem {r7} W.~Rudin, {\it Function Theory on the Unit Ball in ${\bf
C}^n$}, Springer-Verlag, 2008.

\bibitem {r-y} B.~Rynne and M.~Youngson, {\it Linear Functional
Analysis}, 2nd edition, Springer-Verlag, 2008.

\bibitem {ds1} D.~Sarason, {\it Function Theory on the Unit Circle},
Department of Mathematics, Virginia Polytchnic Institute and State
University, 1978.

\bibitem {ds2} D.~Sarason, {\it Complex Function Theory}, 2nd edition,
American Mathematical Society, 2007.

\bibitem {ks} K.~Saxe, {\it Beginning Functional Analysis},
Springer-Verlag, 2002.

\bibitem {sch-w} H.~Schaefer and M.~Wolff, {\it Topological Vector
Spaces}, 2nd edition, Springer-Verlag, 1999.

\bibitem {sch} M.~Schechter, {\it Principles of Functional Analysis},
2nd edition, American Mathematical Society, 2002.

\bibitem {scn} V.~Scheidemann, {\it Introduction to Complex Analysis
in Several Variables}, Birkh\"auser, 2005.

\bibitem {gs} G.~Simmons, {\it Introduction to Topology and Modern
Analysis}, Krieger, 1983.

\bibitem {st1} E.~Stein, {\it Singular Integrals and Differentiability
Properties of Functions}, Princeton University Press, 1970.

\bibitem {st2} E.~Stein, {\it Boundary Behavior of Holomorphic
Functions of Several Complex Variables}, Princeton University Press,
1972.

\bibitem {st3} E.~Stein, {\it Harmonic Analysis: Real-Variable
Methods, Orthogonality, and Oscillatory Integrals}, with the
assistance of T.~Murphy, Princeton University Press, 1993.

\bibitem {st-s1} E.~Stein and R.~Shakarchi, {\it Fourier Analysis: An
Introduction}, Princeton University Press, 2003.

\bibitem {st-s2} E.~Stein and R.~Shakarchi, {\it Complex Analysis},
Princeton University Press, 2003.

\bibitem {st-s3} E.~Stein and R.~Shakarchi, {\it Real Analysis:
Measure Theory, Integration, and Hilbert Spaces}, Princeton University
Press, 2005.

\bibitem {s-w} E.~Stein and G.~Weiss, {\it Introduction to Fourier
Analysis on Euclidean Spaces}, Princeton University Press, 1971.

\bibitem {es1} E.~Stout, {\it The Theory of Uniform Algebras}, Bogden
\& Quigley, 1971.

\bibitem {es2} E.~Stout, {\it Polynomial Convexity}, Birkh\"auser,
2007.

\bibitem {rs1} R.~Strichartz, {\it The Way of Analysis}, Jones and
Bartlett, 1995.

\bibitem {rs2} R.~Strichartz, {\it A Guide to Distribution Theory and
Fourier Transforms}, World Scientific, 2003.

\bibitem {str} K.~Stromberg, {\it Introduction to Classical Real
Analysis}, Wadsworth, 1981.

\bibitem {cs1} C.~Swartz, {\it An Introduction to Functional
Analysis}, Dekker, 1992.

\bibitem {cs2} C.~Swartz, {\it Elementary Functional Analysis}, World
Scientific, 2009.

\bibitem {jt} J.~Taylor, {\it Several Complex Variables with
Connections to Algebraic Geometry and Lie Groups}, American
Mathematical Society, 2002.

\bibitem {mt} M.~Taylor, {\it Measure Theory and Integration},
American Mathematical Society, 2006.

\bibitem {t1} A.~Torchinsky, {\it Real Variables}, Addison-Wesley,
1988.

\bibitem {t2} A.~Torchinsky, {\it Real-Variable Methods in Harmonic
Analysis}, Dover, 2004.

\bibitem {t} F.~Tr\`eves, {\it Topological Vector Spaces,
Distributions, and Kernels}, Dover, 2006.

\bibitem {rw} R.~Walker, {\it The Stone--{\v C}ech Compactification},
Springer-Verlag, 1974.

\bibitem {w-z} R.~Wheeden and A.~Zygmund, {\it Measure and Integral:
An Introduction to Real Analysis}, Dekker, 1977.

\bibitem {w} A.~Wilansky, {\it Modern Methods in Topological Vector
Spaces}, McGraw-Hill, 1978.

\bibitem {y} K.~Yosida, {\it Functional Analysis}, Springer-Verlag,
1995.

\bibitem {ze} A.~Zemanian, {\it Distribution Theory and Transform
Analysis}, 2nd edition, Dover, 1987.

\bibitem {zy} A.~Zygmund, {\it Trigonometric Series}, Volumes I, II,
3rd edition, with a foreword by R.~Fefferman, Cambridge University
Press, 2002.





\end{thebibliography}
\end{document}